\documentclass[reqno,11pt]{amsart}

\usepackage[top=2.0cm,bottom=2.0cm,left=3cm,right=3cm]{geometry}
\usepackage{amsthm,amsmath,amssymb,dsfont}
\usepackage{mathrsfs,amsfonts,functan,extarrows,mathtools}
\usepackage[colorlinks]{hyperref}

\usepackage{marginnote}
\usepackage{xcolor}
\usepackage{stmaryrd}
\usepackage{esint}
\usepackage{graphicx}
\usepackage{bbm}

\usepackage{tikz}
\usetikzlibrary{calc,decorations.pathreplacing}

\newtheorem{theorem}{Theorem}[section]

\newtheorem{definition}{Definition}[section]

\newtheorem{remark}{Remark}[section]
\newtheorem{lemma}{Lemma}[section]

\numberwithin{equation}{section}
\allowdisplaybreaks

\def\p{\partial}
\def\d{\mathrm{d}}
\def\no{\nonumber}
\def\R{\mathbb{R}}
\def\eps{\varepsilon}
\def\div{\mathrm{div}}
\def\u{\mathbf{u}}

\def\exp{\mathrm{exp}}
\def\M{\mathfrak{M}}

\def\B{\mathcal{B}}

\def\u{\mathfrak{u}}
\def\L{\mathcal{L}}

\def\P{\mathcal{P}}
\def\I{\mathcal{I}}

\def\j{\mathfrak{j}}
\def\e{\mathfrak{e}}
\def\h{\mathfrak{h}}

\newcounter{wronumber}

\begin{document}
\title[Decay estimates of pseudo-inverse $\mathcal{L}^{-1}$]
			{Grad-Caflisch pointwise decay estimates revisited}

\author[Ning Jiang]{Ning Jiang}
\address[Ning Jiang]{\newline School of Mathematics and Statistics, Wuhan University, Wuhan, 430072, P. R. China}
\email{njiang@whu.edu.cn}

\author[Yi-Long Luo]{Yi-Long Luo${}^*$}
\address[Yi-Long Luo]
{\newline School of Mathematics, South China University of Technology, Guangzhou, 510641, P. R. China}
\email{luoylmath@scut.edu.cn}

\author[Shaojun Tang]{Shaojun Tang}
\address[Shaojun Tang]
		{\newline Department of Mathematics, Wuhan University of Technology, Wuhan, 430070, P. R. China }
\email{shaojun.tang@whut.edu.cn}
\thanks{${}^*$ Corresponding author \quad \today}

\maketitle

\begin{abstract}
   In the influential paper \cite{Caflish-1980-CPAM} which was the starting point of the employment of Hilbert expansion method to the rigorous justifications of the fluid limits of the Boltzmann equation, Caflisch discovered an elegant and crucial estimate on each expansion term (Proposition 3.1 in \cite{Caflish-1980-CPAM}). The proof essentially relied on an estimate of Grad \cite{Grad-1963}, which was on the pointwise decay properties of $\L^{-1}$, the pseudo-inverse operator of the linearized Boltzmann collision operator $\L$, for the hard potential collision kernel, i.e. the power $0\leq \gamma\leq 1$. Caflisch's arguments need the exponential version of Grad's estimate. However, Grad's original paper was only on the polynomial decay. In this paper, we revisit and provide a full proof of the Caflisch-Grad type decay estimates and the corresponding applications in the compressible Euler limit of the Boltzmann equaiton. The main novelty is that for the case collision kernel power $-\frac{3}{2}<\gamma\leq 1$, the proof of the pointwise estimate does not use any derivatives. So the potential applications of this estimate could be wider than in the Hilbert expansion. For the completeness of the result, we also prove the almost everywhere pointwise estimate using derivatives for the case $-3<\gamma\leq  -\frac{3}{2}$. Furthermore, in the application to fluid limits, $\L^{-1}$ and the derivatives with respect to the parameters (for example, $(t,x)$, this must happen when $\L$ is linearized around local Maxwellian which depends on $(t,x)$) are not commutative. We detailed analyze the estimate of commutators, which was missing in previous literatures of fluid limits of the Boltzmann equation. This estimate is needed in all compressible fluid limits from Boltzmann equation. \\


\end{abstract}





\section{Introduction}

\subsection{The Boltzmann collision operator}

The well-known Boltzmann collision operator $\mathcal{B} (F, F)$ is defined by
\begin{equation}
	\begin{aligned}
		\mathcal{B} (F, F) = \iint_{\mathbb{R}^3 \times \mathbb{S}^2} (F' F_1' - F F_1) b (\omega, v_1 - v) \d \omega \d v_1 \,.
	\end{aligned}
\end{equation}
Here $\omega \in \mathbb{S}^2$ is a unit vector, $\d \omega$ is the rotationally invariant surface ingetral on $\mathbb{S}^2$, while $F_1'$, $F'$, $F_1$ and $F$ are the number density $F(\cdot)$ evaluated at the velocities $v_1'$, $v'$, $v_1$ and $v$ respectively, i.e.
\begin{equation*}
	\begin{aligned}
		F_1' = F(v_1') \,, \ F' = F(v') \,, \ F_1 = F(v_1) \,, \ F = F(v) \,.
	\end{aligned}
\end{equation*}
Here $(v_1', v')$ are the velocities after an elastic binary collision between two molecules with velocities $(v_1, v)$ before the collision, or vice versa.
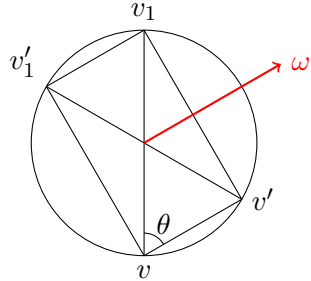
\begin{figure}[h]
	\begin{center}
		\begin{tikzpicture}
			\pgfmathsetmacro\rd{1.5};
			\draw (\rd,0) arc (0: 360: \rd);
			\draw (0, -\rd)--(-30:\rd)node[right]{$v'$}--(90:\rd)--(150:\rd)node[above left]{$v_1'$}--cycle;
			\draw (0, -\rd) node[below]{$v$}--(0,\rd)node[above]{$v_1$};
			\draw (0, -\rd*0.8) arc (90:30:0.2*\rd)node[above]{$\theta$};
			\draw (-30:\rd)--(150:\rd);
			\draw[->,thick,color=red] (0,0)--(30:1.4*\rd) node[right]{$\omega$};
		\end{tikzpicture}
	\end{center}
	\caption{Transformation of velocities.}\label{Fig1}
\end{figure}
Since both momentum and energy are conserved during the elastic collision, $v_1'$ and $v'$ can be expressed in terms of $v_1$ and $v$ as
\begin{equation}\label{Colli-v}
	\begin{aligned}
		v_1' = v_1 - [ (v_1 - v) \cdot \omega ] \omega \,, \ v' = v + [ (v_1 - v) \cdot \omega ] \omega \,,
	\end{aligned}
\end{equation}
where the unit vector $\omega \in \mathbb{S}^2$ is parallel to the deflections $v_1' - v_1$ and $v'-v$, and is therefore perpendicular to the plane of reflection, see Figure \ref{Fig1}. In the collision term $\mathcal{B} (F, F)$, $F_1'F'$ is the gain term, while $F_1F$ is called the loss term.

H. Grad \cite{Grad-1963} suggested the so-called cut-off kernels, i.e. the collision kernel $b (\omega, v_1 - v)$ with the following factored form
\begin{equation}\label{b}
	\begin{aligned}
		b (\omega, v_1 - v) = \beta (\theta) |v_1 - v|^\gamma \,, \ \cos \theta = \frac{(v_1 - v) \cdot \omega}{|v_1 - v|} \,, \ -3 < \gamma \leq 1 \,,
	\end{aligned}
\end{equation}
where $\beta (\theta)$ satisfies the symmetric assumption
\begin{equation}\label{Symmetric-Asmp}
	\begin{aligned}
		\beta (\theta) \sin \theta = \beta (\tfrac{\pi}{2} - \theta) \sin (\tfrac{\pi}{2} - \theta) \,, \ \beta (\theta) \sin \theta = \beta (\pi - \theta) \sin (\pi - \theta) \,,
	\end{aligned}
\end{equation}
and the small deflection cutoff condition
\begin{equation}\label{Cutoff}
	\begin{aligned}
		\int_{\mathbb{S}^2} \beta (\theta) \d \omega = 1 \,, \ 0 \leq \beta (\theta) \leq \beta_0 |\cos \theta|
	\end{aligned}
\end{equation}
for some constant $\beta_0 > 0$. The cases $-3 < \gamma < 0$ and $0 \leq \gamma \leq 1$ are respectively referred to as the ``soft" and ``hard" potential cases. In particular, $\gamma = 0$ is the Maxwell potential case, and $\gamma = 1$ is the hard sphere case, in which $\beta (\theta) = \beta_0 |\cos \theta|$.

The Maxwellians are those number densities $\M$ such that $\mathcal{B}(\M, \M)=0$. Boltzmann's $H\mbox{-}$theorem tells that the Maxwellians must have the form \begin{equation}\label{Maxwellian}
	\begin{aligned}
		\M (v) = \M_{[\rho, \u, T]} (v) = \frac{\rho}{(2 \pi T)^\frac{3}{2}} \exp \Big( - \frac{|v - \u|^2}{2 T} \Big) \,,
	\end{aligned}
\end{equation}
where in the fluid state $(\rho, \u, T)$, $\rho > 0$ is the density, $T > 0$ is the temperature, and $\u \in \R^3$ is the bulk velocity. When $(\rho, \u, T)$ depends on the time $t>0$ and space variable $x$, $\M$ is called the local Maxwellian. Otherwise, it is called the global Maxwellian. In this paper, we consider the local Maxwellians, which are more general.

For the perturbation around $\M$, i.e. $F = \M + \sqrt{\M} f$, we consider the linearized Boltzmann collision operator $\L$ defined by
\begin{equation}\label{Lf}
	\begin{aligned}
		\L f = & - \M^{- \frac{1}{2}} (v) \big[ \mathcal{B} (\M, \sqrt{\M} f) + \mathcal{B} ( \sqrt{\M} f , \M ) \big] \\
		= & \nu (v) f (v) + K f (v) \,,
	\end{aligned}
\end{equation}
where the collision frequency $\nu (v)$ is defined by
\begin{equation}\label{nu}
	\begin{aligned}
		\nu (v) = \iint_{\mathbb{R}^3 \times \mathbb{S}^2} \M (v_1) b (\omega, v_1 - v) \d \omega \d v_1 \overset{\eqref{Cutoff}}{=} \int_{\R^3} |v_1 - v|^\gamma \M (v_1) \d v_1 \,,
	\end{aligned}
\end{equation}
and the operator $K f (v)$ can be decomposed into two parts:
\begin{equation}
	\begin{aligned}\label{K-K1-K2}
		K f (v) = K_1 f (v) - K_2 f (v) \,.
	\end{aligned}
\end{equation}
Here the loss term $K_1 f (v)$ is
\begin{equation}\label{K1}
	\begin{aligned}
		K_1 f(v) = \M^\frac{1}{2} (v) \iint_{\mathbb{R}^3 \times \mathbb{S}^2} f (v_1) \M^\frac{1}{2} (v_1) b (\omega, v_1 - v) \d \omega \d v_1 \,,
	\end{aligned}
\end{equation}
and the gain term $K_2 f (v)$ is
\begin{equation}\label{K2}
	\begin{aligned}
		K_2 f(v) = & \M^\frac{1}{2} (v) \iint_{\mathbb{R}^3 \times \mathbb{S}^2} \big[ \M^{- \frac{1}{2}} (v') f (v') + \M^{- \frac{1}{2}} (v_1') f(v_1') \big] \M (v_1) b (\omega, v_1 - v) \d \omega \d v_1 \\
		= & 2 \M^\frac{1}{2} (v) \iint_{\mathbb{R}^3 \times \mathbb{S}^2} \M^{- \frac{1}{2}} (v') f (v') \M (v_1) b (\omega, v_1 - v) \d \omega \d v_1 \,,
	\end{aligned}
\end{equation}
where the last equality is derived from the symmetric assumption \eqref{Symmetric-Asmp}, see Section 2 of \cite{Grad-1963} for instance.

The null space $\mathrm{Null} (\L)$ of the operator $\L$ is spanned by the basis
\begin{equation*}
	\begin{aligned}
		\psi_0 = \frac{1}{\sqrt{\rho}} \sqrt{\M} \,, \ \psi_i =  \frac{v_i - \u_i}{\sqrt{\rho T}} \sqrt{\M} \,\quad\! (i = 1,2,3) \,, \ \psi_4 = \frac{1}{\sqrt{6 \rho}} \left( \frac{|v - \u|^2}{T} - 3 \right) \sqrt{\M} \,.
	\end{aligned}
\end{equation*}
The basis is orthonormal in the Hilbert space $L^2 : = L^2 (\R^3)$. Let $\mathrm{Null}^\perp (\L)$ be the orthogonal space of $\mathrm{Null} (\L)$ in $L^2$, namely,
\begin{equation*}
	\begin{aligned}
		\mathrm{Null} (\L) \oplus \mathrm{Null}^\perp (\L) = L^2 \,.
	\end{aligned}
\end{equation*}
We define the projection $\P : L^2 \to \mathrm{Null} (\L)$ by
\begin{equation}\label{Projt-P}
	\begin{aligned}
		\P f = \big\{ \tfrac{\rho_f}{\rho} + u_f \cdot \tfrac{v - \u}{T} + \tfrac{\theta_f}{6 T} ( \tfrac{|v - \u|^2}{T} - 3 ) \big\} \sqrt{\M} \,,
	\end{aligned}
\end{equation}
where
\begin{equation*}
	\begin{aligned}
		\rho_f = \int_{\R^3} f \sqrt{\M} \d v \,, \ u_f = \frac{1}{\rho} \int_{\R^3} f (v - \u) \sqrt{\M} \d v \,, \ \theta_f = \frac{1}{\rho} \int_{\R^3} f (|v - \u|^2 - 3 T) \sqrt{\M} \d v \,.
	\end{aligned}
\end{equation*}

It is a classical result that (for the proof, see Proposition 2.3 of \cite{Golse-SRM-05} for instance) $\{ \psi_\alpha \M^{- \frac{1}{2}} : \alpha = 0, 1, \cdots, 4 \}$ forms the the basis of all collision invariants, i.e.,
$$ \langle \mathcal{B} (F, F), \psi_\alpha \M^{- \frac{1}{2}} \rangle = 0 \,, \  \alpha = 0, 1, \cdots, 4 \,.$$
Here the bracket $\langle \cdot, \cdot \rangle$ represents the $L^2$-inner product. Moreover, the Boltzmann $H$-theorem states that $\mathcal{B} (F, F) = 0$ is equivalent to $F = \M$. It thereby infers that
$$ \L f \in \mathrm{Null}^\perp (\L)  $$
for any $f \in L^2$. Furthermore, the restricted operator $\L |_{\mathrm{Null}^\perp (\L)}$ is a one-to-one map from $\mathrm{Null}^\perp (\L)$ onto $\mathrm{Null}^\perp (\L)$, because $K = \L - \nu $ is a Fredholm operator. Thus, one gives the following definition.

\begin{definition}[Pseudo-inverse operator of $\L$]\label{Df-L-1}
	The inverse operator
	\begin{equation*}
		\begin{aligned}
			( \L |_{\mathrm{Null}^\perp (\L)} )^{-1} : \mathrm{Null}^\perp (\L) \to \mathrm{Null}^\perp (\L)
		\end{aligned}
	\end{equation*}
    is called the pseudo-inverse operator of the linearized Boltzmann collision operator $\L$, which is briefly denoted by $\L^{-1}$.
\end{definition}

By \cite{Golse-SRM-05} for instance, the linearized Boltzmann collision operator $\L$ satisfies the well-known hypocoercivity
\begin{equation}\label{Hypo-L}
	\begin{aligned}
		\langle \L f, f \rangle \geq \lambda \| f \|^2_{L^2(\nu)}
	\end{aligned}
\end{equation}
for some constant $\lambda = \lambda (\rho, \u, T) > 0$, where $\| \cdot \|_{L^2 (\nu)} = ( \int_{\R^3} |\cdot |^2 \nu (v) \d v )^\frac{1}{2}$ is the norm of the weighted $L^2$ space. Furthermore, by \cite{LS-2010-KRM}, the collision frequency $\nu (v)$ defined in \eqref{nu} subjects to
\begin{equation}\label{nu1}
	\begin{aligned}
		C_1 (\rho, \u, T) \langle v \rangle^\gamma \leq \nu (v) \leq C_2 (\rho, \u, T) \langle v \rangle^\gamma
	\end{aligned}
\end{equation}
for some positive constants $C_i (\rho, \u, T) > 0$ ($i = 1,2$), where the bracket $\langle v \rangle$ is defined as $\langle v \rangle = 1 + |v|$. For notational simplicity, \eqref{nu1} will be briefly denoted by
\begin{equation}\label{nu-equivalent}
	\begin{aligned}
		\nu (v) \thicksim_{\rho, \u, T} \langle v \rangle^\gamma \,.
	\end{aligned}
\end{equation}
It is noted that the collision frequency $\nu (v)$ is monotone in $|v|$. More specifically,
\begin{equation}
	\begin{aligned}
		\tfrac{\partial \nu (v)}{\partial |v|} \geq 0 \ \textrm{for } 0 \leq \gamma \leq 1 \,, \textrm{and } \tfrac{\partial \nu (v)}{\partial |v|} < 0 \ \textrm{for } - 3 < \gamma < 0 \,.
	\end{aligned}
\end{equation}
The details can be referred to Section 3 of \cite{Grad-1963}.

\subsection{Motivations and decay results of  pseudo-inverse operator}

One of the goals of this paper is to investigate the decay (with respect to $v$) estimates of the pseudo-inverse operator $\L^{-1}$. This decay estimate has been widely applied to hydrodynamic limits of the Boltzmann equations (for examples, \cite{Guo-CPAM06, Guo-Jang-Jiang-2010-CPAM, GHW-2020}). Specifically, in any fluid limits of kinetic equations, once the Hilbert expansion approach is employed, the solution will be formally written as $F_\kappa = \sum_{n = 0}^\infty \kappa^n F_n$, where the kinetic part of $F_n$ is solved from $\L F_n= G_n$. Any estimates of $F_n$ must rely on the decay estimate of $\L^{-1}$. We emphasize that rigorous treatment of the Hilbert expansion approach to the fluid limits of the Boltzmann equation started from Caflisch's work \cite{Caflish-1980-CPAM}.

In Caflisch's paper \cite{Caflish-1980-CPAM}, Page 656, the following result for the decay of $\mathscr{L}^{-1}$ for the hard potential cases $0 \leq \gamma \leq 1$ was stated: \footnote{The symbol $\xi$ in \cite{Caflish-1980-CPAM} corresponds to the velocity $v$ in current work. Moreover, the linearized operator $\mathscr{L}$ in \cite{Caflish-1980-CPAM} equals to $\M^\frac{1}{2} \L$ that we employed.}
\begin{equation*}
	``\textrm{\em Grad has shown in \cite{Grad-1963} that } \mathscr{L}^{-1} \textrm{ \em preserves decay in } \xi, \textrm{\em so that } \Psi_1 \thicksim |\xi|^3 \exp \{ - |\xi - \u|^2 / 2 T \} \,."
\end{equation*}
Here the function $\Psi_1 = \mathscr{L}^{-1} g$ for some determined $g $ given in Page 655 of \cite{Caflish-1980-CPAM}. In our notations, the Caflisch's assertion is equivalent to the following pointwise estimate:
\begin{equation}\label{Caflisch-inq}
	\begin{aligned}
		| \L^{-1} g | \lesssim \langle v \rangle^3 \M^\frac{1}{2} (v) \| \langle v \rangle^3 \M^{- \frac{1}{2}} g \|_{L^\infty}
	\end{aligned}
\end{equation}
for $g \in \mathrm{Null}^\perp (\L)$, provided that $\| \langle v \rangle^3 \M^{- \frac{1}{2}} g \|_{L^\infty} < \infty$. This result has been widely employed in many literatures later. However, this estimate actually was {\em not} implied by Grad's work \cite{Grad-1963}, as stated in \cite{Caflish-1980-CPAM}. Indeed, Grad \cite{Grad-1963} only proved the following weaker inequality than \eqref{Caflisch-inq}:
\begin{equation}\label{Grad-inq}
	\begin{aligned}
		|\L^{-1} g (v)| \lesssim \langle v \rangle^{-m} \| \langle v \rangle^m g \|_{L^\infty} \ (\forall v \in \R^3)
	\end{aligned}
\end{equation}
for any fixed integer $m > 0$, provided that $\| \langle v \rangle^m g \|_{L^\infty} < \infty$. The key difference between the estimates \eqref{Caflisch-inq} and \eqref{Grad-inq} is: the former is an exponential decay, while the latter is a polynomial decay.

In order to illustrate this, we first display the correspondence of notations between Grad's work \cite{Grad-1963} and the current paper as follows:
\begin{center}
\begin{tabular}{|l|c|c|c|c|c|c|}
	\hline
    Notations in Grad's work \cite{Grad-1963} & $L$ & $\hat{f}$ & $K$ & $\xi$ & $\hat{\textrm{g}}$ & $\xi^2$ \\
    Notations in current paper & $\L$ & $f$ & $- K$ & $v$ & $g$ & $|v|^2$ \\
	\hline
\end{tabular}
\end{center}
As shown in (111) and (113), Page 50 of \cite{Grad-1963},
\begin{equation}\label{Grad-1}
	\begin{aligned}
		L [\hat{f}] = \nu \hat{f} - K [\hat{f}] = \hat{\textrm{g}} \,, \ |\hat{f}| = \tfrac{1}{\nu} | \hat{\textrm{g}} + K [\hat{f}] | < \tfrac{1}{\nu_0} ( |\hat{\textrm{g}} | + |K [\hat{f}] | ) \,,
	\end{aligned}
\end{equation}
where $\nu > \nu_0 > 0$ for hard potential is employed. The inequalities (70) in Page 41 and (78)-(79) in Page 43 of \cite{Grad-1963} read
\begin{equation}\label{Grad-2}
	\begin{aligned}
		N_1 [K[\hat{f}]] < k N_2 [\hat{f}] \,, \ N_3^{(r + 1/2)} [K[\hat{f}]] < k N_3^{(r)} [\hat{f}]
	\end{aligned}
\end{equation}
for any $r \geq 0$ and some constant $k > 0$, where the quantities $N_1 [\hat{f}]$, $N_2 [\hat{f}]$ and $N_3^{(r)} [\hat{f}]$ are defined in (75), Page 42 of \cite{Grad-1963}, i.e.,
\begin{equation}
	\begin{aligned}
		N_1 [\hat{f}] = \max_\xi |\hat{f}| \,, \ N_2 [\hat{f}] = \Big[ \int \hat{f}^2 \d \xi \Big]^\frac{1}{2} \,, \ N_3^{(r)} [\hat{f}] = \max_\xi (1 + \xi^2)^r |\hat{f}| \,.
	\end{aligned}
\end{equation}
Following the arguments in Page 50-51 of \cite{Grad-1963}, there holds
\begin{equation}\label{Grad-3}
	\begin{aligned}
		N_3^{(m/2)} [\hat{f}] < \infty
	\end{aligned}
\end{equation}
for hard potential case and any integer $m > 0$ provided that
\begin{equation*}
	\begin{aligned}
		\hat{\textrm{g}} \in \mathrm{Null}^\perp (L) \,, \ \int (\hat{\textrm{g}}^2/\nu) \d \xi < \infty \,, \ N_3^{(m/2)} [\hat{\textrm{g}}] < \infty \,.
	\end{aligned}
\end{equation*}
Specifically, the hypocoercivity \eqref{Hypo-L} and $\int (\hat{\textrm{g}}^2/\nu) \d \xi < \infty$ imply that
\begin{equation*}
	\begin{aligned}
		\{ N_2 [\hat{f}] \}^2 < \tfrac{1}{\nu_0} \int \nu \hat{f}^2 \d \xi < \infty \,.
	\end{aligned}
\end{equation*}
Then, the first inequality in \eqref{Grad-2} implies $N_1 [K[\hat{f}]] < \infty$. Together with \eqref{Grad-1}, $N_1 [\hat{f}] = N_3^{(0)} [\hat{f}] < \infty$ under the condition $N_3^{(0)} [\hat{\textrm{g}}] < \infty$. Moreover, the last inequality of \eqref{Grad-2} means $N_3^{(1/2)} [K[\hat{f}]] < \infty$. Consequently, \eqref{Grad-1} and $N_3^{(1/2)} [\hat{\textrm{g}}] < \infty$ imply $N_3^{(1/2)} [\hat{f}] < \infty$. From repeating the previous procedure $m$ times, one directly derives \eqref{Grad-3}, which means the decay estimate \eqref{Grad-inq}.

The Grad's inequality \eqref{Grad-inq} illustrates the polynomial decay with any fixed order $m > 0$ of the pseudo-inverse operator $\L^{-1}$ for the hard potential cases $0 \leq \gamma \leq 1$. In this paper, we will enhance the polynomial decay \eqref{Grad-inq} to the exponential decay of the pseudo-inverse operator $\L^{-1}$ for the hard potential cases $0 \leq \gamma \leq 1$, and furthermore generalize the exponential decay to part of soft potential cases $- \frac{3}{2} < \gamma < 0$.

We now denote the derivative operator $\partial_v^\alpha$ by
\begin{equation*}
	\begin{aligned}
		\partial_v^\alpha f = \frac{ \partial^{|\alpha|} f }{ \partial v_1^{\alpha_1} \partial v_2^{ \alpha_2 } \partial v_3^{ \alpha_3 } } \,,
	\end{aligned}
\end{equation*}
where the multi-index $ \alpha = ( \alpha_1, \alpha_2, \alpha_3 ) \in \mathbb{N}^3 $ with $|\alpha| = \alpha_1 + \alpha_2 + \alpha_3$. The symbol $\alpha' \leq \alpha$ means that $\alpha_i' \leq \alpha_i$ for $1 \leq i \leq 3$. Moreover, $\alpha' < \alpha$ represents $\alpha' \leq \alpha$ and $|\alpha'| < |\alpha|$. We also define the bracket
\begin{equation*}
	\begin{aligned}
		\langle v \rangle = 1 + |v| \,.
	\end{aligned}
\end{equation*}

We now state the first theorem as follows.

\begin{theorem}[Exponential decay of $\L^{-1}$]\label{Thm1}
	Let $- 3 < \gamma \leq 1$, $0 < q < 1$ and $g \in \mathrm{Null}^\perp (\L)$. The operator $\L^{-1}$ is given in Definition \ref{Df-L-1}.
	\begin{enumerate}
	\item[(I) ] \underline{\bf The cases $- \frac{3}{2} < \gamma \leq 1$:}
	
	Define the positive number $k_\gamma$ by
	\begin{equation*}
		k_\gamma > \left\{
		  \begin{aligned}
		  	\tfrac{3 - \gamma}{2} \,, & \textrm{ if } 0 \leq \gamma \leq 1 \,, \\
		  	\tfrac{3 - 2 \gamma}{2} \,, & \textrm{ if } - \tfrac{3}{2} < \gamma < 0 \,.
		  \end{aligned}
		\right.
	\end{equation*}
    Assume that $\langle v \rangle^{k_\gamma} \M^{- \frac{q}{2}} g \in L^\infty$. Then the following statements hold:
	\begin{enumerate}
		\item[(i)] For the hard potential cases $0 \leq \gamma \leq 1$,
		\begin{equation}\label{hard decay}
			\begin{aligned}
				| \L^{-1} g (v) | \leq C \| \langle v \rangle^{k_\gamma} \M^{- \frac{q}{2}} g \|_{L^\infty} \M^\frac{q}{2} (v) \,, \ \forall v \in \R^3 \,.
			\end{aligned}
		\end{equation}
	    \item[(ii)] For the part of soft potential cases $- \frac{3}{2} < \gamma < 0$,
	    \begin{equation}\label{soft decay}
	    	\begin{aligned}
	    		| \L^{-1} g (v) | \leq C \| \langle v \rangle^{k_\gamma} \M^{- \frac{q}{2}} g \|_{L^\infty} \nu^{-1} (v) \M^\frac{q}{2} (v) \,, \ \forall v \!\quad v \in \R^3 \,.
	    	\end{aligned}
	    \end{equation}
	\end{enumerate}

    \item[(II)] \underline{\bf The cases $-3 < \gamma \leq - \frac{3}{2}$:}

    Assume that
    \begin{equation*}
    	\begin{aligned}
    		\sum_{|\alpha| \leq 2} \| \langle v \rangle^{\gamma (|\alpha| - 3)} \M^{- \frac{q}{2} } \partial^\alpha_v g \|_{L^2} < \infty \,.
    	\end{aligned}
    \end{equation*}
    Then
    \begin{equation}
    	\begin{aligned}
    		| \L^{-1} g (v) | \leq C \sum_{|\alpha| \leq 2} \| \langle v \rangle^{\gamma (|\alpha| - 3)} \M^{- \frac{q}{2} } \partial^\alpha_v g \|_{L^2} \M^\frac{q}{2} (v) \,, \ a.e. \!\quad v \in \R^3 \,.
    	\end{aligned}
    \end{equation}
    \end{enumerate}
    The positive constant $C > 0$ above depends on $\rho$, $\u$ and $T$.
\end{theorem}

\begin{remark}
	Actually, the constant $C = C(\rho, \u, T) > 0$ in Theorem \ref{Thm1} can be regarded as a positive smooth function depending on the variables $(\rho, \u, T) \in \R_+ \times \R^3 \times \R_+$. In Part (I), the function $g (v)$ is only required in weighted $L^\infty$ space. However, $g (v)$ should be assumed in weighted $H^2$ space in Part (II). In this sense, the condition of Part (II) is definitely stronger than that of Part (II).
\end{remark}

Comparing to the hard potential case estimate \eqref{Caflisch-inq} originally stated by Caflisch \cite{Caflish-1980-CPAM}, our estimate \eqref{hard decay} is slightly different. More precisely, our estimate requires $0<q<1$, while the estimate \eqref{Caflisch-inq} basically corresponds to the case $q=1$. Furthermore, besides the hard potential, our result also cover the soft potential cases, i.e. $- 3 < \gamma < 0$.

There were some weighted $L^2$-norms estimates for the pseudo-inverse operator $\L^{-1}$ associated with the global Maxwellian $M (v) = \M_{[1,0,1]} (v)$ by the hypocoercivity of the operator $\L$. In \cite{Ukai-Yang-2006}, Ukai and Yang proved the boundedness of $\L^{-1}$ from $L^2 (\nu^{-1} \d v) \cap \mathrm{Null}^\perp (\L)$ to $L^2 (\nu \d v) \cap \mathrm{Null}^\perp (\L)$ for the hard potential cases $0 \leq \gamma \leq 1$. This result has also been applied in \cite{Liu-Zhao-JDE-2011} for instance. In \cite{Guo-CPAM06}, Guo proved the hypocoercivity of the $\L$ and its $v$-derivatives in the mixed space $H^k_{x,v} (M^{-q} \nu \d v \d x)$ ($k \geq 0$, $0 < q < 1$) for hard and soft potential cases $- 3 < \gamma \leq 1$, which can imply the boundedness of $\L^{-1}$ from  $H^k_{x,v} (M^{-q} \nu^{-1} \d v \d x) \cap \mathrm{Null}^\perp (\L)$ to  $H^k_{x,v} (M^{-q} \nu \d v \d x) \cap \mathrm{Null}^\perp (\L)$. In \cite{Duan-Yu-AM2020}, Duan and Yu proved the similar decay estimate \eqref{hard decay} when the pseudo-inverse of linearized Landau collision operator acted on any polynomial $U (v)$ of $v - \u$ such that $U(v) \M \in \mathrm{Null}^\perp (\L)$. Our Grad-Caflisch type estimate is stronger than the estimates above, since ours is pointwise.

\subsection{Applications to Hilbert expansions of Boltzmann equation with Euler scaling}

In this subsection, the exponential decay estimates in Theorem \ref{Thm1} will be used to control the terms in the Hilbert expansions of the Boltzmann equation with compressible Euler scaling. The scaled Boltzmann equation is
\begin{equation}\label{BE}
	\begin{aligned}
		\partial_t F_\kappa + v \cdot \nabla_x F_\kappa = \tfrac{1}{\kappa} \mathcal{B} (F_\kappa, F_\kappa) \,,
	\end{aligned}
\end{equation}
where $F_\kappa (t,x,v)$  is the density of particles of velocity $v \in \R^3$, at position $x \in \R^3$, time $t \in \R_+$, and $\kappa > 0$ is the Knudsen number, which is assumed to be small. The solution $F_\kappa$ is sought as the Hilbert expansion form
\begin{equation}\label{Expnd-series}
	\begin{aligned}
		F_\kappa = \sum_{n = 0}^\infty \kappa^n F_n \,,
	\end{aligned}
\end{equation}
where $F_0, F_1, \cdots, F_n, \cdots$ are independent of $\kappa$. As shown in \cite{Guo-Jang-Jiang-2010-CPAM, JLT-arXiv-2021} for instance, the expanded terms solve the equations
\begin{equation}\label{H-expnd}
	\begin{aligned}
		\kappa^{-1} : & \qquad 0 = \B (F_0, F_0) \,, \\
		\kappa^0 : & \qquad ( \partial_t + v \cdot \nabla_x ) F_0 = \B (F_0, F_1) + \B (F_1, F_0) \,, \\
		\kappa^1 : & \qquad (\partial_t + v \cdot \nabla_x ) F_1 = \B (F_0, F_2) + \B (F_2, F_0) + \B (F_1, F_1) \,, \\
		& \qquad \qquad \cdots \cdots \\
		\kappa^n : & \qquad (\partial_t + v \cdot \nabla_x ) F_n = \B (F_0, F_{n+1}) + \B (F_{n+1}, F_0) + \sum_{\substack{i+j=n+1 \\ i,j \geq 1}} \B (F_i, F_j) \,, \\
		& \qquad \qquad \cdots \cdots
	\end{aligned}
\end{equation}
From the well-known $H$-theorem, the order of $\kappa^{-1}$ in \eqref{H-expnd} implies that
\begin{equation*}
	\begin{aligned}
		F_0 (t,x,v) = \M (t,x,v) = \M_{[\rho, \u, T]} (v) \,,
	\end{aligned}
\end{equation*}
where $[\rho, \u, T]$ depends on the variables $(t,x)$, which obeys the compressible Euler system
\begin{equation}\label{CES}
	\left\{
	  \begin{aligned}
	  	& \partial_t \rho + \div_x (\rho \u) = 0 \,, \\
	  	& \partial_t (\rho \u) + \div_x (\rho \u \otimes \u) + \nabla_x (\rho T) = 0 \,, \\
	  	& \partial_t [ \rho ( \tfrac{3}{2} T + \tfrac{1}{2} |\u|^2 ) ] + \div_x [ \rho \u ( \tfrac{3}{2} T + \tfrac{1}{2} |\u|^2 ) ] + \div_x (\rho T \u) = 0 \,.
	  \end{aligned}
	\right.
\end{equation}
Moreover, for each $n \geq 1$, $f_n : = \tfrac{F_n}{\sqrt{\M}}$ can be decomposed as the fluid and kinetic parts:
\begin{equation}\label{fn-spilt}
	\begin{aligned}
		f_n : = & \tfrac{F_n}{\sqrt{\M}} = \P f_n + (\I - \P) f_n \\
		= & \big\{ \tfrac{\rho_n}{\rho} + u_n \cdot \tfrac{v - \u}{T} + \tfrac{\theta_n}{6 T} ( \tfrac{|v - \u|^2}{T} - 3 ) \big\} \sqrt{\M} + (\I - \P) f_n \,,
	\end{aligned}
\end{equation}
where the projection $\P$ is defined in \eqref{Projt-P}, the kinetic part is solved by
\begin{equation}\label{kinetic-n}
	\begin{aligned}
		(\I - \P) f_n = \L^{-1} \Big\{ - (\I - \P) \big[ \M^{- \frac{1}{2}} (\partial_t + v \cdot \nabla_x) ( f_{n-1} \sqrt{\M} ) \big] + \sum_{\substack{i+j=n \\ i,j \geq 1}} \Gamma (f_i, f_j) \Big\} \,,
	\end{aligned}
\end{equation}
and the fluid variables $(\rho_n, u_n, \theta_n)$ are solved by the linearized compressible Euler system
\begin{equation}\label{CELS}
	\left\{
	\begin{aligned}
		& \partial_t \rho_n + \div_x ( \rho u_n + \rho_n \u ) = 0 \,, \\
		& \rho (\partial_t u_n + u_n \cdot \nabla_x \u + \u \cdot \nabla_x u_n) - \tfrac{\nabla_x (\rho T)}{\rho} \rho_n + \nabla_x \Big( \tfrac{\rho \theta_n + 3 T \rho_n}{3} \Big) = \mathcal{F}_u^\perp ((\I-\P)f_n) \,, \\
		& \rho \Big( \partial_t \theta_n + \u \cdot \nabla_x \theta_n + \tfrac{2}{3} (\theta_n \div_x \u + 3 T \div_x u_n) + 3 u_n \cdot \nabla_x T \Big) = \mathcal{G}_\theta^\perp ((\I-\P)f_n) \,.
	\end{aligned}
    \right.
\end{equation}
Here the nonlinear collision operator $\Gamma (f_i, f_j)$ denotes
\begin{equation}\label{Gamma-fi-fj}
	\begin{aligned}
		\Gamma (f_i, f_j) = \tfrac{1}{\sqrt{\M}} \B (f_i \sqrt{\M}, f_j \sqrt{\M}) \,.
	\end{aligned}
\end{equation}
The source terms $\mathcal{F}_u^\perp ((\I - \P)f_n) \in \R^3$ and $\mathcal{G}_\theta^\perp ((\I - \P)f_n) \in \R$ are defined by
\begin{equation}\label{FG-fn}
	\begin{aligned}
		& \mathcal{F}_{u,i}^\perp ((\I-\P)f_n) = - \sum_{j=1}^3 \partial_{x_j} \langle T \mathscr{A}_{ij}, (\I-\P) f_n \rangle \ (i = 1,2,3) \,, \\
		& \mathcal{G}_\theta^\perp ((\I-\P)f_n) = - \div_x \big( 2 T^\frac{3}{2} \langle \mathscr{B}, (\I-\P) f_n \rangle + 2 T \u \cdot \langle \mathscr{A}, (\I-\P) f_n \rangle  \big) - 2 \u \cdot \mathcal{F}_u^\perp ((\I-\P)f_n) \,,
	\end{aligned}
\end{equation}
where $\mathscr{A} \in \R^{3 \times 3}$ and $\mathscr{B} \in \R^3$ are the Burnett functions belonging to $\mathrm{Null}^\perp (\L)$ with entries
\begin{equation}\label{Burnett}
	\begin{aligned}
		\mathscr{A}_{ij} & = \Big\{ \tfrac{(v_i - \u_i) (v_j - \u_i)}{T} - \delta_{ij} \tfrac{|v - \u|^2}{3 T} \Big\} \sqrt{\M} \ (1 \leq i,j \leq 3) \,, \\
		\mathscr{B}_i & = \tfrac{v_i - \u_i}{2 \sqrt{T}} \Big( \tfrac{|v - \u|^2}{T} - 5 \Big) \sqrt{\M} \ (1 \leq i \leq 3) \,.
	\end{aligned}
\end{equation}
Once $f_n$ is rewritten as the following abstract and intuitive form:
\begin{equation}
	\begin{aligned}
		f_n = Kinetic(n) + Fluid (n) \,,
	\end{aligned}
\end{equation}
where $Kinetic (n) = (\I - \P) f_n$ and $Fluid (n) = \P f_n$, one can summarize the orders to solve $f_n$ ($n \geq 1$) as the following Figure \ref{Fig3}.
\begin{figure}[h]
	\begin{center}
		\begin{tikzpicture}
			\pgfmathsetmacro\rd{2};
			\draw (\rd, 0.5) node[above]{$Fluid(0)$};
			\draw (\rd,-1) node[above]{$Fluid(1)$};
			\draw (\rd,-2.5) node[above]{$Fluid(2)$};
			\draw (\rd,-4) node[above]{$Fluid(n)$};
			\draw (-\rd,-1) node[above]{$Kinetic(1)$};
			\draw (-\rd,-2.5) node[above]{$Kinetic(2)$};
			\draw (-\rd,-4) node[above]{$Kinetic(n)$};
			\draw[->,thick] (-0.9,-0.7)--(1.1, -0.7);
			\draw[->,thick] (-0.9,-2.1)--(1.1, -2.1);
			\draw[->,thick] (-0.9,-3.7)--(1.1, -3.7);
			\draw[->,thick] (1.2,0.7)--(-1,-0.4);
			\draw[->,thick] (1.2,-0.9)--(-1.1,-1.9);
			\draw[->,thick,dashed] (1.2,-2.3)--(-1.1,-3.5);
			\draw[->,thick,dashed] (1.2,-3.8)--(-1.1,-4.9);
		\end{tikzpicture}
	\end{center}
	\caption{The orders of solving $f_n$ ($n \geq 1$): $Fluid(0)$ is $F_0/\sqrt{\M} = \sqrt{\M}$ determined by the solution $(\rho, \u, T)$ to the compressible Euler system \eqref{CES}.}\label{Fig3}
\end{figure}
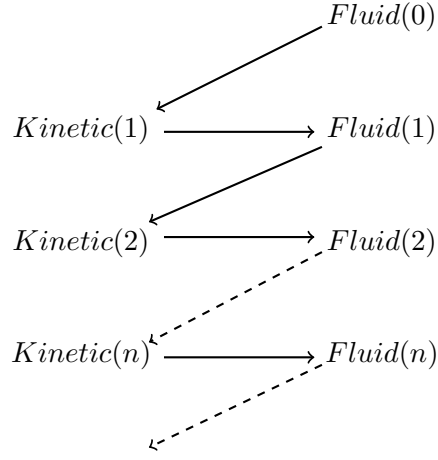

In this paper, the following two hypotheses will be imposed while controlling the expanded terms $F_n (t,x,v)$ ($n \geq 1$):

{\bf (H1)} From Schochet's work \cite{Schochet-1986-CMP}, one can assume $(\rho, \u, T)$ is a bounded smooth solution to the compressible Euler system \eqref{CES}, and $\rho, T$ have both positive lower and upper bounds. Moreover, all derivatives of $(\rho, \u, T)$ are also bounded. Namely, for any large $m \geq 0$,
\begin{equation}\label{e-m}
	\begin{aligned}
		\inf_{t,x} \rho > 0 \,, \ \inf_{t,x} T > 0 \,, \mathfrak{e}_m = \tfrac{1}{\inf_{t,x} \rho} + \tfrac{1}{\inf_{t,x} T} + \sum_{0 \leq |\beta| \leq m} \sup_{t,x} |\partial_{t,x}^{\beta} (\rho, \u, T)| < \infty \,,
	\end{aligned}
\end{equation}
where the derivative operator $\partial^\beta_{t,x}$ is defined by
\begin{equation*}
	\begin{aligned}
		\partial^\beta_{t,x} = \tfrac{\partial^{|\beta|}}{\partial t^{\beta_0} \partial x_1^{\beta_1} \partial x_2^{\beta_2} \partial x_3^{\beta_2}}
	\end{aligned}
\end{equation*}
for $\beta = (\beta_0, \beta_1, \beta_2, \beta_3) \in \mathbb{N}^4$ with $|\beta| = \beta_0 + \beta_1 + \beta_2 + \beta_3$.

{\bf (H2)} Based on (H1) and Theorem \ref{Thm2} below, the solution $(\rho_n, u_n, \theta_n)$ to the linearized compressible Euler system \eqref{CELS} can also be assumed to be bounded and smooth, {\em provided that the source terms $\mathcal{F}_u^\perp ((\I-\P)f_n)$ and $\mathcal{G}_\theta^\perp ((\I-\P)f_n)$ are bounded and smooth}. Furthermore, all derivatives of $(\rho_n, u_n, \theta_n)$ are also bounded. More precisely, for any $m \geq 0$ and $n \geq 1$,
\begin{equation}\label{h-m}
	\begin{aligned}
		\h_m^n : = \sum_{0 \leq |\beta| \leq m} \sup_{t,x} |\partial_{t,x}^\beta (\rho_n, u_n, \theta_n)| < \infty \,.
	\end{aligned}
\end{equation}

Remark that, for any $m \geq 0$ and $n \geq 1$,
\begin{equation*}
	\begin{aligned}
		\e_m \leq \e_{m+1} \,, \ \h_m^n \leq \h_{m+1}^n \,.
	\end{aligned}
\end{equation*}

The goal of this subsection is to control the expanded terms $F_n = f_n \sqrt{\M}$ for each $n \geq 1$. For the hard potential cases $0 \leq \gamma \leq 1$, Caflisch stated the following inequality in Proposition 3.1, Page 656 of \cite{Caflish-1980-CPAM}:

{\em
Let $(\rho, \u ,T)$ be a smooth solution of the compressible equations \eqref{CES}, and form the Maxwellian $F_0 = \M$ as in \eqref{Maxwellian}. Then the terms $F_1, \cdots, F_6$ of the Hilbert expansion\footnote{In \cite{Caflish-1980-CPAM}, the expansion was truncated after the sixth term, due to showing the convergence of the expanded series \eqref{Expnd-series} by controlling the associated remainder system uniformly in $\kappa > 0$.} are smooth in $(t,x)$ and have decay given by
\begin{equation}\label{Caflisch-inq-F}
	\begin{aligned}
		F_i (t,x,v) \leq c |v|^{3 i} \M (t,x,v) \,,
	\end{aligned}
\end{equation}
where $c$ is a constant independent of $v, x, t$.
}

As we mentioned in the last subsection, in Grad \cite{Grad-1963} the inequality \eqref{Grad-inq} was shown. Based on it the following bound of $F_n$ could be derived:
\begin{equation}\label{Caflisch-inq-FC}
	\begin{aligned}
		|F_i (t,x,v)| \leq c \langle v \rangle^{-m} \M^\frac{1}{2} (t,x,v) \ \textrm{for any integer } m > 0 \,,
	\end{aligned}
\end{equation}
which is weaker than \eqref{Caflisch-inq-F}.

We use Theorem \ref{Thm1} to prove some new estimates of the terms (and their derivatives over $(t,x)$) in the Hilbert expansion \eqref{Expnd-series}. This is the second main result of this paper.

First, in the Hilbert expansion framework, the linearized operator $\L$ depends on the parameters $(t,x)$, and so $\L^{-1}$ does. We will also derive the decay estimates of the $(t,x)$-derivatives of the pseudo-inverse operator $\L^{-1}$. Then the following decay theorem is first displayed.

\begin{theorem}[Exponential decay of parameter derivatives of $\L^{-1}$]\label{Thm2}
	Let $- 3 < \gamma \leq 1$, $0 < q < 1$ and $g \in \mathrm{Null}^\perp (\L)$. For any fixed integer $m \geq 0$, let $\beta = (\beta_0, \beta_1, \beta_2, \beta_3) \in \mathbb{N}^4$ with $|\beta| = m$. Let $(\rho, \u, T)$ be a solution to \eqref{CES} satisfying hypothesis (H1).
	\begin{itemize}
	\item[(I)] \underline{\bf The cases $- \frac{3}{2} < \gamma \leq 1$:}
	
	Assume that
	\begin{equation}\label{Xm-g}
		\begin{aligned}
			\mathcal{W}^{m, \infty}_{q_0, \cdots, q_m} (g) : = \sum_{i=0}^m \sum_{|\beta'| = i} \| \langle v \rangle^{k_\gamma} \M^{- \frac{q_i}{2}} (v) \partial_{t,x}^{\beta'} g \|_{L^\infty} < \infty
		\end{aligned}
	\end{equation}
    for any fixed $0 < q_m < q_{m-1} < \cdots < q_2 < q_1 < q_0 = q < 1$. Then there is a positive constant $C = C (\e_m) > 0$ such that
    \begin{equation}\label{L-inverse-Dtx}
    	\begin{aligned}
    		| \partial_{t,x}^\beta [ \L^{-1} g (t,x,v) ] | \leq C \mathcal{W}^{m, \infty}_{q_0, \cdots, q_m} (g) \Theta_\gamma (v) \M^\frac{q_m}{2} (t,x,v) \ \forall t,x,v \,,
    	\end{aligned}
    \end{equation}
    where $\e_m$ is mentioned in \eqref{e-m}, and $\Theta_\gamma (v)$ is defined by
    \begin{equation}\label{Theta-gamma}
    	\begin{aligned}
    		\Theta_\gamma (v) = 1 \textrm{ if } 0 \leq \gamma \leq 1 \,, \textrm{ and } \Theta_\gamma (v) = \nu^{-1} (v) \textrm{ if } - \tfrac{3}{2} < \gamma < 0 \,.
    	\end{aligned}
    \end{equation}

    \item[(II)] \underline{\bf The cases $- 3 < \gamma \leq - \frac{3}{2}$:}

    Assume that
    \begin{equation}\label{Wm2-g}
    	\begin{aligned}
    		\mathcal{W}_q^{m,2} (g) = \sum_{i=0}^m \sum_{|\beta'| = i} \sum_{|\alpha| \leq 2} \| \langle v \rangle^{\gamma (|\alpha| - 3)} \M^{- \frac{q}{2} } (v) \partial_v^\alpha \partial_{t,x}^{\beta'} g \|_{L^2} < \infty \,.
    	\end{aligned}
    \end{equation}
    Then there is a constant $C = C (\mathfrak{e}_m) > 0$ such that
    \begin{equation}\label{L-inverse-Dtx-VSP}
    	\begin{aligned}
    		| \partial_{t,x}^\beta [ \L^{-1} g (t,x,v) ] | \leq C \mathcal{W}_q^{m,2} (g) \M^\frac{q}{2} (t,x,v) \ \forall \, t,x,v \,.
    	\end{aligned}
    \end{equation}
    \end{itemize}
\end{theorem}

Then, based on Theorem \ref{Thm2}, the following new estimates are proved.

\begin{theorem}\label{Thm3}
	Let $0 < q < 1$, $- 3 < \gamma \leq 1$ and integer $m \geq 0$ be arbitrary. Let $(\rho, \u, T)$ be a bounded and smooth solution of the compressible system \eqref{CES}, and form the Maxwellian $F_0 = \M$. Assume that the hypotheses (H1)-(H2) hold. Then the terms $F_1 (t,x,v), \cdots, F_n (t,x,v), \cdots$ of the Hilbert expansion are bounded and smooth in $(t,x)$, and satisfy the decay
	\begin{equation}\label{C-bnd}
		\begin{aligned}
			\sum_{|\beta| = m} |\partial_{t,x}^\beta F_n (t,x,v)| \leq C \M^\frac{1+q}{2} (t,x,v) \ \ \forall n \geq 1 \,,
		\end{aligned}
	\end{equation}
    where the constant $C = C (\e_{m+n}, \h_{m+n-1}^1, \h_{m+n-2}^2, \cdots, \h_{m+1}^{n-1}, \h^n_m) > 0$ is independent of the variables $t,x,v$. Here $\e_{m+n}$ and $\h^i_{m+n-i}$ $(1 \leq i \leq n)$ are defined in \eqref{e-m} and \eqref{h-m}, respectively.
\end{theorem}

We note that, in \cite{Guo-CPAM06}, based on the boundedness of $\L^{-1}$ from  $H^k_{x,v} (M^{-q} \nu^{-1} \d v \d x) \cap \mathrm{Null}^\perp (\L)$ to  $H^k_{x,v} (M^{-q} \nu \d v \d x) \cap \mathrm{Null}^\perp (\L)$ and the Sobolev embedding theory, Guo proved the similar bounds \eqref{C-bnd} of the expanded terms $F_n (t,x,v)$ when proving the incompressible Navier-Stokes limit of the Boltzmann equation. Note that Guo's work is about the global Maxwellians, which is independent of $(t,x)$. This makes the estimates much simpler.

\begin{remark}
	It is seemly to be confused that the bounds of $\sum_{|\beta| = m} |\partial_{t,x}^\beta F_n (t,x,v)| $ depends on the quantity $\h^n_m = \sum_{0 \leq |\beta| \leq m} \sup_{t,x} |\partial_{t,x}^\beta (\rho_n, u_n, \theta_n)|$ associated with the fluid variables $(\rho_n, u_n, \theta_n)$ of $F_n (t,x,v)$. However, $(\rho_n, u_n, \theta_n)$ is solved by the linearized compressible Euler system \eqref{CELS} with the source terms $\mathcal{F}_u^\perp ((\I - \P)f_n) $ and $\mathcal{G}_\theta^\perp ((\I - \P)f_n) $ given in \eqref{FG-fn}. They subject to the bounds \eqref{Fu-1}-\eqref{Gtheta-1} and \eqref{Fu-Gtheta-2} below, i.e.,
	\begin{equation*}
		\begin{aligned}
			\Big| \partial_{t,x}^\beta  \mathcal{F}_u^\perp ((\I-\P)f_1) \Big| + \Big| \partial_{t,x}^\beta \mathcal{G}_\theta^\perp ((\I-\P)f_1) \Big| \leq C (\e_{m+2}) < \infty
		\end{aligned}
	\end{equation*}
    and
    \begin{equation*}
    	\begin{aligned}
    		|\partial_{t,x}^\beta \mathcal{F}_u^\perp ((\I-\P)f_n)| + |\partial_{t,x}^\beta \mathcal{G}_\theta^\perp & ((\I-\P) f_n)| \\
    		& \leq C (\e_{m+n+1}) C (\h_{m+n+1-\mathfrak{p}}^\mathfrak{p}; 1 \leq \mathfrak{p} \leq n - 1) < \infty  \ (\forall n \geq 2)
    	\end{aligned}
    \end{equation*}
    for any $|\beta| = m \geq 0$. By the existence theory of the system \eqref{CELS} (see \cite{GHW-2020, JLT-arXiv-2021}, for instance), the all quantities $\h^n_m$ $(n \geq 1, m \geq 0)$ can be bounded by the quantity $\e_{m'} = \tfrac{1}{\inf_{t,x} \rho} + \tfrac{1}{\inf_{t,x} T} + \sum_{0 \leq |\beta| \leq m'} \sup_{t,x} |\partial_{t,x}^{\beta} (\rho, \u, T)|$ for sufficiently large integer $m' = m' (n,m) \geq 0$, which corresponds to the solution $(\rho, \u, T)$ of the compressible Euler system \eqref{CES}. Nevertheless, the well-posedness of the system \eqref{CELS} is not the goal of this paper. We thereby represent the upper bounds of $F_n$ depending on the quantities $\h_{m+n-1}^1, \h_{m+n-2}^2, \cdots, \h_{m+1}^{n-1}, \h^n_m$.
\end{remark}

\begin{remark}
	Let $m = 0$ in the bound \eqref{C-bnd}. There holds
	\begin{equation}\label{C-bnd-0}
		\begin{aligned}
			| F_n (t,x,v)| \leq C \M^\frac{1+q}{2} (t,x,v) \ \ \forall n \geq 1
		\end{aligned}
	\end{equation}
	for $0 < q < 1$ and $- \frac{3}{2} < \gamma \leq 1$.
	For the hard potential cases $0 \leq \gamma \leq 1$, the bound \eqref{C-bnd-0} is stronger than the bound \eqref{Caflisch-inq-FC}, but weaker than the bound \eqref{Caflisch-inq-F}. However, the estimate \eqref{Caflisch-inq-F} is hard to be verified. The closest results to Caflisch's inequality that we can achieve is estimates of $F_1 (t,x,v)$ with respect to the hard sphere case $\gamma = 1$. More precisely, we compare the bound of $F_1 (t,x,v)$, about which Caflisch claimed that
	\begin{equation}\label{Caflisch-F1}
		\begin{aligned}
			| F_1 (t,x,v) | \leq c |v|^3 \M (t,x,v) \,.
		\end{aligned}
	\end{equation}
	About our results, by the equality \eqref{f1-1} and \eqref{f1-2} below, the kinetic part of $f_1 = F_1 / \sqrt{\M}$ can be expressed as
	\begin{equation}
		\begin{aligned}
			(\I - \P) f_1 = - \nabla_x \u : \widehat{\mathscr{A}} (v) + \tfrac{\nabla_x T}{\sqrt{T}} \cdot \widehat{\mathscr{B}} (v) \,,
		\end{aligned}
	\end{equation}
	where $\widehat{\mathscr{A}} (v) = \L^{-1} \mathscr{A} (v)$ and $\widehat{\mathscr{B}} (v) = \L^{-1} \mathscr{B} (v)$. Following the equation (23) in Section 3 of \cite{BGL1}, Theorem 2.1 of \cite{DG-1994} or Section 2.3.3 of \cite{BGP-2000}, it can be shown that there exist two scalar valued functions $\mathfrak{a}, \mathfrak{b} : [0, \infty) \to \R$ such that
	\begin{equation}
		\begin{aligned}
			\widehat{\mathscr{A}} (v) = \mathfrak{a} (|v|) \mathscr{A} (v) \quad \textrm{and} \quad \widehat{\mathscr{B}} (v) = \mathfrak{b} (|v|) \mathscr{B} (v) \,.
		\end{aligned}
	\end{equation}
	Moreover, Proposition 6.5 of \cite{Golse-SRM-05} shown that, for a hard sphere gas $(\gamma = 1)$, the functions $\mathfrak{a}$ and $\mathfrak{b}$ satisfy the growth estimate
	\begin{equation}
		\begin{aligned}
			| \mathfrak{a} (|v|) | + | \mathfrak{b} (|v|) | \leq C \langle v \rangle \ \forall v \in \R^3 \,.
		\end{aligned}
	\end{equation}
	It thereby holds that
	\begin{equation*}
		\begin{aligned}
			| (\I - \P) f_1 | \leq C (\e_1) \langle v \rangle^4 \M^\frac{1}{2} \,.
		\end{aligned}
	\end{equation*}
	Together with $\P f_1 = \big\{ \tfrac{\rho_1}{\rho} + u_1 \cdot \tfrac{v - \u}{T} + \tfrac{\theta_1}{6 T} ( \tfrac{|v - \u|^2}{T} - 3 ) \big\} \sqrt{\M}$, one has
	\begin{equation*}
		\begin{aligned}
			|f_1 (t,x,v)| \leq C(\e_1, \h_0^1) \langle v \rangle^4 \M^\frac{1}{2} (t,x,v) \,,
		\end{aligned}
	\end{equation*}	
	which means that
	\begin{equation}
		\begin{aligned}
			|F_1 (t,x,v)| \leq C(\e_1, \h_0^1) \langle v \rangle^4 \M (t,x,v) \,.
		\end{aligned}
	\end{equation}
	However, this is also weaker than the Caflisch's inequality \eqref{Caflisch-F1} with respect to the case $n=1$ and hard sphere case $\gamma = 1$.
\end{remark}

\subsection{Sketch of proofs}

The first part of this paper is to derive the decay estimates of $f = \L^{-1} g$ for $f, g \in \mathrm{Null}^\perp (\L)$, that is, Theorem \ref{Thm1}. Note that $g = \L f = \nu f + K f = \nu f + K_1 f - K_2 f$. For preparation of the work, with general local Maxwellian $\M$ with fluid variable $(\rho, \u, T)$, we initially calculate the integral kernel of the operators $K_1$ and $K_2$ in Lemma \ref{Lmm-k1-k2}, and estimate their bounds in Lemma \ref{Lmm-k1k2-bnd}, which also shows that the integral kernels are both in $L^1$ for $-3 < \gamma \leq 1$. Based on these results, we will justify the conclusions given in Theorem \ref{Thm1} separately for the cases $- \frac{3}{2} < \gamma \leq 1$ and $- 3 < \gamma \leq - \frac{3}{2}$.

(I) The case $- \frac{3}{2} < \gamma \leq 1$ in Theorem \ref{Thm1}: The first step is to prove that $\M^{- \frac{q}{2}} (v) |K_i f (v)|$ can be bounded the quantity $\| \M^{- \frac{q}{2}} f \|_{L^2}$ for $0 < q < 1$ and $- \frac{3}{2} < \gamma \leq 1$, see Lemma \ref{Lmm-K1K2-L2}. Intuitively, since $\M^{- \frac{q}{2}} (v) |K_i f (v)| \leq \int_{\R^3} | \mathbbm{k}_i (v, v_1) | \M^{- \frac{q}{2}} (v) \M^\frac{q}{2} (v_1) | \M^{- \frac{q}{2}} (v_1) f(v_1) | \d v_1$, the $L^2$-integrability of $\mathbbm{k}_i (v, v_1)  \M^{- \frac{q}{2}} (v) \M^\frac{q}{2} (v_1)$ uniformly in $v$ is required. By the spirit of Lemma \ref{Lmm-k1k2-bnd}, the constraint $- \frac{3}{2} < \gamma \leq 1$ is necessary.
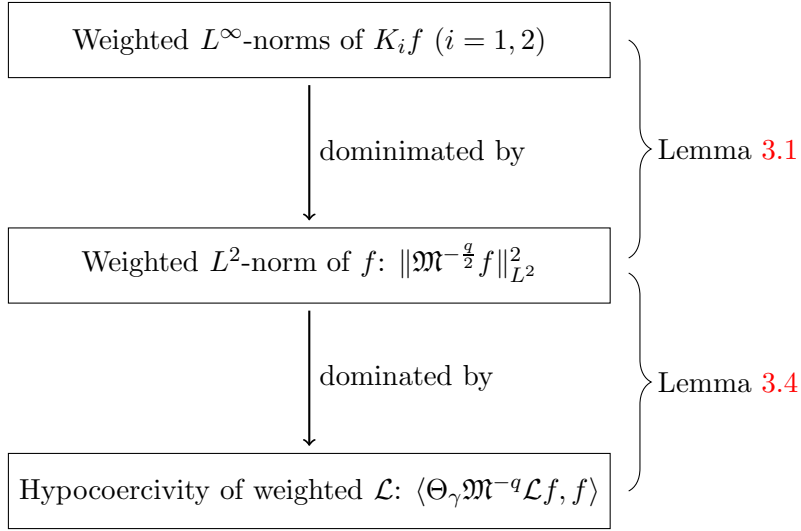
\begin{figure}[h]
	\begin{center}
		\begin{tikzpicture}
			\pgfmathsetmacro\rd{4};
			\draw (-\rd, 0)--(\rd, 0)--(\rd, -1)--(-\rd, -1)--(-\rd,0)--cycle;
			\draw (-\rd, -3)--(\rd, -3)--(\rd, -4)--(-\rd, -4)--(-\rd,-3)--cycle;
			\draw (-\rd, -6)--(\rd, -6)--(\rd, -7)--(-\rd, -7)--(-\rd,-6)--cycle;
			\draw[->, thick] (0,-1.1)--(0,-2.9);
			\draw[->,thick] (0,-4.1)--(0,-5.9);
			\draw (0,-2) node[right]{dominated by};
			\draw (0,-5) node[right]{dominated by};
			\draw (0,-1) node[above=1.5mm]{Weighted $L^\infty$-norms of $K_i f$ ($i = 1,2$)};
			\draw (0,-4) node[above=1.5mm]{Weighted $L^2$-norm of $f$: $\| \M^{- \frac{q}{2}} f \|^2_{L^2}$};
			\draw (0, -7) node[above=1.5mm]{Hypocoercivity of weighted $\L$: $\langle \Theta_\gamma \M^{-q} \L f, f \rangle$};
			\draw[decorate,decoration={brace,raise=7pt, amplitude=0.3cm}] (\rd,-0.5) -- ($(0,-3.4)!(\rd,-0.5)!(\rd,-3.4)$);
			\draw (\rd, -1.95) node[right=5mm]{Lemma \ref{Lmm-K1K2-L2}};
			\draw[decorate,decoration={brace,raise=7pt, amplitude=0.3cm}] (\rd,-3.6) -- ($(0,-6.5)!(\rd,-3.6)!(\rd,-6.5)$);
			\draw (\rd, -5.05) node[right=5mm]{Lemma \ref{Lmm-wL-Hypo}};
		\end{tikzpicture}
	\end{center}
	\caption{Outlines of proof of decay estimates for the pseudo-inverse $\L$.}\label{Fig4}
\end{figure}

The second step is to derive the hypocoercivity of the weighted linearized Boltzmann operator $\L$. To be more precise, one should prove that the quantity $\| \M^{- \frac{q}{2}} f \|^2_{L^2}$ can be bounded by $\langle \Theta_\gamma \M^{- \frac{q}{2}} \L f , \M^{- \frac{q}{2}} f \rangle + \| f \|^2_{L^2 (\nu)}$, as in Lemma \ref{Lmm-wL-Hypo}, where $\Theta_\gamma$ is given in \eqref{Theta-gamma}. The key ingredients of proof is to derive the low velocities estimates about $\langle \Theta_\gamma \M^{- \frac{q}{2}} K^{1 - \chi} f , \M^{- \frac{q}{2}} f \rangle$ in Lemma \ref{Lmm-K1K2-Hypo1} and the high velocities estimates on $\langle \Theta_\gamma \M^{- \frac{q}{2}} K^{\chi} f , \M^{- \frac{q}{2}} f \rangle$ in Lemma \ref{Lmm-K1K2-Hypo2}. The low velocities estimates obtain the upper bound $\| f \|^2_{L^2 (\nu)}$. On the other hand, the high velocities estimates gain the upper bound $l (r) \| \M^{- \frac{q}{2}} f \|^2_{L^2} + \| f \|^2_{L^2 (\nu)}$ for some factor $l (r) \to 0$ as $ r \to + \infty$, which is important to obtain the hypocoercivity of the weighted operator $\L$. Here $r > 0$ characterizes the magnitude of the velocity $v$. Remark that the quantity $\| f \|^2_{L^2(\nu)}$ can be bounded by $\langle \L f, f \rangle$ via the hypocoercivity \eqref{Hypo-L} of $\L$. Collecting the low and high velocities estimates before reduces to the hypocoercivity of the weighted linearized Boltzmann operator $\L$. Consequently, the proof of Theorem \ref{Thm1} can be completed in Subsection \ref{Subsec-Thm1}. Moreover, the key ingredients of the proof for decay estimates of $\L^{-1}$ can be intuitively sketched as the Figure \ref{Fig4} above.

When carrying out the high velocities estimates of $K^\chi_2 f$, there are different ways to deal with the integral factor $I (V_\shortparallel, \zeta_\perp) : = \int_{V_\perp \perp V_\shortparallel} \tfrac{\chi (|V|)}{|V|^{1 - \gamma}} \exp \big( - \tfrac{|V_\perp + \zeta_\perp|^2}{2 T} \big) \d V_ \perp$ of $K^\chi_2 f$ in \eqref{K2-chi-bnd} for the hard potential cases $0 \leq \gamma \leq 1$ and the soft potential cases $-3 < \gamma < 0$.
\begin{itemize}
	\item Hard potential cases $0 \leq \gamma \leq 1$: We dominate $I (V_\shortparallel, \zeta_\perp)$ by $\frac{2 \pi T}{r^{1 - \gamma}}$. If $0 \leq \gamma < 1$, the small factor $l (r) = \frac{2 \pi T}{r^{1 - \gamma}}$ is exactly what we required for sufficiently large $r > 0$. If $\gamma = 1$, by $\nu (v) \thicksim \langle v \rangle$, $\| \M^{- \frac{q}{2}} f \|^2_{L^2}$ can be bounded by $\| f \|^2_{L^2 (\nu)} + \tfrac{1}{1 + r} \| \M^{- \frac{q}{2}} f \|^2_{L^2 (\nu)}$. Here the small factor $l (r) = \tfrac{1}{1 + r}$ is required. For details, see Case 1 of proof of Lemma \ref{Lmm-K1K2-Hypo2}.
	
	\item Soft potential cases $- \frac{3}{2} < \gamma < 0$: We first decompose $I (V_\shortparallel, \zeta_\perp)$ into two parts
	\begin{equation*}
		\begin{aligned}
			I (V_\shortparallel, \zeta_\perp) = \Big\{ \underbrace{ \int_{|V_\perp'| > \tau |\zeta_\perp|} }_{I_1 (V_\shortparallel, \zeta_\perp)} + \underbrace{ \int_{|V_\perp'| > \tau |\zeta_\perp|} }_{I_2 (V_\shortparallel, \zeta_\perp)} \Big\} \tfrac{\chi (\sqrt{|V_\shortparallel|^2 + |V_\perp' - \zeta_\perp|^2})}{(|V_\shortparallel|^2 + |V_\perp' - \zeta_\perp|^2)^{\frac{1 - \gamma}{2}}} \exp \big( - \tfrac{|V_\perp'|^2}{2 T} \big) \d V_\perp'
		\end{aligned}
	\end{equation*}
    for some $\tau \in (0, 1)$. For the quantity $I_1 (V_\shortparallel, \zeta_\perp)$, we use the bound $\tfrac{\chi(\sqrt{|V_\shortparallel|^2 + |V_\perp' - \zeta_\perp|^2})}{( |V_\shortparallel|^2 + |V_\perp' - \zeta_\perp|^2)^{\frac{1 - \gamma}{2}}} \leq \tfrac{1}{r^{1-\gamma}}$ to get the small factor $l (r) = \frac{1}{r^{1 - \gamma}}$ as required (see Case 2(a) of proof of Lemma \ref{Lmm-K1K2-Hypo2}). For the quantity $I_2 (V_\shortparallel, \zeta_\perp)$, the inequality $$\tfrac{\chi(\sqrt{|V_\shortparallel|^2 + |V_\perp' - \zeta_\perp|^2})}{(|V_\shortparallel|^2 + |V_\perp' - \zeta_\perp|^2)^{\frac{1 - \gamma}{2}}} \leq \tfrac{ 2^{1 - \gamma} }{(1 + |v - \u|^2 + |v - \u + V_\shortparallel|^2)^{\frac{1 - \gamma}{2}}}$$  in Case 2(b) will be first employed. Then we apply the high-low velocities estimates arguments for the quantity $|v - \u| + |v - \u + V_\shortparallel|$. As in Case 2(b-1) of proof of Lemma \ref{Lmm-K1K2-Hypo2}, the high velocities estimates gain the bound $\tfrac{1}{1 + r} \| \M^{- \frac{q}{2}} f \|^2_{L^2}$ with the small factor $l (r) = \tfrac{1}{1 + r}$ required. On the other hand, as in Case 2(b-2) of proof of Lemma \ref{Lmm-K1K2-Hypo2}, the low velocities estimate gets a bound $\| f \|^2_{L^2 (\nu)}$.
\end{itemize}

(II) The case $- 3 < \gamma \leq - \frac{3}{2}$ in Theorem \ref{Thm1}: As illustrated in the previous case, the weighted $L^2$-norms of $f$ cannot dominate the weighted $L^\infty$ norms of $K f$ for $- 3 < \gamma \leq - \frac{3}{2}$, due to the integrability of $\mathbbm{k}_i (v, v_1)  \M^{- \frac{q}{2}} (v) \M^\frac{q}{2} (v_1)$. Inspired by \cite{Strain-Guo-2008-ARMA}, we can first derive the weighted $H^2$ bound of $\L^{-1} g$ in terms of that of $g$. Then the Sobolev embedding $H^2 \hookrightarrow L^\infty$ can cover the $L^\infty$ estimate of $\L^{-1} g$ with exponential decay weight. The key ingredient of proof is to obtain the inequality \eqref{L2-derivative} in Lemma \ref{Lmm-L-L2} below, namely, for any $0 < q < 1$, $l \in \R$ and $N \geq 0$,
\begin{equation*}
	\begin{aligned}
		\sum_{|\alpha| \leq N} \| \langle v \rangle^{\gamma (|\alpha| - l + \frac{1}{2})} \M^{ - \frac{q}{2} } \partial^\alpha_v f \|^2_{L^2} \leq C \sum_{|\alpha| \leq N} \| \langle v \rangle^{\gamma (|\alpha| - l - \frac{1}{2})} \M^{- \frac{q}{2} } \partial^\alpha_v g \|^2_{L^2} \,.
	\end{aligned}
\end{equation*}
Here $f = \L^{-1} g$. We remark that Strain-Guo \cite{Strain-Guo-2008-ARMA} proved the similar estimate with some sufficiently small $q > 0$ and the constant state $(\rho, \u, T) = (1,0,1)$. In this paper, we refine the proof in \cite{Strain-Guo-2008-ARMA} such that the above estimate holds for any $q \in (0, 1)$ and constant state $(\rho, \u, T)$ with $\rho, T > 0$.

The second part of this paper is to apply the decay estimates of $\L^{-1}$ to the Hilbert expansion of the scaled Boltzmann equation. Because the fluid variables $(\rho, \u, T)$, depending on $(t,x)$, subject to the compressible Euler equations, the linearized operator $\L$ thereby involves the parameters $t,x$. Consequently, we first justify the decay of $(t,x)$-derivatives of $\L^{-1}$, i.e., Theorem \ref{Thm2}. Denote by $f = \L^{-1} g$, i.e., $\L f = g$. We will finish our proof (bounding $\partial_{ t, x }^\beta f$) starting from $ \partial_{ t, x }^\beta ( \L f ) = \partial_{ t, x }^\beta g $, which is equivalent to
\begin{equation*}
	\begin{aligned}
		\L \{ (\I - \P) \partial_{t,x}^\beta f \} = \partial_{t,x}^\beta g - [\partial_{t,x}^\beta , \L ] f \,.
	\end{aligned}
\end{equation*}
Theorem \ref{Thm1} can establish the pointwise estimates of $(\I - \P) \partial_{t,x}^\beta f$ in terms of $\partial_{t,x}^\beta g$ and $[\partial_{t,x}^\beta , \L ] f$. So the key point is to control the commutator $[\partial_{t,x}^\beta , \L ] f$ by the induction arguments for $|\beta| = m$. For the case $ - \frac{3}{2} < \gamma \leq 1$, we can obtain the commutator estimate \eqref{f-beta-j+1-2} below, i.e.,
\begin{equation*}
	\begin{aligned}
		\| \langle v \rangle^{k_\gamma} \M^{- \frac{q_{\j + 1}}{2}} [ \partial_{t,x}^\beta \,, \L ] f (v_1) \d v_1 \|_{L^\infty} \leq C (\e_{\j + 1}) \mathcal{W}^{\j + 1, \infty}_{q_0,\cdots,q_{\j + 1}} (g) \,.
	\end{aligned}
\end{equation*}
For the case $- 3 < \gamma \leq - \frac{3}{2} $, the bound \eqref{Finish9} below are established, hence,
\begin{equation*}
	\begin{aligned}
		& | \langle \langle v \rangle^{  2 \gamma l} \M^{ - q } \partial_v^\alpha [ \partial_{t,x}^\beta , \L ] f , \partial_v^\alpha \partial_{t,x}^\beta f \rangle | \\
		\lesssim & \mathfrak{e}_m^m \| \langle v \rangle^{ \gamma ( l + \frac{1}{2} ) } \M^{ - \frac{q}{2} } \partial_v^\alpha \partial_{t,x}^\beta f \|_{L^2} \\
		& \times \sum_{ 0 \neq \beta' \leq \beta } \Big\{ \| \partial_{t,x}^{\beta - \beta'} f \|_{ L^2 ( \nu ) } + \sum_{|\alpha'| \leq |\alpha|} \| \langle v \rangle^{ \gamma ( l + \frac{1}{2} ) } \M^{ - \frac{q}{2} } \partial_v^{ \alpha' } \partial_{t,x}^{ \beta - \beta' } f \|_{ L^2 } \Big\} \,.
	\end{aligned}
\end{equation*}
While controlling the fluid part $ \P \partial_{ t, x }^\beta f $, the equation \eqref{P-f-beta-j+1-0} below, that is,
\begin{equation*}
	\begin{aligned}
		\langle \partial_{t,x}^\beta f, \phi \rangle = - \sum_{0 \neq \tilde{\beta}' \leq \tilde{\beta}} C_{\tilde{\beta}}^{\tilde{\beta}'} \langle \partial_{t,x}^{\tilde{\beta} + \bar{\beta} - \tilde{\beta}'} f , \partial_{t,x}^{\tilde{\beta}'} \phi \rangle - \sum_{\tilde{\beta}'' \leq \tilde{\beta}} C_{\tilde{\beta}}^{\tilde{\beta}''} \langle \partial_{t,x}^{\tilde{\beta} - \tilde{\beta}''} f , \partial_{t,x}^{\tilde{\beta}'' + \bar{\beta}} \phi \rangle \,.
	\end{aligned}
\end{equation*}
is very important in the induction arguments for $|\beta| = m$. Here $\phi = \phi (t,x,v)$ is given in \eqref{phi} below.

Finally, in Theorem \ref{Thm3}, we prove the decay bounds of the Hilbert expansions $F_n (t,x,v)$ $(n \geq 1)$ and their $(t,x)$-derivatives. Denote by $f_n = F_n / \sqrt{\M}$, which satisfies $(\I - \P) f_n = \L^{-1} g_n$, see \eqref{f2-I-P} below. The kernel of proving Theorem \ref{Thm3} is to verify that uniformly in $t,x$,
\begin{equation*}
	\begin{aligned}
		\mathcal{W}^{m, \infty}_{q_0, \cdots, q_m} (g_n) < \infty \ ( - \frac{3}{2} < \gamma \leq 1 ) \textrm{ and } \mathcal{W}_q^{m, 2} (g_n) < \infty \ ( - 3 < \gamma \leq - \frac{3}{2} ) \,.
	\end{aligned}
\end{equation*}
Here the quantities $\mathcal{W}^{m, \infty}_{q_0, \cdots, q_m} (g_n)$ and $\mathcal{W}_q^{m, 2} (g_n)$ are defined in Theorem \ref{Thm2}. The ideas are that, based on the decay estimates in Theorem \ref{Thm2}, the induction arguments for $n \geq 1$ can achieve our goal.

\subsection{Organization of this paper}

In the next section, the expressions of the integral kernels of $\mathbbm{k}_1 (v,v_1)$ and $\mathbbm{k}_2 (v,v_1)$ are derived in Lemma \ref{Lmm-k1-k2}, and their bounds are verified in Lemma \ref{Lmm-k1k2-bnd}. In Section \ref{Sec-3}, the pointwise decay estimates of $\L^{-1}$ in Theorem \ref{Thm1} are given. Section \ref{Sec-4} is aim at justifying decay estimates of $(t,x)$-derivatives of $\L^{-1}$ in Theorem \ref{Thm2}, and deriving the bounds of the expanded terms $F_n (t,x,v)$ ($n \geq 1$) in Theorem \ref{Thm3}. In Section \ref{Sec:K-LHVE-K}, the details of proof for low-high velocities estimates of $K$ (the case $ - \frac{3}{2} < \gamma \leq 1 $) are derived, i.e., proving Lemma \ref{Lmm-K1K2-Hypo1} and \ref{Lmm-K1K2-Hypo2}. In Section \ref{Sec:VSP}, the weighted $L^2$ estimates of $K f$ for the case $- 3 < \gamma \leq - \frac{3}{2}$ are obtained, namely, verifying Lemma \ref{Lmm-VSP}. In Section \ref{Sec:Lmm4.1}, we justify the bound of weighted $L^2$ norm of the mixed derivative $\partial_v^\alpha \partial_{ t, x }^\beta f$, hence, proving Lemma \ref{Lmm-SVP-Der}.

\section{Estimates of the integral kernel of the operator $K$}\label{Sec-2}

This section is aimed at estimating the integral kernel of the operator $K$. Recalling the decomposition $K f = K_1 f - K_2 f$ in \eqref{K-K1-K2}, one respectively rewrites the operators $K_1 f$ and $K_2 f$ in \eqref{K1} and \eqref{K2} as
\begin{equation}\label{K1K2-integ-kernel-forms}
	\begin{aligned}
		K_1 f (v) = \int_{\R^3} \mathbbm{k}_1 (v,v_1) f(v_1) \d v_1 \,, \ K_2 f (v) = \int_{\R^3} \mathbbm{k}_2 (v,v_1) f(v_1) \d v_1 \,.
	\end{aligned}
\end{equation}
Then the integral kernels $\mathbbm{k}_1 (v,v_1)$ and $\mathbbm{k}_2 (v,v_1)$ can be explicitly expressed as follows.

\begin{lemma}\label{Lmm-k1-k2}
	Let $-3 < \gamma \leq 1$, and the collision kernel $b$ in \eqref{b} subject to the conditions \eqref{Symmetric-Asmp}-\eqref{Cutoff}. Then the integral kernel $\mathbbm{k}_1 (v,v_1)$ is
	\begin{equation}\label{k1}
		\begin{aligned}
			\mathbbm{k}_1 (v,v_1) = \M^\frac{1}{2} (v) \M^\frac{1}{2} (v_1) |v_1 - v|^\gamma \,,
		\end{aligned}
	\end{equation}
    and the integral kernel $\mathbbm{k}_2 (v,v_1)$ can be expressed as follows:
    \begin{enumerate}
    	\item if $\gamma = 1$,
    	\begin{equation}\label{k2-1}
    		\begin{aligned}
    			\mathbbm{k}_2 (v,v_1) = \tfrac{4 \beta_0 \rho}{(2 \pi T)^\frac{1}{2}} | v_1 - v |^{-1} \exp \big[ - \tfrac{|v_1 - v|^2}{8 T} - \tfrac{(|v_1|^2 - |v|^2)^2}{8 T |v_1 - v|^2} + \tfrac{|v_1|^2 - |v|^2}{2 T |v_1 - v|^2} \u \cdot (v_1 - v) \big] \,;
    		\end{aligned}
    	\end{equation}
        \item if $- 3 < \gamma < 1$,
        \begin{equation}\label{k2-2}
        	\begin{aligned}
        		\mathbbm{k}_2 (v,v_1) = \tfrac{4 \rho}{(2 \pi T)^\frac{3}{2}} |v_1 - v|^{-2} \exp \Big[ - \tfrac{|v_1 - v|^2}{8 T} - \tfrac{(|v_1|^2 - |v|^2)^2}{8 T |v_1 - v|^2} + \tfrac{|v_1|^2 - |v|^2}{2 T |v_1 - v|^2} \u \cdot (v_1 - v) \\
        		\times \int_{y \perp (v_1 - v)} \exp ( - \tfrac{|y + \zeta_2 - \u|^2}{2 T} ) ( |v_1 - v|^2 + |y|^2 )^\frac{\gamma}{2} \beta ( \arccos \tfrac{|v_1 - v|}{\sqrt{|v_1 - v|^2 + |y|^2}} ) \d y \Big] \,,
        	\end{aligned}
        \end{equation}
        where $\zeta_2 = \frac{(v_1 - v) \wedge (v_1 + v) \wedge (v_1 - v)}{2 |v_1 - v|^2}$.
    \end{enumerate}	
\end{lemma}

\begin{proof}
	The proof for the case $[\rho, \u, T] = [1,0,1]$ can be found in \cite{Grad-1963} and \cite{LS-2010-KRM}, and the arguments of the general case $[\rho, \u, T]$ are similar to that therein. However, for convenience for readers, one sketches the outline of proof for the general case $[\rho, \u, T]$ here, in which there is an additional quantity $\tfrac{|v_1|^2 - |v|^2}{2 T |v_1 - v|^2} \u \cdot (v_1 - v)$ in the expression of $\mathbbm{k}_2 (v,v_1)$.
	
	The verification of $\mathbbm{k}_1 (v,v_1)$ is trivial. The goal is to focus on the integral kernel $\mathbbm{k}_2 (v,v_1)$. Let
	\begin{equation*}
		\begin{aligned}
			v_1 - v = \xi + y \,, \xi = [\omega \cdot (v_1 - v)] \omega \,, y = (v_1 - v) - [\omega \cdot (v_1 - v)] \omega = \omega \wedge (v_1 - v) \wedge \omega \,,
		\end{aligned}
	\end{equation*}
	where $(v_1, v, \omega)$ is given in Figure \ref{Fig1}. From Section 2 of \cite{Grad-1963}, one has
	\begin{equation}
		\begin{aligned}
			\d \omega \d v_1 = \d \omega \d (v_1 - v) = \frac{2}{|\xi|^2} \d \xi \d y \,.
		\end{aligned}
	\end{equation}
    Moreover, the relations in Figure \ref{Fig1} show that
    \begin{equation}
    	\begin{aligned}
    		v' = \xi + v \,, \ v_1 = v + \xi + y \,.
    	\end{aligned}
    \end{equation}
	Then, from \eqref{K2},
	\begin{equation}
		\begin{aligned}
			K_2 f (v) = 4 \M^\frac{1}{2} (v) \int_{\R^3} \int_{y \perp \xi} \M^{- \frac{1}{2}} (\xi + v) \M (v+ \xi + y) f (v + \xi) \\
			\times | \xi |^{-2} \beta (\arccos \tfrac{|\xi|}{|\xi + y|} ) |\xi + y|^\gamma \d y \d \xi \,,
		\end{aligned}
	\end{equation}
	where the relation $\cos \theta = \frac{(v_1 - v) \cdot \omega}{|v_1 - v|} = \frac{|\xi|}{|\xi + y|}$ is used. Since $v' = v + \xi$, one has $\d v' = \d \xi$ for any fixed $v$. Then there holds
	\begin{equation}
		\begin{aligned}
			K_2 f (v) = 4 \M^\frac{1}{2} (v) \int_{\R^3} \int_{y \perp (v'-v)} \M^{- \frac{1}{2}} (v') \M (v' + y) f (v') | v'-v |^{-2} \\
			\times \beta (\arccos \tfrac{|v'-v|}{|v'-v + y|} ) |v'-v + y|^\gamma \d y \d v' \,,
		\end{aligned}
	\end{equation}
	which means that
	\begin{equation}
		\begin{aligned}
			\mathbbm{k}_2 (v,v') = 4 \M^\frac{1}{2} (v) \M^{- \frac{1}{2}} (v') & | v'-v |^{-2} \\ & \times \int_{y \perp (v'-v)} \M (v' + y) \beta (\arccos \tfrac{|v'-v|}{|v'-v + y|} ) |v'-v + y|^\gamma \d y \,.
		\end{aligned}
	\end{equation}
    One first splits
    \begin{equation*}
    	\begin{aligned}
    		\tfrac{v'+v}{2} = \zeta_1 + \zeta_2 \,, \ \zeta_1 \parallel (v' - v) \,, \ \zeta_2 \parallel y \,.
    	\end{aligned}
    \end{equation*}
    More specifically,
    \begin{equation*}
    	\begin{aligned}
    		\zeta_1 = \frac{|v'|^2 - |v|^2}{2 |v' - v|^2} (v' - v) \,, \ \zeta_2 = \frac{(v_1 - v) \wedge (v_1 + v) \wedge (v_1 - v)}{2 |v_1 - v|^2} \,.
    	\end{aligned}
    \end{equation*}
	By the facts $y \cdot (v'-v) = y \cdot \zeta_1 = \zeta_1 \cdot \zeta_2 = 0$, a direct calculation implies
	\begin{equation}
		\begin{aligned}
			& - \tfrac{|v - \u|^2}{4 T} + \tfrac{|v' - \u|^2}{4 T} - \tfrac{|v'+y-\u|^2}{2 T} \\
			= & - \tfrac{|v' - v|^2}{8 T} - \tfrac{|y + \frac{v'+v}{2} - \u|^2}{2 T} \\
			= & - \tfrac{|v' - v|^2}{8 T} - \tfrac{(|v'|^2 - |v|^2)^2}{8 T |v' - v|^2} + \tfrac{|v'|^2 - |v|^2}{2 T |v' - v|^2} \u \cdot (v' - v) - \tfrac{|y + \zeta_2 - \u|^2}{2 T} \,.
		\end{aligned}
	\end{equation}
	Thus,
	\begin{equation*}
		\begin{aligned}
			\M^\frac{1}{2} (v) \M^{- \frac{1}{2}} (v') & \M (v' + y) \\
			= & \tfrac{\rho}{(2 \pi T)^\frac{3}{2}} \exp \big[ - \tfrac{|v' - v|^2}{8 T} - \tfrac{(|v'|^2 - |v|^2)^2}{8 T |v' - v|^2} + \tfrac{|v'|^2 - |v|^2}{2 T |v' - v|^2} \u \cdot (v' - v) - \tfrac{|y + \zeta_2 - \u|^2}{2 T} \big] \,.
		\end{aligned}
	\end{equation*}
	It therefore infers that
	\begin{equation}
		\begin{aligned}
			\mathbbm{k}_2 (v,v') = \tfrac{4 \rho}{(2 \pi T)^\frac{3}{2}} |v' - v|^{-2} \exp \big[ - \tfrac{|v' - v|^2}{8 T} - \tfrac{(|v'|^2 - |v|^2)^2}{8 T |v' - v|^2} + \tfrac{|v'|^2 - |v|^2}{2 T |v' - v|^2} \u \cdot (v' - v) \big] \\
			\times \int_{y \perp (v'-v)} \exp \big[ - \tfrac{|y + \zeta_2 - \u|^2}{2 T} \big] \beta (\arccos \tfrac{|v'-v|}{|v'-v + y|} ) |v'-v + y|^\gamma \d y \,.
		\end{aligned}
	\end{equation}

	If $\gamma = 1$, $\beta (\theta) = \beta_0 |\cos \theta|$. Then
	\begin{equation}
		\begin{aligned}
			\int_{y \perp (v'-v)} \exp \big[ - \tfrac{|y + \zeta_2 - \u|^2}{2 T} \big] \beta (\arccos \tfrac{|v'-v|}{|v'-v + y|} ) |v'-v + y|^\gamma \d y \\
			= \beta_0 |v' - v| \int_{y \perp (v'-v)} \exp \big[ - \tfrac{|y + \zeta_2 - \u|^2}{2 T} \big] \d y = 2 \pi T \beta_0 |v' - v| \,.
		\end{aligned}
	\end{equation}
	Namely, \eqref{k2-1} holds.
	
	If $- 3 < \gamma < 1$, $y \perp (v' -v)$ implies that $|v'-v+y| = \sqrt{|v'-v|^2 + |y|^2}$. Then \eqref{k2-2} holds, and the proof of Lemma \ref{Lmm-k1-k2} is completed.
\end{proof}

Next, the bounds of the integral kernels $\mathbbm{k}_1 (v, v_1)$ and $\mathbbm{k}_2 (v, v_1)$ will be derived as following lemma.

\begin{lemma}\label{Lmm-k1k2-bnd}
	Let $- 3 < \gamma \leq  1$. Set $b_0 \in [0, 1 - \gamma]$ if $- 1 < \gamma \leq 1$, and $b_0 \in (- \gamma - 1, 2)$ if $- 3 < \gamma \leq -1$. Then, for any fixed $\epsilon \in (0,1)$,
	\begin{equation}\label{k1k2-bnd}
		\begin{aligned}
			\mathbbm{k}_1 (v, v_1) = & \tfrac{\rho}{(2 \pi T)^\frac{3}{2}} |v_1 - v|^\gamma \exp \big( - \tfrac{|v - \u|^2 + |v_1 - \u|^2}{4 T} \big) \leq \mathbf{k}_1 (v_1 - v) \,, \\
			\mathbbm{k}_2 (v, v_1) \leq & \tfrac{c_1 \rho}{(2 \pi T)^\frac{3}{2}} (1 + T) e^\frac{|\u|^2}{2 \epsilon T} \tfrac{1}{|v_1 - v|^{2 - b_0 - \gamma}} \exp \big[ - \tfrac{|v_1 - v|^2}{8 T} - \tfrac{(1-\epsilon)(|v_1|^2 - |v|^2)^2}{8 T |v_1 - v|^2} \big] \leq \mathbf{k}_2 (v_1 - v) \,,
		\end{aligned}
	\end{equation}
	where $c_1 > 0$ is a constant, and
	\begin{equation*}
		\begin{aligned}
			\mathbf{k}_1 (v_1 - v) = & \tfrac{\rho}{(2 \pi T)^\frac{3}{2}} |v_1 - v|^\gamma \exp \big( - \tfrac{ |v_1 - v|^2}{8 T} \big) \,, \\
			\mathbf{k}_2 (v_1 - v) = & \tfrac{c_1 \rho}{(2 \pi T)^\frac{3}{2}} (1 + T) e^\frac{|\u|^2}{2 \epsilon T} \tfrac{1}{|v_1 - v|^{2 - b_0 - \gamma}} \exp \big[ - \tfrac{|v_1 - v|^2}{8 T} \big] \,.
		\end{aligned}
	\end{equation*}
	Moreover, $\mathbf{k}_1 (v)$ and $\mathbf{k}_2 (v)$ are both in $L^1$.	
\end{lemma}

\begin{proof}
	The first inequality on $\mathbbm{k}_1 (v, v_1)$ in \eqref{k1k2-bnd} can be directly implied by the elementary inequality $|v - \u|^2 + |v_1 - \u|^2 \geq \frac{1}{2} |v_1 - v|^2$.
	
	Now we focus on the second inequality in \eqref{k1k2-bnd}. One first claims that if $b_0 < 2$ and $0 \leq b_0 \leq 1 - \gamma$,
	\begin{equation}\label{ClaimA}
		\begin{aligned}
			A := \int_{y \perp (v_1 - v)} \exp ( - \tfrac{|y + \zeta_2 - \u|^2}{2 T} ) ( |v_1 - v|^2 + |y|^2 )^\frac{\gamma}{2} \beta ( \arccos \tfrac{|v_1 - v|}{\sqrt{|v_1 - v|^2 + |y|^2}} ) \d y \\
			\leq 2 \pi c_0 \beta_0 (T + \tfrac{1}{2 - b_0}) |v_1 - v|^{b_0 + \gamma} \,.
		\end{aligned}
	\end{equation}
	Indeed, by \eqref{Cutoff}, one has
	\begin{equation}
		\begin{aligned}
			\beta (\arccos \tfrac{|v_1 -v|}{\sqrt{|v_1 -v|^2 + |y|^2}}) \leq \beta_0|v_1-v| \big( |v_1-v|^2 + |y|^2 \big)^{- \frac{1}{2}} \,.
		\end{aligned}
	\end{equation}
	Then, by Lemma \ref{Lmm-f-bnd} below,
	\begin{equation}
		\begin{aligned}
			A \leq & \beta_0 |v_1 - v| \int_{y \perp (v_1 - v)} \exp ( - \tfrac{|y + \zeta_2 - \u|^2}{2 T} ) ( |v_1 - v|^2 + |y|^2 )^\frac{\gamma - 1}{2} \d y \\
			\leq & c_0 \beta_0 |v_1 - v|^{b_0 + \gamma} \int_{y \perp (v_1 - v)} |y|^{-b_0} \exp ( - \tfrac{|y + \zeta_2 - \u|^2}{2 T} )  \d y \\
			\leq & c_0 \beta_0 |v_1 - v|^{b_0 + \gamma} \Big( \int_{\substack{y \perp (v_1 - v) \\ |y| \leq 1}} |y|^{-b_0} \d y + \int_{\substack{y \perp (v_1 - v) \\ |y| > 1}} \exp ( - \tfrac{|y + \zeta_2 - \u|^2}{2 T} )  \d y \Big) \\
			\leq & 2 \pi c_0 \beta_0 (T + \tfrac{1}{2 - b_0}) |v_1 - v|^{b_0 + \gamma} \,,
		\end{aligned}
	\end{equation}
	provided that $b_0 < 2$ and $0 \leq b_0 \leq 1 - \gamma$. The claim \eqref{ClaimA} holds.
	
	From \eqref{k2-1}-\eqref{k2-2} and \eqref{ClaimA}, it is derived that
	\begin{equation}
		\begin{aligned}
			\mathbbm{k}_2 (v,v_1) \leq \tfrac{c_1 \rho}{(2 \pi T)^\frac{3}{2}} (1 + T) & \tfrac{1}{|v_1 - v|^{2 - b_0 - \gamma}} \\
			\times & \exp \big[ - \tfrac{|v_1 - v|^2}{8 T} - \tfrac{(|v_1|^2 - |v|^2)^2}{8 T |v_1 - v|^2} + \tfrac{|v_1|^2 - |v|^2}{2 T |v_1 - v|^2} \u \cdot (v_1 - v) \big] \,,
		\end{aligned}
	\end{equation}
	where $c_1 = 8 \pi \beta_0 \max \{ 1, \tfrac{c_0 (3-b_0)}{2 - b_0} \} > 0$. For any $\epsilon \in (0,1)$, straightforward calculations reduce to
	\begin{equation}\label{Cross-u}
		\begin{aligned}
			& - \tfrac{ \epsilon (|v_1|^2 - |v|^2)^2}{8 T |v_1 - v|^2} + \tfrac{|v_1|^2 - |v|^2}{2 T |v_1 - v|^2} \u \cdot (v_1 - v) \\
			= &  - \tfrac{\epsilon}{8 T} \big( \tfrac{|v_1|^2 - |v|^2}{|v_1 - v|} + \tfrac{2}{\epsilon} \u \cdot \tfrac{v_1 - v}{|v_1 - v|} \big)^2 + \tfrac{1}{2 \epsilon T} (\u \cdot \tfrac{v_1 - v}{|v_1 - v|})^2 \leq \tfrac{|\u|^2}{2 \epsilon T} \,.
		\end{aligned}
	\end{equation}
	It thereby infers that
	\begin{equation}
		\begin{aligned}
			\mathbbm{k}_2 (v,v_1) \leq \tfrac{c_1 \rho}{(2 \pi T)^\frac{3}{2}} (1 + T) e^\frac{|\u|^2}{2 \epsilon T} \tfrac{1}{|v_1 - v|^{2 - b_0 - \gamma}} \exp \big[ - \tfrac{|v_1 - v|^2}{8 T} - \tfrac{(1 - \epsilon)(|v_1|^2 - |v|^2)^2}{8 T |v_1 - v|^2} \big]
		\end{aligned}
	\end{equation}
	for any fixed $\epsilon \in (0,1)$, hence the second inequality in \eqref{k1k2-bnd} holds.
	
	Finally, the $L^1$-integrability of $\mathbf{k}_1 (v)$ and $\mathbf{k}_2 (v)$ on $v \in \R^3$ will be justified. One asserts that, for any fixed $a > 0$,
	\begin{equation}\label{Assert-integrable}
		\begin{aligned}
			\textrm{the integral } \int_{\R^3} |v|^k e^{- a |v|^2} \d v \textrm{ converges, if and only if } k > - 3 \,.
		\end{aligned}
	\end{equation}
	Indeed, one has
	\begin{equation}
		\begin{aligned}
			\int_{\R^3} |v|^k e^{- a |v|^2} \d v = \Big\{ \int_{|v| \leq 1} + \int_{|v| > 1} \Big\} |v|^k e^{- a |v|^2} \d v \,.
		\end{aligned}
	\end{equation}
	The last part is obviously convergent for any $k \in \R$. The first part can be bounded by
	\begin{equation}
		\begin{aligned}
			e^{-a} \int_{|v| \leq 1} |v|^k \d v \leq \int_{|v| \leq 1} |v|^k e^{- a |v|^2} \d v \leq \int_{|v| \leq 1} |v|^k \d v \,.
		\end{aligned}
	\end{equation}
	It is easy to verify that $\int_{|v| \leq 1} |v|^k \d v$ converges if and only if $k > - 3$. The assertion \eqref{Assert-integrable} holds. Therefore, $\mathbf{k}_1 (v)$ and $\mathbf{k}_2 (v)$ are both in $L^1$, if and only if $\gamma > - 3$ and $b_0 + \gamma - 2 > -3$. Together with the constraints $b_0 < 2$ and $0 \leq b_0 \leq 1 - \gamma$, $b_0$ exactly subjects to the conditions given in Lemma \ref{Lmm-k1k2-bnd}. Then the proof is completed.	
\end{proof}

\begin{lemma}\label{Lmm-f-bnd}
	Let $- 3 < \gamma \leq 1$. For any $0 \leq b_0 \leq 1 - \gamma$ and $a > 0$, the function
	\begin{equation*}
		\begin{aligned}
			f (x) = x^{b_0} (x^2 + a^2)^\frac{\gamma - 1}{2}
		\end{aligned}
	\end{equation*}
    is bounded on $(0, + \infty)$, and
    \begin{equation*}
    	\begin{aligned}
    		0 \leq f (x) \leq c_0 a^{b_0 + \gamma - 1} \,,
    	\end{aligned}
    \end{equation*}
    where the constant $c_0 > 0$ is independent of $a > 0$.
\end{lemma}

\begin{proof}
	If $\gamma = 1$, $b_0 = 0$ and $f (x) \equiv 1$. The conclusion holds. If $- 3 < \gamma < 1$, a direct calculation shows
	\begin{equation*}
		\begin{aligned}
			f'(x) = \tfrac{(b_0 + \gamma - 1) x^2 + b_0 a^2}{x (x^2 + a^2)} f (x) \,.
		\end{aligned}
	\end{equation*}
	Note that $f(x) > 0$ for $x > 0$.
	
	(1) If $b_0 + \gamma - 1 = 0$, then $ f'(x) \geq 0$ for $x > 0$. Namely, $f (x)$ is increasing on $(0, + \infty)$. It thereby infers that $0 = \lim_{x \to 0+} f (x) \leq f (x) \leq \lim_{x \to + \infty} f(x) = 1 = a^{b_0 + \gamma - 1}$. The conclusion holds.
	
	(2) If $b_0 = 0$, then $f' (x) < 0$ for $x > 0$. Hence $f(x)$ is decreasing on $x > 0$. Then there holds that $0 = \lim_{x \to + \infty} f(x) \leq f(x) \leq \lim_{x \to 0+} f(x) = a^{\gamma - 1} = a^{b_0 + \gamma - 1}$. The conclusion holds.
	
	(3) If $0 < b_0 < 1 - \gamma$, then $f(x)$ is increasing on $(0, \sqrt{\frac{b_0 a^2}{1-\gamma - b_0}})$, and decreasing on $( \sqrt{\frac{b_0 a^2}{1-\gamma - b_0}}, + \infty )$. Observe that $\lim_{x \to 0+} f(x) = \lim_{x \to + \infty} f(x) = 0$ and $f(\sqrt{\frac{b_0 a^2}{1-\gamma - b_0}}) = c_0 a^{b_0 + \gamma - 1}$, where $c_0 = \frac{b_0^\frac{b_0}{2} (1 - \gamma)^\frac{\gamma-1}{2}}{(1 - \gamma - b_0)^\frac{b_0 + \gamma - 1}{2}} > 0$. Then $0 \leq f (x) \leq c_0 a^{b_0 + \gamma - 1}$, and the proof of Lemma \ref{Lmm-f-bnd} is finished.	
\end{proof}

\section{Pointwise decay estimates of $\L^{-1}$: Proof of Theorem \ref{Thm1}}\label{Sec-3}

In this section, the pointwise decay estimates of the pseudo-inverse operator $\L^{-1}$ with $-3 < \gamma \leq 1$. For the cases $- \frac{3}{2} < \gamma \leq 1$, the proof will be derived by three steps: (1) to bound $|K f|$ by some weighted $L^2$ norms of $f$; (2) to control the weighted $L^2$ norms of $f$ by employing the hypocoercivity of some weighted-$\L$ operators; (3) to finish the pointwise decay estimates of the pseudo-inverse operator $\L^{-1}$. Remark that the restriction $\gamma > - \tfrac{3}{2}$ will be required in the first steps. For the cases $- 3 < \gamma \leq - \frac{3}{2}$, we mainly aim at proving that the quantity $\| \M^{- q_0} f \|_{H^2}$ can be bounded by the quantity $\sum_{|\beta| \leq 2} \| \langle v \rangle^{\gamma (|\beta| - 3)} \M^{- q_0} \partial^\beta_v g \|_{L^2}$. Then the Sobolev embedding theory infers the required pointwise decay estimate of $\L^{-1}$.

Throughout this section, the notation $A \lesssim_\alpha B$ will be employed to denote by the inequality $A \leq C(\alpha) B$ for some constant $C(\alpha) > 0$ depending on the parameter $\alpha$.

\subsection{Pointwise controls of $|K f|$}

This subsection is aimed at dominating the quantities $|K_1 f|$ and $|K_2 f|$ by some proper weighted $L^2$ norms of $f$, where $K_1 f$ and $K_2 f$ are respectively defined in \eqref{K1} and \eqref{K2} before. More precisely, the following results about the hard potential case and part of soft potential case will be given.

\begin{lemma}\label{Lmm-K1K2-L2}
	For any fixed $- \frac{3}{2} < \gamma \leq 1$ and $0 < q < 1$, there hold
	\begin{equation}\label{K1-L2Bnd-H}
		\begin{aligned}
			\M^{- \frac{q}{2}} (v) |K_1 f (v)| \lesssim_{\rho, \u, T} \| f \|_{L^2} \,, \ \forall \, v \in \R^3 \,,
		\end{aligned}
	\end{equation}
	and
	\begin{equation}\label{K2-L2Bnd-H}
		\begin{aligned}
			\M^{- \frac{q}{2}} (v) |K_2 f (v)| \lesssim_{\rho, \u, T} \| \M^{- \frac{q}{2}} f \|_{L^2} \,, \ \forall \, v \in \R^3 \,.
		\end{aligned}
	\end{equation}	
    Furthermore, by the decomposition $K f = K_1 f - K_2 f$, there holds
    \begin{equation}
    	\begin{aligned}
    		\M^{- \frac{q}{2}} (v) |K f (v)| \lesssim_{\rho, \u, T} \| \M^{- \frac{q}{2}} f \|_{L^2} + \| f \|_{L^2} \,, \ \forall \, v \in \R^3 \,.
    	\end{aligned}
    \end{equation}
\end{lemma}

\begin{remark}
	As indicated in the proof of Lemma \ref{Lmm-K1K2-L2} below, the condition $- \frac{3}{2} < \gamma \leq 1$ is actually due to the requirement of the square integrability of the integral kernels $\mathbbm{k}_1 (v,v_1)$ and $\mathbbm{k}_2 (v,v_1)$, or corresponding to the bounds $\mathbf{k}_1 (v_1 - v)$ and $\mathbf{k}_2 (v_1 - v)$ given in Lemma \ref{Lmm-k1k2-bnd}.
\end{remark}

\begin{proof}[Proof of Lemma \ref{Lmm-K1K2-L2}]
	Recalling the definition of $K_1 f$ in \eqref{K1} and Lemma \ref{Lmm-k1k2-bnd}, one has
	\begin{equation}
		\begin{aligned}
			\M^{- \frac{q}{2}} |K_1 f| & \lesssim_{\rho, \u, T} \exp ( \tfrac{q |v - \u|^2}{4 T} ) \int_{\R^3} |v_1 - v|^\gamma \exp ( - \tfrac{|v - \u|^2 + |v_1 - \u|^2}{4 T} ) |f (v_1)| \d v_1 \\
			& \lesssim_{\rho, \u, T} \underbrace{ \exp ( \tfrac{q |v - \u|^2}{4 T} ) \Big( \int_{\R^3} |v_1 - v|^{2 \gamma} \exp ( - \tfrac{|v - \u|^2 + |v_1 - \u|^2}{2 T} ) \d v_1 \Big)^\frac{1}{2} }_{B} \| f \|_{L^2} \,,
		\end{aligned}
	\end{equation}
	where the last inequality is derived from the H\"older inequality. Since $|v - \u|^2 + |v_1 - \u|^2 \geq \frac{1}{2} |v_1 - v |$ and $0 < q < 1$, one derives from \eqref{Assert-integrable} that if $2 \gamma > - 3$,
	\begin{equation}\label{B2}
		\begin{aligned}
			B^2 \leq & \int_{\R^3} |v_1 - v|^{2 \gamma} \exp \big[ - \tfrac{(1-q) (|v - \u|^2 + |v_1 - \u|^2)}{2 T} \big] \d v_1 \\
			\leq & \int_{\R^3} |v_1 - v|^{2 \gamma} \exp \big[ - \tfrac{(1-q) |v_1 - v|^2 }{4 T} \big] \d v_1 \\
			\leq & \int_{|v| \leq 1} |v|^{2 \gamma} \d v + \int_{\R^3} \exp \big[ - \tfrac{(1-q) |v|^2 }{4 T} \big] \d v = \tfrac{4 \pi }{2 \gamma + 3} + ( \tfrac{4 \pi T}{1 - q} )^\frac{3}{2} < \infty \,.
		\end{aligned}
	\end{equation}
	It therefore implies the estimate \eqref{K1-L2Bnd-H}, provided that $2 \gamma > - 3$ is further assumed.
	
	Next, one tries to verify the bound \eqref{K2-L2Bnd-H}. By Lemma \ref{Lmm-k1k2-bnd} and the H\"older inequality, there hold
	\begin{equation}
		\begin{aligned}
			\M^{- \frac{q}{2}} |K_2 f| \lesssim_{\rho, \u, T} \exp & \big( \tfrac{q |v - \u|^2}{4 T} \big) \\
			\times & \int_{\R^3} \tfrac{1}{|v_1 - v|^{2 - b_0 - \gamma}} \exp \big[ - \tfrac{|v_1 - v|^2}{8 T} - \tfrac{(1 - \epsilon) (|v_1|^2 - |v|^2)^2}{8 T} \big] |f (v_1)| \d v_1 \\
			\lesssim_{\rho, \u, T} D \| & \M^{- \frac{q}{2}} f \|_{L^2} \,,
		\end{aligned}
	\end{equation}
	where $\epsilon \in (0,1) $ is arbitrary, $b_0$ is assumed in Lemma \ref{Lmm-k1k2-bnd}, and $D$ is expressed by
	\begin{equation}
		\begin{aligned}
			D : = \Big( \int_{\R^3} \tfrac{1}{|v_1 - v|^{2( 2 - b_0 - \gamma )}} \exp \big[ \tfrac{q ( |v - \u|^2 - |v_1 - \u|^2 )}{2 T} - \tfrac{|v_1 - v|^2}{4 T} - \tfrac{(1 - \epsilon) (|v_1|^2 - |v|^2)^2}{4 T} \big] \d v_1 \Big)^\frac{1}{2} \,.
		\end{aligned}
	\end{equation}
	Note that
	\begin{equation}
		\begin{aligned}
			|v - \u|^2 - |v_1 - \u|^2 = - |v_1 - v|^2 - 2 (v_1 - v) \cdot (v - \u) \,,
		\end{aligned}
	\end{equation}
	and
	\begin{equation}
		\begin{aligned}
			|v_1|^2 - |v|^2 = |v_1 - v|^2 + 2 v \cdot (v_1 - v) \,.
		\end{aligned}
	\end{equation}
	Then there holds
	\begin{equation}
		\begin{aligned}
			& \tfrac{q ( |v - \u|^2 - |v_1 - \u|^2 )}{2 T} - \tfrac{|v_1 - v|^2}{4 T} - \tfrac{(1 - \epsilon) (|v_1|^2 - |v|^2)^2}{4 T} \\
			= & - \tfrac{(2 q + 2 - \epsilon) |v_1 - v|^2}{4 T} - \tfrac{1}{T} \Big\{ (q + 1 -\epsilon) \tfrac{v \cdot (v_1 - v)}{|v_1 - v|} |v_1 - v| + \tfrac{(1 - \epsilon) [v \cdot (v_1 - v)]^2}{|v_1 - v|^2} - q \u \cdot (v_1 - v) \Big\} \\
			= & - (1 - \tfrac{q^2}{1 - \epsilon} - 4 \eps_0) \tfrac{ |v_1 - v|^2}{4 T} + \tfrac{q |\u|^2}{4 \eps_0 T} \\
			& \qquad - \tfrac{1-\epsilon}{T} \left| \tfrac{v \cdot (v_1 - v)}{|v_1 - v|} + \tfrac{q+1-\epsilon}{2 (1 - \epsilon)} |v_1 - v| \right|^2 - \tfrac{\eps_0}{T} |v_1 - v - \tfrac{q}{2 \eps_0} \u|^2 \\
			\leq & - (1 - \tfrac{q^2}{1 - \epsilon} - 4 \eps_0) \tfrac{ |v_1 - v|^2}{4 T} + \tfrac{q |\u|^2}{4 \eps_0 T}
		\end{aligned}
	\end{equation}
	for some small $\eps_0 > 0$ to be determined. Then, for any fixed $0 < q < 1$, one will find the proper $\epsilon \in (0,1)$ and $\eps_0 > 0$ such that
	\begin{equation}
		\begin{aligned}
			1 - \tfrac{q^2}{1 - \epsilon} - 4 \eps_0 > 0 \,.
		\end{aligned}
	\end{equation}
	Actually, by taking $\epsilon = \frac{1}{2} (1 - q^2) > 0$ and $\eps_0 = \frac{1 - q^2}{8 (1 + q^2)} > 0$, one has
	\begin{equation}
		\begin{aligned}
			1 - \tfrac{q^2}{1 - \epsilon} - 4 \eps_0 = \tfrac{1 - q^2}{2 (1 + q^2)} > 0 \,.
		\end{aligned}
	\end{equation}
	Namely,
	\begin{equation}
		\begin{aligned}
			\tfrac{q ( |v - \u|^2 - |v_1 - \u|^2 )}{2 T} - \tfrac{|v_1 - v|^2}{4 T} - \tfrac{(1 - \epsilon) (|v_1|^2 - |v|^2)^2}{4 T} \leq - \tfrac{ (1 - q^2) |v_1 - v|^2}{8 (1 + q^2) T} + \tfrac{2 q (1 + q^2) |\u|^2}{(1	- q^2) T} \,.
		\end{aligned}
	\end{equation}
	Then, by \eqref{Assert-integrable}, if $- 2( 2 - b_0 - \gamma ) > - 3$,
	\begin{equation}\label{D2}
		\begin{aligned}
			D^2 \leq & \exp \big[ \tfrac{2 q (1 + q^2) |\u|^2}{(1	- q^2) T} \big] \int_{\R^3} \tfrac{1}{|v_1 - v|^{2( 2 - b_0 - \gamma )}} \exp \big[ - \tfrac{ (1 - q^2) |v_1 - v|^2}{8 (1 + q^2) T} \big] \d v_1 \\
			\leq & \exp \big[ \tfrac{2 q (1 + q^2) |\u|^2}{(1	- q^2) T} \big] \Big\{ \int_{|v| \leq 1}  \tfrac{1}{| v|^{2( 2 - b_0 - \gamma )}} \d v + \int_{\R^3} \exp \big[ - \tfrac{ (1 - q^2) | v|^2}{8 (1 + q^2) T} \big] \d v \Big\} \\
			= & \big[ \tfrac{4 \pi}{2 (b_0 + \gamma - 2) + 3} + ( \tfrac{8 \pi T (1 + q^2)}{1 - q^2} )^\frac{3}{2} \big] \exp \big[ \tfrac{2 q (1 + q^2) |\u|^2}{(1	- q^2) T} \big] < \infty \,.
		\end{aligned}
	\end{equation}
	It therefore infers that the bound \eqref{K2-L2Bnd-H} holds, provided that $- 2( 2 - b_0 - \gamma ) > - 3$ is further assumed.
	
	Recall that the constraints of $b_0$ in Lemma \ref{Lmm-k1k2-bnd} is that
	\begin{equation}
		\begin{aligned}
			\Phi_1 (\gamma) : = \{ b_0 | \, b_0 \in [0, 1 - \gamma] \textrm{ if } - 1 < \gamma \leq 1, \textrm{ and } b_0 \in (- \gamma - 1, 2) \textrm{ if }- 3 < \gamma \leq -1 \} \,.
		\end{aligned}
	\end{equation}
	The restrictions given in \eqref{B2} and \eqref{D2} can be further characterized by a set
	\begin{equation}
		\begin{aligned}
			\Phi_2 (\gamma) : = \{ b_0 | \, - 2( 2 - b_0 - \gamma ) > - 3 \,, 2 \gamma > - 3 \}.
		\end{aligned}
	\end{equation}
	One only requires to accurately find the set
	\begin{equation}
		\begin{aligned}
			\Xi = \{ \gamma | \, \Phi_1 (\gamma) \cap \Phi_2 (\gamma) \neq \emptyset, - 3 < \gamma \leq 1 \} \,,
		\end{aligned}
	\end{equation}
	because the conclusions of Lemma \ref{Lmm-K1K2-L2} hold for all $\gamma \in \Xi$. If $- 1 < \gamma \leq 1$, it is easy to know that $b_0 =1 - \gamma \in \Phi_1 (\gamma) \cap \Phi_2 (\gamma)$, hence, $(-1,1] \subseteq \Xi$. If $-3 < \gamma \leq - 1$, the set $\Phi_1 (\gamma)$ reads $- \gamma - 1 < b_0 < 2$. However, the set $\Phi_2 (\gamma)$ indicates that $\gamma > - \frac{3}{2}$ and $b_0 > - \gamma + \frac{1}{2}$. As shown in Figure \ref{Fig2}, 	
	\begin{figure}[h]
		\begin{center}
			\begin{tikzpicture}
				\pgfmathsetmacro\rd{0.1};
				\draw (\rd,0) arc (0: 360: \rd);
				\draw (3,0) arc (0:360:\rd);
				\draw (-3,0) arc (0:360:\rd);
				\draw (-5, 0)--(-3 - 2*\rd, 0);
				\draw (-3,0)--(-\rd,0);
				\draw (\rd,0)--(3-2*\rd, 0);
				\draw[->] (3,0)--(5,0) node[right]{$b_0$};
				\draw (-3 - \rd, \rd)node[below=2mm]{$-\gamma-1$}--(-3 - \rd, 1)--(3 - \rd, 1)--(3 - \rd, \rd)node[below=2mm]{2};
				\draw (0,\rd)node[below=1mm]{$-\gamma + \frac{1}{2}$}--(0,1.5)--(4.5, 1.5);
				\filldraw[yellow, draw=black] (\rd, 0) arc (0:90:\rd)--(0,1)--(3 - \rd, 1)--(3 - \rd, \rd) arc (90:180:\rd)--(\rd, 0)--cycle;
				\draw (1.5,1)node[below=2mm]{$\Phi_1 (\gamma) \cap \Phi_2 (\gamma)$};
			\end{tikzpicture}
		\end{center}
		\caption{The set $\Phi_1 (\gamma) \cap \Phi_2 (\gamma)$ with $- 3 < \gamma \leq - 1$.}\label{Fig2}
	\end{figure}
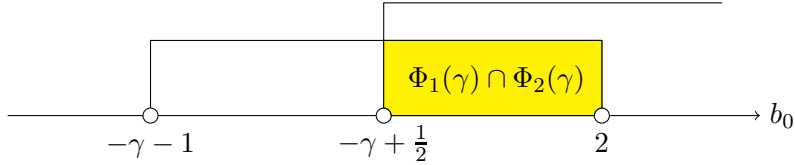
	$\Phi_1 (\gamma) \cap \Phi_2 (\gamma) \neq \emptyset$ is equivalent to $- \gamma + \frac{1}{2} < 2$ and $-3 < \gamma \leq - 1$, namely, $- \frac{3}{2} < \gamma \leq - 1$. It thereby means that $(- \frac{3}{2}, -1] \subseteq \Xi$ and $(-3, - \frac{3}{2}] \subseteq \R - \Xi$, the complement of $\Xi$ in $\R$. Consequently, one concludes that $\Xi = (- \frac{3}{2}, 1 ]$, and the proof of Lemma \ref{Lmm-K1K2-L2} is completed.
\end{proof}

\subsection{Hypocoercivity of the weighted linearized Boltzmann operator $\L$}

In this subsection, the hypocoercivity of the weighted linearized Boltzmann operator will be derived, which is the modified version of the hypocoercivity of the linearized operator $\L$ in \eqref{Hypo-L}. These results are applied to deal with the weighted $L^2$ norms of $f$.

First, for any fixed $r > 0$, the smooth monotone cutoff function $\chi (s) \in C^\infty (0,+ \infty)$ is introduced by
\begin{equation}\label{Chi}
	\begin{aligned}
		\chi (s) \equiv 0 \textrm{ for } 0 < s \leq r \,, \ \chi (s) \equiv 1 \textrm{ for } s \geq 2 r \,, \ 0 \leq \chi (s) \leq 1 \textrm{ for all } s > 0 \,.
	\end{aligned}
\end{equation}
Split the operator $K f$ as follows:
\begin{equation}\label{K-Lambda}
	\begin{aligned}
		K f = K^{1 - \chi} f + K^\chi f \,, \ K^\Lambda f = K_1^\Lambda f - K_2^\Lambda f \textrm{ for } \Lambda = 1-\chi \textrm{ and } \chi \,,
	\end{aligned}
\end{equation}
where
\begin{equation}\label{K1-Lambda}
	\begin{aligned}
		K_1^\Lambda f (v) = \M^\frac{1}{2} (v) \iint_{\mathbb{R}^3 \times \mathbb{S}^2} \Lambda (|v_1 - v|) \M^\frac{1}{2} (v_1) f (v_1) b (\omega, v_1 - v) \d \omega \d v_1 \,,
	\end{aligned}
\end{equation}
and
\begin{equation}\label{K2-Lambda}
	\begin{aligned}
		K_2^\Lambda f (v) = 2 \M^\frac{1}{2} (v) \iint_{\mathbb{R}^3 \times \mathbb{S}^2} \Lambda (|v_1 - v|) \M^{- \frac{1}{2}} (v') \M (v_1) f (v') b (\omega, v_1 -v) \d \omega \d v_1 \,.
	\end{aligned}
\end{equation}

\begin{lemma}[Low velocities estimates]\label{Lmm-K1K2-Hypo1}
	Let $0 < q < 1$, $- \frac{3}{2} < \gamma \leq 1$ and any fixed $r > 0$. Then there holds
		\begin{equation}\label{K-1-chi-bnd}
			\begin{aligned}
				\big| \langle \Theta_\gamma \M^{- \frac{q}{2}} K^{1 - \chi} f , \M^{- \frac{q}{2}} f \rangle \big| \lesssim_{\rho, \u, T, r} \| f \|^2_{L^2 (\nu)} \,.
			\end{aligned}
		\end{equation}
	where $\Theta_\gamma = \Theta_\gamma (v)$ is defined in \eqref{Theta-gamma}, i.e.,
	\begin{equation}
		\begin{aligned}
			\Theta_\gamma (v) = 1 \textrm{ if } 0 \leq \gamma \leq 1 \,, \textrm{ and } \Theta_\gamma (v) = \nu^{-1} (v) \textrm{ if } - \tfrac{3}{2} < \gamma < 0 \,.
		\end{aligned}
	\end{equation}
\end{lemma}

Remark that the lemma still holds for the cases $- 3 < \gamma \leq - \frac{3}{2}$. But it is not useful to prove our main results for $- 3 < \gamma \leq - \frac{3}{2}$. The details of proof of Lemma \ref{Lmm-K1K2-Hypo1} will be given in Section \ref{Sec:K-LHVE-K} later.

\begin{lemma}[High velocities estimates]\label{Lmm-K1K2-Hypo2}
	Let $0 < q < 1$, $- \frac{3}{2} < \gamma \leq 1$ and any fixed $r > 0$.
	\begin{enumerate}
		\item There holds
		\begin{equation}\label{K1-chi}
			\begin{aligned}
				\big| \langle \Theta_\gamma \M^{- \frac{q}{2}} K_1^\chi f, \M^{- \frac{q}{2}} f \rangle \big| \lesssim_{\rho, \u, T} \exp \big[ - \tfrac{(1-q) r^2}{32 T} \big] \| \M^{- \frac{q}{2}} f \|^2_{L^2} \,,
			\end{aligned}
		\end{equation}
		where $\Theta_\gamma = \Theta_\gamma (v)$ is given in \eqref{Theta-gamma}.
		
		\item If $\gamma = 1$, then
		\begin{equation}\label{gamma-HS}
			\begin{aligned}
				\big| \langle \M^{- \frac{q}{2}} K_2^\chi f, \M^{- \frac{q}{2}} f \rangle \big| \lesssim_{\rho, \u, T} \frac{1}{1 + r} \| \M^{- \frac{q}{2}} f \|^2_{L^2(\nu)} + \exp (\tfrac{q r^2}{2 T}) \| f \|^2_{L^2 (\nu)} \,.
			\end{aligned}
		\end{equation}
	    If $0 \leq \gamma < 1$, then
	    \begin{equation}\label{gamma-HP}
	    	\begin{aligned}
	    		\big| \langle \M^{- \frac{q}{2}} K_2^\chi f, \M^{- \frac{q}{2}} f \rangle \big| \lesssim_{\rho, \u, T} \frac{1}{r^{1 - \gamma}} \| \M^{- \frac{q}{2}} f \|^2_{L^2} \,.
	    	\end{aligned}
	    \end{equation}
        If $- \frac{3}{2} < \gamma < 0$, then
        \begin{equation}\label{gamma-SP}
        	\begin{aligned}
        		| \langle \nu^{-1} \M^{- \frac{q}{2}} K_2^\chi f , \M^{- \frac{q}{2}} f \rangle | \lesssim_{\rho, \u, T} \frac{1}{1 + r} \| \M^{- \frac{q}{2}} f \|^2_{L^2} + (1 + r)^{- \gamma} \exp (\tfrac{2 q r^2}{T}) \| f \|^2_{L^2 (\nu)} \,.
        	\end{aligned}
        \end{equation}
	\end{enumerate}
\end{lemma}

Remark that the lemma still holds for the cases $- 3 < \gamma \leq - \frac{3}{2}$. But it is not useful to prove our main results for $- 3 < \gamma \leq - \frac{3}{2}$. The proof of Lemma \ref{Lmm-K1K2-Hypo2} will be given in Section \ref{Sec:K-LHVE-K} later.

Based on Lemma \ref{Lmm-K1K2-Hypo1} and Lemma \ref{Lmm-K1K2-Hypo2}, the following hypocoercivity of the weighted linearized Boltzmann collision operator can be derived.

\begin{lemma}\label{Lmm-wL-Hypo}
	Let $0 < q < 1$ and $- \frac{3}{2} < \gamma \leq 1$. Then there is a constant $C = C(\rho, \u, T) > 0$ such that
	\begin{enumerate}
		\item if $0 \leq \gamma \leq 1$,
	\begin{equation}\label{Hypo-HP}
		\begin{aligned}
			\langle \M^{- \frac{q}{2}} \L f , \M^{- \frac{q}{2}} f \rangle \geq \tfrac{1}{2} \| \M^{- \frac{q}{2}} f \|^2_{L^2 (\nu)} - C \| f \|^2_{L^2 (\nu)} \,;
		\end{aligned}
	\end{equation}

    \item if $- \frac{3}{2} < \gamma < 0$,
    \begin{equation}\label{Hypo-SP}
    	\begin{aligned}
    		\langle \nu^{-1} \M^{- \frac{q}{2}} \L f , \M^{- \frac{q}{2}} f \rangle \geq \tfrac{1}{2} \| \M^{- \frac{q}{2}} f \|^2_{L^2} - C \| f \|^2_{L^2 (\nu)} \,.
    	\end{aligned}
    \end{equation}
    \end{enumerate}
\end{lemma}

\begin{proof}
	
	(1) The hard potential cases for $0 \leq \gamma \leq 1$.
	
	Since $\L f = \nu f + K^{1 - \chi} f + K_1^\chi f - K_2^\chi f$ by \eqref{Lf} and \eqref{K-Lambda}, one derives from Lemma \ref{Lmm-K1K2-Hypo1}-\ref{Lmm-K1K2-Hypo2} and $\nu (v) \geq \nu_0 > 0$ for $0 \leq \gamma \leq 1$ that
	\begin{equation}\label{Sq-1}
		\begin{aligned}
			\langle \M^{- \frac{q}{2}} \L f , \M^{- \frac{q}{2}} f \rangle \geq & \| \M^{- \frac{q}{2}} f \|^2_{L^2 (\nu)} - \big| \langle \M^{- \frac{q}{2}} K^{1 - \chi} f , \M^{- \frac{q}{2}} f \rangle \big| \\
			& - \big| \langle \M^{- \frac{q}{2}} K_1^\chi f , \M^{- \frac{q}{2}} f \rangle \big| - \big| \langle \M^{- \frac{q}{2}} K_2^\chi f , \M^{- \frac{q}{2}} f \rangle \big| \\
			\geq & \Big\{ 1 - C(\rho, \u, T) \big[ \exp ( - \tfrac{(1 - q) r^2}{32 T} ) + p_\gamma (r) \big] \Big\} \| \M^{- \frac{q}{2}} f \|^2_{L^2 (\nu)} \\
			& - C (\rho, \u, T, r) \| f \|^2_{L^2 (\nu)} \,,
		\end{aligned}
	\end{equation}
	where
	\begin{equation*}
		\begin{aligned}
			p_\gamma (r) = \left\{
			  \begin{aligned}
			  	\tfrac{1}{r^{1 - \gamma}} \,, & \ \textrm{ if } 0 \leq \gamma < 1 \,, \\
			  	\tfrac{1}{r+1} \,, & \ \textrm{ if } \gamma = 1 \,,
			  \end{aligned}
			\right.
		\end{aligned}
	\end{equation*}
	satisfying $p_\gamma (r) \to 0$ as $r \to + \infty$. Then, one can take a sufficiently large $r = r (\rho, \u, T) > 0$ such that
	\begin{equation*}
		\begin{aligned}
			1 - C(\rho, \u, T) \big[ \exp ( - \tfrac{(1 - q) r^2}{32 T} ) + p_\gamma (r) \big] \geq \tfrac{1}{2} \,.
		\end{aligned}
	\end{equation*}
	It thereby concludes the inequality \eqref{Hypo-HP} in Lemma \ref{Lmm-wL-Hypo} from \eqref{Sq-1}.
	
	\vspace*{3mm}
	
	(2) The soft potential cases for $- \frac{3}{2} < \gamma < 0$.
	
	Note that $\nu^{-1} \L f = f + \nu^{-1} K^{1-\chi} f + \nu^{-1} K_1^\chi f - \nu^{-1} K_2^\chi f$. Then the inequalities \eqref{K-1-chi-bnd} in Lemma \ref{Lmm-K1K2-Hypo1}, \eqref{K1-chi} and \eqref{gamma-SP} in Lemma \ref{Lmm-K1K2-Hypo2} tell us that
	\begin{equation}
		\begin{aligned}
			\langle \nu^{-1} \M^{- \frac{q}{2}} \L f , \M^{- \frac{q}{2}} f \rangle \geq & \| \M^{- \frac{q}{2}} f \|^2_{L^2} - \big| \langle \nu^{-1} \M^{- \frac{q}{2}} K^{1 - \chi} f , \M^{- \frac{q}{2}} f \rangle \big| \\
			& - \big| \langle \nu^{-1} \M^{- \frac{q}{2}} K_1^\chi f , \M^{- \frac{q}{2}} f \rangle \big| - \big| \langle \nu^{-1} \M^{- \frac{q}{2}} K_2^\chi f , \M^{- \frac{q}{2}} f \rangle \big| \\
			\geq & \Big\{ 1 - C(\rho, \u, T) \big[ \exp ( - \tfrac{(1 - q) r^2}{32 T} ) + \tfrac{1}{1 + r} \big] \Big\} \| \M^{- \frac{q}{2}} f \|^2_{L^2} \\
			& - C (\rho, \u, T, r) \| f \|^2_{L^2 (\nu)} \,.
		\end{aligned}
	\end{equation}
    Because $\exp ( - \tfrac{(1 - q) r^2}{32 T} ) + \tfrac{1}{1 + r} \to 0$ as $r \to + \infty$, there is a large enough $r = r (\rho, \u, T) > 0$ such that
    \begin{equation*}
    	\begin{aligned}
    		1 - C(\rho, \u, T) \big[ \exp ( - \tfrac{(1 - q) r^2}{32 T} ) + \tfrac{1}{1 + r} \big] \geq \tfrac{1}{2} \,.
    	\end{aligned}
    \end{equation*}
    Consequently, the inequality \eqref{Hypo-SP} holds, and the proof of Lemma \ref{Lmm-wL-Hypo} is finished.
\end{proof}

Next we derive the weighted $L^2$ estimates for the operator $K f$ with $- 3 < \gamma \leq - \frac{3}{2}$.

\begin{lemma}\label{Lmm-VSP}
	Let $\alpha \in \mathbb{N}^3$ with $|\alpha| \geq 0$, $- 3 < \gamma \leq - \frac{3}{2}$, $l \in \R$ and $0 < q < 1$. Then for any $\eta > 0$ there is a small $C_\eta > 0$ such that
	\begin{equation}\label{L2-fg1}
		\begin{aligned}
			\langle \langle v \rangle^{2 \gamma l} \M^{- q} \partial^\alpha_v (\nu f), \partial^\alpha_v f \rangle \geq & \| \langle v \rangle^{\gamma (l + \frac{1}{2} )} \M^{- \frac{q}{2}} \partial^\alpha_v f \|^2_{L^2} \\
			& - \eta \sum_{\alpha' \leq \alpha} \| \langle v \rangle^{\gamma (l + \frac{1}{2} )} \M^{- \frac{q}{2}} \partial^{\alpha'}_v f \|^2_{L^2} - C_\eta \| f \|^2_{L^2 (\nu)} \,,
		\end{aligned}
	\end{equation}
	and
	\begin{equation}\label{L2-fg2}
		\begin{aligned}
			& \big| \langle \langle v \rangle^{2 \gamma l} \M^{- q} \partial^\alpha_v (K f), g \rangle \big| \\
			\leq & \Big\{ \eta \sum_{|\alpha'| \leq |\alpha|} \| \langle v \rangle^{\gamma (l + \frac{1}{2} )} \M^{- \frac{q}{2}} \partial^{\alpha'}_v f \|_{L^2} + C_\eta \| f \|_{L^2 (\nu)} \Big\} \| \langle v \rangle^{\gamma (l + \frac{1}{2} )} \M^{- \frac{q}{2}} g \|_{L^2} \,.
		\end{aligned}
	\end{equation}
\end{lemma}

The proof of Lemma \ref{Lmm-VSP} will be given in Section \ref{Sec:VSP} later.

\subsection{Decay of the pseudo-inverse $\L^{-1}$: Proof of Theorem \ref{Thm1}}\label{Subsec-Thm1}

In this subsection, one aims at justifying the decay estimates of the pseudo-inverse operator $\L^{-1}$ by using the results obtained in the previous subsection, i.e., proving Theorem \ref{Thm1}.

\subsubsection{\texttt{Proof of Part (I) of Theorem \ref{Thm1}: The cases $- \tfrac{3}{2} < \gamma \leq 1$}}
	We will separately proof the results in Part (I) of Theorem \ref{Thm1} for the hard potential cases $0 \leq \gamma \leq 1$ and the part of soft potential cases $- \frac{3}{2} < \gamma < 0$.
	
	\vspace*{3mm}
	
	{\em (1) The hard potential cases $0 \leq \gamma \leq 1$.}
	
	\vspace*{3mm}
	
	Let $f = \L^{-1} g \in \mathrm{Null}^\perp (\L)$, i.e., $\L f = g$. Recalling the definition of $\L$ in \eqref{Lf}, namely, $\L f = \nu f + K f$, one has
	\begin{equation*}
		\begin{aligned}
			f = \nu^{-1} g - \nu^{- 1} K f \,.
		\end{aligned}
	\end{equation*}
	Because $\nu^{-1} \leq \nu_0^{-1}$ for $0 \leq \gamma \leq 1$, there hold
	\begin{equation}\label{Sq-2}
		\begin{aligned}
			\M^{- \frac{q}{2}} (v) |f (v)| & \leq \nu^{-1} (v) \M^{- \frac{q}{2}} (v) |g (v)| + \nu^{-1} (v) \M^{- \frac{q}{2}} (v) |K f (v)| \\
			& \lesssim_{\rho, \u, T} \M^{- \frac{q}{2}} (v) |g (v)| + \M^{- \frac{q}{2}} (v) |K f (v)| \\
			& \lesssim_{\rho, \u, T} \M^{- \frac{q}{2}} (v) |g (v)| + \| \M^{- \frac{q}{2}} f \|^2_{L^2} + \| f \|_{L^2} \\
			& \lesssim_{\rho, \u, T} \M^{- \frac{q}{2}} (v) |g (v)| + \| \M^{- \frac{q}{2}} f \|^2_{L^2 (\nu)} \,,
		\end{aligned}
	\end{equation}
	where the last second inequality is derived from Lemma \ref{Lmm-K1K2-L2}, and the last one is implied by the facts $\nu (v) \geq \nu_0 > 0$ for $0 \leq \gamma \leq 1$ and $\M^{- \frac{q}{2}} (v) \geq C(\rho, T) > 0$ for $0 < q < 1$. Moreover, Lemma \ref{Lmm-wL-Hypo} and the hypocoercivity of $\L$ in \eqref{Hypo-L} yield that for $0 \leq \gamma \leq 1$ and $0 < q < 1$,
	\begin{equation}\label{Sq-3}
		\begin{aligned}
			\| \M^{- \frac{q}{2}} f \|^2_{L^2 (\nu)} & \lesssim_{\rho, \u, T} \langle \M^{- \frac{q}{2}} \L f, \M^{- \frac{q}{2}} f \rangle + \| f \|^2_{L^2 (\nu)} \\
			& \lesssim_{\rho, \u, T} \langle \M^{- \frac{q}{2}} \L f, \M^{- \frac{q}{2}} f \rangle + \langle \L f, f \rangle \\
			& \lesssim_{\rho, \u, T} \int_{\R^3} \M^{- \frac{q}{2}} (v) |g (v)| \cdot |\M^{- \frac{q}{2}} (v) f (v)| \d v \\
			& \lesssim_{\rho, \u, T} \| \M^{- \frac{q}{2}} f \|_{L^2 (\nu)} \| \langle v \rangle^k \M^{- \frac{q}{2}} g \|_{L^\infty} \Big( \int_{\R^3} \tfrac{1}{\langle v \rangle^{2 k + \gamma}} \d v \Big)^\frac{1}{2}
		\end{aligned}
	\end{equation}
	for any fixed $k > \frac{3 - \gamma}{2}$, where $\nu (v) \thicksim_{\rho, \u, T} \langle v \rangle^\gamma$ in \eqref{nu-equivalent} has been used. Noticing that the integral $\int_{\R^3} \tfrac{1}{\langle v \rangle^{2 k + \gamma}} \d v$ converges under $k > \frac{3 - \gamma}{2}$, there holds
	\begin{equation}
		\begin{aligned}
			\| \M^{- \frac{q}{2}} f \|_{L^2 (\nu)} \lesssim_{\rho, \u, T} \| \langle v \rangle^k \M^{- \frac{q}{2}} g \|_{L^\infty} \,.
		\end{aligned}
	\end{equation}
	Consequently, the inequalities \eqref{Sq-2} and \eqref{Sq-3} conclude the results in Part (I)-(i) of Theorem \ref{Thm1}.
	
	\vspace*{3mm}
	
	{\em (2) The part of soft potential cases $- \frac{3}{2} < \gamma < 0$.}
	
	\vspace*{2mm}
	
	Noticing that $\nu f = \L f - K f = g - K f$, where $f = \L^{-1} g \in \mathrm{Null}^\perp (\L)$, one derives from the similar arguments in \eqref{Sq-2} that
	\begin{equation}\label{Sq-4}
		\begin{aligned}
			\nu (v) \M^{- \frac{q}{2}} (v) |f (v)| \lesssim_{\rho, \u, T} \M^{- \frac{q}{2}} (v) |g (v)| + \| \M^{- \frac{q}{2}} f \|_{L^2} \,.
		\end{aligned}
	\end{equation}
	Moreover, the inequality \eqref{Hypo-SP} in Lemma \ref{Lmm-wL-Hypo} and the hypocoercivity of $\L$ in \eqref{Hypo-L} imply that
	\begin{equation*}
		\begin{aligned}
			\| \M^{- \frac{q}{2}} f \|^2_{L^2} \lesssim_{\rho, \u, T} & \langle \nu^{-1} \M^{- \frac{q}{2}} \L f , \M^{- \frac{q}{2}} f \rangle + \| f \|^2_{L^2 (\nu)} \\
			\lesssim_{\rho, \u, T} & \langle \nu^{-1} \M^{- \frac{q}{2}} \L f , \M^{- \frac{q}{2}} f \rangle + \langle \L f , f \rangle \\
			\lesssim_{\rho, \u, T} & \int_{\R^3} \nu^{-1} (v) \M^{- \frac{q}{2}} (v) |g (v)| \, |\M^{- \frac{q}{2}} (v) f (v)| \d v \\
			\lesssim_{\rho, \u, T} & \| \langle v \rangle^k \M^{- \frac{q}{2}} g \|_{L^\infty} \| \M^{- \frac{q}{2}} f \|_{L^2} \Big( \int_{\R^3} \tfrac{\d v}{\langle v \rangle^{2 k + 2 \gamma}} \Big)^\frac{1}{2} \,,
		\end{aligned}
	\end{equation*}
	which means that
	\begin{equation}\label{Sq-5}
		\begin{aligned}
			\| \M^{- \frac{q}{2}} f \|_{L^2} \lesssim_{\rho, \u, T} \| \langle v \rangle^k \M^{- \frac{q}{2}} g \|_{L^\infty} \Big( \int_{\R^3} \tfrac{\d v}{ \langle v \rangle^{2 k + 2 \gamma}} \Big)^\frac{1}{2} \,.
		\end{aligned}
	\end{equation}
	Since the integral $\int_{\R^3} \tfrac{\d v}{\langle v \rangle^{2 k + 2 \gamma}}$ converges provided that $k > \tfrac{3 - 2 \gamma}{2}$. As a consequence, the bounds \eqref{Sq-4} and \eqref{Sq-5} imply the results in Part (I)-(ii) of Theorem \ref{Thm1}.

\subsubsection{\texttt{Proof of Part (II) of Theorem \ref{Thm1}: The cases $- 3 < \gamma \leq - \tfrac{3}{2}$}}

Before proving our main results, we first introduce the following lemma.

\begin{lemma}\label{Lmm-L-L2}
	Let $- 3 < \gamma \leq - \frac{3}{2}$, $N \geq 0$, $l \geq 0$, $0 < q < 1$. Denote by $f = \L^{-1} g$, where $g$ is given in Theorem \ref{Thm1}. Then there is a positive constant $C > 0$ such that
	\begin{equation}\label{L2-derivative}
		\begin{aligned}
			\sum_{|\alpha| \leq N} \| \langle v \rangle^{\gamma (|\alpha| - l + \frac{1}{2})} \M^{ - \frac{q}{2} } \partial^\alpha_v f \|^2_{L^2} \leq C \sum_{|\alpha| \leq N} \| \langle v \rangle^{\gamma (|\alpha| - l - \frac{1}{2})} \M^{- \frac{q}{2} } \partial^\alpha_v g \|^2_{L^2} \,.
		\end{aligned}
	\end{equation}
\end{lemma}

\begin{proof}
	Observe that $\nu (v) \thicksim \langle v \rangle^\gamma$ for $- 3 < \gamma \leq - \tfrac{3}{2}$. Then by the similar arguments of Lemma 2 of \cite{Strain-Guo-2008-ARMA}, the following estimates hold: for any $\alpha \in \mathbb{N}^3$ with $|\alpha| \geq 0$, $l \in \R$, $\eta \in (0,1)$ and $0 < q < 1$,
	\begin{equation}\label{L2-1}
		\begin{aligned}
			\langle \langle v \rangle^{2 \gamma l} \M^{- q} \partial^\alpha_v (\nu f), \partial^\alpha_v f \rangle \geq & \| \langle v \rangle^{\gamma (l + \frac{1}{2} )} \M^{- \frac{q}{2} } \partial^\alpha_v f \|^2_{L^2} \\
			& - \eta \sum_{\alpha_1 \leq \alpha} \| \langle v \rangle^{\gamma (l + \frac{1}{2} )} \M^{- \frac{q}{2} } \partial^{\alpha_1}_v f \|^2_{L^2} - C_\eta \| f \|^2_{L^2 (\nu)} \,,
		\end{aligned}
	\end{equation}
	and
	\begin{equation}\label{L2-2}
		\begin{aligned}
			& \big| \langle \langle v \rangle^{2 \gamma l} \M^{- q} \partial^\alpha_v (K f), \partial^\alpha_v f \big| \\
			\leq & \Big\{ \eta \sum_{\alpha' \leq \alpha} \| \langle v \rangle^{\gamma (l + \frac{1}{2} )} \M^{- \frac{q}{2} } \partial^{\alpha_1}_v f \|_{L^2} + C_\eta \| f \|_{L^2 ( \nu ) } \Big\} \| \langle v \rangle^{\gamma (l + \frac{1}{2} )} \M^{- \frac{q}{2} } \partial^\beta_v f \|_{L^2} \,.
		\end{aligned}
	\end{equation}
	Recall that $\L f = \nu f + K f$. By replacing $l$ in \eqref{L2-1}-\eqref{L2-2} by $|\alpha| - l$, it infers
	\begin{equation}\label{L2-3}
		\begin{aligned}
			& \langle \langle v \rangle^{2 \gamma ( |\alpha| - l )} \M^{- q } \partial^\alpha_v (\L f), \partial^\alpha_v f \rangle \\
			\geq & (1 - 2 \eta) \| \langle v \rangle^{\gamma ( |\alpha| - l + \frac{1}{2} )} \M^{- \frac{q}{2} } \partial^\alpha_v f \|^2_{L^2} - 2 \eta \sum_{\alpha_1 < \alpha} \| \langle v \rangle^{\gamma ( |\alpha| - l + \frac{1}{2} )} \M^{- \frac{q}{2} } \partial^{\alpha_1}_v f \|^2_{L^2} - C_\eta \| f \|^2_{L^2 (\nu)}
		\end{aligned}
	\end{equation}
	for any small $\eta \in (0, \frac{1}{2})$. Recall that $\L f = g$. The inequality \eqref{L2-3} indicates
	\begin{equation}\label{L2-4}
		\begin{aligned}
			& \tfrac{1}{2} (1 - 2 \eta) \| \langle v \rangle^{\gamma ( |\alpha| - l + \frac{1}{2} )} \M^{- \frac{q}{2} } \partial^\alpha_v f \|^2_{L^2} \\
			\leq & 2 \eta \sum_{\alpha_1 < \alpha} \| \langle v \rangle^{\gamma ( |\alpha| - l + \frac{1}{2} )} \M^{- \frac{q}{2} } \partial^{\alpha_1}_v f \|^2_{L^2} + C_\eta \| f \|^2_{L^2 (\nu)} + C_\eta \| \langle v \rangle^{\gamma ( |\alpha| - l - \frac{1}{2} )} \M^{- \frac{q}{2} } \partial^\alpha_v g \|^2_{L^2} \,.
		\end{aligned}
	\end{equation}
	From summing up \eqref{L2-4} over $|\alpha| \leq N$ and taking $\eta > 0$ sufficiently small, it follows
	\begin{equation}
		\begin{aligned}
			\sum_{|\alpha| \leq N} \| \langle v \rangle^{\gamma ( |\alpha| - l + \frac{1}{2} )} \M^{- \frac{q}{2} } \partial^\alpha_v f \|^2_{L^2} \leq C \sum_{|\alpha| \leq N} \| \langle v \rangle^{\gamma ( |\alpha| - l - \frac{1}{2} )} \M^{- \frac{q}{2} } \partial^\alpha_v g \|^2_{L^2} + C \| f \|^2_{L^2 (\nu)} \,.
		\end{aligned}
	\end{equation}
	
	Recalling that \eqref{Hypo-L}, i.e., $\langle \L f , f \rangle \geq \lambda \| f \|_{L^2 (\nu)}^2$, one has
	\begin{equation*}
		\begin{aligned}
			\| f \|_{L^2 (\nu)}^2 \leq \lambda \langle \L f , f \rangle = \lambda \langle g , f \rangle \leq C \| f \|_{L^2 (\nu)} \| \langle v \rangle^{- \frac{\gamma}{2} } g \|_{L^2} \,,
		\end{aligned}
	\end{equation*}
	which means that
	\begin{equation}\label{L2-5}
		\begin{aligned}
			\| f \|_{L^2 (\nu)}^2 \leq C \| \langle v \rangle^{- \frac{\gamma}{2} } g \|^2_{L^2} \,.
		\end{aligned}
	\end{equation}
	Observe that $\langle v \rangle^{- \frac{\gamma}{2} } \leq C \langle v \rangle^{\gamma (- l - \frac{1}{2})} \M^{- \frac{q}{2} } (v)$ for any $l \geq 0$ and $0 < q < 1$. It thereby follows
	\begin{equation}\label{L2-6}
		\begin{aligned}
			\| \langle v \rangle^{- \frac{\gamma}{2} } g \|_{L^2} \leq C \| \langle v \rangle^{\gamma (- l - \frac{1}{2})} \M^{- \frac{q}{2} } g \|_{L^2} \,.
		\end{aligned}
	\end{equation}
	Then the inequalities \eqref{L2-4}-\eqref{L2-5}-\eqref{L2-6} conclude the estimate \eqref{L2-derivative}. The proof of Lemma \ref{Lmm-L-L2} is therefore finished.
\end{proof}

Now we focus on the proof of Part (II) of Theorem \ref{Thm1}. One claims that
\begin{equation}\label{Claim-H2}
	\begin{aligned}
		\| \M^{- \frac{q}{2} } f \|^2_{H^2} \leq C \sum_{|\alpha| \leq 2} \| \langle v \rangle^{\gamma (|\alpha| - l + \frac{1}{2})} \M^{- \frac{q}{2} } \partial^\alpha_v f \|^2_{L^2}
	\end{aligned}
\end{equation}
for $- 3 < \gamma \leq - \frac{3}{2}$ and $l \geq \frac{5}{2}$. Indeed,
\begin{equation*}
	\begin{aligned}
		\| \M^{- \frac{q}{2} } f \|^2_{H^2} = \int_{\R^3} |\M^{- \frac{q}{2} } f|^2 \d v + \sum_{|\alpha| = 1} \int_{\R^3} |\partial^\alpha_v ( \M^{- \frac{q}{2} } f )|^2 \d v + \sum_{|\alpha| = 2} \int_{\R^3} |\partial^\alpha_v ( \M^{- \frac{q}{2} } f )|^2 \d v \,.
	\end{aligned}
\end{equation*}
It is easy to see
\begin{equation*}
	\begin{aligned}
		| \partial^\alpha_v \M^{- \frac{q}{2} } (v) | \leq C \langle v \rangle^{|\alpha|} \M^{- \frac{q}{2} } (v)
	\end{aligned}
\end{equation*}
for any $|\alpha| \geq 0$.

For $|\alpha| = 2$, it follows
\begin{equation*}
	\begin{aligned}
		|\partial^\alpha_v ( \M^{- \frac{q}{2} } f )| \leq & | \partial^\alpha_v \M^{- \frac{q}{2} } | | f| + \sum_{\alpha_1 \leq |\alpha|, |\alpha_1| = 1} C_\alpha^{\alpha_1} |\partial^{\alpha - \alpha_1}_v \M^{- \frac{q}{2} } | |\partial^{\alpha_1}_v f| + \M^{- \frac{q}{2} } |\partial^\alpha_v f| \\
		\leq & C \langle v \rangle^2 | \M^{- \frac{q}{2} } f | + C \sum_{|\alpha_1| = 1} \langle v \rangle |\M^{- \frac{q}{2} } \partial^{\alpha_1}_v f| + C | \M^{- \frac{q}{2} } \partial^\alpha_v f| \,.
	\end{aligned}
\end{equation*}
For $|\alpha| = 1$, it infers
\begin{equation*}
	\begin{aligned}
		|\partial^\alpha_v ( \M^{- \frac{q}{2} } f )| \leq | \partial^\alpha_v \M^{- \frac{q}{2} } f | + | \M^{- \frac{q}{2} } \partial^\alpha_v f | \leq C \langle v \rangle | \M^{- \frac{q}{2} } f | + | \M^{- \frac{q}{2} } \partial^\alpha_v f | \,.
	\end{aligned}
\end{equation*}
Then one has
\begin{equation*}
	\begin{aligned}
		\| \M^{- \frac{q}{2} } f \|^2_{H^2} \leq C \| \langle v \rangle^2 \M^{- \frac{q}{2} } f \|^2_{L^2} + C \sum_{|\alpha| = 1} \| \langle v \rangle \M^{- \frac{q}{2} } \partial^\alpha_v f \|^2_{L^2} + C \sum_{|\alpha| = 2} \| \M^{- \frac{q}{2} } \partial^\alpha_v f \|^2_{L^2} \,.
	\end{aligned}
\end{equation*}
Notice that
\begin{equation*}
	\begin{aligned}
		\| \langle v \rangle^2 \M^{- \frac{q}{2} } f \|^2_{L^2} \leq \| \langle v \rangle^{\gamma (0 -l + \frac{1}{2})} \M^{- \frac{q}{2} } f \|^2_{L^2}
	\end{aligned}
\end{equation*}
for $l \geq \frac{1}{2} - \frac{2}{\gamma}$ with $- 3 < \gamma \leq - \frac{3}{2}$. Moreover, for $|\alpha| = 1$ one has
\begin{equation*}
	\begin{aligned}
		\| \langle v \rangle \M^{- \frac{q}{2} } \partial^\alpha_v f \|^2_{L^2} \leq \| \langle v \rangle^{\gamma (1 - l + \frac{1}{2})} \M^{- \frac{q}{2} } \partial^\alpha_v f \|^2_{L^2}
	\end{aligned}
\end{equation*}
for $l \geq \frac{3}{2} - \frac{1}{\gamma}$ with $- 3 < \gamma \leq - \frac{3}{2}$. It further holds that for $|\alpha| = 2$
\begin{equation*}
	\begin{aligned}
		\| \M^{- \frac{q}{2} } \partial^\alpha_v f \|^2_{L^2} \leq \| \langle v \rangle^{\gamma (2 - l + \frac{1}{2})} \M^{- \frac{q}{2} } \partial^\alpha_v f \|^2_{L^2}
	\end{aligned}
\end{equation*}
for $l \geq \frac{5}{2}$. In summary, there holds
\begin{equation*}
	\begin{aligned}
		\| \M^{- \frac{q}{2} } f \|^2_{H^2} \leq C \sum_{|\alpha| \leq 2} \| \langle v \rangle^{\gamma (|\alpha| - l + \frac{1}{2})} \M^{- \frac{q}{2} } \partial^\alpha_v f \|^2_{L^2}
	\end{aligned}
\end{equation*}
for $l \geq \max \{ \frac{5}{2}, \frac{3}{2} - \frac{1}{\gamma}, \frac{1}{2} - \frac{2}{\gamma} \}$ with $- 3 < \gamma \leq - \frac{3}{2}$. Observe that
\begin{equation*}
	\begin{aligned}
		\max \{ \tfrac{5}{2}, \tfrac{3}{2} - \tfrac{1}{\gamma}, \tfrac{1}{2} - \tfrac{2}{\gamma} \} = \tfrac{5}{2}
	\end{aligned}
\end{equation*}
for $- 3 < \gamma \leq - \frac{3}{2}$. Thus the claim \eqref{Claim-H2} holds.

Now we take $N = 2$ in \eqref{L2-derivative} of Lemma \ref{Lmm-L-L2}, and choose $l = \frac{5}{2}$ in both \eqref{L2-derivative} and \eqref{Claim-H2}. It therefore follows
\begin{equation}\label{H2-bnd}
	\begin{aligned}
		\| \M^{- \frac{q}{2} } f \|^2_{H^2} \leq C \sum_{|\alpha| \leq 2} \| \langle v \rangle^{\gamma (|\alpha| - 3)} \M^{- \frac{q}{2} } \partial^\alpha_v g \|^2_{L^2}
	\end{aligned}
\end{equation}
for $- 3 < \gamma \leq - \frac{3}{2}$ and $0 < q < 1$. Furthermore, by the Sobolev embedding $H^2 \hookrightarrow L^\infty$,
\begin{equation}\label{Embedding}
	\begin{aligned}
		\| \M^{- \frac{q}{2} } f \|_{L^\infty} \leq C \| \M^{- \frac{q}{2} } f \|_{H^2} \,.
	\end{aligned}
\end{equation}
Consequently, the bounds \eqref{H2-bnd} and \eqref{Embedding} conclude the Part (II) of Theorem \ref{Thm1}.

\section{Applications to Hilbert expansion: Proof of Theorem \ref{Thm2}-\ref{Thm3}}\label{Sec-4}

In this section, one mainly aims at proving Theorem \ref{Thm2}-\ref{Thm3}, namely, justifying the decay estimates of $(t,x)$-derivatives of $\L^{-1}$ in Theorem \ref{Thm2}, and deriving the bounds of the expanded terms $F_n$ in Theorem \ref{Thm3}.

\subsection{Decay estimates for $\partial_{t,x}^\beta \L^{-1}$: Proof of Theorem \ref{Thm2}}

In this subsection, one devotes to prove the bounds \eqref{L-inverse-Dtx} and \eqref{L-inverse-Dtx-VSP} for $m \geq 0$. Note that the case $m = 0$ can be directly resulted from Theorem \ref{Thm2}. We thereby focus on the cases $m \geq 1$. For notational simplicity, we still denote by $\M (v) = \M (t,x,v)$ in the following.

\subsubsection{\texttt{Proof of Part (I) of Theorem \ref{Thm2}: The cases $- \frac{3}{2} < \gamma \leq 1$}}

We will employ the induction arguments for the integer $m \geq 1$ to achieve our goal.

\textbf{\emph{Step 1. $m = 1$.}}

For $|\beta| = 1$, $\partial_{t,x}^\beta (\L f) = \partial_{t,x}^\beta g$. Then one has
\begin{equation*}
	\begin{aligned}
		\L \{ (\I - \P) \partial_{t,x}^\beta f \} = \L \{ \partial_{t,x}^\beta f \} = \partial_{t,x}^\beta g - [\partial_{t,x}^\beta , \L ] f : = g_1 \,,
	\end{aligned}
\end{equation*}
where $[A, B]$ represents the commutator operator via meaning $AB - BA$. By Theorem \ref{Thm1}, there holds
\begin{equation}\label{f-beta-perp-0}
	\begin{aligned}
		|(\I - \P) \partial_{t,x}^\beta f (t,x,v) | \leq C \| \langle v \rangle^{k_\gamma} \M^{- \frac{q_1}{2}} g_1 \|_{L^\infty} \Theta_\gamma (v) \M^\frac{q_1}{2} (v)
	\end{aligned}
\end{equation}
for any $0 < q_1 < q_0 = q < 1$, provided that $\| \langle v \rangle^{k_\gamma} \M^{- \frac{q_1}{2}} g_1 \|_{L^\infty} < \infty$. Noticing that
\begin{equation}\label{f-beta-perp-1}
	\begin{aligned}
		\| \langle v \rangle^{k_\gamma} \M^{- \frac{q_1}{2}} g_1 \|_{L^\infty} \leq \| \langle v \rangle^{k_\gamma} \M^{- \frac{q_1}{2}} \partial_{t,x}^\beta g \|_{L^\infty} + \| \langle v \rangle^{k_\gamma} \M^{- \frac{q_1}{2}} [\partial_{t,x}^\beta , \L ] f \|_{L^\infty} \,,
	\end{aligned}
\end{equation}
it suffices to show $\| \langle v \rangle^{k_\gamma} \M^{- \frac{q_1}{2}} [\partial_{t,x}^\beta , \L ] f \|_{L^\infty} < \infty$. By definition of $\L$ in \eqref{Lf}-\eqref{K2} and the expressions of $K_1 f, K_2 f$ in \eqref{K1K2-integ-kernel-forms}, one has
\begin{equation}\label{Comm}
	\begin{aligned}
		[\partial_{t,x}^\beta , \L ] f = (\partial_{t,x}^\beta \nu) f + \int_{\R^3} \partial_{t,x}^\beta \mathbbm{k}_1 (v, v_1) f (v_1) \d v_1 - \int_{\R^3} \partial_{t,x}^\beta \mathbbm{k}_2 (v, v_1) f (v_1) \d v_1
	\end{aligned}
\end{equation}
for $|\beta| = 1$, where the integral kernels $\mathbbm{k}_2 (v, v_1)$ and $\mathbbm{k}_2 (v, v_1)$ are given in Lemma \ref{Lmm-k1-k2}.

A straightforward calculation implies
\begin{equation}\label{nu-beta-1}
	\begin{aligned}
		\partial_{t,x}^\beta \M = \big\{ \tfrac{\partial_{t,x}^\beta \rho}{\rho} + \partial_{t,x}^\beta \u \cdot \tfrac{v - \u}{T} + \tfrac{3 \partial_{t,x}^\beta T}{6 T} ( \tfrac{|v - \u|^2}{T} - 3 ) \big\} \M
	\end{aligned}
\end{equation}
for $|\beta| = 1$. Then, by \eqref{nu},
\begin{equation}\label{nu-beta-2}
	\begin{aligned}
		\partial_{t,x}^\beta \nu = \int_{\R^3} \big\{ \tfrac{\partial_{t,x}^\beta \rho}{\rho} + \partial_{t,x}^\beta \u \cdot \tfrac{v_1 - \u}{T} + \tfrac{3 \partial_{t,x}^\beta T}{6 T} ( \tfrac{|v_1 - \u|^2}{T} - 3 ) \big\} |v_1 - v|^\gamma \M (v_1) \d v_1 \,.
	\end{aligned}
\end{equation}
Since $|\big\{ \tfrac{\partial_{t,x}^\beta \rho}{\rho} + \partial_{t,x}^\beta \u \cdot \tfrac{v_1 - \u}{T} + \tfrac{3 \partial_{t,x}^\beta T}{6 T} ( \tfrac{|v_1 - \u|^2}{T} - 3 ) \big\}| \M^\frac{1}{2} (v_1) \leq C (\mathfrak{e}_1) < \infty$ uniformly in $t,x,v$, where the constant $C (\mathfrak{e}_1) > 0$ depends on the quantity $\mathfrak{e}_1 > 0$ defined in \eqref{e-m}, one has
\begin{equation*}
	\begin{aligned}
		|\partial_{t,x}^\beta \nu| \leq C \int_{\R^3} |v_1 - v|^\gamma \M^\frac{1}{2} (v_1) \d v_1 \leq C (\mathfrak{e}_1) \nu \,,
	\end{aligned}
\end{equation*}
where the last inequality is derived from Lemma 2.3 of \cite{LS-2010-KRM} and \eqref{nu-equivalent}. Furthermore, Theorem \ref{Thm1} shows us that
\begin{equation}\label{L-d1-1}
	\begin{aligned}
		|f (t,x,v)| \leq C (\mathfrak{e}_0) \| \langle v \rangle^{k_\gamma} \M^{- \frac{q_0}{2}} g \|_{L^\infty} \Theta_\gamma (v) \M^\frac{q_0}{2} (v)
	\end{aligned}
\end{equation}
for $- \frac{3}{2} < \gamma \leq 1$ and any fixed $0 < q_0 < 1$. Here $\Theta_\gamma (v)$ is given in \eqref{Theta-gamma}. Then, one sees
\begin{equation}\label{Comm-1}
	\begin{aligned}
		\| \langle v \rangle^{k_\gamma} \M^{- \frac{q_1}{2}} (\partial_{t,x}^\beta \nu) f \|_{L^\infty} \leq & C (\mathfrak{e}_1) \| \langle v \rangle^{k_\gamma + |\gamma|} \M^\frac{q_0 - q_1}{2} (v) \|_{L^\infty} \| \langle v \rangle^{k_\gamma} \M^{- \frac{q_0}{2}} g \|_{L^\infty} \\
		\leq & C (\mathfrak{e}_1) \| \langle v \rangle^{k_\gamma} \M^{- \frac{q_0}{2}} g \|_{L^\infty} < \infty
	\end{aligned}
\end{equation}
for any $0 < q_1 < q_0 < 1$.

Then we estimate the quantity $\| \langle v \rangle^{k_\gamma} \M^{- \frac{q_1}{2}} \int_{\R^3} \partial_{t,x}^\beta \mathbbm{k}_1 (v, v_1) f (v_1) \d v_1 \|_{L^\infty}$. Recalling the expression of $\mathbbm{k}_1 (v, v_1)$ in \eqref{k1} of Lemma \ref{Lmm-k1-k2}, it infers from a direct calculation that, for $|\beta| = 1$,
\begin{equation}\label{k1-beta-1}
	\begin{aligned}
		\partial_{t,x}^\beta \mathbbm{k}_1 (v, v_1) =  \big\{ \tfrac{\partial_{t,x}^\beta \rho}{\rho} + \partial_{t,x}^\beta \u \cdot \tfrac{(v_1 - \u) + (v - \u)}{T} + \tfrac{\partial_{t,x}^\beta T}{2 T} ( \tfrac{|v_1 - \u|^2  + |v - \u|^2}{2 T} - 3 ) \big\} \mathbbm{k}_1 (v, v_1) \,.
	\end{aligned}
\end{equation}
Then, by introducing a polynomial
\begin{equation}\label{pn}
	\begin{aligned}
		p_n (\varsigma) = 1 + \varsigma + \cdots + \varsigma^n \,, \ \varsigma \in \R
	\end{aligned}
\end{equation}
for any fixed integer $n \geq 0$, there holds
\begin{equation}\label{L-d1-2}
	\begin{aligned}
		|\partial_{t,x}^\beta \mathbbm{k}_1 (v, v_1)| \leq C (\mathfrak{e}_1) [ p_2 (|v_1 - \u|) + p_2 (|v - \u|) ] \mathbbm{k}_1 (v, v_1) \\
		\leq C (\mathfrak{e}_1) \M^{- \frac{1+3 q_0}{4}} (v_1) \M^{- \frac{1 - q_0}{4}} (v) \mathbbm{k}_1 (v, v_1) \,,
	\end{aligned}
\end{equation}
where the last inequality is derived from the bounds
$$[ p_2 (|v_1 - \u|) + p_2 (|v - \u|) ] \M^{ \frac{1+3 q_0}{4}} (v_1) \M^{ \frac{1 - q_0}{4}} (v) \leq C (\mathfrak{e}_0) < \infty$$
uniformly in $t,x,v$. The bounds \eqref{L-d1-1} and \eqref{L-d1-2} thereby imply that
\begin{equation*}
	\begin{aligned}
		& \langle v \rangle^{k_\gamma} \M^{- \frac{q_1}{2}} (v) \Big| \int_{\R^3} \partial_{t,x}^\beta \mathbbm{k}_1 (v, v_1) f (v_1) \d v_1 \Big| \\
		\leq &  C(\mathfrak{e}_1) \mathcal{W}^{1,\infty}_{q_0, q_1} (g) \langle v \rangle^{k_\gamma} \Theta_\gamma (v) \M^\frac{q_0 - q_1}{2} (v) \int_{\R^3} |v_1 - v|^\gamma \M^\frac{1-q_0}{4} (v_1) \M^\frac{1 - q_0}{4} (v) \d v_1 \\
		\leq & C(\mathfrak{e}_1) \mathcal{W}^{1,\infty}_{q_0, q_1} (g) \langle v \rangle^{k_\gamma} \Theta_\gamma (v) \nu (v) \M^\frac{1 + q_0 - 2 q_1}{4} (v) \leq C(\mathfrak{e}_1) \mathcal{W}^{1,\infty}_{q_0, q_1} (g) \,,
	\end{aligned}
\end{equation*}
where the last second inequality is derived from Lemma 2.3 of \cite{LS-2010-KRM} and \eqref{nu-equivalent}. This infers that
\begin{equation}\label{Comm-2}
	\begin{aligned}
		\| \langle v \rangle^{k_\gamma} \M^{- \frac{q_1}{2}} \int_{\R^3} \partial_{t,x}^\beta \mathbbm{k}_1 (v, v_1) f (v_1) \d v_1 \|_{L^\infty} \leq C(\mathfrak{e}_1) \mathcal{W}^{1,\infty}_{q_0, q_1} (g) \,.
	\end{aligned}
\end{equation}

Then, one focuses on the quantity $\| \langle v \rangle^{k_\gamma} \M^{- \frac{q_1}{2}} \int_{\R^3} \partial_{t,x}^\beta \mathbbm{k}_2 (v, v_1) f (v_1) \d v_1 \|_{L^\infty}$. Recalling the expression of $\mathbbm{k}_2 (v, v_1)$ in \eqref{k2-1} and \eqref{k2-2}, it infers from a direct calculation that
\begin{equation}\label{k2-tilde+hat}
	\begin{aligned}
		\partial_{t,x}^\beta \mathbbm{k}_2 (v, v_1) = \partial_{t,x}^\beta \widetilde{\mathbbm{k}}_2 (v, v_1) \widehat{\mathbbm{k}}_2 (v,v_1) + \widetilde{\mathbbm{k}}_2 (v, v_1) \partial_{t,x}^\beta \widehat{\mathbbm{k}}_2 (v,v_1) \,,
	\end{aligned}
\end{equation}
where $\mathbbm{k}_2 (v,v_1) = \widetilde{\mathbbm{k}}_2 (v, v_1) \widehat{\mathbbm{k}}_2 (v,v_1)$,
\begin{equation}\label{k2-tilde}
	\begin{aligned}
		\widetilde{\mathbbm{k}}_2 (v, v_1) = \tfrac{4 \rho}{(2 \pi T)^\frac{3}{2}} |v_1 - v|^{-2} \exp \Big[ - \tfrac{|v_1 - v|^2}{8 T} - \tfrac{(|v_1|^2 - |v|^2)^2}{8 T |v_1 - v|^2} + \tfrac{|v_1|^2 - |v|^2}{2 T |v_1 - v|^2} \u \cdot (v_1 - v) \Big] \,,
	\end{aligned}
\end{equation}
\begin{equation}\label{k2-hat}
	\begin{aligned}
		\widehat{\mathbbm{k}}_2 (v,v_1) = \int_{y \perp (v_1 - v)} \exp ( - \tfrac{|y + \zeta_2 - \u|^2}{2 T} ) ( |v_1 - v|^2 + |y|^2 )^\frac{\gamma}{2} \beta ( \arccos \tfrac{|v_1 - v|}{\sqrt{|v_1 - v|^2 + |y|^2}} ) \d y \,,
	\end{aligned}
\end{equation}
and
\begin{equation}\label{k2-tilde-beta-1}
	\begin{aligned}
		\partial_{t,x}^\beta \widetilde{\mathbbm{k}}_2 (v, v_1) = & \Big[ \tfrac{T^\frac{3}{2}}{\rho} \partial_{t,x}^\beta ( \tfrac{\rho}{T^\frac{3}{2}} ) - \partial_{t,x}^\beta (\tfrac{1}{T}) \big( \tfrac{|v_1 - v|^2}{8} + \tfrac{(|v_1|^2 - |v|^2)^2}{8 |v_1 - v|^2} \big) \\
		& \qquad \qquad \qquad \qquad \qquad + \partial_{t,x}^\beta (\tfrac{\u}{T}) \cdot (v_1 - v) \tfrac{|v_1|^2 - |v|^2}{2 |v_1 - v|^2}  \Big] \widetilde{\mathbbm{k}}_2 (v, v_1) \,,
	\end{aligned}
\end{equation}
\begin{equation}\label{k2-hat-beta-1}
	\begin{aligned}
		\partial_{t,x}^\beta \widehat{\mathbbm{k}}_2 (v, v_1) = & \int_{y \perp (v_1 - v)} \big[ \partial_{t,x}^\beta \u \cdot \tfrac{y + \zeta_2 - \u}{T} + \tfrac{\partial_{t,x}^\beta T}{2 T} \tfrac{|y + \zeta_2 - \u|^2}{T} \big] \exp ( - \tfrac{|y + \zeta_2 - \u|^2}{2 T} ) \\
		& \qquad \qquad \qquad \qquad \times ( |v_1 - v|^2 + |y|^2 )^\frac{\gamma}{2} \beta ( \arccos \tfrac{|v_1 - v|}{\sqrt{|v_1 - v|^2 + |y|^2}} ) \d y \,.
	\end{aligned}
\end{equation}
Observe that
\begin{equation*}
	\begin{aligned}
		\big| \partial_{t,x}^\beta \u \cdot \tfrac{y + \zeta_2 - \u}{T} + \tfrac{\partial_{t,x}^\beta T}{2 T} \tfrac{|y + \zeta_2 - \u|^2}{T} \big| \exp ( - \tfrac{|y + \zeta_2 - \u|^2}{4 T} ) \leq C (\mathfrak{e}_1)
	\end{aligned}
\end{equation*}
uniformly in $y,t,x$. Then the same arguments in \eqref{ClaimA} imply that
\begin{equation}\label{k2-hat-beta}
	\begin{aligned}
		| \widehat{\mathbbm{k}}_2 (v, v_1)| + |\partial_{t,x}^\beta \widehat{\mathbbm{k}}_2 (v, v_1)| \leq C (\mathfrak{e}_0) |v_1 - v|^{b_0 + \gamma} \,,
	\end{aligned}
\end{equation}
where $b_0$ is given in Lemma \ref{Lmm-k1k2-bnd}. By the relation \eqref{Cross-u}, one has
\begin{equation}\label{k2-tilde-1}
	\begin{aligned}
		\widetilde{\mathbbm{k}}_2 (v, v_1) \leq C (\mathfrak{e}_0) |v_1 - v|^{-2} \exp \Big[ - \tfrac{|v_1 - v|^2}{8 T} - \tfrac{(1 - \epsilon)(|v_1|^2 - |v|^2)^2}{8 T |v_1 - v|^2} \Big]
	\end{aligned}
\end{equation}
for any fixed $\epsilon \in (0,1)$. Moreover, there hold
\begin{equation*}
	\begin{aligned}
		\Big| \tfrac{T^\frac{3}{2}}{\rho} \partial_{t,x}^\beta ( \tfrac{\rho}{T^\frac{3}{2}} ) - \partial_{t,x}^\beta (\tfrac{1}{T}) \big( \tfrac{|v_1 - v|^2}{8} + \tfrac{(|v_1|^2 - |v|^2)^2}{8 |v_1 - v|^2} \big) + \partial_{t,x}^\beta (\tfrac{\u}{T}) \cdot (v_1 - v) \tfrac{|v_1|^2 - |v|^2}{2 |v_1 - v|^2}  \Big| \\
		\leq C (\mathfrak{e}_1) (1 + |v_1 - v|^2) (1 + \tfrac{(|v_1|^2 - |v|^2)^2}{ |v_1 - v|^2})
	\end{aligned}
\end{equation*}
and
\begin{equation*}
	\begin{aligned}
		(1 + |v_1 - v|^2) (1 + \tfrac{(|v_1|^2 - |v|^2)^2}{ |v_1 - v|^2}) \exp \Big[ - \tfrac{\epsilon' |v_1 - v|^2}{8 T} - \tfrac{ \epsilon (|v_1|^2 - |v|^2)^2}{8 T |v_1 - v|^2} \Big] \leq C (\mathfrak{e}_0) < \infty
	\end{aligned}
\end{equation*}
for any fixed $\epsilon \in (0,\frac{1}{2})$ and $\epsilon' \in (0,1)$. Consequently, together with \eqref{k2-tilde-1},
\begin{equation}\label{k2-tilde-beta}
	\begin{aligned}
		| \widetilde{\mathbbm{k}}_2 (v, v_1) | + | \partial_{t,x}^\beta \widetilde{\mathbbm{k}}_2 (v, v_1) | \leq C(\mathfrak{e}_1) |v_1 - v|^{-2} \exp ( - \tfrac{(1 - \epsilon' )|v_1 - v|^2}{8 T} - \tfrac{ (1 - 2 \epsilon ) (|v_1|^2 - |v|^2)^2}{8 T |v_1 - v|^2} )  \,.
	\end{aligned}
\end{equation}
From the relations \eqref{k2-tilde+hat}, \eqref{k2-hat-beta} and \eqref{k2-tilde-beta}, it infers that, for any fixed $\epsilon \in (0,\frac{1}{2})$ and $\epsilon' \in (0,1)$,
\begin{equation}\label{k2-beta}
	\begin{aligned}
		|\partial_{t,x}^\beta \mathbbm{k}_2 (v,v_1)| \leq C(\mathfrak{e}_1) |v_1 - v|^{b_0 + \gamma - 2} \exp (- \tfrac{(1 - \epsilon')|v_1 - v|^2}{8 T} - \tfrac{ (1 - 2 \epsilon ) (|v_1|^2 - |v|^2)^2}{8 T |v_1 - v|^2}) \,.
	\end{aligned}
\end{equation}
Together with \eqref{L-d1-1} and \eqref{k2-beta}, one sees that, for any fixed $\epsilon \in (0,\tfrac{1}{2})$, $\epsilon' \in (0,1)$ and $0 < q_1 < q_0' < q_0 < 1$,
\begin{equation}\label{k2-beta-1}
	\begin{aligned}
		& \langle v \rangle^{k_\gamma} \M^{- \frac{q_1}{2}} (v) | \int_{\R^3} \partial_{t,x}^\beta \mathbbm{k}_2 (v, v_1) f (v_1) \d v_1 | \\
		\leq & C (\mathfrak{e}_1) \| \langle v \rangle^{k_\gamma} \M^{- \frac{q_0}{2}} g \|_{L^\infty}  \langle v \rangle^{k_\gamma} \M^{- \frac{q_1}{2}} (v) \\
		& \times \int_{\R^3} \Theta_\gamma (v_1) \M^\frac{q_0}{2} (v_1) |v_1 - v|^{b_0 + \gamma - 2} \exp (- \tfrac{(1 - \epsilon')|v_1 - v|^2}{8 T} - \tfrac{ (1 - 2 \epsilon ) (|v_1|^2 - |v|^2)^2}{8 T |v_1 - v|^2}) \d v_1 \\
		\leq & C (\mathfrak{e}_1) \mathcal{W}^{1,\infty}_{q_0,q_1} (g) \langle v \rangle^{k_\gamma} \M^{ \frac{q_0' - q_1}{2}} (v) \int_{\R^3} \Theta_\gamma (v_1) \M^\frac{q_0 - q_0'}{2} (v_1) |v_1 - v|^{b_0 + \gamma - 2} \\
		& \qquad \qquad \qquad \qquad \times \exp (\tfrac{q_0' (|v - \u|^2 - |v_1 - \u|^2)}{4 T} - \tfrac{(1 - \epsilon')|v_1 - v|^2}{8 T} - \tfrac{ (1 - 2 \epsilon ) (|v_1|^2 - |v|^2)^2}{8 T |v_1 - v|^2}) \d v_1 \\
		\leq & C (\mathfrak{e}_1) \mathcal{W}^{1,\infty}_{q_0,q_1} (g) \int_{\R^3} |v_1 - v|^{b_0 + \gamma - 2} \\
		& \qquad \qquad \qquad \qquad \times \exp (\tfrac{q_0' (|v - \u|^2 - |v_1 - \u|^2)}{4 T} - \tfrac{(1 - \epsilon')|v_1 - v|^2}{8 T} - \tfrac{ (1 - 2 \epsilon ) (|v_1|^2 - |v|^2)^2}{8 T |v_1 - v|^2}) \d v_1 \,,
	\end{aligned}
\end{equation}
where the last inequality is derived from the uniform bounds $\langle v \rangle^{k_\gamma} \M^{ \frac{q_0' - q_1}{2}} (v) \leq C (\mathfrak{e}_0)$ and $\Theta_\gamma (v_1) \M^\frac{q_0 - q_0'}{2} (v_1) \leq C (\mathfrak{e}_0)$.

Noticing that $|v_1|^2 - |v|^2 = |v_1 - v|^2 + 2 v \cdot (v_1 - v)$, one has
\begin{equation*}
	\begin{aligned}
		|v - \u|^2 - |v_1 - \u|^2 = & - (v_1 - v) \cdot (v_1 + v - 2 \u) \\
		= & - |v_1 - v|^2 - 2 (v - \u) \cdot (v_1 - v) \\
		= & - (|v_1|^2 - |v|^2) + 2 \u \cdot (v_1 - v) \,.
	\end{aligned}
\end{equation*}
Then, for any $\epsilon \in (0,\tfrac{1}{2})$, $\epsilon' \in (0,1)$ and $0 < q_1 < q_0' < q_0 < 1$, a direct calculation implies
\begin{equation*}
	\begin{aligned}
		& \tfrac{q_0' (|v - \u|^2 - |v_1 - \u|^2)}{4 T} - \tfrac{(1 - \epsilon')|v_1 - v|^2}{8 T} - \tfrac{ (1 - 2 \epsilon ) (|v_1|^2 - |v|^2)^2}{8 T |v_1 - v|^2} \\
		= & - \tfrac{1 - 2 \epsilon}{8 T} \Big[ (\tfrac{|v_1|^2 - |v|^2}{|v_1 - v|})^2 + 2 \cdot\tfrac{q_0'|v_1 - v|}{1-2\epsilon} \cdot \tfrac{|v_1|^2 - |v|^2}{|v_1 - v|} \Big] - \tfrac{(1 - \epsilon')|v_1 - v|^2}{8 T} + \tfrac{q_0' \u \cdot (v_1 - v)}{2 T} \\
		= & - \tfrac{1 - 2 \epsilon}{8 T} \Big( \tfrac{|v_1|^2 - |v|^2}{|v_1 - v|} + \tfrac{q_0'|v_1 - v|}{1-2\epsilon} \Big)^2 - (1 - \epsilon' - \tfrac{q_0'^2}{1 - 2 \epsilon}) \tfrac{|v_1 - v|^2}{8 T} + \tfrac{q_0' \u \cdot (v_1 - v)}{2 T} \\
		= & - \tfrac{1 - 2 \epsilon}{8 T} \Big( \tfrac{|v_1|^2 - |v|^2}{|v_1 - v|} + \tfrac{q_0'|v_1 - v|}{1-2\epsilon} \Big)^2 - \tfrac{1}{2} (1 - \epsilon' - \tfrac{q_0'^2}{1 - 2 \epsilon}) \tfrac{|v_1 - v|^2}{8 T} \\
		& - \tfrac{1 - \epsilon' - \tfrac{q_0'^2}{1 - 2 \epsilon}}{16 T} \big| v_1 - v - 4 q_0' (1 - \epsilon' - \tfrac{q_0'^2}{1 - 2 \epsilon})^{-1} \u \big|^2 + q_0'^2 (1 - \epsilon' - \tfrac{q_0'^2}{1 - 2 \epsilon})^{-1} \tfrac{|\u|^2}{T} \,.
	\end{aligned}
\end{equation*}
By taking $\epsilon = \frac{q - q_0'^2}{4} \in (0, \frac{1}{2})$ and $\epsilon' =\tfrac{1 - q_0'^2}{2 (1 + q_0'^2)}\in (0, \frac{1}{2})$ such that
\begin{equation*}
	\begin{aligned}
		1 - \epsilon' - \tfrac{q_0'^2}{1 - 2 \epsilon} = \tfrac{1 - q_0'^2}{2 (1 + q_0'^2)} > 0 \,,
	\end{aligned}
\end{equation*}
it infers that
\begin{equation}\label{k2-beta-2}
	\begin{aligned}
		\tfrac{q_0' (|v - \u|^2 - |v_1 - \u|^2)}{4 T} - \tfrac{(1 - \epsilon')|v_1 - v|^2}{8 T} - \tfrac{ (1 - 2 \epsilon ) (|v_1|^2 - |v|^2)^2}{8 T |v_1 - v|^2} \\
		\leq - \tfrac{1 - q_0'^2}{4 (1 + q_0'^2)} \tfrac{|v_1 - v|^2}{8 T} + \tfrac{2 q_0'^2 (1 + q_0'^2)}{1 - q_0'^2} \tfrac{|\u|^2}{T} \,.
	\end{aligned}
\end{equation}
Collecting the estimates \eqref{k2-beta-1} and \eqref{k2-beta-2}, one has
\begin{equation*}
	\begin{aligned}
		& \langle v \rangle^{k_\gamma} \M^{- \frac{q_1}{2}} (v) | \int_{\R^3} \partial_{t,x}^\beta \mathbbm{k}_2 (v, v_1) f (v_1) \d v_1 | \\
		\leq & C (\mathfrak{e}_1) \mathcal{W}^{1,\infty}_{q_0,q_1} (g) \int_{\R^3} |v_1 - v|^{b_0 + \gamma - 2} \exp \big[ - \tfrac{1 - q_0'^2}{4 (1 + q_0'^2)} \tfrac{|v_1 - v|^2}{8 T} \big] \d v_1 \\
		= & C (\mathfrak{e}_1) \mathcal{W}^{1,\infty}_{q_0,q_1} (g) \int_{\R^3} |v_1 |^{b_0 + \gamma - 2} \exp \big[ - \tfrac{1 - q_0'^2}{4 (1 + q_0'^2)} \tfrac{|v_1 |^2}{8 T} \big] \d v_1 \leq C (\mathfrak{e}_1) \mathcal{W}^{1,\infty}_{q_0,q_1} (g) \,,
	\end{aligned}
\end{equation*}
where the fact $\exp \big[ \tfrac{2 q_0'^2 (1 + q_0'^2)}{1 - q_0'^2} \tfrac{|\u|^2}{T} \big] \leq C (\mathfrak{e}_0)$ has been used in the first inequality, and the last inequality is followed from the same arguments of proving $\mathbf{k}_2 (v) \in L^1$ in Lemma \ref{Lmm-k1k2-bnd}. Namely, for any fixed $0 < q_1 < q_0 < 1$,
\begin{equation}\label{Comm-3}
	\begin{aligned}
		\| \langle v \rangle^{k_\gamma} \M^{- \frac{q_1}{2}} \int_{\R^3} \partial_{t,x}^\beta \mathbbm{k}_2 (v, v_1) f (v_1) \d v_1 \|_{L^\infty} \leq C (\mathfrak{e}_1) \mathcal{W}^{1,\infty}_{q_0,q_1} (g) \,.
	\end{aligned}
\end{equation}
Then, \eqref{Comm}, \eqref{Comm-1}, \eqref{Comm-2} and \eqref{Comm-3} conclude that
\begin{equation}\label{f-beta-perp-2}
	\begin{aligned}
		\| \langle v \rangle^{k_\gamma} \M^{- \frac{q_1}{2}} [\partial_{t,x}^\beta , \L ] f \|_{L^\infty} \leq C (\mathfrak{e}_1) \mathcal{W}^{1,\infty}_{q_0,q_1} (g) \,.
	\end{aligned}
\end{equation}
Consequently, the relations \eqref{f-beta-perp-0}, \eqref{f-beta-perp-1} and \eqref{f-beta-perp-2} imply that, for $|\beta| = 1$ and any $0 < q_1 < q_0 < 1$,
\begin{equation}\label{f-perp-1}
	\begin{aligned}
		|(\I - \P) \partial_{t,x}^\beta f (t,x,v) | \leq C (\mathfrak{e}_1) \mathcal{W}^{1,\infty}_{q_0,q_1} (g) \Theta_\gamma (v) \M^\frac{q_1}{2} (v) \,.
	\end{aligned}
\end{equation}

Next we control $\P \partial_{t,x}^\beta f (t,x,v)$. Recalling the definition of $\P$ in \eqref{Projt-P}, we have
\begin{equation}\label{P-f-beta}
	\begin{aligned}
		\P \partial_{t,x}^\beta f = \big\{ \tfrac{ \rho_\beta}{\rho} + \u_\beta \cdot \tfrac{v - \u}{T} + \tfrac{ \theta_\beta }{6 T} ( \tfrac{|v - \u|^2}{T} - 3 ) \big\} \sqrt{\M} \,,
	\end{aligned}
\end{equation}
where
\begin{equation}\label{P-f-beta-fluid}
	\begin{aligned}
		\rho_\beta = \langle \partial_{t,x}^\beta f , \sqrt{\M} \rangle \,, \ u_\beta = \tfrac{1}{\rho} \langle  \partial_{t,x}^\beta f , \tfrac{v - \u}{T} \sqrt{\M} \rangle \,, \ \theta_\beta = \tfrac{1}{\rho} \langle \partial_{t,x}^\beta f , (|v - \u|^2 - 3 T) \sqrt{\M} \rangle \,.
	\end{aligned}
\end{equation}
It suffices to control the fluid variables $(\rho_\beta, u_\beta, \theta_\beta)$. For convenience of notations, we denote by
\begin{equation}\label{phi}
	\begin{aligned}
		\phi (v) = \sqrt{\M} \,, \tfrac{v - \u}{\rho T} \sqrt{\M} \,, \tfrac{1}{\rho} (|v - \u|^2 - 3 T) \sqrt{\M} \,.
	\end{aligned}
\end{equation}
Then the quantity $\langle \partial_{t,x}^\beta f , \phi \rangle$ can be transformed to
\begin{equation}\label{f-beta-phi-beta}
	\begin{aligned}
		\langle \partial_{t,x}^\beta f , \phi \rangle = \partial_{t,x}^\beta \langle f , \phi \rangle - \langle f , \partial_{t,x}^\beta \phi \rangle = - \langle f , \partial_{t,x}^\beta \phi \rangle \,,
	\end{aligned}
\end{equation}
since $f \in \mathrm{Null}^\perp (\L)$ and $\phi \in \mathrm{Null} (\L)$. It is easy to see that
\begin{equation*}
	\begin{aligned}
		|\partial_{t,x}^\beta \phi | \leq C(\mathfrak{e}_1) p_4 (|v - \u|) \M^\frac{1}{2} (v) \,,
	\end{aligned}
\end{equation*}
where the polynomial $p_4$ is given in \eqref{pn}. Together with \eqref{L-d1-1}, there holds
\begin{equation*}
	\begin{aligned}
		|\langle \partial_{t,x}^\beta f , \phi \rangle| = & | \langle f , \partial_{t,x}^\beta \phi \rangle | \leq \int_{\R^3} |f (t,x,v)| |\partial_{t,x}^\beta \phi (v)| \d v \\
		\leq & C (\mathfrak{e}_1) \int_{\R^3} \| \langle v \rangle^{k_\gamma} \M^{- \frac{q_0}{2}} g \|_{L^\infty} \Theta_\gamma (v) \M^\frac{q_0}{2} (v) \cdot p_4 (|v - \u|) \M^\frac{1}{2} (v) \d v \\
		\leq & C (\mathfrak{e}_1) \mathcal{W}^{1,\infty}_{q_0,q_1} (g) \,,
	\end{aligned}
\end{equation*}
where $\int_{\R^3} \Theta_\gamma (v) \M^\frac{q_0}{2} (v) \cdot p_4 (|v - \u|) \M^\frac{1}{2} (v) \d v \leq C (\mathfrak{e}_0)$ has been used. As a result,
\begin{equation*}
	\begin{aligned}
		|(\rho_\beta, u_\beta, \theta_\beta)| \leq C (\mathfrak{e}_1) \mathcal{W}^{1,\infty}_{q_0,q_1} (g) \,.
	\end{aligned}
\end{equation*}
It thereby infers that
\begin{equation}\label{f-null-1}
	\begin{aligned}
		| \P \partial_{t,x}^\beta f | = & \big| \big\{ \tfrac{ \rho_\beta}{\rho} + \u_\beta \cdot \tfrac{v - \u}{T} + \tfrac{ \theta_\beta }{6 T} ( \tfrac{|v - \u|^2}{T} - 3 ) \big\} \big| \sqrt{\M} \\
		\leq & C (\mathfrak{e}_0) |(\rho_\beta, u_\beta, \theta_\beta)| p_2 (|v-\u|) \M^\frac{1}{2} (v) \leq C (\mathfrak{e}_1) \mathcal{W}^{1,\infty}_{q_0,q_1} (g) \Theta_\gamma (v) \M^\frac{q_1}{2} (v)
	\end{aligned}
\end{equation}
for any $0 < q_1 < 1$, where $p_2 (|v - \u|) \Theta_\gamma^{-1} (v) \M^\frac{1-q_1}{2} (v) \leq C (\mathfrak{e}_0)$ is utilized. Therefore, the bounds \eqref{f-perp-1} and \eqref{f-null-1} conclude the results in Theorem \ref{Thm2} with $m = 1$.

\textbf{\emph{Step 2. Induction assumptions for $m \leq \mathfrak{j}$ with any fixed $\mathfrak{j} \geq 1$.}}

Assume that the results in Theorem \ref{Thm2} hold for all $m = 1, \cdots, \mathfrak{j}$. That is, for any $|\beta| = m \in \{ 1, \cdots, \mathfrak{j} \}$ and $0 < q_m < q_{m-1} < \cdots < q_1 < q_0 < 1$,
\begin{equation}\label{Assume-f}
	\begin{aligned}
		|\partial_{t,x}^\beta f (t,x,v)| \leq C (\mathfrak{e}_m) \mathcal{W}^{m,\infty}_{q_0, \cdots , q_m} (g) \Theta_\gamma (v) \M^\frac{q_m}{2} (v) \,.
	\end{aligned}
\end{equation}

\textbf{\emph{Step 3. $m = \mathfrak{j} + 1$.}}

For $|\beta| = \mathfrak{j} + 1$, one has
\begin{equation*}
	\begin{aligned}
		\L \{ (\I - \P) \partial_{t,x}^\beta f \} = \L \{ \partial_{t,x}^\beta f \} = \partial_{t,x}^\beta g - [\partial_{t,x}^\beta , \L ] f : = g_{\mathfrak{j} + 1} \,,
	\end{aligned}
\end{equation*}
where
\begin{equation}\label{Comm-beta-0}
	\begin{aligned}
		[\partial_{t,x}^\beta \,, \L ] f = [ \partial_{t,x}^\beta, \nu] f + \int_{\R^3} [ \partial_{t,x}^\beta \,, \mathbbm{k}_1 (v, v_1) ] f (v_1) \d v_1 - \int_{\R^3} [ \partial_{t,x}^\beta \,, \mathbbm{k}_2 (v, v_1) ] f (v_1) \d v_1 \,.
	\end{aligned}
\end{equation}
Theorem \ref{Thm1} tells us that
\begin{equation}\label{f-beta-j+1-0}
	\begin{aligned}
		| (\I - \P) \partial_{t,x}^\beta f (t,x,v) | \leq C \| \langle v \rangle^{k_\gamma} \M^{- \frac{q_{\j + 1}}{2}} g_{\j + 1} \|_{L^\infty} \Theta_\gamma (v) \M^\frac{q_{\j+1}}{2} (v)
	\end{aligned}
\end{equation}
for any $0 < q_{\j + 1} < q_\j < \cdots < q_1 < q_0 < 1$, where $\Theta_\gamma (v)$ is given in \eqref{Theta-gamma}. Observe that
\begin{equation}\label{f-beta-j+1-1}
	\begin{aligned}
		\| \langle v \rangle^{k_\gamma} \M^{- \frac{q_{\j + 1}}{2}} g_{\j + 1} \|_{L^\infty} \leq \mathcal{W}^{\j+1,\infty}_{q_0,\cdots,q_{\j+1}} (g) + \| \langle v \rangle^{k_\gamma} \M^{- \frac{q_{\j + 1}}{2}} [\partial_{t,x}^\beta \,, \L ] f \|_{L^\infty} \,.
	\end{aligned}
\end{equation}

We first consider the quantity $\| \langle v \rangle^{k_\gamma} \M^{- \frac{q_{\j + 1}}{2}} [ \partial_{t,x}^\beta, \nu] f \|_{L^\infty}$ for $|\beta| = \j + 1$. Notice that
\begin{equation*}
	\begin{aligned}
		[ \partial_{t,x}^\beta, \nu] f = \sum_{0 \neq \beta' \leq \beta} C_\beta^{\beta'} \partial_{t,x}^{\beta'} \nu \partial_{t,x}^{\beta - \beta'} f \,,
	\end{aligned}
\end{equation*}
where $C_\beta^{\beta'}$ is the constant depending only on $\beta$ and $\beta'$. Observe that, by \eqref{nu-beta-1} and \eqref{nu-beta-2},
\begin{equation*}
	\begin{aligned}
		|\partial_{t,x}^{\beta'} \nu | \leq & C (\mathfrak{e}_{|\beta'|}) \int_{\R^3} p_{2 |\beta'|} (|v_1 - \u|) |v_1 - v|^\gamma \M (v_1) \d v_1 \\
		\leq & C (\mathfrak{e}_{\j + 1}) \int_{\R^3} |v_1 - v|^\gamma \M^\frac{1}{2} (v_1) \d v_1 \leq C (\mathfrak{e}_{\j + 1}) \nu \,,
	\end{aligned}
\end{equation*}
where $p_{2 |\beta'|} (|v_1 - \u|) \M^\frac{1}{2} (v_1) \leq C (\mathfrak{e}_0)$ is utilized in the last second inequality. For $|\beta| = \j + 1$ and $\beta' \neq 0$, one has $0 \leq |\beta - \beta'| \leq \j$. Together with \eqref{L-d1-1} and \eqref{Assume-f}, we have
\begin{equation*}
	\begin{aligned}
		| \partial_{t,x}^{\beta - \beta'} f | \leq C (\mathfrak{e}_{|\beta - \beta'|}) \mathcal{W}^{|\beta - \beta'|, \infty}_{q_0, \cdots, q_{|\beta - \beta'|}} (g) \Theta_\gamma (v) \M^\frac{q_{|\beta - \beta'|}}{2} (v) \,.
	\end{aligned}
\end{equation*}
Due to $q_{\j + 1} < q_\j \leq q_{|\beta - \beta'|}$ for $0 \leq |\beta - \beta' | \leq \j$, there holds
$$\langle v \rangle^{k_\gamma} \M^{- \frac{q_{\j + 1}}{2}} (v) \nu (v) \Theta_\gamma (v) \M^\frac{q_{|\beta - \beta'|}}{2} (v) \leq C (\mathfrak{e}_0) $$
uniformly in $t,x,v$. Consequently,
\begin{equation}\label{Comm-beta-1}
	\begin{aligned}
		\| \langle v \rangle^{k_\gamma} \M^{- \frac{q_{\j + 1}}{2}} [ \partial_{t,x}^\beta, \nu] f \|_{L^\infty} & \leq C (\mathfrak{e}_{\j + 1}) \sum_{0 \neq \beta' \leq \beta} C (\mathfrak{e}_{|\beta - \beta'|}) \mathcal{W}^{|\beta - \beta'|, \infty}_{q_0, \cdots, q_{|\beta - \beta'|}} (g) \\
		& \times \| \langle v \rangle^{k_\gamma} \M^{- \frac{q_{\j + 1}}{2}} (v) \nu (v) \Theta_\gamma (v) \M^\frac{q_{|\beta - \beta'|}}{2} (v) \|_{L^\infty} \\
		& \leq C (\mathfrak{e}_{\j + 1}) \mathcal{W}^{\j, \infty}_{q_0, \cdots, q_{\j}} (g) \,.
	\end{aligned}
\end{equation}

We then consider the quantity $\| \langle v \rangle^{k_\gamma} \M^{- \frac{q_{\j + 1}}{2}} \int_{\R^3} [ \partial_{t,x}^\beta \,, \mathbbm{k}_1 (v, v_1) ] f (v_1) \d v_1 \|_{L^\infty}$ for $|\beta| = \j + 1$. Observe that
\begin{equation*}
	\begin{aligned}
		[ \partial_{t,x}^\beta \,, \mathbbm{k}_1 (v, v_1) ] f (v_1) = \sum_{0 \neq \beta' \leq \beta} C_\beta^{\beta'} \partial_{t,x}^{\beta'} \mathbbm{k}_1 (v, v_1) \partial_{t,x}^{\beta - \beta'} f (v_1) \,.
	\end{aligned}
\end{equation*}
Repeating the arguments of \eqref{k1-beta-1}, one easily knows that
\begin{equation*}
	\begin{aligned}
		|\partial_{t,x}^{\beta'} \mathbbm{k}_1 (v, v_1)| \leq C (\mathfrak{e}_{|\beta'|}) [ p_{2 |\beta'|} (|v - \u|) + p_{2 |\beta'|} (|v_1 - \u|)] \mathbbm{k}_1 (v, v_1) \,,
	\end{aligned}
\end{equation*}
where $p_{2 |\beta'|}$ is given in \eqref{pn}. Due to $0 \leq |\beta - \beta'| \leq \j$, $1 \leq |\beta'| \leq \j + 1$ and $0 < q_{\j + 1} < q_\j \leq q_{|\beta - \beta'|}$, one has $\M^\frac{q_{|\beta - \beta'|}}{2} (v_1) \leq \M^\frac{q_\j}{2} (v_1)$ and
\begin{equation*}
	\begin{aligned}
		 p_{2 |\beta'|} (|v - \u|) + p_{2 |\beta'|} (|v_1 - \u|) \leq & p_{2 \j + 2} (|v - \u|) + p_{2 \j + 2} (|v_1 - \u|) \\
		 \leq & C (\e_0) \M^{- \frac{1 + 3 q_\j}{4}} (v_1) \M^{- \frac{1 - q_\j}{4}} (v) \,.
	\end{aligned}
\end{equation*}
Then, combining with \eqref{L-d1-1} and \eqref{Assume-f},
\begin{equation*}
	\begin{aligned}
		& \int_{\R^3} \big| [ \partial_{t,x}^\beta \,, \mathbbm{k}_1 (v, v_1) ] f (v_1) \big| \d v_1 \\
		\leq & \sum_{0 \neq \beta' \leq \beta} C (\e_{|\beta'|}) \int_{\R^3} [ p_{2 |\beta'|} (|v - \u|) + p_{2 |\beta'|} (|v_1 - \u|)] \mathbbm{k}_1 (v, v_1) \\
		& \qquad \qquad \times C (\e_{|\beta - \beta'|}) \mathcal{W}^{|\beta - \beta'|, \infty}_{q_0, \cdots, q_{|\beta - \beta'|}} (g) \Theta_\gamma (v_1) \M^\frac{q_{|\beta - \beta'|}}{2} (v_1) \d v_1 \\
		\leq & C (\e_{\j + 1}) \mathcal{W}^{\j + 1, \infty}_{q_0, \cdots, q_{\j + 1}} (g) \M^\frac{1 + q_\j}{4} (v) \int_{\R^3} \Theta_\gamma (v_1) |v_1 - v|^\gamma \M^\frac{1-q_\j}{4} (v_1) \d v_1 \\
		\leq &  C (\e_{\j + 1}) \mathcal{W}^{\j + 1, \infty}_{q_0, \cdots, q_{\j + 1}} (g) \M^\frac{1 + q_\j}{4} (v) \int_{\R^3} |v_1 - v|^\gamma \M^\frac{1-q_\j}{8} (v_1) \d v_1 \\
		\leq &  C (\e_{\j + 1}) \mathcal{W}^{\j + 1, \infty}_{q_0, \cdots, q_{\j + 1}} (g) \nu (v) \M^\frac{1 + q_\j}{4} (v) \,,
	\end{aligned}
\end{equation*}
which implies that
\begin{equation}\label{Comm-beta-2}
	\begin{aligned}
		& \| \langle v \rangle^{k_\gamma} \M^{- \frac{q_{\j + 1}}{2}} \int_{\R^3} [ \partial_{t,x}^\beta \,, \mathbbm{k}_1 (v, v_1) ] f (v_1) \d v_1 \|_{L^\infty} \\
		\leq & C (\e_{\j + 1}) \mathcal{W}^{\j + 1, \infty}_{q_0, \cdots, q_{\j + 1}} (g) \| \langle v \rangle^{k_\gamma} \nu (v) \M^\frac{1 + q_\j - 2 q_{\j + 1}}{4} (v) \|_{L^\infty} \leq C (\e_{\j + 1}) \mathcal{W}^{\j + 1, \infty}_{q_0, \cdots, q_{\j + 1}} (g)
	\end{aligned}
\end{equation}
for all $0 < q_{\j + 1} < q_\j < 1$.

One now focuses on the quantity $\| \langle v \rangle^{k_\gamma} \M^{- \frac{q_{\j + 1}}{2}} \int_{\R^3} [ \partial_{t,x}^\beta \,, \mathbbm{k}_2 (v, v_1) ] f (v_1) \d v_1 \|_{L^\infty}$ for $|\beta| = \j + 1$. Observe that
\begin{equation}\label{Comm-k2-beta}
	\begin{aligned}
		[ \partial_{t,x}^\beta \,, \mathbbm{k}_2 (v, v_1) ] f (v_1) = \sum_{0 \neq \beta' \leq \beta} C_\beta^{\beta'} \partial_{t,x}^{\beta'} \mathbbm{k}_2 (v, v_1) \partial_{t,x}^{\beta - \beta'} f (v_1) \,.
	\end{aligned}
\end{equation}
Recall that $\mathbbm{k}_2 (v,v_1) = \widetilde{\mathbbm{k}}_2 (v, v_1) \widehat{\mathbbm{k}}_2 (v,v_1)$, where $\widetilde{\mathbbm{k}}_2 (v, v_1)$ and $\widehat{\mathbbm{k}}_2 (v,v_1)$ are given in \eqref{k2-tilde} and \eqref{k2-hat}, respectively. Then
\begin{equation*}
	\begin{aligned}
		\partial_{t,x}^{\beta'} \mathbbm{k}_2 (v, v_1) = \sum_{\beta'' \leq \beta'} C_{\beta'}^{\beta''} \partial_{t,x}^{\beta''} \widetilde{\mathbbm{k}}_2 (v, v_1) \partial_{t,x}^{\beta' - \beta''} \widehat{\mathbbm{k}}_2 (v,v_1) \,.
	\end{aligned}
\end{equation*}
Repeating the processes of taking derivatives as in \eqref{k2-tilde-beta-1} and \eqref{k2-hat-beta-1}, it infers that
\begin{equation*}
	\begin{aligned}
		|\partial_{t,x}^{\beta''} \widetilde{\mathbbm{k}}_2 (v, v_1)| \leq C (\e_{|\beta''|}) \big\{ p_{2 |\beta''|} ( |v_1 - v|) + p_{2 |\beta''|} ( \big| \tfrac{|v_1|^2 - |v|^2}{|v_1 - v|} \big|) \big\} \widetilde{\mathbbm{k}}_2 (v, v_1) \,,
	\end{aligned}
\end{equation*}
and
\begin{equation*}
	\begin{aligned}
		|\partial_{t,x}^{\beta' - \beta''} \widehat{\mathbbm{k}}_2 (v, v_1)| \leq C (\e_{|\beta' - \beta''|}) \int_{y \perp (v_1 - v)} p_{2 |\beta' - \beta''|} ( |y + \zeta_2 - \u| ) \exp ( - \tfrac{|y + \zeta_2 - \u|^2}{2 T} ) \\
		\times ( |v_1 - v|^2 + |y|^2 )^\frac{\gamma}{2} \beta ( \arccos \tfrac{|v_1 - v|}{\sqrt{|v_1 - v|^2 + |y|^2}} ) \d y \,.
	\end{aligned}
\end{equation*}
Following the similar derivations of \eqref{k2-tilde-beta}, there holds
\begin{equation*}
	\begin{aligned}
		|\partial_{t,x}^{\beta''} \widetilde{\mathbbm{k}}_2 (v, v_1)| \leq C (\e_{\j + 1}) |v_1 - v|^{-2} \exp ( - \tfrac{(1 - \epsilon' )|v_1 - v|^2}{8 T} - \tfrac{ (1 - 2 \epsilon ) (|v_1|^2 - |v|^2)^2}{8 T |v_1 - v|^2} )
	\end{aligned}
\end{equation*}
for any $\beta'' \leq \beta'$, $\epsilon \in (0, \frac{1}{2})$ and $\epsilon' \in (0,1)$. By the similar arguments of \eqref{k2-hat-beta}, it is implied that
\begin{equation*}
	\begin{aligned}
		|\partial_{t,x}^{\beta' - \beta''} \widehat{\mathbbm{k}}_2 (v, v_1)| \leq C (\e_{\j + 1}) |v_1 - v|^{b_0 + \gamma}
	\end{aligned}
\end{equation*}
for all $\beta'' \leq \beta'$, where $b_0$ is given in Lemma \ref{Lmm-k1k2-bnd}. Therefore,
\begin{equation*}
	\begin{aligned}
		| \partial_{t,x}^{\beta'} \mathbbm{k}_2 (v, v_1) | \leq C (\e_{\j + 1}) |v_1 - v|^{b_0 + \gamma -2} \exp ( - \tfrac{(1 - \epsilon' )|v_1 - v|^2}{8 T} - \tfrac{ (1 - 2 \epsilon ) (|v_1|^2 - |v|^2)^2}{8 T |v_1 - v|^2} )
	\end{aligned}
\end{equation*}
for any $0 \neq \beta' \leq \beta$, $\epsilon \in (0, \frac{1}{2})$ and $\epsilon' \in (0,1)$. Moreover, from \eqref{L-d1-1} and \eqref{Assume-f}, it is derived that
\begin{equation*}
	\begin{aligned}
		|\partial_{t,x}^{\beta - \beta'} f (v_1)| \leq C (\e_\j) \mathcal{W}^{\j, \infty}_{q_0,\cdots,q_\j} (g) \Theta_\gamma (v_1) \M^\frac{q_{|\beta - \beta'|}}{2} (v_1)
	\end{aligned}
\end{equation*}
for all $0 \neq \beta' \leq \beta$. Together with \eqref{Comm-k2-beta}, it follows that
\begin{equation*}
	\begin{aligned}
		| [ \partial_{t,x}^\beta \,, \mathbbm{k}_2 (v, v_1) ] f (v_1) | \leq C (\e_{\j + 1}) \mathcal{W}^{\j, \infty}_{q_0,\cdots,q_\j} (g) |v_1 - v|^{b_0 + \gamma -2} \Theta_\gamma (v_1) \M^\frac{q_\j}{2} (v_1) \\
		\times \exp ( - \tfrac{(1 - \epsilon' )|v_1 - v|^2}{8 T} - \tfrac{ (1 - 2 \epsilon ) (|v_1|^2 - |v|^2)^2}{8 T |v_1 - v|^2} )
	\end{aligned}
\end{equation*}
for any $\epsilon \in (0, \frac{1}{2})$ and $\epsilon' \in (0,1)$.

Following the same arguments from \eqref{k2-beta} to \eqref{Comm-3}, it is derived that
\begin{equation}\label{Comm-beta-3}
	\begin{aligned}
		\| \langle v \rangle^{k_\gamma} \M^{- \frac{q_{\j + 1}}{2}} \int_{\R^3} [ \partial_{t,x}^\beta \,, \mathbbm{k}_2 (v, v_1) ] f (v_1) \d v_1 \|_{L^\infty} \leq C (\e_{\j + 1}) \mathcal{W}^{\j, \infty}_{q_0,\cdots,q_\j} (g) \,.
	\end{aligned}
\end{equation}
Collecting the relations \eqref{Comm-beta-0}, \eqref{Comm-beta-1}, \eqref{Comm-beta-2} and \eqref{Comm-beta-3}, one has
\begin{equation}\label{f-beta-j+1-2}
	\begin{aligned}
		\| \langle v \rangle^{k_\gamma} \M^{- \frac{q_{\j + 1}}{2}} [ \partial_{t,x}^\beta \,, \L ] f (v_1) \d v_1 \|_{L^\infty} \leq C (\e_{\j + 1}) \mathcal{W}^{\j + 1, \infty}_{q_0,\cdots,q_{\j + 1}} (g) \,.
	\end{aligned}
\end{equation}
Therefore, the relations \eqref{f-beta-j+1-0}, \eqref{f-beta-j+1-1} and \eqref{f-beta-j+1-2} imply that, for $|\beta| = \j + 1$ and $0 < q_{\j + 1} < q_\j < \cdots < q_1 < q_0 < 1$,
\begin{equation}\label{f-perp-j+1}
	\begin{aligned}
		| (\I - \P) \partial_{t,x}^\beta f (t,x,v) | \leq C (\e_{\j + 1}) \mathcal{W}^{\j + 1, \infty}_{q_0,\cdots,q_{\j + 1}} (g) \Theta_\gamma (v) \M^\frac{q_{\j+1}}{2} (v) \,.
	\end{aligned}
\end{equation}

At the end, $\P \partial_{t,x}^\beta f (t,x,v)$ should be estimated. The expression of $\P \partial_{t,x}^\beta f (t,x,v)$ is given in \eqref{P-f-beta}-\eqref{P-f-beta-fluid} for $|\beta| = \j + 1$. The main goal is to control the quantity $\langle \partial_{t,x}^\beta f, \phi \rangle$, where $\phi = \phi (t,x,v)$ is given in \eqref{phi}. Splitting $\beta = \tilde{\beta} + \bar{\beta}$ with $|\tilde{\beta}| = \j$ and $| \bar{\beta} | = 1 $, one has $\partial_{t,x}^\beta = \partial_{t,x}^{\tilde{\beta}} \partial_{t,x}^{\bar{\beta}}$. Moreover, the equality \eqref{f-beta-phi-beta} reads $\langle \partial_{t,x}^{\bar{\beta}} f , \phi \rangle = - \langle f , \partial_{t,x}^{\bar{\beta}} \phi \rangle$, which implies that
\begin{equation*}
	\begin{aligned}
		\partial_{t,x}^{\tilde{\beta}} \langle \partial_{t,x}^{\bar{\beta}} f , \phi \rangle = - \partial_{t,x}^{\tilde{\beta}} \langle f , \partial_{t,x}^{\bar{\beta}} \phi \rangle \,.
	\end{aligned}
\end{equation*}
Then, by $\beta = \tilde{\beta} + \bar{\beta}$, one has
\begin{equation}\label{P-f-beta-j+1-0}
	\begin{aligned}
		\langle \partial_{t,x}^\beta f, \phi \rangle = - \sum_{0 \neq \tilde{\beta}' \leq \tilde{\beta}} C_{\tilde{\beta}}^{\tilde{\beta}'} \langle \partial_{t,x}^{\tilde{\beta} + \bar{\beta} - \tilde{\beta}'} f , \partial_{t,x}^{\tilde{\beta}'} \phi \rangle - \sum_{\tilde{\beta}'' \leq \tilde{\beta}} C_{\tilde{\beta}}^{\tilde{\beta}''} \langle \partial_{t,x}^{\tilde{\beta} - \tilde{\beta}''} f , \partial_{t,x}^{\tilde{\beta}'' + \bar{\beta}} \phi \rangle \,.
	\end{aligned}
\end{equation}
Noticing that $0 \neq \tilde{\beta}' \leq \tilde{\beta}$, $\tilde{\beta}'' \leq \tilde{\beta}$, $|\tilde{\beta}| = \j$ and $| \bar{\beta} | = 1$, one knows that
\begin{equation*}
	\begin{aligned}
		1 \leq |\tilde{\beta} + \bar{\beta} - \tilde{\beta}'| \leq \j \,, 0 \leq |\tilde{\beta} - \tilde{\beta}''| \leq \j \,, 1 \leq |\tilde{\beta}'| \leq \j \,, 1 \leq |\tilde{\beta}'' + \bar{\beta}| \leq \j + 1 \,.
	\end{aligned}
\end{equation*}
As a result, the bounds \eqref{L-d1-1} and \eqref{Assume-f} reduce to
\begin{equation}\label{P-f-beta-j+1-1}
	\begin{aligned}
		\big| \partial_{t,x}^{\tilde{\beta} + \bar{\beta} - \tilde{\beta}'} f (t,x,v) \big| + \big| \partial_{t,x}^{\tilde{\beta} - \tilde{\beta}''} f (t,x,v) \big| \leq C (\mathfrak{e}_\j) \mathcal{W}^{\j,\infty}_{q_0, \cdots , q_\j} (g) \Theta_\gamma (v) \M^\frac{q_\j}{2} (v) \,.
	\end{aligned}
\end{equation}
By the expressions of $\phi$ in \eqref{phi}, it is easy to know that
\begin{equation}\label{P-f-beta-j+1-2}
	\begin{aligned}
		\big| \partial_{t,x}^{\tilde{\beta}'} \phi (t,x,v) \big| + \big| \partial_{t,x}^{\tilde{\beta}'' + \bar{\beta}} \phi (t,x,v) \big| \leq C (\e_{\j + 1}) p_{2 \j + 4} (|v - \u|) \M^\frac{1}{2} (v) \,,
	\end{aligned}
\end{equation}
where the polynomial $p_{2 \j + 4} (|v - \u|)$ is given in \eqref{pn}. Collecting the relations \eqref{P-f-beta-j+1-0}, \eqref{P-f-beta-j+1-1} and \eqref{P-f-beta-j+1-2}, there holds
\begin{equation*}
	\begin{aligned}
		| \langle \partial_{t,x}^\beta f, \phi \rangle | \leq & C (\e_{\j + 1}) \mathcal{W}^{\j,\infty}_{q_0, \cdots , q_\j} (g) \int_{\R^3} p_{2 \j + 4} (|v - \u|) \Theta_\gamma (v) \M^\frac{1 + q_\j}{2} (v) \d v \\
		\leq & C (\e_{\j + 1}) \mathcal{W}^{\j,\infty}_{q_0, \cdots , q_\j} (g) \,,
	\end{aligned}
\end{equation*}
which means that
\begin{equation*}
	\begin{aligned}
		|(\rho_\beta, u_\beta, \theta_\beta)| \leq C (\e_{\j + 1}) \mathcal{W}^{\j,\infty}_{q_0, \cdots , q_\j} (g) \,.
	\end{aligned}
\end{equation*}
Then, as similar as in \eqref{f-null-1}, one has
\begin{equation}\label{f-null-j+1}
	\begin{aligned}
		| \P \partial_{t,x}^\beta f (t,x,v) | \leq C (\e_{\j + 1}) \mathcal{W}^{\j + 1,\infty}_{q_0, \cdots , q_{\j + 1}} (g) \Theta_\gamma (v) \M^\frac{q_{\j + 1}}{2} (v)
	\end{aligned}
\end{equation}
for any $0 < q_{\j + 1} < \cdots < q_0 < 1$ and $|\beta| = \j + 1$. The bounds \eqref{f-perp-j+1} and \eqref{f-null-j+1} thereby imply the bound \eqref{L-inverse-Dtx} in Theorem \ref{Thm2} for the case $m = \j + 1$. Therefore, the Induction Principle finally concludes the results in Part (I) of Theorem \ref{Thm2}.

\subsubsection{\texttt{Proof of Part (II) of Theorem \ref{Thm2}: The cases $- 3 < \gamma \leq - \frac{3}{2}$.}}

Before justifying the results of this subsection, we first give the following lemma.

\begin{lemma}\label{Lmm-SVP-Der}
	Let $- 3 < \gamma \leq - \frac{3}{2}$, $0 < q < 1$, $l \in \R$, $\alpha \in \mathbb{N}^3$ and $\beta \in \mathbb{N}^4$ with $|\beta| = m \geq 1$. Assume $g (t,x,v) \in \mathrm{Null}^\perp (\L)$ satisfies
	\begin{equation*}
		\begin{aligned}
			\widetilde{\mathcal{W}}_{\alpha, \beta}^{\gamma, q, l} (g) : = \sum_{\alpha' \leq \alpha \,, \, \beta' \leq \beta} \| \langle v \rangle^{\gamma (l - \frac{1}{2})} \M^{- \frac{q}{2} } (v) \partial_v^{\alpha'} \partial_{t,x}^{\beta'} g \|^2_{L^2} < \infty \,.
		\end{aligned}
	\end{equation*}
    Then the solution $f (t,x,v) \in \mathrm{Null}^\perp (\L)$ to $\L f = g$ enjoys the bound
    \begin{equation}\label{L2-Dtxv}
    	\begin{aligned}
    		\| \langle v \rangle^{ \gamma ( l + \frac{1}{2} ) } \M^{ - \frac{q}{2} } \partial_v^\alpha \partial_{t,x}^\beta f \|^2_{L^2} \lesssim ( \mathfrak{e}_m )^{2 m^3} \widetilde{\mathcal{W}}_{\alpha, \beta}^{\gamma, q, l} (g) \,,
    	\end{aligned}
    \end{equation}
    where the quantity $\mathfrak{e}_m$ is defined in \eqref{e-m}.
\end{lemma}

The proof of Lemma \ref{Lmm-SVP-Der} will be given in Section \ref{Sec:Lmm4.1} later.

We first take $l$ in \eqref{L2-Dtxv} as $|\alpha| - l$ with $l \geq \frac{5}{2}$ and sum up the resultant inequality for $|\alpha| \leq 2$. Then one obtains
\begin{equation*}
	\begin{aligned}
		\sum_{|\alpha| \leq 2} \| \langle v \rangle^{ \gamma ( l + \frac{1}{2} ) } \M^{ - \frac{q}{2} } \partial_v^\alpha \partial_{t,x}^\beta f \|^2_{L^2} \lesssim ( \mathfrak{e}_m )^{2 m^3} \sum_{|\alpha| \leq 2} \widetilde{\mathcal{W}}_{\alpha, \beta}^{\gamma, q, |\alpha| - l} (g) = ( \mathfrak{e}_m )^{2 m^3} \mathcal{W}_q^{m,2} (g)
	\end{aligned}
\end{equation*}
for $l \geq \frac{5}{2}$. By \eqref{Claim-H2} and \eqref{Embedding}, one knows that
\begin{equation*}
	\begin{aligned}
		\| \M^{ - \frac{q}{2} } \partial_{t,x}^\beta f \|_{L^\infty}^2 \lesssim \| \M^{ - \frac{q}{2} } \partial_{t,x}^\beta f \|^2_{H^2} \lesssim \sum_{|\alpha| \leq 2} \| \langle v \rangle^{ \gamma ( l + \frac{1}{2} ) } \M^{ - \frac{q}{2} } \partial_v^\alpha \partial_{t,x}^\beta f \|^2_{L^2}
	\end{aligned}
\end{equation*}
for $l \geq \frac{5}{2}$. Then
\begin{equation}\label{f-Dtxv-bnd}
	\begin{aligned}
		\| \M^{ - \frac{q}{2} } \partial_{t,x}^\beta f \|_{L^\infty}^2 \lesssim ( \mathfrak{e}_m )^{2 m^3} \sum_{|\alpha| \leq 2} \widetilde{\mathcal{W}}_{\alpha, \beta}^{\gamma, q, |\alpha| - l} (g) \,.
	\end{aligned}
\end{equation}
Note that
\begin{equation*}
	\begin{aligned}
		\mathcal{W}_q^{m,2} (g) = \sum_{|\alpha| \leq 2} \widetilde{\mathcal{W}}_{\alpha, \beta}^{\gamma, q, |\alpha| - \frac{5}{2}} (g) \,,
	\end{aligned}
\end{equation*}
where the quantity $\mathcal{W}_q^{m,2} (g)$ is defined in \eqref{Wm2-g}. It thereby follows the bound \eqref{L-inverse-Dtx-VSP} from taking $l = \frac{5}{2}$ in \eqref{f-Dtxv-bnd}. Then the proof of Part (II) of Theorem \ref{Thm2} is completed.

\subsection{Bounds for $F_n (t,x,v)$: Proof of Theorem \ref{Thm3}}

For any $0 \leq q < 1$ and any fixed integer $m \geq 0$, we first define a set $\mathbf{S}_{\infty m} (q,1)$ as
\begin{equation}\label{S-infty-m-q1}
	\begin{aligned}
		\mathbf{S}_{\infty m} (q,1) = \big\{ (q_{ij})_{\substack{i \geq 1 \\ 0 \leq j \leq m}} \big| \, q < q_{ij} < 1 \,, q_{ij} < q_{i,j-1} (\forall i \geq 1, 0 \leq j \leq m) \,, q_{i+1,0} < q_{im} (\forall i \geq 1) \big\} \,.
	\end{aligned}
\end{equation}
It is easy to see that $\mathbf{S}_{\infty m} (q,1) \neq \emptyset$ for any $0 \leq q < 1$ and integer $m \geq 0$.

Before estimating the bounds of $F_n$, one first derive the following assertion on the functions $f_n = F_n/\sqrt{\M}$ ($n \geq 1$):

\emph{ For any $m \geq 0$, $n \geq 1$ and any $\beta \in \mathbb{N}^4$ with $|\beta| = m$, the following statements hold:}
\begin{itemize}
\item \emph{\underline{Case $- \frac{3}{2} < \gamma \leq 1$:} For any $(q_{ij})_{\substack{i \geq 1 \\ 0 \leq j \leq m}} \in \mathbf{S}_{\infty m} (q,1)$,}
\begin{equation}\label{Assert-fn}
	\begin{aligned}
		\big| \partial_{t,x}^\beta f_n (t,x,v) \big| \leq C (\e_{m+n}) C (\h^i_{m+n-i}; 1 \leq i \leq n) \Theta_\gamma (v) \M^\frac{q_{nm}}{2} (t,x,v) \,;
	\end{aligned}
\end{equation}

\item \emph{\underline{Case $- 3 < \gamma \leq - \frac{3}{2}$:} For any $q \in [0,1)$, $q_0 = q < q_n < q_{n-1} < 1$ $(\forall \, n \geq 1)$, }
\begin{equation}\label{Assert-fn-1}
	\begin{aligned}
		\sum_{ |\alpha| \leq 2 } \| \langle v \rangle^{ \gamma ( |\alpha| - 2 ) } \M^{ - \frac{q_n}{2} } \partial_v^\alpha \partial_{ t, x }^\beta f_n (t,x, \cdot) \|^2_{L^2} \leq C (\e_{m+n}) C (\h^i_{m+n-i}; 1 \leq i \leq n)
	\end{aligned}
\end{equation}
\emph{uniformly in $t,x$.} \emph{Here $\e_{m+n}$ and $\h^i_{m+n-i}$ are given in \eqref{e-m} and \eqref{h-m}, respectively.
}
\end{itemize}

Once the assertions \eqref{Assert-fn} and \eqref{Assert-fn-1} hold, one can conclude the results in Theorem \ref{Thm3}. Indeed, for the case $- \frac{3}{2} < \gamma \leq 1$,
\begin{equation*}
	\begin{aligned}
		| \partial^\beta_{t,x} F_n | = & |\partial_{t,x}^\beta (f_n \sqrt{\M})| \leq \sum_{\beta' \leq \beta} C_\beta^{\beta'} |\partial_{t,x}^{\beta'} f_n \partial_{t,x}^{\beta - \beta'} \sqrt{\M}| \\
		\leq & C (\e_{m+n}) C (\h^i_{m+n-i}; 1 \leq i \leq n) \Theta_\gamma (v) \M^\frac{q_{nm}}{2} (t,x,v) \sum_{\beta' \leq \beta} | \partial_{t,x}^{\beta - \beta'} \sqrt{\M}| \\
		\leq & C (\e_{m+n}) C (\h^i_{m+n-i}; 1 \leq i \leq n) \langle v \rangle^{2m + |\gamma|} \M^\frac{1 + q_{nm}}{2} (t,x,v) \,,
	\end{aligned}
\end{equation*}
where one has utilized the bounds $\Theta_\gamma (v) \leq C \langle v \rangle^{|\gamma|}$ and
\begin{equation*}
	\begin{aligned}
		\sum_{\beta' \leq \beta} | \partial_{t,x}^{\beta - \beta'} \sqrt{\M}| \leq \sum_{\beta' \leq \beta} C (\e_{|\beta - \beta'|}) \langle v \rangle^{2 |\beta - \beta'|} \M^\frac{1}{2} \leq C (\e_m) \langle v \rangle^{2m} \,.
	\end{aligned}
\end{equation*}
From the definition of $\mathbf{S}_{\infty m} (q,1)$ in \eqref{S-infty-m-q1}, one knows that $0 < q < q_{nm} < 1$. Then, by the bound $\langle v \rangle^{2m + |\gamma|} \M^\frac{q_{nm} - q}{2} (v) \leq C (\e_0)$, there holds
\begin{equation}\label{Fn}
	\begin{aligned}
		| \partial^\beta_{t,x} F_n | \leq C (\e_{m+n}) C (\h^i_{m+n-i}; 1 \leq i \leq n) \M^\frac{1 + q}{2} (t,x,v)
	\end{aligned}
\end{equation}
for $- \frac{3}{2} < \gamma \leq 1$.

For the case $- 3 < \gamma \leq - \frac{3}{2}$, together \eqref{Claim-H2} and \eqref{Embedding}, the bound \eqref{Assert-fn-1} implies
\begin{equation*}
	\begin{aligned}
		| \partial_{ t, x }^\beta f_n (t,x,v) | \leq C (\e_{m+n}) C (\h^i_{m+n-i}; 1 \leq i \leq n) \M^\frac{q_n}{2} (t,x,v)
	\end{aligned}
\end{equation*}
for $0 < q < q_n < 1$. Following the similar arguments in the case $- \frac{3}{2} < \gamma \leq 1$ above, one sees that the bound \eqref{Fn} holds for $- 3 < \gamma \leq - \frac{3}{2}$. Namely, the bound \eqref{C-bnd} in Theorem \ref{Thm3} holds for $- 3 < \gamma \leq 1$.

It remains to prove the assertion \eqref{Assert-fn}. As in \eqref{fn-spilt}, the functions $f_n = F_n/\sqrt{\M}$ ($n \geq 1$) are all split into kinetic part $(\I - \P) f_n$ and fluid part $\P f_n$. Based on Figure \ref{Fig3} associated with the orders to solve $f_n$, the assertion \eqref{Assert-fn} can be justified by the induction arguments for $n \geq 1$.

\textbf{\emph{Step 1. $n = 1$.}}

We first estimate the kinetic part $(\I - \P) f_1$ in \eqref{kinetic-n}, and then control the fluid part $\P f_1$ associated with the fluid variables $(\rho_1, u_1, \theta_1)$ in \eqref{CELS}.

\underline{\texttt{Case 1.1: $- \frac{3}{2} < \gamma \leq 1$.}} Observe that, from \eqref{kinetic-n},
\begin{equation}\label{f1-1}
	\begin{aligned}
		(\I - \P) f_1 = \L^{-1} g_1 \,, \ g_1 = - (\I - \P) \big[ \M^{- \frac{1}{2}} (\partial_t + v \cdot \nabla_x ) \M \big] \,.
	\end{aligned}
\end{equation}
By \eqref{nu-beta-1}, one has
\begin{equation*}
	\begin{aligned}
		\partial \M = \big\{ \tfrac{\partial \rho}{\rho} + \partial \u \cdot \tfrac{v - \u}{T} + \tfrac{3 \partial T}{6 T} ( \tfrac{|v - \u|^2}{T} - 3 ) \big\} \M
	\end{aligned}
\end{equation*}
for $\partial = \partial_t, \nabla_x$. Then, together with the definition $\P$ in \eqref{Projt-P}, $g_1$ can be simplified as
\begin{equation}\label{f1-2}
	\begin{aligned}
		g_1 = - \nabla_x \u : \mathscr{A} (v) + \tfrac{\nabla_x T}{\sqrt{T}} \cdot \mathscr{B} (v) \,,
	\end{aligned}
\end{equation}
where the Burnett functions $\mathscr{A} (v)$ and $\mathscr{B} (v)$ are given in \eqref{Burnett}.

We first claim that
\begin{equation}\label{Claim-W-m}
	\begin{aligned}
		\mathcal{W}_{q_{10}, \cdots, q_{1m}}^{m, \infty} (g_1) = \sum_{i-0}^m \sum_{|\beta'| = i} \| \langle v \rangle^{k_\gamma} \M^{- \frac{q_{1i}}{2}} (v) \partial_{t,x}^{\beta'} g_1 \|_{L^\infty} \leq C(\e_{m+1})
	\end{aligned}
\end{equation}
for any $m \geq 0$ and $0 < q_{1m} < \cdots < q_{10} < 1$. Indeed, for $\beta' \in \mathbb{N}^4$ with $|\beta'| = i$,
\begin{equation*}
	\begin{aligned}
		|\partial_{t,x}^{\beta'} g_1| \leq C (\e_{i+1}) p_{2i + 3} (|v-\u|) \M^\frac{1}{2} (v) ,,
	\end{aligned}
\end{equation*}
where $p_{2i + 3}$ is given in \eqref{pn}. Then, one has
\begin{equation*}
	\begin{aligned}
		\| \langle v \rangle^{k_\gamma} \M^{- \frac{q_{1i}}{2}} (v) \partial_{t,x}^{\beta'} g_1 \|_{L^\infty} \leq C (\e_{i+1}) \| \langle v \rangle^{k_\gamma} p_{2i + 3} (|v-\u|) \M^{1 - \frac{q_{1i}}{2}} (v) \partial_{t,x}^{\beta'} \|_{L^\infty} \leq C (\e_{i+1}) < \infty \,,
	\end{aligned}
\end{equation*}
which means that
\begin{equation*}
	\begin{aligned}
		\mathcal{W}_{q_{10}, \cdots, q_{1m}}^{m, \infty} (g_1) \leq \sum_{i=0}^m \sum_{|\beta'| = i} C (\e_{i+1}) \leq C (\e_{m+1}) < \infty \,.
	\end{aligned}
\end{equation*}
Namely, the claim \eqref{Claim-W-m} holds. Based on the claim \eqref{Claim-W-m}, Theorem \ref{Thm1}-\ref{Thm2} imply that, for $|\beta| = m \geq 0$ and $0 < q_{1m} < \cdots < q_{10} < 1$,
\begin{equation}\label{f1-perp}
	\begin{aligned}
		|\partial_{t,x}^\beta (\I - \P) f_1| \leq C (\e_m) \mathcal{W}_{q_{10}, \cdots, q_{1m}}^{m, \infty} (g_1) \Theta_\gamma (v) \M^\frac{q_{1m}}{2} (v) \leq C (\e_{m+1}) \Theta_\gamma (v) \M^\frac{q_{1m}}{2} (v) \,.
	\end{aligned}
\end{equation}

We now estimate the fluid part $\P f_1$, which, by \eqref{fn-spilt}, can be expressed by
\begin{equation}\label{Pf1}
	\begin{aligned}
		\P f_1 = \big\{ \tfrac{\rho_1}{\rho} + u_1 \cdot \tfrac{v - \u}{T} + \tfrac{\theta_1}{6 T} ( \tfrac{|v - \u|^2}{T} - 3 ) \big\} \sqrt{\M} \,.
	\end{aligned}
\end{equation}
Here $(\rho_1, u_1, \theta_1)$ are solved by \eqref{CELS} with source terms $ \mathcal{F}_u^\perp ((\I-\P)f_1)$ and $\mathcal{G}_\theta^\perp ((\I-\P)f_1)$ defined in \eqref{FG-fn}. For any $\beta \in \mathbb{N}^4$ with $|\beta| = m \geq 0$,
\begin{equation}\label{Fu-1}
	\begin{aligned}
		\Big| \partial_{t,x}^\beta  \mathcal{F}_u^\perp ((\I-\P)f_1) \Big| = & \Big| \partial_{t,x}^\beta \big\{ - \div_x \langle T \mathscr{A} , (\I - \P) f_1 \rangle \big\} \Big| \\
		\leq & C (\e_{m+1}) \sum_{i=0}^{m+1} \sum_{|\beta'| = i} \langle p_{2i + 2} (|v - \u|) , | \partial_{t,x}^{\beta'} (\I - \P) f_1 | \rangle \\
		\leq & C (\e_{m+1}) \sum_{i=0}^{m+1} C (\e_{i+1}) \langle p_{2i + 2} (|v - \u|) , \Theta_\gamma (v) \M^\frac{q_{1i}}{2} (v) \rangle \\
		\leq & C (\e_{m+2}) < \infty \,.
	\end{aligned}
\end{equation}
Similarly, for any $\beta \in \mathbb{N}^4$ with $|\beta| = m \geq 0$,
\begin{equation}\label{Gtheta-1}
	\begin{aligned}
		\Big| \partial_{t,x}^\beta \mathcal{G}_\theta^\perp ((\I-\P)f_1) \Big| \leq C (\e_{m+2}) < \infty \,.
	\end{aligned}
\end{equation}
The bounds \eqref{Fu-1} and \eqref{Gtheta-1} yield that the source terms $ \mathcal{F}_u^\perp ((\I-\P)f_1)$ and $\mathcal{G}_\theta^\perp ((\I-\P)f_1)$ are both smooth and bounded provided that $(\rho, \u, T)$ is sufficiently smooth and bounded. Then, by the hypothesis (H2), we can assume $\h_m^1 < \infty$ for sufficiently large $m \geq 0$, where $\h_m^1$ is defined in \eqref{h-m}. Consequently,
\begin{equation}\label{f1-null}
	\begin{aligned}
		\big| \partial_{t,x}^\beta \P f_1 \big| \leq C (\h_m^1) C (\e_m) p_{2m+2} (|v - \u|) \M^\frac{1}{2} (v) \leq C (\h_m^1) C (\e_m) \Theta_\gamma (v) \M^\frac{q_{1m}}{2} (v)
	\end{aligned}
\end{equation}
for any $|\beta| = m \geq 0$, where the uniform bound $p_{2m+2} (|v - \u|) \Theta_\gamma^{-1} (v) \M^\frac{1 - q_{1m}}{2} (v) \leq C (\e_0) $ has been used. Then the bounds \eqref{f1-perp} and \eqref{f1-null} reduce to
\begin{equation}\label{f1-beta}
	\begin{aligned}
		|\partial_{t,x}^\beta f_1| \leq |\partial_{t,x}^\beta (\I - \P) f_1| + |\partial_{t,x}^\beta \P f_1| \leq C (\e_{m+1}) C (\h_m^1) \Theta_\gamma (v) \M^\frac{q_{1m}}{2} (v)
	\end{aligned}
\end{equation}
for any $|\beta| = m \geq 0$ and $0 < q_{1m} < \cdots < q_{10} < 1$.

\underline{\texttt{Case 1.2: $- 3 < \gamma \leq - \frac{3}{2}$.}} By \eqref{L2-Dtxv} in Lemma \ref{Lmm-SVP-Der} with $l = |\alpha| - \frac{5}{2}$ and the equation \eqref{f1-1}, one has
\begin{equation*}
	\begin{aligned}
		& \sum_{ |\alpha| \leq 2 } \| \langle v \rangle^{\gamma ( |\alpha| - 2 )} \M^{ - \frac{q_1}{2} } \partial_v^\alpha \partial_{ t, x }^\beta (\I - \P) f_1 \|^2_{L^2} \lesssim \mathfrak{e}_m^{2 m^3} \sum_{ |\alpha| \leq 2 }  \widetilde{\mathcal{W}}_{\alpha, \beta}^{\gamma, q_1 , |\alpha| - \frac{5}{2}} (g_1) = \mathfrak{e}_m^{2 m^3} \mathcal{W}_{q_1}^{m,2} (g_1) \,.
	\end{aligned}
\end{equation*}
From the expression of $g_1$ in \eqref{f1-2}, it easily follows that
\begin{equation*}
	\begin{aligned}
		\sum_{\alpha' \leq \alpha \,, \beta' \leq \beta} | \partial_v^{\alpha'} \partial_{ t, x }^{\beta'} g_1 | \lesssim \mathfrak{e}_{m+1}^2 \langle v \rangle^{ 3 + |\alpha| + 3m } \M^\frac{1}{2} (v) \,,
	\end{aligned}
\end{equation*}
which means that
\begin{equation*}
	\begin{aligned}
		\mathcal{W}_{q_1}^{m,2} (g_1) = & \sum_{ |\alpha| \leq 2 }  \widetilde{\mathcal{W}}_{\alpha, \beta}^{\gamma, q_1 , |\alpha| - \frac{5}{2}} (g_1) = \sum_{ |\alpha| \leq 2 } \sum_{\alpha' \leq \alpha \,, \beta' \leq \beta} \| \langle v \rangle^{ \gamma ( |\alpha| - 3 ) } \M^{ - \frac{q_1}{2} } \partial_v^{\alpha'} \partial_{ t, x }^{\beta'} g_1 \|^2_{L^2} \\
		\leq & C \mathfrak{e}_{m+1}^4 \| \langle v \rangle^{ \gamma ( |\alpha| - 3 ) + 3 + |\alpha| + 3m } \M^\frac{1}{2} (v) \|^2_{L^2} \leq C ( \mathfrak{e}_{m+1} ) < \infty \,.
	\end{aligned}
\end{equation*}
It thereby holds
\begin{equation}\label{f1-bnd0}
	\begin{aligned}
		\sum_{ |\alpha| \leq 2 } \| \langle v \rangle^{\gamma ( |\alpha| - 2 )} \M^{ - \frac{q_1}{2} } \partial_v^\alpha \partial_{ t, x }^\beta (\I - \P) f_1 \|^2_{L^2} \leq C ( \mathfrak{e}_{m+1} ) \,.
	\end{aligned}
\end{equation}

By \eqref{Pf1}, one easily derives
\begin{equation*}
	\begin{aligned}
		\sum_{ |\alpha| \leq 2 } | \partial_v^\alpha \partial_{ t, x }^\beta (\I - \P) f_1 | \leq C ( \mathfrak{h}_m^1 ) \langle v \rangle^{2 + |\alpha| + 2 m} \M^\frac{1}{2} (v) \,,
	\end{aligned}
\end{equation*}
which means that
\begin{equation}\label{f1-bnd1}
	\begin{aligned}
		\sum_{ |\alpha| \leq 2 } \| \langle v \rangle^{\gamma ( |\alpha| - 2 )} \M^{ - \frac{q_1}{2} } \partial_v^\alpha \partial_{ t, x }^\beta \P f_1 \|^2_{L^2} \leq C ( \mathfrak{h}_m^1 ) \| \langle v \rangle^{ \gamma ( |\alpha| - 2 ) + 2 + |\alpha| + 2 m} \M^\frac{1}{2} \|^2_{L^2} \leq C ( \mathfrak{h}_m^1 ) \,.
	\end{aligned}
\end{equation}
As a result, \eqref{f1-bnd0} and \eqref{f1-bnd1} show that
\begin{equation}
	\begin{aligned}
		\sum_{ |\alpha| \leq 2 } \| \langle v \rangle^{\gamma ( |\alpha| - 2 )} \M^{ - \frac{q_1}{2} } \partial_v^\alpha \partial_{ t, x }^\beta f_1 (t,x, \cdot) \|^2_{L^2} \leq C ( \mathfrak{e}_{m+1} ) C ( \mathfrak{h}_m^1 )
	\end{aligned}
\end{equation}
uniformly in $t,x$.

\textbf{\emph{Step 2. Induction hypotheses for $1 \leq \mathfrak{k} \leq n-1$.}}

Assume that the assertion \eqref{Assert-fn} and \eqref{Assert-fn-1} hold for all $\mathfrak{k} = 1, \cdots, n-1$. Namely, for any $|\beta| = m \geq 0$, $\mathfrak{k} \in \{ 1, \cdots, n-1 \}$, the following statements hold:
\begin{itemize}
	\item \underline{Case $- \frac{3}{2} < \gamma \leq 1$:} For any $(q_{ij})_{\substack{i \geq 1 \\ 0 \leq j \leq m}} \in \mathbf{S}_{\infty m} (q,1)$,
	\begin{equation}\label{Assume-fk-beta}
		\begin{aligned}
			| \partial_{t,x}^\beta f_\mathfrak{k} (t,x,v) | \leq C (\e_{m + \mathfrak{k}}) C(\h^i_{m+ \mathfrak{k} - i}; 1 \leq i \leq \mathfrak{k} ) \Theta_\gamma (v) \M^\frac{q_{\mathfrak{k} m}}{2} (t,x,v) \,.
		\end{aligned}
	\end{equation}
	
	\item \underline{Case $- 3 < \gamma \leq - \frac{3}{2}$:} For any $q \in [0,1)$, $q_0 = q < q_n < q_{n-1} < 1$ $(\forall \, n \geq 1)$,
	\begin{equation}\label{Assume-fk-beta-1}
		\begin{aligned}
			\sum_{ |\alpha| \leq 2 } \| \langle v \rangle^{ \gamma ( |\alpha| - 2 ) } \M^{ - \frac{ q_\mathfrak{k} }{2} } \partial_v^\alpha \partial_{ t, x }^\beta f_\mathfrak{k} (t,x, \cdot) \|^2_{L^2} \leq C (\e_{m + \mathfrak{k}}) C(\h^i_{m+ \mathfrak{k} - i}; 1 \leq i \leq \mathfrak{k} )
		\end{aligned}
	\end{equation}
	holds uniformly in $t,x$.
\end{itemize}

\textbf{\emph{Step 3. General $n \geq 2$.}}

\underline{\texttt{Case 3.1: $- \frac{3}{2} < \gamma \leq 1$.}} From \eqref{kinetic-n}, one has
\begin{equation}\label{f2-I-P}
	\begin{aligned}
		(\I - \P) f_n = \L^{-1} g_n \,, \ g_n = - (\I - \P) \underbrace{\big[ \M^{- \frac{1}{2}} (\partial_t + v \cdot \nabla_x) (f_{n-1} \sqrt{\M}) \big] }_{G_n} + \sum_{\substack{i+j=n \\ i,j \geq 1}} \Gamma (f_i, f_j) \,.
	\end{aligned}
\end{equation}
It is thereby derived from Theorem \ref{Thm1}-\ref{Thm2} that, for any $|\beta| = m \geq 0$ and $0 < q < q_{2m} < \cdots < q_{20} < q_{1m}$,
\begin{equation}\label{f2-I-P-1}
	\begin{aligned}
		|\partial_{t,x}^\beta (\I - \P) f_n| \leq C(\e_m) \mathcal{W}^{m, \infty}_{q_{n0},\cdots, q_{nm}} (g_n) \Theta_\gamma (v) \M^\frac{q_{nm}}{2} (v)
	\end{aligned}
\end{equation}
It remains to show that $\mathcal{W}^{m, \infty}_{q_{n0},\cdots, q_{nm}} (g_n) = \sum_{\mathfrak{r} = 0}^m \sum_{|\beta'| = \mathfrak{r}} \| \langle v \rangle^{k_\gamma} \M^{- \frac{q_{n \mathfrak{r}}}{2}} (v) \partial_{t,x}^{\beta'} g_n \|_{L^\infty} < \infty$.
A straightforward calculation implies that
\begin{equation}\label{G2-1}
	\begin{aligned}
		G_n = (\partial_t + v \cdot \nabla_x) f_{n-1} - \big[ \tfrac{(\partial_t + v \cdot \nabla_x) \rho}{\rho} + (\partial_t + v \cdot \nabla_x) \u \cdot \tfrac{v - \u}{T} + \tfrac{(\partial_t + v \cdot \nabla_x) T}{6 T} ( \tfrac{|v - \u|^2}{T} - 3 ) \big] f_{n-1} \,.	
	\end{aligned}
\end{equation}
Note that
\begin{equation}\label{G2-1-I-P}
	\begin{aligned}
		| \partial_{t,x}^{\beta'} (\I - \P) G_n | \leq | \partial_{t,x}^{\beta'} G_n | + | \partial_{t,x}^{\beta'} \P G_n | \,.
	\end{aligned}
\end{equation}
By \eqref{G2-1}, it is easy to know that
\begin{equation}\label{G2-1-beta}
	\begin{aligned}
	    | \partial_{t,x}^{\beta'} G_n | \leq & \langle v \rangle \sum_{|\beta''| = m + 1} |\partial_{t,x}^{\beta''} f_{n-1}| + C (\e_{m+1}) p_3 (|v|) \sum_{\beta'' \leq \beta'} |\partial_{t,x}^{\beta''} f_{n-1}| \\
	    \leq & C (\e_{m+1}) p_3 (|v|) \sum_{|\beta''| \leq m+1} |\partial_{t,x}^{\beta''} f_{n-1}| \\
	    \leq & C (\e_{m+n}) C (\h^\mathfrak{p}_{m+n-\mathfrak{p}}; 1 \leq \mathfrak{p} \leq n-1) p_3 (|v|) \sum_{\mathfrak{r} = 0}^{m+1} \Theta_\gamma (v) \M^\frac{q_{n-1,\mathfrak{r}}}{2} (v) \\
	    \leq & C (\e_{m+n}) C (\h^\mathfrak{p}_{m+n-\mathfrak{p}}; 1 \leq \mathfrak{p} \leq n-1) \Theta_\gamma (v) \M^\frac{q_{n-1,m+1}}{2} (v) \,,
	\end{aligned}
\end{equation}
where $q_{n-1,m+1} \in (q_{n0}, q_{n-1,m})$, the last second inequality is derived from \eqref{Assume-fk-beta}, and the last one is implied by the relation $0 < q < q_{1,m+1} < q_{1m} < \cdots < q_{10} < 1$.

By the definition of $\P f$ in \eqref{Projt-P}, one has
\begin{equation*}
	\begin{aligned}
		\P G_n = \big\{ \tfrac{\rho_G}{\rho} + u_G \cdot \tfrac{v - \u}{T} + \tfrac{\theta_G}{6 T} ( \tfrac{|v-\u|^2}{T} - 3 ) \big\} \sqrt{\M} \,,
	\end{aligned}
\end{equation*}
where $\rho_G = \langle G_n, \sqrt{\M} \rangle $, $u_G = \langle G_n, \tfrac{1}{\rho} (v - \u) \sqrt{\M} \rangle$ and $\theta_G = \langle G_n, \tfrac{1}{\rho} (|v - \u|^2 - 3 T) \sqrt{\M} \rangle$. Then there holds
\begin{equation*}
	\begin{aligned}
		| \partial_{t,x}^{\beta'} \P G_n | \leq C (\e_{|\beta'|}) p_{2 |\beta'| + 2} (|v|) \sum_{\beta'' \leq \beta'} \big| \partial_{t,x}^{\beta''} ( \rho_G, u_G, \theta_G ) \big| \M^\frac{1}{2} (v) \,.
	\end{aligned}
\end{equation*}
Consider the quantities $\partial_{t,x}^{\beta''} \langle G_n , \phi \rangle $, where $\phi$ is given in \eqref{phi}. To be precise, one has
\begin{equation*}
	\begin{aligned}
		\partial_{t,x}^{\beta''} \langle G_n , \phi \rangle = \sum_{\beta^* \leq \beta''} C_{\beta''}^{\beta^*} \langle \partial_{t,x}^{\beta^*} G_n , \partial_{t,x}^{\beta'' - \beta^*} \phi \rangle \,.
	\end{aligned}
\end{equation*}
Noting that $|\beta^*| \leq |\beta''| \leq |\beta'| \leq m$ and $|\beta'' - \beta^*| \leq |\beta''| \leq m$, one has
$$| \partial_{t,x}^{\beta'' - \beta^*} \phi | \leq C (\e_{|\beta'' - \beta^*|}) p_{2 |\beta'' - \beta^*| + 3} (|v|) \M^\frac{1}{2} (v) \leq C (\e_m) p_{2m + 3} (|v|) \M^\frac{1}{2} (v) \,. $$
Together with \eqref{G2-1-beta}, we have
\begin{equation*}
	\begin{aligned}
		\big| \partial_{t,x}^{\beta''} \langle G_n , \phi \rangle \big| \leq & C (\e_{m+n}) C (\h^\mathfrak{p}_{m+n-\mathfrak{p}}; 1 \leq \mathfrak{p} \leq n-1) \\
		& \qquad \times \langle \Theta_\gamma (v) \M^\frac{q_{n-1,m+1}}{2} (v), p_{2m + 3} (|v|) \M^\frac{1}{2} (v) \rangle \\
		\leq & C (\e_{m+n}) C (\h^\mathfrak{p}_{m+n-\mathfrak{p}}; 1 \leq \mathfrak{p} \leq n-1) \,.
	\end{aligned}
\end{equation*}
This means that
\begin{equation*}
	\begin{aligned}
		\sum_{\beta'' \leq \beta'} \big| \partial_{t,x}^{\beta''} ( \rho_G, u_G, \theta_G ) \big| \leq C (\e_{m+n}) C (\h^\mathfrak{p}_{m+n-\mathfrak{p}}; 1 \leq \mathfrak{p} \leq n-1) < \infty \,.
	\end{aligned}
\end{equation*}
Consequently, it infers that
\begin{equation}\label{G2-P-1-beta}
	\begin{aligned}
		| \partial_{t,x}^{\beta'} \P G_n | \leq & C (\e_{|\beta'|}) p_{2 |\beta'| + 2} (|v|) C (\e_{m+n}) C (\h^\mathfrak{p}_{m+n-\mathfrak{p}}; 1 \leq \mathfrak{p} \leq n-1) \M^\frac{1}{2} (v) \\
		\leq & C (\e_{m+n}) C (\h^\mathfrak{p}_{m+n-\mathfrak{p}}; 1 \leq \mathfrak{p} \leq n-1) \Theta_\gamma (v) \M^\frac{q_{n-1,m+1}}{2} (v) \,.
	\end{aligned}
\end{equation}
Therefore, the bounds \eqref{G2-1-I-P}, \eqref{G2-1-beta} and \eqref{G2-P-1-beta} imply that, for any $0 \leq i \leq m$,
\begin{equation}\label{G2-1-I-P-bnd}
	\begin{aligned}
		\sum_{|\beta'| = \mathfrak{r}} | \partial_{t,x}^{\beta'} (\I - \P) G_n | \leq C (\e_{m+n}) C (\h^\mathfrak{p}_{m+n-\mathfrak{p}}; 1 \leq \mathfrak{p} \leq n-1) \Theta_\gamma (v) \M^\frac{q_{n-1,m+1}}{2} (v) \,.
	\end{aligned}
\end{equation}

Next, we will control the quantities $ \partial_{t,x}^{\beta'} \Gamma (f_i, f_j)$ for $| \beta' | = \mathfrak{r}$ with $0 \leq \mathfrak{r} \leq m$ and $i+j = n$ with $i,j \geq 1$. Recalling the definition of $\Gamma (f_i, f_j)$ as in \eqref{Gamma-fi-fj}, one has
\begin{equation}\label{Gamma-split}
	\begin{aligned}
		\Gamma (f_i, f_j) = & \M^{- \frac{1}{2}} (v) \B (f_i \sqrt{\M}, f_j \sqrt{\M} ) \\
		= & \M^{- \frac{1}{2}} (v) \iint_{\R^3 \times \mathbb{S}^2} \big[ f_i (v') \M^\frac{1}{2} (v') f_j (v_1') \M^\frac{1}{2} (v'_1) - f_i (v) \M^\frac{1}{2} (v) f_j (v_1) \M^\frac{1}{2} (v_1) \big] \\
		& \qquad \qquad \qquad \qquad \qquad \qquad \times b (\omega, v_1 - v) \d \omega \d v_1 \\
		= & \underbrace{ \iint_{\R^3 \times \mathbb{S}^2} f_i (v') f_j (v_1') \M^\frac{1}{2} (v_1) b (\omega, v_1 - v) \d \omega \d v_1 }_{\Gamma_{gain} (f_i, f_j)} \\
		& - \underbrace{ \iint_{\R^3 \times \mathbb{S}^2} f_i (v) f_j (v_1) \M^\frac{1}{2} (v_1) b (\omega, v_1 - v) \d \omega \d v_1 }_{\Gamma_{loss} (f_i, f_j)} \,.
	\end{aligned}
\end{equation}
where the equality $\M^\frac{1}{2} (v') \M^\frac{1}{2} (v_1') = \M^\frac{1}{2} (v) \M^\frac{1}{2} (v_1)$ has been used. Then, for any $|\beta'| = \mathfrak{r}$ with $0 \leq \mathfrak{r} \leq m$ , and any $i+j=n$ with $i,j \geq 1$,
\begin{equation*}
	\begin{aligned}
		\partial_{t,x}^{\beta'} \Gamma_{gain} (f_i, f_j) = \sum_{\beta'' \leq \beta'} \sum_{\beta^* \leq \beta''} C_{\beta'}^{\beta''} C_{\beta''}^{\beta^*} \iint_{\R^3 \times \mathbb{S}^2} \partial_{t,x}^{\beta^*} f_i (v') \partial_{t,x}^{\beta'' - \beta^*} f_j (v_1') \partial_{t,x}^{\beta' - \beta''} \M^\frac{1}{2} (v_1) \\
		\times b(\omega, v_1 - v) \d \omega \d v_1 \,.
	\end{aligned}
\end{equation*}
Observe that
\begin{equation*}
	\begin{aligned}
		| \partial_{t,x}^{\beta' - \beta''} \M^\frac{1}{2} (v_1) | \leq C (\e_m) p_{2 |\beta' - \beta''|} (|v_1|) \M^\frac{1}{2} (v_1) \leq C (\e_m) p_{2 m} (|v_1|) \M^\frac{1}{2} (v_1) \,,
	\end{aligned}
\end{equation*}
where $p_{2m} $ is given in \eqref{pn}. Together with \eqref{f1-beta} and \eqref{Assume-fk-beta}, it follows that
\begin{equation*}
	\begin{aligned}
		& |\partial_{t,x}^{\beta^*} f_i (v')| \leq C (\e_{m+i}) C (\h^\mathfrak{p}_{m+i - \mathfrak{p}}; 1 \leq \mathfrak{p} \leq i) \Theta_\gamma (v') \M^\frac{q_{im}}{2} (v') \,, \\
		& | \partial_{t,x}^{\beta'' - \beta^*} f_j (v_1') | \leq C (\e_{m+j}) C (\h^\mathfrak{p}_{m+j - \mathfrak{p}}; 1 \leq \mathfrak{p} \leq j) \Theta_\gamma (v'_1) \M^\frac{q_{jm}}{2} (v'_1)
	\end{aligned}
\end{equation*}
for any $i + j = n$ with $i,j \geq 1$. Noticing that $q_{im} \geq q_{n-1, m} > q > 0$ and $q_{jm} \geq q_{n-1,m} > q > 0$, one therefore obtains that, for any fixed $\varpi \in ( 0, \tfrac{q_{n-1,m} - q_{n-1, m+1}}{2} )$,
\begin{equation}\label{Gamma-gain}
	\begin{aligned}
		\big| \partial_{t,x}^{\beta'} \Gamma_{gain} (f_i, f_j) \big| \leq & C (\e_{m+n-1}) C (\h_{m+n-1-\mathfrak{p}}^\mathfrak{p}; 1 \leq \mathfrak{p} \leq n - 1) \iint_{\R^3 \times \mathbb{S}^2} \Theta_\gamma (v') \Theta_\gamma (v_1') \\
		& \qquad \times \big[ \M (v') \M (v_1') \big]^\frac{q_{n-1,m}}{2} p_{2 m} (|v_1|) \M^\frac{1}{2} (v_1) b(\omega, v_1 - v) \d \omega \d v_1 \\
		\leq & C (\e_{m+n-1}) C (\h_{m+n-1-\mathfrak{p}}^\mathfrak{p}; 1 \leq \mathfrak{p} \leq n - 1) \\
		& \qquad \times \iint_{\R^3 \times \mathbb{S}^2} \big[ \M (v) \M (v_1) \big]^{\frac{q_{n-1,m}}{2} - \varpi} b(\omega, v_1 - v) \d \omega \d v_1 \\
		\leq & C (\e_{m+n-1}) C (\h_{m+n-1-\mathfrak{p}}^\mathfrak{p}; 1 \leq \mathfrak{p} \leq n - 1) \nu (v) \M^{\frac{q_{n-1,m}}{2} - \varpi} (v) \\
		\leq & C (\e_{m+n-1}) C (\h_{m+n-1-\mathfrak{p}}^\mathfrak{p}; 1 \leq \mathfrak{p} \leq n - 1) \M^{\frac{q_{n-1,m+1}}{2} } (v) \,,
	\end{aligned}
\end{equation}
where the last second inequality is derived from the relations $\M (v') \M (v_1') = \M (v) \M (v_1)$, Lemma 2.3 of \cite{LS-2010-KRM} and \eqref{nu-equivalent}, the last inequality is implied by the bound
$$\nu (v) \M^{ \frac{q_{n-1,m} - q_{n-1,m+1}}{2} - \varpi } (v) \leq C (\e_0) < \infty \,.$$
Moreover, by similar arguments of \eqref{Gamma-gain},
\begin{equation}\label{Gamma-loss}
	\begin{aligned}
		\big| \partial_{t,x}^{\beta'} \Gamma_{loss} (f_i, f_j) \big| \leq & C (\e_{m+n-1}) C (\h_{m+n-1-\mathfrak{p}}^\mathfrak{p}; 1 \leq \mathfrak{p} \leq n - 1) \iint_{\R^3 \times \mathbb{S}^2} \Theta_\gamma (v) \Theta_\gamma (v_1) \\
		& \qquad \times \big[ \M (v) \M (v_1) \big]^\frac{q_{n-1,m}}{2} p_{2 m} (|v_1|) \M^\frac{1}{2} (v_1) b(\omega, v_1 - v) \d \omega \d v_1 \\
		\leq & C (\e_{m+n-1}) C (\h_{m+n-1-\mathfrak{p}}^\mathfrak{p}; 1 \leq \mathfrak{p} \leq n - 1) \M^{\frac{q_{n-1,m+1}}{2} } (v) \,.
	\end{aligned}
\end{equation}
Consequently, it follow from \eqref{Gamma-gain} and \eqref{Gamma-loss} that, for any $q_{n-1,m+1} \in (q_{n0}, q_{n-1,m})$, any $\beta' \in \mathbb{N}^4$ with $0 \leq |\beta'| = \mathfrak{r} \leq m$ and any $i+j = n$ with $i,j \geq 1$,
\begin{equation}\label{Gamma-beta-bnd}
	\begin{aligned}
		\big| \partial_{t,x}^{\beta'} \Gamma (f_i, f_j) \big| \leq C (\e_{m+n-1}) C (\h_{m+n-1-\mathfrak{p}}^\mathfrak{p}; 1 \leq \mathfrak{p} \leq n - 1) \M^{\frac{q_{n-1,m+1}}{2} } (v) \,.
	\end{aligned}
\end{equation}
It is thereby derived from \eqref{f2-I-P}, \eqref{G2-1-I-P-bnd} and \eqref{Gamma-beta-bnd} that, for any $0 < q < q_{nm} < \cdots < q_{n0} < q_{n-1,m+1} < q_{n-1,m} $,
\begin{equation}\label{f2-I-P-2}
	\begin{aligned}
		\mathcal{W}^{m, \infty}_{q_{n0},\cdots, q_{nm}} (g_n) = & \sum_{\mathfrak{r} = 0}^m \sum_{|\beta'| = \mathfrak{r}} \| \langle v \rangle^{k_\gamma} \M^{- \frac{q_{n \mathfrak{r}}}{2}} (v) \partial_{t,x}^{\beta'} g_n \|_{L^\infty} \\
		\leq & \sum_{\mathfrak{r} = 0}^m \sum_{|\beta'| = \mathfrak{r}} \Big\{ \| \langle v \rangle^{k_\gamma} \M^{- \frac{q_{n \mathfrak{r}}}{2}} (v) \partial_{t,x}^{\beta'} (\I - \P) G_n \|_{L^\infty} \\
		& \qquad \qquad + \sum_{\substack{i+j=n \\ i,j \geq 1}} \| \langle v \rangle^{k_\gamma} \M^{- \frac{q_{n \mathfrak{r}}}{2}} (v) \partial_{t,x}^{\beta'} \Gamma (f_i, f_j) \|_{L^\infty} \Big\} \\
		\leq & C (\e_{m+n}) C (\h_{m+n-\mathfrak{p}}^\mathfrak{p}; 1 \leq \mathfrak{p} \leq n - 1) \| \langle v \rangle^{k_\gamma} \Theta_\gamma (v) \M^\frac{q_{n-1,m+1} - q_{n0}}{2} \|_{L^\infty} \\
		\leq & C (\e_{m+n}) C (\h_{m+n-\mathfrak{p}}^\mathfrak{p}; 1 \leq \mathfrak{p} \leq n - 1) < \infty \,.
	\end{aligned}
\end{equation}
Therefore, \eqref{f2-I-P-1} and \eqref{f2-I-P-2} show us that, for any $|\beta| = m \geq 0$,
\begin{equation}\label{f2-perp}
	\begin{aligned}
		|\partial_{t,x}^\beta (\I - \P) f_n| \leq C (\e_{m+n}) C (\h_{m+n-\mathfrak{p}}^\mathfrak{p}; 1 \leq \mathfrak{p} \leq n - 1) \Theta_\gamma (v) \M^\frac{q_{nm}}{2} (v) \,.
	\end{aligned}
\end{equation}

We now estimate the fluid part $\P f_n$, which, by the expression of $\P f$ in \eqref{Projt-P}, satisfies
\begin{equation*}
	\begin{aligned}
		\P f_n = \big\{ \tfrac{\rho_n}{\rho} + u_n \cdot \tfrac{v - \u}{T} + \tfrac{\theta_n}{6 T} ( \tfrac{|v - \u|^2}{T} ) \big\} \sqrt{\M} \,.
	\end{aligned}
\end{equation*}
Here $(\rho_n, u_n, \theta_n)$ is determined by the linearized compressible Euler system with the source terms $ \mathcal{F}_u^\perp ((\I-\P)f_n)$ and $\mathcal{G}_\theta^\perp ((\I-\P)f_n)$ defined in \eqref{FG-fn}. Following the same arguments in \eqref{Fu-1} and \eqref{Gtheta-1}, one has
\begin{equation}\label{Fu-Gtheta-2}
	\begin{aligned}
		|\partial_{t,x}^\beta \mathcal{F}_u^\perp ((\I-\P)f_n)| + |\partial_{t,x}^\beta \mathcal{G}_\theta^\perp & ((\I-\P) f_n)| \\
		& \leq C (\e_{m+n+1}) C (\h_{m+n+1-\mathfrak{p}}^\mathfrak{p}; 1 \leq \mathfrak{p} \leq n - 1) < \infty
	\end{aligned}
\end{equation}
for any $\beta \in \mathbb{N}^4$ with $|\beta | = m \geq 0$. The bound \eqref{Fu-Gtheta-2} yields that the source terms $ \mathcal{F}_u^\perp ((\I-\P)f_n)$ and $\mathcal{G}_\theta^\perp ((\I-\P)f_n)$ are both smooth and bounded provided that $(\rho, \u, T)$ is sufficiently smooth and bounded. Then, by the hypothesis (H2), one can assume $\h_m^n < \infty$ for large $m \geq 0$, where $\h_m^2$ is given in \eqref{h-m}. Then, it follows from the same arguments in \eqref{f1-null} that
\begin{equation}\label{f2-null}
	\begin{aligned}
		|\partial_{t,x}^\beta \P f_n| \leq C (\e_m) C(\h_m^n) \Theta_\gamma (v) \M^\frac{q_{nm}}{2} (v)
	\end{aligned}
\end{equation}
for any $\beta \in \mathbb{N}^4$ with $|\beta| = m \geq 0$. Collecting the bounds \eqref{f2-perp} and \eqref{f2-null}, one has
\begin{equation}\label{fn-beta}
	\begin{aligned}
		|\partial_{t,x}^\beta f_n| \leq C (\e_{m+n}) C(\h_{m+n-\mathfrak{p}}^\mathfrak{p}; 1 \leq \mathfrak{p} \leq n) \Theta_\gamma (v) \M^\frac{q_{nm}}{2} (v)
	\end{aligned}
\end{equation}
for any $\beta \in \mathbb{N}^4$ with $|\beta| = m \geq 0$ and $0 < q < q_{nm} < \cdots < q_{n0} < q_{n-1,m} $.

By the Induction Principle, \eqref{fn-beta} concludes the assertion \eqref{Assert-fn}.

\underline{\texttt{Case 3.2: $ - 3 < \gamma \leq - \frac{3}{2}$.}} We first investigate the bound of $(\I - \P) f_n$. By \eqref{L2-Dtxv} in Lemma \ref{Lmm-SVP-Der},
\begin{equation*}
	\begin{aligned}
		& \sum_{ |\alpha| \leq 2 } \| \langle v \rangle^{ \gamma ( |\alpha| - 2 ) } \M^{ - \frac{q_n}{2} } \partial_v^\alpha \partial_{ t, x }^\beta ( \I- \P ) f_n \|^2_{L^2} \\
		\lesssim & (\mathfrak{e}_m)^{2 m^3} \sum_{ |\alpha| \leq 2 } \widetilde{\mathcal{W}}_{\alpha, \beta}^{\gamma, q_n, |\alpha| - \frac{5}{2}} (g_n) = (\mathfrak{e}_m)^{2 m^3} {\mathcal{W}}_{q_n}^{m, 2} (g_n) \,,
	\end{aligned}
\end{equation*}
where $g_n$ is given in \eqref{f2-I-P}. It remains to show that
\begin{equation}\label{W-0}
	\begin{aligned}
		{\mathcal{W}}_{q_n}^{m, 2} (g_n) \leq C ( \mathfrak{e}_{m+n} ) C ( \h_{ m + n - \mathfrak{p}}^\mathfrak{p}; 1 \leq \mathfrak{p} \leq n) \,.
	\end{aligned}
\end{equation}
Recall that
\begin{equation}\label{W-1}
	\begin{aligned}
		{\mathcal{W}}_{q_n}^{m, 2} (g_n) = \sum_{ |\alpha| \leq 2 } \sum_{ \alpha' \leq \alpha, \beta' \leq \beta } \| \langle v \rangle^{ \gamma ( |\alpha| - 3 ) } \M^{ - \frac{q_n}{2} } \partial_v^{\alpha'} \partial_{ t, x }^{\beta'} g_n \|^2_{L^2} \,.
	\end{aligned}
\end{equation}
By the induction hypothesis \eqref{Assume-fk-beta-1} in Step 2, one has
\begin{equation*}
	\begin{aligned}
		\sum_{ |\alpha| \leq 2 } \| \langle v \rangle^{ \gamma ( |\alpha| - 2 ) } \M^{ - \frac{q}{2} } \partial_v^\alpha \partial_{ t, x }^\beta f_\mathfrak{k} (t,x, \cdot) \|^2_{L^2} \leq C (\e_{m + \mathfrak{k}}) C(\h^i_{m+ \mathfrak{k} - i}; 1 \leq i \leq \mathfrak{k} )
	\end{aligned}
\end{equation*}
for $\mathfrak{k} \in \{ 1, \cdots, n-1 \}$. By using \eqref{Claim-H2} and \eqref{Embedding}, it holds
\begin{equation*}
	\begin{aligned}
		\| \M^{ - \frac{q}{2} } f_\mathfrak{k} (t,x, \cdot) \|_{L^\infty} \leq & C \sum_{ |\alpha| \leq 2 } \| \langle v \rangle^{ \gamma ( |\alpha| - 2 ) } \M^{ - \frac{q}{2} } \partial_v^\alpha \partial_{ t, x }^\beta f_\mathfrak{k} (t,x, \cdot) \|^2_{L^2} \\
		\leq & C (\e_{m + \mathfrak{k}}) C(\h^i_{m+ \mathfrak{k} - i}; 1 \leq i \leq \mathfrak{k} ) \,,
	\end{aligned}
\end{equation*}
which means that
\begin{equation*}
	\begin{aligned}
		| f_\mathfrak{k} (t, x, v) | \leq C (\e_{m + \mathfrak{k}}) C(\h^i_{m+ \mathfrak{k} - i}; 1 \leq i \leq \mathfrak{k} ) \M^\frac{ q_\mathfrak{k} }{2} (t,x,v)
	\end{aligned}
\end{equation*}
for all $\mathfrak{k} \in \{ 1, \cdots, n-1 \}$. It then follows from the similar arguments in \eqref{G2-1-I-P-bnd} and \eqref{Gamma-beta-bnd} that
\begin{equation}\label{W-2}
	\begin{aligned}
		\sum_{ \alpha' \leq \alpha , \beta' \leq \beta } | \partial_v^{\alpha'} \partial_{ t, x }^{\beta'} g_n (t,x,v) | \leq C ( \mathfrak{e}_{m+n} ) C ( \h_{ m + n - \mathfrak{p}}^\mathfrak{p}; 1 \leq \mathfrak{p} \leq n-1) \langle v \rangle^{N_0} \M^\frac{q_{n-1}}{2} (t,x,v)
	\end{aligned}
\end{equation}
for some large constant $N_0 > 0$. Then \eqref{W-1} and \eqref{W-2} indicate that
\begin{equation*}
	\begin{aligned}
		{\mathcal{W}}_{q_n}^{m, 2} (g_n) \leq C ( \mathfrak{e}_{m+n} ) C ( \h_{ m + n - \mathfrak{p}}^\mathfrak{p}; 1 \leq \mathfrak{p} \leq n-1) \sum_{ |\alpha| \leq 2 } \| \langle v \rangle^{ \gamma ( |\alpha| - 3 ) + N_0} \M^\frac{ q_{n-1} - q_n }{2} \|^2_{L^2} \,.
	\end{aligned}
\end{equation*}
Note that $q_{n-1} - q_n > 0$. One easily has
\begin{equation*}
	\begin{aligned}
		\| \langle v \rangle^{ \gamma ( |\alpha| - 3 ) + N_0} \M^\frac{q_{n-1}}{2} \M^\frac{ q_{n-1} - q_n }{2} \|^2_{L^2} \leq C ( \mathfrak{e}_0 ) \,.
	\end{aligned}
\end{equation*}
One therefore concludes that the assertion \eqref{W-0} holds. Namely,
\begin{equation}\label{N-1}
	\begin{aligned}
		\sum_{ |\alpha| \leq 2 } \| \langle v \rangle^{ \gamma ( |\alpha| - 2 ) } \M^{ - \frac{q_n}{2} } \partial_v^\alpha \partial_{ t, x }^\beta ( \I- \P ) f_n \|^2_{L^2} \leq C ( \mathfrak{e}_{m+n} ) C ( \h_{ m + n - \mathfrak{p}}^\mathfrak{p}; 1 \leq \mathfrak{p} \leq n-1) \,.
	\end{aligned}
\end{equation}

Furthermore, by the similar arguments in \eqref{f2-null}, one has
\begin{equation*}
	\begin{aligned}
		| \partial_v^{\alpha} \partial_{ t, x }^{\beta} \P f_n (t,x,v) | \leq C (\mathfrak{e}_m) C (\mathfrak{h}_m^n) \langle v \rangle^{ 2 + |\alpha| + 2 |\beta| } \M^\frac{1}{2} (t,x,v) \,,
	\end{aligned}
\end{equation*}
which follows
\begin{equation}\label{N-2}
	\begin{aligned}
		& \sum_{ |\alpha| \leq 2 } \| \langle v \rangle^{ \gamma ( |\alpha| - 2 ) } \M^{ - \frac{q_n}{2} } \partial_v^\alpha \partial_{ t, x }^\beta \P f_n \|^2_{L^2} \\
		\leq & C (\mathfrak{e}_m) C (\mathfrak{h}_m^n) \sum_{ |\alpha| \leq 2 } \| \langle v \rangle^{ \gamma ( |\alpha| - 2 ) + 2 + |\alpha| + 2 |\beta| } \M^{ \frac{1 - q_n}{2} } \|^2_{L^2} \\
		\leq & C (\mathfrak{e}_m) C (\mathfrak{h}_m^n) \,.
	\end{aligned}
\end{equation}
The bounds \eqref{N-1} and \eqref{N-2} show that
\begin{equation*}
	\begin{aligned}
		\sum_{ |\alpha| \leq 2 } \| \langle v \rangle^{ \gamma ( |\alpha| - 2 ) } \M^{ - \frac{q_n}{2} } \partial_v^\alpha \partial_{ t, x }^\beta ( \I- \P ) f_n \|^2_{L^2} \leq C ( \mathfrak{e}_{m+n} ) C ( \h_{ m + n - \mathfrak{p}}^\mathfrak{p}; 1 \leq \mathfrak{p} \leq n) \,.
	\end{aligned}
\end{equation*}
Together with \eqref{Assume-fk-beta-1}, the Induction Principle tells us that the assertion \eqref{Assert-fn-1} holds. Therefore, the proof of Theorem \ref{Thm3} is completed.

\section{Low-High velocities estimates of $K$ for $- \frac{3}{2} < \gamma \leq 1$: Proof of Lemma \ref{Lmm-K1K2-Hypo1} and \ref{Lmm-K1K2-Hypo2}}\label{Sec:K-LHVE-K}

In this section, we will separately derive the low and high velocities estimates of the operator $K$ for the cases $- \frac{3}{2} < \gamma \leq 1$, hence, prove Lemma \ref{Lmm-K1K2-Hypo1} and \ref{Lmm-K1K2-Hypo2}.

\subsection{Low velocities estimates for $K$: Proof of Lemma \ref{Lmm-K1K2-Hypo1}}

	Observe that $K^{1-\chi} f = K_1^{1 - \chi} f - K_2^{1 - \chi} f$. It suffices to prove that
	\begin{equation}\label{K12}
		\begin{aligned}
			\big| \langle \Theta_\gamma \M^{- \frac{q}{2}} K_i^{1 - \chi} f , \M^{- \frac{q}{2}} f \rangle \big| \lesssim_{\rho, \u, T, r} \| f \|^2_{L^2 (\nu)}
		\end{aligned}
	\end{equation}
	for $i = 1,2$.
	
	\vspace*{3mm}
	
	{\em Step 1. Control the quantity $\big| \langle \Theta_\gamma \M^{- \frac{q}{2}} K_1^{1 - \chi} f , \M^{- \frac{q}{2}} f \rangle \big|$.}
	
	\vspace*{3mm}
	
	By the definition of $\chi (s)$ in \eqref{Chi}, one knows that $1 - \chi (|v_1 - v|) \equiv 0$ if $|v_1 - v| \geq 2 r$. Recalling the expression of $K_1^{1-\chi} f$ in \eqref{K1-Lambda}, one has
	\begin{equation}\label{K1-1-chi-1}
		\begin{aligned}
			\big| \langle \Theta_\gamma \M^{- \frac{q}{2}} K_1^{1 - \chi} f , \M^{- \frac{q}{2}} f \rangle \big| \leq \iiint_{\mathbb{R}^3 \times \R^3 \times \mathbb{S}^2} [1 - & \chi (|v_1 - v|) ] \Theta_\gamma (v) \M^\frac{1}{2} (v) \M^\frac{1}{2} (v_1) \M^{-q} (v) \\
			& \times | f (v_1) f(v) | b (\omega, v_1 - v) \d \omega \d v_1 \d v \,.
		\end{aligned}
	\end{equation}
	Note that, for $|v_1 - v| \leq 2 r$,
	\begin{equation*}
		\begin{aligned}
			|v - \u|^2 = & |v_1 - v|^2 - 2 (v_1 - v) \cdot (v_1 - \u) + |v_1 - \u|^2 \\
			\leq & (1 + \eps_1) |v_1 - \u|^2 + (1 + \eps_1^{-1}) |v_1 - v|^2 \\
			\leq & (1 + \eps_1) |v_1 - \u|^2 + 4 (1 + \eps_1^{-1}) r^2 \,,
		\end{aligned}
	\end{equation*}
	where the small constant $\eps_1 > 0$ is to be determined. Then, for the undetermined $0 < q' \leq q$,
	\begin{equation*}
		\begin{aligned}
			\M^{-q'} (v) = & \big[ \tfrac{\rho}{(2 \pi T)^\frac{3}{2}} \big]^{-q'} \exp (\tfrac{q' |v-\u|^2}{2 T}) \\
			\leq & \big[ \tfrac{\rho}{(2 \pi T)^\frac{3}{2}} \big]^{-q'} \exp \big[ \tfrac{q' (1 + \eps_1) |v_1 - \u|^2 + 4 q' (1 + \eps_1^{-1}) r^2}{2 T} \big] \\
			= & \big[ \tfrac{\rho}{(2 \pi T)^\frac{3}{2}} \big]^{q' \eps_1} \exp \big[ \tfrac{2 (1 + \eps_1^{-1}) r^2}{T} \big] \M^{- q' (1 + \eps_1)} (v_1) \,.
		\end{aligned}
	\end{equation*}
	It therefore infers that, if $|v_1 - v| \leq 2 r$,
	\begin{equation}\label{K1-1-chi-2}
		\begin{aligned}
			\M^\frac{1}{2} (v) \M^\frac{1}{2} (v_1) \M^{-q} (v) \leq \big[ \tfrac{\rho}{(2 \pi T)^\frac{3}{2}} \big]^{q' \eps_1} \exp \big[ \tfrac{2 (1 + \eps_1^{-1}) r^2}{T} \big] \M^{\frac{1}{2} - (q - q')} (v) \M^{\frac{1}{2} - q' (1 + \eps_1)} (v_1)
		\end{aligned}
	\end{equation}
	for some constants $\eps_1 > 0$ and $0 < q' < q$ to be determined. For any fixed $0 < q < 1$, one should find a pair of $(\eps_1, q')$ such that
	\begin{equation}\label{K1-1-chi-3}
		\begin{aligned}
			\tfrac{1}{2} - (q - q') > 0 \,, \ \tfrac{1}{2} - q' (1 + \eps_1) > 0 \,, 0 < q' < q \,, \eps_1 > 0 \,,
		\end{aligned}
	\end{equation}
	which is equivalent to $0 < \eps_1 < \frac{1}{2 q'} - 1$ and $ \max \{ 0, q - \frac{1}{2} \} < q' < \min \{ q, \frac{1}{2} \}$. Then, for any fixed $0 < q < 1$, $q'$ exists, and $\eps_1$ can be taken as $\eps_1 = \tfrac{1}{2} (\frac{1}{2 q'} - 1) > 0$. Namely, the required pair of $(\eps_1, q')$ exists for $0 < q < 1$. Remark that if $q = 1$, the needed pair of $(\eps_1, q')$ does not exist.
	
	From \eqref{Cutoff}, \eqref{K1-1-chi-1}, \eqref{K1-1-chi-2}, \eqref{K1-1-chi-3} and the H\"older inequality, it is deduced that
	\begin{equation}
		\begin{aligned}
			\big| \langle & \Theta_\gamma \M^{- \frac{q}{2}} K_1^{1 - \chi} f , \M^{- \frac{q}{2}} f \rangle \big| \\
			& \lesssim_{\rho, \u, T, r} \iint_{\R^3 \times \R^3} \Theta_\gamma (v) |f (v) f (v_1)| \M^{\frac{1}{2} - (q - q')} (v) \M^{\frac{1}{2} - q' (1 + \eps_1)} (v_1) |v_1 - v|^\gamma \d v_1 \d v \\
			& \lesssim_{\rho, \u, T, r} \Big( \int_{\R^3}  | \Theta_\gamma (v) f (v) |^2 \M^{\frac{1}{2} - (q - q')} (v) [ \int_{\R^3} \M^{\frac{1}{2} - q' (1 + \eps_1)} (v_1) |v_1 - v|^\gamma \d v_1 ] \d v \Big)^\frac{1}{2} \\
			& \qquad \times \Big( \int_{\R^3} | f (v_1)|^2 \M^{\frac{1}{2} - q' (1 + \eps_1)} (v_1) [ \int_{\R^3} \M^{\frac{1}{2} - (q - q')} (v) |v_1 - v|^\gamma \d v ] \d v_1 \Big)^\frac{1}{2} \,.
		\end{aligned}
	\end{equation}
	Observe that $|\Theta_\gamma (v)|^2 \M^{\frac{1}{2} - (q - q')} (v)$ is uniformly bounded in $v \in \R^3$. Moreover, together with \eqref{nu-equivalent}, Lemma 2.3 of \cite{LS-2010-KRM} implies that
	\begin{equation}
		\begin{aligned}
			\int_{\R^3} \M^{\frac{1}{2} - q' (1 + \eps_1)} (v_1) |v_1 - v|^\gamma \d v_1 \thicksim_{\rho, \u, T} \nu (v) \,, \\
			\int_{\R^3} \M^{\frac{1}{2} - (q - q')} (v) |v_1 - v|^\gamma \d v \thicksim_{\rho, \u, T} \nu (v_1) \,.
		\end{aligned}
	\end{equation}
	Then there holds
	\begin{equation}
		\begin{aligned}
			\big| \langle \Theta_\gamma \M^{- \frac{q}{2}} K_1^{1 - \chi} f , & \M^{- \frac{q}{2}} f \rangle \big| \\
			& \lesssim_{\rho, \u, T, r} \Big( \int_{\R^3}  | f (v) |^2 \nu (v) \d v \Big)^\frac{1}{2} \Big( \int_{\R^3}  | f (v_1) |^2 \nu (v_1) \d v_1 \Big)^\frac{1}{2} = \| f \|^2_{L^2 (\nu)} \,.
		\end{aligned}
	\end{equation}
	Namely, The case $i = 1$ in \eqref{K12} holds.
	
	\vspace*{3mm}
	
	{\em Step 2. Control the quantity $\big| \langle \Theta_\gamma \M^{- \frac{q}{2}} K_2^{1 - \chi} f , \M^{- \frac{q}{2}} f \rangle \big|$.}
	
	\vspace*{3mm}
	
	By \eqref{K2-Lambda} and the fact $ \M^\frac{1}{2} (v) \M^\frac{1}{2} (v_1) = \M^\frac{1}{2} (v') \M^\frac{1}{2} (v_1') $, one has
	\begin{equation}\label{S0}
		\begin{aligned}
			\langle \Theta_\gamma \M^{- \frac{q}{2}} K_2^{1 - \chi} f , \M^{- \frac{q}{2}} f \rangle = 2 \iiint_{\mathbb{R}^3 \times \R^3 \times \mathbb{S}^2} & [1 - \chi (|v_1 - v|)] \Theta_\gamma (v) f (v) f (v') \\
			\times \M^\frac{1}{2} & (v_1') \M^\frac{1}{2} (v_1) \M^{-q} (v) b (\omega, v_1 - v) \d \omega \d v_1 \d v \,.
		\end{aligned}
	\end{equation}
	Since $|v_1' - v| \leq |v_1 - v| \leq 2 r$ (see Figure \ref{Fig1}), there holds that $|v_1' - v| \geq |v - \u| - |v_1' - v| \geq |v - \u| - 2 r$, i.e., $|v-\u| \leq |v_1' - \u| + 2 r$. Then, for some small undetermined $\eps_2 > 0$,
	\begin{equation*}
		\begin{aligned}
			|v-\u|^2 \leq (|v_1' - \u| + 2 r)^2 \leq (1 + \eps_2) | v_1' - \u |^2 + 4 (1 + \eps_2^{-1}) r^2 \,,
		\end{aligned}
	\end{equation*}
	Namely,
	\begin{equation}\label{S1}
		\begin{aligned}
			|v_1' - \u|^2 \geq \tfrac{1}{1 + \eps_2} |v - \u|^2 - \tfrac{4}{\eps_2} r^2 \,.
		\end{aligned}
	\end{equation}
	It immediately infers that
	\begin{equation}\label{S2}
		\begin{aligned}
			\M^\frac{1}{2} (v_1') \leq \big[ \tfrac{\rho}{(2 \pi T)^\frac{3}{2}} \big]^\frac{1}{2} \exp (\tfrac{r^2}{\eps_2 T}) \exp \big[ - \tfrac{|v - \u|^2}{4 (1 + \eps_2) T} \big] \lesssim_{\rho, T} \exp (\tfrac{r^2}{\eps_2 T}) \M^\frac{1}{2 (1 + \eps_2)} (v) \,.
		\end{aligned}
	\end{equation}
	Similarly in \eqref{S1}, there holds
	\begin{equation}\label{S3}
		\begin{aligned}
			|v - \u|^2 \leq (1 + \delta) |v_1 - \u|^2 + 4 (1 + \delta^{-1}) r^2
		\end{aligned}
	\end{equation}
	for some small $\delta > 0$ to be determined. Let $\sigma \in (0,1)$ be an any fixed small parameter. Then, by \eqref{S3},
	\begin{equation*}
		\begin{aligned}
			q |v - \u|^2 = & (1 - \sigma) q |v - \u|^2 + \sigma q |v - \u|^2 \\
			\leq & (1 - \sigma) q |v - \u|^2 + \sigma (1 + \delta) q |v_1 - \u|^2 + 4 \sigma q (1 + \delta^{-1}) r^2 \,.
		\end{aligned}
	\end{equation*}
	which immediately implies that
	\begin{equation}\label{S5}
		\begin{aligned}
			\M^{-q} (v) \lesssim_{\rho, T} \exp \big[ \tfrac{2 \sigma q (1 + \delta^{-1}) r^2}{T} \big] \M^{- (1 - \sigma) q} (v) \M^{- \sigma (1 + \delta) q} (v_1) \,.
		\end{aligned}
	\end{equation}
	Then, \eqref{S2} and \eqref{S5} reduce to
	\begin{equation}\label{S6}
		\begin{aligned}
			\M^\frac{1}{2} (v_1') \M^\frac{1}{2} (v_1) \M^{-q} (v) \lesssim_{\rho, \u, T, r} \M^{\frac{1}{2 (1 + \eps_2)}- (1 - \sigma) q} (v) \M^{\frac{1}{2}- \sigma (1 + \delta) q} (v_1) \,.
		\end{aligned}
	\end{equation}
	For any fixed $0 < q < 1$, one then wants to find some small parameters $\eps_2$, $\sigma$, $\delta > 0$ such that
	\begin{equation}\label{S4}
		\begin{aligned}
			s_1 : =\tfrac{1}{2 (1 + \eps_2)}- (1 - \sigma) q > 0 \,, \ s_2 : = \tfrac{1}{2}- \sigma (1 + \delta) q > 0 \,.
		\end{aligned}
	\end{equation}
	which is equivalent to
	\begin{equation*}
		\begin{aligned}
			0 < \eps_2 < \tfrac{1}{2 (1 - \sigma) q} - 1 \,, \ q - \tfrac{1}{2} < \sigma q < \tfrac{1}{2 (1 + \delta)} \,.
		\end{aligned}
	\end{equation*}
	If $\frac{1}{2} < q < 1$, one first takes $\delta \in (0, \tfrac{2 (1 - q)}{2 q - 1})$ such that $q - \tfrac{1}{2} < \tfrac{1}{2 (1 + \delta)}$; then takes $\sigma \in ( \tfrac{2 q - 1}{2 q}, \tfrac{1}{2 q (1 + \delta)} )$ such that $q - \tfrac{1}{2} < \sigma q < \tfrac{1}{2 (1 + \delta)}$, which implies $\tfrac{1}{2 (1 - \sigma) q} - 1 > 0$; and finally takes $\eps_2 \in (0, \tfrac{1}{2 (1 - \sigma) q} - 1)$. If $0 < q \leq \frac{1}{2}$, one can first choose any $\delta > 0$ and $0 < \sigma < \tfrac{1}{2 q (1 + \delta)}$ such that $q - \tfrac{1}{2} \leq 0 < \sigma q < \tfrac{1}{2 (1 + \delta)}$ and $\tfrac{1}{2 (1 - \sigma) q} - 1 > 0$; then choose $\eps_2 \in (0, \tfrac{1}{2 (1 - \sigma) q} - 1)$. In summary, there is at least one group of $(\eps_2, \sigma, \delta)$ such that \eqref{S4} holds.
	
	Let $s = \min \{ s_1, s_2 \} > 0$. Based on \eqref{S4}, one derives from plugging \eqref{S6} into \eqref{S0} that
	\begin{equation}\label{S7}
		\begin{aligned}
			\big| \langle \Theta_\gamma \M^{- \frac{q}{2}} K_2^{1 - \chi} f , \M^{- \frac{q}{2}} f \rangle \big| \lesssim_{\rho, \u, T, r} & \iiint_{\mathbb{R}^3 \times \R^3 \times \mathbb{S}^2} [1 - \chi (|v_1 - v|)] \Theta_\gamma (v) | f (v) f (v') | \\
			& \qquad \times \M^\frac{1}{2} (v_1') \M^\frac{1}{2} (v_1) \M^{-q} (v) b (\omega, v_1 - v) \d \omega \d v_1 \d v \\
			\lesssim_{\rho, \u, T, r} & \iiint_{\mathbb{R}^3 \times \R^3 \times \mathbb{S}^2} \Theta_\gamma (v) | f (v) f (v') | \d \mu (\omega, v_1, v) \\
			\lesssim_{\rho, \u, T, r} & \Big( \iiint_{\mathbb{R}^3 \times \R^3 \times \mathbb{S}^2} | \Theta_\gamma (v) f (v) |^2 \d \mu (\omega, v_1, v) \Big)^\frac{1}{2}  \\
			& \qquad \times \Big( \iiint_{\mathbb{R}^3 \times \R^3 \times \mathbb{S}^2} | f (v') |^2 \d \mu (\omega, v_1, v) \Big)^\frac{1}{2} \,,
		\end{aligned}
	\end{equation}
	where the nonnegative measure $\d \mu (\omega, v_1, v) = \M^s (v) \M^s (v_1) b (\omega, v_1 - v) \d \omega \d v_1 \d v$. Notice that $\d \mu (\omega, v_1, v) = \d \mu (\omega, v_1', v')$ under the change of variables $(v, v_1) \mapsto (v', v_1')$, see Subsection 2.5 of \cite{LS-2010-KRM} for instance. Then
	\begin{equation}\label{S8}
		\begin{aligned}
			\iiint_{\mathbb{R}^3 \times \R^3 \times \mathbb{S}^2} | f (v') |^2 \d \mu (\omega, v_1, v) = \iiint_{\mathbb{R}^3 \times \R^3 \times \mathbb{S}^2} | f (v) |^2 \d \mu (\omega, v_1, v) \\
			= \int_{\R^3} |f (v)|^2 \M^s (v) [\iint_{\R^3 \times \mathbb{S}^2} \M^s (v_1) b (\omega, v_1 - v) \d \omega \d v_1 ] \d v \\
			= \int_{\R^3} |f (v)|^2 \M^s (v) [\int_{\R^3} \M^s (v_1) |v_1 - v|^\gamma \d v_1 ] \d v \\
			\lesssim_{\rho, \u, T} \int_{\R^3} |f (v)|^2 \nu (v) \d v = \| f \|^2_{L^2 (\nu)} \,,
		\end{aligned}
	\end{equation}
	where the fact $\int_{\R^3} \M^s (v_1) |v_1 - v|^\gamma \d v_1 \thicksim_{\rho, \u, T} \nu (v)$ has been used, which is derived from \eqref{nu-equivalent} and Lemma 2.3 of \cite{LS-2010-KRM}. Moreover, together with the fact $\sup_{v \in \R^3} |\Theta_\gamma (v)|^2 \M^s (v) < \infty$, a similar argument yields that
	\begin{equation}\label{S9}
		\begin{aligned}
			\iiint_{\mathbb{R}^3 \times \R^3 \times \mathbb{S}^2} | \Theta_\gamma (v) f (v) |^2 \d \mu (\omega, v_1, v) \lesssim_{\rho, \u, T} \| f \|^2_{L^2 (\nu)} \,.
		\end{aligned}
	\end{equation}
	Consequently, \eqref{S7}, \eqref{S8} and \eqref{S9} imply that
	\begin{equation*}
		\begin{aligned}
			\big| \langle \Theta_\gamma \M^{- \frac{q}{2}} K_2^{1 - \chi} f , \M^{- \frac{q}{2}} f \rangle \big| \lesssim_{\rho, \u, T, r} \| f \|^2_{L^2 (\nu)} \,.
		\end{aligned}
	\end{equation*}
	Namely, the inequality \eqref{K12} holds for $i = 2$. Therefore the proof of Lemma \ref{Lmm-K1K2-Hypo1} is finished.	

\subsection{High velocities estimates for $K$: Proof of Lemma \ref{Lmm-K1K2-Hypo2}}

We now start our proof.

	(1) One first focuses on the quantity $\big| \langle \Theta_\gamma \M^{- \frac{q}{2}} K_1^\chi f, \M^{- \frac{q}{2}} f \rangle \big|$. By the definitions of $\chi (s)$ in \eqref{Chi} and $K_1^\chi f$ in \eqref{K1-Lambda},
	\begin{equation}\label{T0}
		\begin{aligned}
			& \big| \langle \Theta_\gamma \M^{- \frac{q}{2}} K_1^\chi f, \M^{- \frac{q}{2}} f \rangle \big| \\
			\leq & \iiint_{\mathbb{R}^3 \times \R^3 \times \mathbb{S}^2} \M^{-q} (v) \M^\frac{1}{2} (v) \M^\frac{1}{2} (v_1) \chi (|v_1 - v|) \Theta_\gamma (v) |f (v_1) f (v)| b (\omega, v_1 - v) \d \omega \d v_1 \d v \\
			\leq & \iint_{|v_1 - v| \geq r} |f (v) \M^{- \frac{q}{2}} (v)| \, |f(v_1) \M^{- \frac{q}{2}} (v_1)| \Theta_\gamma (v) \M^\frac{1-q}{2} (v) \M^\frac{1+q}{2} (v_1) |v_1 - v|^\gamma \d v_1 \d v \\
			\lesssim&_{\rho, \u, T} \Big( \iint_{|v_1 - v| \geq r} |f (v) \M^{- \frac{q}{2}} (v)|^2 \M^\frac{3(1-q)}{8} (v) \M^\frac{1+q}{2} (v_1) |v_1 - v|^\gamma \d v_1 \d v \Big)^\frac{1}{2} \\
			& \quad \times \Big( \iint_{|v_1 - v| \geq r} |f (v_1) \M^{- \frac{q}{2}} (v_1)|^2 \M^\frac{3(1-q)}{8} (v) \M^\frac{1+q}{2} (v_1) |v_1 - v|^\gamma \d v_1 \d v \Big)^\frac{1}{2} \,,
		\end{aligned}
	\end{equation}
	where the cutoff condition \eqref{Cutoff} and the uniform bound $\sup_{v \in \R^3} \Theta_\gamma (v) \M^\frac{1-q}{8} (v) \leq C(\rho, \u, T) < \infty$ have also been utilized. Since $|v_1 - \u|^2 + |v - \u|^2 \geq \frac{1}{2} |v_1 - v|^2 \geq \frac{1}{2} r^2$, one has
	\begin{equation}
		\begin{aligned}
			\big[ \M (v) \M (v_1) \big]^\frac{1-q}{4} = & \big[ \tfrac{\rho}{(2 \pi T)^\frac{3}{2}} \big]^\frac{1-q}{2} \exp \big[ - \tfrac{(1-q) (|v-\u|^2 + |v_1 - \u|^2)}{8 T} \big] \\
			\leq & \big[ \tfrac{\rho}{(2 \pi T)^\frac{3}{2}} \big]^\frac{1-q}{2} \exp \big[ - \tfrac{(1-q) r^2}{16 T} \big] \,.
		\end{aligned}
	\end{equation}
	It thereby infers that
		\begin{align}\label{T1}
			\no \iint&_{|v_1 - v| \geq r} |f (v) \M^{- \frac{q}{2}} (v)|^2 \M^\frac{3(1-q)}{8} (v) \M^\frac{1+q}{2} (v_1) |v_1 - v|^\gamma \d v_1 \d v \\
			\no \lesssim_{\rho, \u, T} & \exp \big[ - \tfrac{(1-q) r^2}{16 T} \big] \int_{\R^3} |f (v) \M^{- \frac{q}{2}} (v)|^2 \M^\frac{1-q}{8} (v) \big[ \int_{\R^3} \M^\frac{1+3q}{4} (v_1) |v_1 - v|^\gamma \d_1 \big] \d v \\
			\no \lesssim_{\rho, \u, T} & \exp \big[ - \tfrac{(1-q) r^2}{16 T} \big] \int_{\R^3} |f (v) \M^{- \frac{q}{2}} (v)|^2 \M^\frac{1-q}{8} (v) \nu (v) \d v \\
			\lesssim_{\rho, \u, T} & \exp \big[ - \tfrac{(1-q) r^2}{16 T} \big] \| \M^{- \frac{q}{2}} f \|^2_{L^2} \,,
		\end{align}
	where the last second inequality is implied by the fact $\int_{\R^3} \M^\frac{1+3q}{4} (v_1) |v_1 - v|^\gamma \d_1 \thicksim_{\rho, \u, T} \nu (v)$ derived from Lemma 2.3 of \cite{LS-2010-KRM} and \eqref{nu-equivalent}, and the last inequality is resulted from the uniform bound $\sup_{v \in \R^3} \M^\frac{1-q}{8} (v) \nu (v) \leq C(\rho, \u, T) < \infty$. Moreover, there similarly hold
	\begin{equation}\label{T2}
		\begin{aligned}
			\iint_{|v_1 - v| \geq r} |f (v_1) \M^{- \frac{q}{2}} (v_1)|^2 \M^\frac{3(1-q)}{8} (v) \M^\frac{1+q}{2} (v_1) |v_1 - v|^\gamma \d v_1 \d v \\
			\leq \int_{\R^3} |f (v_1) \M^{- \frac{q}{2}} (v_1)|^2 \M^\frac{1+q}{2} (v_1) \big[ \int_{\R^3} \M^\frac{3(1-q)}{8} (v) |v_1 - v|^\gamma \d v \big] \d v_1 \\
			\lesssim_{\rho, \u, T} \int_{\R^3} |f (v_1) \M^{- \frac{q}{2}} (v_1)|^2 \M^\frac{1+q}{2} (v_1) \nu (v_1) \d v_1 \lesssim_{\rho, \u, T} \| \M^{- \frac{q}{2}} f \|^2_{L^2} \,.
		\end{aligned}
	\end{equation}
	Therefore, the inequality \eqref{K1-chi} is concluded by plugging \eqref{T1} and \eqref{T2} into \eqref{T0}.
	
	\vspace*{3mm}
	
	(2) Recalling the definition of $K_2^\chi f$ in \eqref{K2-Lambda}, one has
	\begin{equation}\label{T3}
		\begin{aligned}
			|K_2^\chi f (v)| \leq 2 \iint_{\R^3 \times \mathbb{S}^2} \chi (|v_1 - v|) \M^\frac{1}{2} (v) \M^{- \frac{1}{2}} (v') \M (v_1) |f (v')| b (\omega, v_1 - v) \d \omega \d v_1 \,.
		\end{aligned}
	\end{equation}
	Let
	\begin{equation*}
		\begin{aligned}
			V = v_1 - v \,, \ V_\shortparallel = (V \cdot \omega) \omega \,,  V_\perp = V - V_\shortparallel \,, \omega = \frac{v' - v}{|v' - v|} \,,
		\end{aligned}
	\end{equation*}
	whose intuitively geometric relations are illustrated in Figure \ref{Fig1} before. As shown in Section 2 of \cite{Grad-1963}, one has
	\begin{equation}\label{T4}
		\begin{aligned}
			\d \omega \d v_1 = \d \omega \d V = \frac{2 \d V_\perp \d V_\shortparallel}{|V_\shortparallel|^2} \,.
		\end{aligned}
	\end{equation}
	Then, \eqref{T3} and \eqref{T4} reduce to
	\begin{equation*}
		\begin{aligned}
			|K_2^\chi f (v)| \leq 4 \int_{\R^3} \int_{V_\perp \perp V_\shortparallel} \chi (|V|) \M^\frac{1}{2} (v) \M^{- \frac{1}{2}} (v + V_\shortparallel) \M (v + V) |f (v + V_\shortparallel)| b (\omega , V) \frac{\d V_\perp \d V_\shortparallel}{|V_\shortparallel|^2} \,.
		\end{aligned}
	\end{equation*}
	Denote by
	\begin{equation*}
		\begin{aligned}
			\eta = v + V_\shortparallel \,, \ \zeta = \tfrac{1}{2} (v + \eta) \,,
		\end{aligned}
	\end{equation*}
	which leads to $\zeta - v = \tfrac{1}{2} V_\shortparallel$ and $v + V = \eta + V_\perp$. It is easy to see that $ V_\perp \cdot \zeta = V_\perp \cdot v = V_\perp \cdot \eta $. Then, direct calculations imply
	\begin{equation*}
		\begin{aligned}
			- \tfrac{|v-\u|^2}{2} + \tfrac{|\eta - \u|^2}{2} - |\eta + V_\perp - \u |^2 = - \tfrac{|v-\u|^2}{2} - \tfrac{|\eta - \u|^2}{2} - 2 V_\perp \cdot (\eta - \u) - |V_\perp|^2 \\
			= - \tfrac{|v-\u|^2}{2} - \tfrac{|\eta - \u|^2}{2} - 2 V_\perp \cdot (\zeta - \u) - |V_\perp|^2 \,,
		\end{aligned}
	\end{equation*}
	and
	\begin{equation*}
		\begin{aligned}
			|\eta - \u |^2 = & |2 \zeta - v - \u|^2 = |2 (\zeta - \u) - (v - \u)|^2 \\
			= & 4 |\zeta - \u|^2 - 4 (\zeta - \u) \cdot (v - \u) + |v - \u|^2 \,.
		\end{aligned}
	\end{equation*}
	It thereby infers that
	\begin{equation*}
		\begin{aligned}
			- \tfrac{|v-\u|^2}{2} + \tfrac{|\eta - \u|^2}{2} - |\eta + V_\perp - \u |^2 = & - |\zeta - v|^2 - |V_\perp + \zeta - \u|^2 \\
			= & - \tfrac{1}{4} |V_\shortparallel|^2 - |\zeta_\shortparallel|^2 - |V_\perp + \zeta_\perp|^2 \,,
		\end{aligned}
	\end{equation*}
	where $\zeta_\shortparallel = [(\zeta - \u)\cdot \omega] \omega$ and $\zeta_\perp = (\zeta - \u) - \zeta_\shortparallel$. As a result,
	\begin{equation*}
		\begin{aligned}
			|K_2^\chi f (v)| \leq \tfrac{4 \rho}{(2 \pi T)^\frac{3}{2}} \int_{\R^3} \int_{V_\perp \perp V_\shortparallel} \tfrac{\chi (|V|)}{|V_\shortparallel |^2} |f (v + V_\shortparallel)| b (\omega, V) \exp \big[ - \tfrac{|V_\shortparallel |^2}{8 T} - \tfrac{|\zeta_\shortparallel |^2}{2 T} - \tfrac{|V_\perp + \zeta_\perp |^2}{2 T} \big] \d V_\perp \d V_\shortparallel \,.
		\end{aligned}
	\end{equation*}
	By \eqref{Cutoff}, one has
	\begin{equation*}
		\begin{aligned}
			b (\omega, V) \leq \beta_0 |\cos \theta| |V|^\gamma = \beta_0 \tfrac{|\omega \cdot V|}{|V|} |V|^\gamma = \tfrac{\beta_0 |V_\shortparallel|}{|V|^{1 - \gamma}} \,.
		\end{aligned}
	\end{equation*}
	It thereby implies that
	\begin{equation}\label{K2-chi-bnd}
		\begin{aligned}
			|K_2^\chi f (v)| \leq \tfrac{4 \beta_0 \rho}{(2 \pi T)^\frac{3}{2}} \int_{\R^3} |V_\shortparallel|^{-1} & |f (v + V_\shortparallel)| \exp \big[ - \tfrac{|V_\shortparallel|^2}{8 T} - \tfrac{|\zeta_\shortparallel|^2}{2 T} \big] \\
			& \times \big[ \int_{V_\perp \perp V_\shortparallel} \tfrac{\chi (|V|)}{|V|^{1 - \gamma}} \exp \big( - \tfrac{|V_\perp + \zeta_\perp|^2}{2 T} \big) \d V_ \perp \big] \d V_\shortparallel \,.
		\end{aligned}
	\end{equation}
	
	\vspace*{3mm}
	
	\textbf{\emph{\underline{Case 1. $0 \leq \gamma \leq 1$ for the hard potential.}}}
	
	\vspace*{3mm}
	
	Since $1 - \gamma \geq 0$, the definition of $\chi (s)$ in \eqref{Chi} implies
	\begin{equation*}
		\begin{aligned}
			\int_{V_\perp \perp V_\shortparallel} \tfrac{\chi (|V|)}{|V|^{1 - \gamma}} \exp \big( - \tfrac{|V_\perp + \zeta_\perp|^2}{2 T} \big) \d V_ \perp \leq \frac{1}{r^{1 - \gamma}} \int_{V_\perp \perp V_\shortparallel} \exp \big( - \tfrac{|V_\perp + \zeta_\perp|^2}{2 T} \big) \d V_ \perp = \frac{2 \pi T}{r^{1 - \gamma}} \,.
		\end{aligned}
	\end{equation*}
	Therefore, one has
	\begin{equation*}
		\begin{aligned}
			|K_2^\chi f (v)| \leq \tfrac{4 \beta_0 \rho}{\sqrt{2 \pi T}} \tfrac{1}{r^{1 - \gamma}} \int_{\R^3} |V_\shortparallel|^{-1} |f (v + V_\shortparallel)| \exp \big[ - \tfrac{|V_\shortparallel|^2}{8 T} - \tfrac{|\zeta_\shortparallel|^2}{2 T} \big] \d V_\shortparallel \,,
		\end{aligned}
	\end{equation*}
	which reduces to
	\begin{equation}
		\begin{aligned}
			\big| \langle & \M^{- \frac{q}{2}} K_2^\chi f , \M^{ - \frac{q}{2}} f \rangle \big| \\
			\lesssim_{\rho, \u, T} & \tfrac{1}{r^{1 - \gamma}} \iint_{\R^3 \times \R^3} \tfrac{1}{|V_\shortparallel|} |f (v + V_\shortparallel) f (v)| \M^{-q} (v) \exp \big[ - \tfrac{|V_\shortparallel|^2}{8 T} - \tfrac{|\zeta_\shortparallel|^2}{2 T} \big] \d V_\shortparallel \d v \\
			\lesssim_{\rho, \u, T} & \tfrac{1}{r^{1 - \gamma}} \iint_{\R^3 \times \R^3} \tfrac{1}{|V_\shortparallel|} |\M^{- \frac{q}{2}} (v + V_\shortparallel) f (v + V_\shortparallel) | \, | \M^{- \frac{q}{2}} (v) f (v)| \\
			& \qquad \quad \times \exp \big( \tfrac{q |v - \u|^2}{4 T} - \tfrac{q |v + V_\shortparallel -\u|^2}{4 T} - \tfrac{|V_\shortparallel|^2}{8 T} - \tfrac{|\zeta_\shortparallel|^2}{2 T} \big) \d V_\shortparallel \d v \,.
		\end{aligned}
	\end{equation}
	By $\zeta = v + \frac{1}{2} V_\shortparallel$ and $V_\shortparallel \cdot (\zeta - \u) = V_\shortparallel \cdot \zeta_\shortparallel$, a simple calculation implies
	\begin{equation}\label{Qexpand-1}
		\begin{aligned}
			& \tfrac{q |v - \u|^2}{4 T} - \tfrac{q |v + V_\shortparallel -\u|^2}{4 T} - \tfrac{|V_\shortparallel|^2}{8 T} - \tfrac{|\zeta_\shortparallel|^2}{2 T} \\
			= & - (1 - q) \big( \tfrac{|V_\shortparallel|^2}{8 T} + \tfrac{|\zeta_\shortparallel|^2}{2 T} \big) + \tfrac{q}{2 T} \big[ - \tfrac{|V_\shortparallel|^2}{2} - V_\shortparallel \cdot (\zeta - \u - \tfrac{1}{2} V_\shortparallel ) - \tfrac{|V_\shortparallel|^2}{4} - |\zeta_\shortparallel|^2 \big] \\
			= & - (1 - q) \big( \tfrac{|V_\shortparallel|^2}{8 T} + \tfrac{|\zeta_\shortparallel|^2}{2 T} \big) + \tfrac{q}{2 T} \big( - V_\shortparallel \cdot \zeta_\shortparallel - \tfrac{|V_\shortparallel|^2}{4} - |\zeta_\shortparallel|^2 \big) \\
			= & - (1 - q) \big( \tfrac{|V_\shortparallel|^2}{8 T} + \tfrac{|\zeta_\shortparallel|^2}{2 T} \big) - \tfrac{q |\zeta_\shortparallel + \tfrac{1}{2} V_\shortparallel|^2}{2 T} \leq - \tfrac{(1 - q) |V_\shortparallel|^2}{8 T} \,.
		\end{aligned}
	\end{equation}
	Then, one has
	\begin{equation*}
		\begin{aligned}
			\big| \langle \M^{- \frac{q}{2}} K_2^\chi f , \M^{ - \frac{q}{2}} f \rangle \big| \lesssim_{\rho, \u, T} & \tfrac{1}{r^{1 - \gamma}} \iint_{\R^3 \times \R^3} \tfrac{1}{|V_\shortparallel|} |\M^{- \frac{q}{2}} (v + V_\shortparallel) f (v + V_\shortparallel) | \, | \M^{- \frac{q}{2}} (v) f (v)| \\
			& \qquad \qquad \qquad \qquad \qquad \qquad \times \exp \big[ - \tfrac{(1 - q) |V_\shortparallel|^2}{8 T} \big] \d V_\shortparallel d v \\
			\lesssim_{\rho, \u, T} & \tfrac{1}{r^{1 - \gamma}} \| \M^{- \frac{q}{2}} f \|^2_{L^2} \int_{\R^3} \tfrac{1}{|V_\shortparallel|} \exp \big[ - \tfrac{(1 - q) |V_\shortparallel|^2}{8 T} \big] \d V_\shortparallel \\
			\lesssim_{\rho, \u, T} & \tfrac{1}{r^{1 - \gamma}} \| \M^{- \frac{q}{2}} f \|^2_{L^2} \,,
		\end{aligned}
	\end{equation*}
	where the convergence $\int_{\R^3} \tfrac{1}{|V_\shortparallel|} \exp \big[ - \tfrac{(1 - q) |V_\shortparallel|^2}{8 T} \big] \d V_\shortparallel \leq C(T) < \infty$ has been utilized.
	
	If $0 \leq \gamma < 1$, the inequality \eqref{gamma-HP} holds. If $\gamma = 1$, it infers from $\nu (v) \thicksim_{\rho, \u, T} \langle v \rangle$ in \eqref{nu-equivalent} that
	\begin{equation*}
		\begin{aligned}
			\| \M^{- \frac{q}{2}} f \|^2_{L^2} \lesssim_{\rho, T} \int_{\R^3} |\M^{- \frac{q}{2}} (v) f (v)|^2 \tfrac{\nu (v)}{\langle v \rangle} \d v & = \Big\{ \int_{|v| \leq r} + \int_{|v| > r}  \Big\} |\M^{- \frac{q}{2}} (v) f (v)|^2 \tfrac{\nu (v)}{\langle v \rangle} \d v \\
			& \lesssim_{\rho, \u, T} \exp (\tfrac{q r^2}{2 T}) \| f \|^2_{L^2 (\nu)} + \tfrac{1}{r + 1} \| \M^{- \frac{q}{2}} f \|^2_{L^2 (\nu)} \,,
		\end{aligned}
	\end{equation*}
	which concludes the inequality \eqref{gamma-HS}.
	
	\vspace*{3mm}
	
	\emph{\underline{\textbf{ Case 2. $-3 < \gamma < 0$ for the soft potential.}}}
	
	\vspace*{3mm}
	
	By changing variables $V_\perp \to  V_\perp' = V_\perp - \zeta_\perp$, one has
	\begin{equation}
		\begin{aligned}
			\int_{V_\perp \perp V_\shortparallel} \tfrac{\chi (|V|)}{|V|^{1 - \gamma}} \exp \big( - \tfrac{|V_\perp + \zeta_\perp|^2}{2 T} \big) \d V_ \perp = \int_{V_\perp' \perp V_\shortparallel} \tfrac{\chi (\sqrt{|V_\shortparallel|^2 + |V_\perp' - \zeta_\perp|^2})}{(|V_\shortparallel|^2 + |V_\perp' - \zeta_\perp|^2)^{\frac{1 - \gamma}{2}}} \exp \big( - \tfrac{|V_\perp'|^2}{2 T} \big) \d V_\perp' \,.
		\end{aligned}
	\end{equation}
	Then the bound \eqref{K2-chi-bnd} of $K_2^\chi f (v)$ can be split into two parts:
	\begin{equation}\label{K2chi-W1+W2}
		\begin{aligned}
			|K_2^\chi f (v)| \leq W_1 (v) + W_2 (v) \,,
		\end{aligned}
	\end{equation}
	where
	\begin{equation}
		\begin{aligned}
			W_1 (v) = \tfrac{4 \beta_0 \rho}{(2 \pi T)^\frac{3}{2}} \int_{\R^3} |V_\shortparallel|^{-1} & |f (v + V_\shortparallel)| \exp \big[ - \tfrac{|V_\shortparallel|^2}{8 T} - \tfrac{|\zeta_\shortparallel|^2}{2 T} \big] \\
			& \times \big[ \int_{|V_\perp'| > \tau |\zeta_\perp|} \tfrac{\chi (\sqrt{|V_\shortparallel|^2 + |V_\perp' - \zeta_\perp|^2})}{(|V_\shortparallel|^2 + |V_\perp' - \zeta_\perp|^2)^{\frac{1 - \gamma}{2}}} \exp \big( - \tfrac{|V_\perp'|^2}{2 T} \big) \d V_\perp' \big] \d V_\shortparallel \,,
		\end{aligned}
	\end{equation}
	and
	\begin{equation}\label{W2}
		\begin{aligned}
			W_2 (v) = \tfrac{4 \beta_0 \rho}{(2 \pi T)^\frac{3}{2}} \int_{\R^3} |V_\shortparallel|^{-1} & |f (v + V_\shortparallel)| \exp \big[ - \tfrac{|V_\shortparallel|^2}{8 T} - \tfrac{|\zeta_\shortparallel|^2}{2 T} \big] \\
			& \times \big[ \int_{|V_\perp'| \leq \tau |\zeta_\perp|} \tfrac{\chi (\sqrt{|V_\shortparallel|^2 + |V_\perp' - \zeta_\perp|^2})}{(|V_\shortparallel|^2 + |V_\perp' - \zeta_\perp|^2)^{\frac{1 - \gamma}{2}}} \exp \big( - \tfrac{|V_\perp'|^2}{2 T} \big) \d V_\perp' \big] \d V_\shortparallel \,.
		\end{aligned}
	\end{equation}
	Here $\tau \in (0,1)$ is an any fixed constant.
	
	\vspace*{2mm}
	
	{\em Case 2(a). Estimates of $\langle \nu^{-1} \M^{- \frac{q}{2}} W_1 , \M^{- \frac{q}{2}} f \rangle$}.
	
	\vspace*{2mm}
	
	Notice that, by \eqref{Chi},
	\begin{equation*}
		\begin{aligned}
			\tfrac{\chi(\sqrt{|V_\shortparallel|^2 + |V_\perp' - \zeta_\perp|^2})}{( |V_\shortparallel|^2 + |V_\perp' - \zeta_\perp|^2)^{\frac{1 - \gamma}{2}}} \leq \tfrac{1}{r^{1-\gamma}} \,.
		\end{aligned}
	\end{equation*}
	Moreover, under $|V_\perp'| > \tau |\zeta_\perp|$ and for any fixed $0 < q' < 1$,
	\begin{equation*}
		\begin{aligned}
			\tfrac{\chi (\sqrt{|V_\shortparallel|^2 + |V_\perp' - \zeta_\perp|^2})}{(|V_\shortparallel|^2 + |V_\perp' - \zeta_\perp|^2)^{\frac{1 - \gamma}{2}}} \exp \big( - \tfrac{|V_\perp'|^2}{2 T} \big) \leq & \tfrac{1}{r^{1-\gamma}} \exp \big( - \tfrac{|V_\perp'|^2}{2 T} \big) \leq \tfrac{1}{r^{1-\gamma}} \exp \big( - \tfrac{q'|V_\perp'|^2}{2 T} - \tfrac{(1 - q') \tau^2 |\zeta_\perp|^2}{2 T} \big) \\
			= & \tfrac{1}{r^{1-\gamma}} \exp \big[ \tfrac{(1 - q') \tau^2 |\zeta_\shortparallel|^2}{2 T} - \tfrac{q'|V_\perp'|^2}{2 T} - \tfrac{(1 - q') \tau^2 (|\zeta_\perp|^2 + |\zeta_\shortparallel|^2)}{2 T} \big] \,.
		\end{aligned}
	\end{equation*}
	Furthermore, $|\zeta_\perp|^2 + |\zeta_\shortparallel|^2 = |\zeta_\perp + \zeta_\shortparallel|^2 = |\zeta - \u|^2 = |v - \u + \tfrac{1}{2} V_\shortparallel|^2$. It thereby infers that
	\begin{equation}
		\begin{aligned}
			W_1 (v) \lesssim_{\rho, T} \tfrac{1}{r^{1-\gamma}} \int_{\R^3} |V_\shortparallel|^{-1} |f (v+V_\shortparallel)| \exp \big[ - \tfrac{|V_\shortparallel|^2}{8 T} - \tfrac{[1 - (1 - q') \tau^2]|\zeta_\shortparallel|^2}{2 T} \big] \\
			\times \Big\{ \int_{|V_\perp'| > \tau |\zeta_\perp|} \exp \big( - \tfrac{q'|V_\perp'|^2}{2 T} - \tfrac{(1 - q') \tau^2 |v - \u + \tfrac{1}{2} V_\shortparallel|^2 }{2 T} \big) \d V_\perp' \Big\} \d V_\shortparallel \,.
		\end{aligned}
	\end{equation}
	One expands
	\begin{equation}\label{Qexpand-2}
		\begin{aligned}
			|v - \u + \tfrac{1}{2} V_\shortparallel|^2 = & |v - \u|^2 + \tfrac{1}{4} |V_\shortparallel|^2 + (v - \u) \cdot V_\shortparallel \\
			= & \tfrac{3}{4} |v - \u|^2 + \tfrac{1}{4} |v - \u + V_\shortparallel|^2 + \tfrac{1}{2} (v - \u) \cdot V_\shortparallel \\
			\geq & \tfrac{3}{4} |v - \u|^2 + \tfrac{1}{4} |v - \u + V_\shortparallel|^2 - \tfrac{1}{4} |v - \u|^2 - \tfrac{1}{4} |V_\shortparallel|^2 \\
			= & \tfrac{1}{2} |v - \u|^2 + \tfrac{1}{4} |v - \u + V_\shortparallel|^2 - \tfrac{1}{4} |V_\shortparallel|^2 \,.
		\end{aligned}
	\end{equation}
	Denote by $\tau_* : = (1 - q') \tau^2 \in (0, 1)$. The arbitrariness of $q', \tau \in (0,1)$ guarantees that $\tau_*$ can be arbitrary in $(0, 1)$. Observing that $\int_{|V_\perp'| > \tau |\zeta_\perp|} \exp \big( - \tfrac{q'|V_\perp'|^2}{2 T} \big) \d V_\perp' \leq \tfrac{2 \pi T}{q'} < \infty$, there hold
	\begin{equation}
		\begin{aligned}
			W_1 (v) \lesssim_{\rho, T} & \tfrac{1}{r^{1-\gamma}} \exp \big[ - \tfrac{\tau_* |v - \u|^2}{4 T} \big] \int_{\R^3} |V_\shortparallel|^{-1} |f (v+V_\shortparallel)| \\
			& \qquad \times \exp \big[ - \tfrac{( 1 - \tau_* ) |V_\shortparallel|^2}{8 T} - \tfrac{( 1 - \tau_* ) |\zeta_\shortparallel|^2}{2 T} - \tfrac{\tau_* |v - \u + V_\shortparallel|^2}{8 T} \big] \d V_\shortparallel \,.
		\end{aligned}
	\end{equation}
	Then, by \eqref{nu-equivalent} and $\sup_{v \in \R^3} (1 + |v|)^{- \gamma} \exp \big( - \tfrac{\tau_* |v - \u|^2}{4 T} \big) \leq C(\u, T) < \infty$, one has
	\begin{equation}
		\begin{aligned}
			| \langle \nu^{-1} & \M^{- \frac{q}{2}} W_1 , \M^{- \frac{q}{2}} f \rangle | \\
			\lesssim_{\rho, \u, T} & \tfrac{1}{r^{1-\gamma}} \iint_{\R^3 \times \R^3} (1 + |v|)^{- \gamma} \exp \big( - \tfrac{\tau_* |v - \u|^2}{4 T} \big) |V_\shortparallel|^{-1} | \M^{- \frac{q}{2}} (v + V_\shortparallel) f (v+V_\shortparallel)| \\
			& \times | \M^{- \frac{q}{2}} (v) f (v) | \exp \big[ \tfrac{q |v - \u|^2}{4 T} - \tfrac{(1 - \tau_*) |V_\shortparallel|^2}{8 T} - \tfrac{(1 - \tau_*) |\zeta_\shortparallel|^2}{2 T} - \tfrac{(q + \tau_*/2) |v - \u + V_\shortparallel|^2}{4 T} \big] \d V_\shortparallel \d v \\
			\lesssim_{\rho, \u, T} & \tfrac{1}{r^{1-\gamma}} \iint_{\R^3 \times \R^3} |V_\shortparallel|^{-1} | \M^{- \frac{q}{2}} (v) f (v) \M^{- \frac{q}{2}} (v + V_\shortparallel) f (v+V_\shortparallel)| \\
			& \qquad \qquad \times \exp \big[ \tfrac{q |v - \u|^2}{4 T} - \tfrac{q |v - \u + V_\shortparallel|^2}{4 T} - \tfrac{(1 - \tau_*) |V_\shortparallel|^2}{8 T} - \tfrac{(1 - \tau_*) |\zeta_\shortparallel|^2}{2 T} \big] \d V_\shortparallel \d v \,.
		\end{aligned}
	\end{equation}
	By the same arguments in \eqref{Qexpand-1}, it is implied that
	\begin{equation}
		\begin{aligned}
			& \tfrac{q |v - \u|^2}{4 T} - \tfrac{q |v - \u + V_\shortparallel|^2}{4 T} - \tfrac{(1 - \tau_*) |V_\shortparallel|^2}{8 T} - \tfrac{(1 - \tau_*) |\zeta_\shortparallel|^2}{2 T} \\
			= & - (1 - \tau_* - q) \big( \tfrac{ |V_\shortparallel|^2}{8 T} + \tfrac{ |\zeta_\shortparallel|^2}{2 T} \big) - \tfrac{q |\zeta_\shortparallel + \frac{1}{2} V_\shortparallel|^2}{2 T} \\
			\leq & - (1 - \tau_* - q) \big( \tfrac{ |V_\shortparallel|^2}{8 T} + \tfrac{ |\zeta_\shortparallel|^2}{2 T} \big) \,.
		\end{aligned}
	\end{equation}
	By taking a $\tau_*$ such that $0 < \tau_* < 1 -q$, i.e., $1 - \tau_* - q > 0$, one has
	\begin{equation*}
		\begin{aligned}
			\tfrac{q |v - \u|^2}{4 T} - \tfrac{q |v - \u + V_\shortparallel|^2}{4 T} - \tfrac{(1 - \tau_*) |V_\shortparallel|^2}{8 T} - \tfrac{(1 - \tau_*) |\zeta_\shortparallel|^2}{2 T} \leq - \tfrac{(1 - \tau_* - q) |V_\shortparallel|^2}{8 T} \,.
		\end{aligned}
	\end{equation*}
	Then, by the inequality $\iint_{\R^3 \times \R^3} |g (v) g (v + V_\shortparallel) w (V_\shortparallel)| \d V_\shortparallel \d v \leq \int_{\R^3} |g (v)|^2 \d v \cdot \int_{\R^3} |w (V_\shortparallel)| \d V_\shortparallel$ and the convergence $\int_{\R^3} |V_\shortparallel|^{-1} \exp \big[ - \tfrac{(1 - \tau_* - q) |V_\shortparallel|^2}{8 T} \big] \d V_\shortparallel \leq 4 \pi + (\tfrac{8 \pi T}{1 - \tau_* - q})^\frac{3}{2} < \infty$, one easily obtains
	\begin{equation}\label{W1-bnd}
		\begin{aligned}
			| \langle \nu^{-1} & \M^{- \frac{q}{2}} W_1 , \M^{- \frac{q}{2}} f \rangle | \\
			\lesssim_{\rho, \u, T} & \tfrac{1}{r^{1-\gamma}} \iint_{\R^3 \times \R^3} | \M^{- \frac{q}{2}} (v) f (v) \M^{- \frac{q}{2}} (v + V_\shortparallel) f (v+V_\shortparallel)| \\
			& \qquad \qquad \qquad \qquad \times |V_\shortparallel|^{-1} \exp \big[ - \tfrac{(1 - \tau_* - q) |V_\shortparallel|^2}{8 T} \big] \d V_\shortparallel \d v \\
			\lesssim_{\rho, \u, T} & \tfrac{1}{r^{1-\gamma}} \int_{\R^3} |\M^{- \frac{q}{2}} (v) f (v)| \d v \cdot \int_{\R^3} |V_\shortparallel|^{-1} \exp \big[ - \tfrac{(1 - \tau_* - q) |V_\shortparallel|^2}{8 T} \big] \d V_\shortparallel \\
			\lesssim_{\rho, \u, T} & \tfrac{1}{r^{1-\gamma}} \| \M^{- \frac{q}{2}} f \|^2_{L^2} \,.
		\end{aligned}
	\end{equation}
	
	\vspace*{2mm}
	
	{\em Case 2(b). Estimates of $\langle \nu^{-1} \M^{- \frac{q}{2}} W_2 , \M^{- \frac{q}{2}} f \rangle$}.
	
	\vspace*{2mm}
	
	Without loss of generality, we assume $r > 1$ in \eqref{Chi}, because our later arguments will require the positive constant $r$ large enough. Then
	\begin{equation}
		\begin{aligned}
			\tfrac{\chi(\sqrt{|V_\shortparallel|^2 + |V_\perp' - \zeta_\perp|^2})}{(|V_\shortparallel|^2 + |V_\perp' - \zeta_\perp|^2)^{\frac{1 - \gamma}{2}}} \leq \tfrac{ 2^{1 - \gamma} }{(1 + |V_\shortparallel|^2 + |V_\perp' - \zeta_\perp|^2)^{\frac{1 - \gamma}{2}}}\,.
		\end{aligned}
	\end{equation}
	As a result, there holds
	\begin{equation}
		\begin{aligned}
			\int_{|V_\perp'| \leq \tau |\zeta_\perp|} \tfrac{\chi (\sqrt{|V_\shortparallel|^2 + |V_\perp' - \zeta_\perp|^2})}{(|V_\shortparallel|^2 + |V_\perp' - \zeta_\perp|^2)^{\frac{1 - \gamma}{2}}} & \exp \big( - \tfrac{|V_\perp'|^2}{2 T} \big) \d V_\perp' \\
			& \leq 2^{1 - \gamma} \int_{|V_\perp'| \leq \tau |\zeta_\perp|} \tfrac{ \exp \big( - \tfrac{|V_\perp'|^2}{2 T} \big) \d V_\perp' }{(1 + |V_\shortparallel|^2 + |V_\perp' - \zeta_\perp|^2)^{\frac{1 - \gamma}{2}}} \,.
		\end{aligned}
	\end{equation}
	Under the restriction $|V_\perp'| \leq \tau |\zeta_\perp|$, $|V_\perp' - \zeta_\perp| \geq |\zeta_\perp| - |V_\perp'| \geq (1 - \tau) |\zeta_\perp|$ ($0 < \tau < 1$). Moreover, $|\zeta_\perp|^2 + |\zeta_\shortparallel|^2 = |\zeta_\perp + \zeta_\shortparallel|^2 = |\zeta - \u|^2 = |v - \u + \tfrac{1}{2} V_\shortparallel|^2$. Thus, for any fixed $0 < q'' < 1$,
	\begin{equation}
		\begin{aligned}
			\int_{|V_\perp'| \leq \tau |\zeta_\perp|} & \tfrac{\chi (\sqrt{|V_\shortparallel|^2 + |V_\perp' - \zeta_\perp|^2})}{(|V_\shortparallel|^2 + |V_\perp' - \zeta_\perp|^2)^{\frac{1 - \gamma}{2}}} \exp \big( - \tfrac{|V_\perp'|^2}{2 T} \big) \d V_\perp' \\
			& \lesssim_{\rho, \u, T} \int_{|V_\perp'| \leq \tau |\zeta_\perp|} \tfrac{ \exp \big( - \tfrac{|V_\perp'|^2}{2 T} + \tfrac{q'' |\zeta_\shortparallel|^2}{2 T} \big) \d V_\perp' }{(1 + |V_\shortparallel|^2 + | \zeta_\perp |^2 + |\zeta_\shortparallel|^2)^{\frac{1 - \gamma}{2}}} \lesssim_{\rho, \u, T} \tfrac{ \exp \big( \tfrac{q'' |\zeta_\shortparallel|^2}{2 T} \big) }{(1 + |V_\shortparallel|^2 + |v - \u + \tfrac{1}{2} V_\shortparallel|^2)^{\frac{1 - \gamma}{2}}} \,,
		\end{aligned}
	\end{equation}
	where the bound $\int_{|V_\perp'| \leq \tau |\zeta_\perp|} \exp \big( - \tfrac{|V_\perp'|^2}{2 T} \big) \d V_\perp' \leq 2 \pi T < \infty$ has been used. Together with the definition of $W_2 (v)$ in \eqref{W2} and the inequality \eqref{Qexpand-2}, one has
	\begin{equation}\label{W-W1-W2}
		\begin{aligned}
			W_2 (v) & \lesssim_{\rho, \u, T} \int_{\R^3} \tfrac{|f (v + V_\shortparallel)| \exp \big[ - \tfrac{|V_\shortparallel|^2}{8 T} - \tfrac{(1 - q'') |\zeta_\shortparallel|^2}{2 T} \big] }{ |V_\shortparallel| ( 1 + |v - \u|^2 + |v - \u + V_\shortparallel|^2 + |V_\shortparallel|^2 )^\frac{1-\gamma}{2}} \d V_\shortparallel \\
			& \lesssim_{\rho, \u, T} \int_{\R^3} \tfrac{|f (v + V_\shortparallel)| \exp \big[ - \tfrac{|V_\shortparallel|^2}{8 T} - \tfrac{(1 - q'') |\zeta_\shortparallel|^2}{2 T} \big] }{ |V_\shortparallel| ( 1 + |v - \u| + |v - \u + V_\shortparallel| )^{1-\gamma}} \d V_\shortparallel = W_2^\phi (v) + W_2^{1 - \phi} (v)
		\end{aligned}
	\end{equation}
	for any fixed $0 < q'' < 1$, where
	\begin{equation}\label{W2-phi}
		\begin{aligned}
			W_2^\phi (v) = \int_{\R^3} \tfrac{ \phi (v, V_\shortparallel) |f (v + V_\shortparallel)| \exp \big[ - \tfrac{|V_\shortparallel|^2}{8 T} - \tfrac{(1 - q'') |\zeta_\shortparallel|^2}{2 T} \big] }{ |V_\shortparallel| ( 1 + |v - \u| + |v - \u + V_\shortparallel| )^{1-\gamma}} \d V_\shortparallel \,,
		\end{aligned}
	\end{equation}
	and
	\begin{equation}\label{W2-1-phi}
		\begin{aligned}
			W_2^{1 - \phi} (v) = \int_{\R^3} \tfrac{ [1 - \phi (v, V_\shortparallel) ] |f (v + V_\shortparallel)| \exp \big[ - \tfrac{|V_\shortparallel|^2}{8 T} - \tfrac{(1 - q'') |\zeta_\shortparallel|^2}{2 T} \big] }{ |V_\shortparallel| ( 1 + |v - \u| + |v - \u + V_\shortparallel| )^{1-\gamma}} \d V_\shortparallel \,.
		\end{aligned}
	\end{equation}
	Here the smooth cutoff function is defined by
	\begin{equation*}
		\begin{aligned}
			\phi (v, V_\shortparallel) = \chi (|v - \u| + |v - \u + V_\shortparallel|) \,,
		\end{aligned}
	\end{equation*}
	where $\chi$ is introduced in \eqref{Chi}.
	
	\vspace*{2mm}
	
	\underline{Case 2(b-1). Estimates of $\langle \nu^{-1} \M^{- \frac{q}{2}} W_2^\phi , \M^{- \frac{q}{2}} f \rangle$.}
	
	\vspace*{2mm}
	
	By \eqref{nu-equivalent} and $(1 + |v|)^{- \gamma} \leq (1 + |\u|)^{- \gamma} (1 + |v - \u|)^{- \gamma}$ for $- 3 < \gamma < 0$,
	\begin{equation}\label{W2-phi-1}
		\begin{aligned}
			\nu^{-1} \M^{- \frac{q}{2}} (v) \exp \big[ - \tfrac{|V_\shortparallel|^2}{8 T} - \tfrac{(1 - q'') |\zeta_\shortparallel|^2}{2 T} \big] \lesssim_{\rho, \u, T} (1 + |v-\u|)^{- \gamma} \exp \big( \tfrac{q |v - \u|^2}{4 T}  - \tfrac{|V_\shortparallel|^2}{8 T} - \tfrac{(1 - q'') |\zeta_\shortparallel|^2}{2 T} \big) \,.
		\end{aligned}
	\end{equation}
	Recalling the relations $V_\shortparallel \cdot \zeta_\perp =0$ and $\zeta - \u = \zeta_\shortparallel + \zeta_\perp = v - \u + \tfrac{1}{2} V_\shortparallel$, a direct calculation implies that
	\begin{equation*}
		\begin{aligned}
			& |v - \u|^2 + |v - \u + V_\shortparallel - \zeta_\perp|^2 \\
			= & |v - \u|^2 + |v - \u + V_\shortparallel|^2 + |\zeta_\perp|^2 - 2 (v - \u + V_\shortparallel) \cdot \zeta_\perp \\
			= & |v - \u|^2 + |v - \u + V_\shortparallel|^2 + |\zeta_\perp|^2 - 2 (v - \u ) \cdot \zeta_\perp \\
			= & |v - \u - \zeta_\perp|^2 + |v - \u + V_\shortparallel|^2 \\
			= & |\zeta - \u - \zeta_\perp - \tfrac{1}{2} V_\shortparallel|^2 + |v - \u + V_\shortparallel|^2 \\
			= & |\zeta_\shortparallel - \tfrac{1}{2} V_\shortparallel|^2 + |v - \u + V_\shortparallel|^2
		\end{aligned}
	\end{equation*}
	which means that
	\begin{equation}\label{W2-phi-2}
		\begin{aligned}
			|v - \u|^2 \leq |v - \u + V_\shortparallel|^2 + |\zeta_\shortparallel - \tfrac{1}{2} V_\shortparallel|^2 \,.
		\end{aligned}
	\end{equation}
	Moreover, $|\zeta_\shortparallel - \tfrac{1}{2} V_\shortparallel|^2 \leq 2 |\zeta_\shortparallel|^2 + \tfrac{1}{2} |V_\shortparallel|^2$, namely,
	\begin{equation}\label{W2-phi-3}
		\begin{aligned}
			|\zeta_\shortparallel|^2 \geq \tfrac{1}{2} |\zeta_\shortparallel - \tfrac{1}{2} V_\shortparallel|^2 - \tfrac{1}{4} |V_\shortparallel|^2 \,.
		\end{aligned}
	\end{equation}
	Then, it is deduced from \eqref{W2-phi-1}, \eqref{W2-phi-2} and \eqref{W2-phi-3} that
	\begin{equation}\label{W2-phi-4}
		\begin{aligned}
			\nu^{-1} & \M^{- \frac{q}{2}} (v) \exp \big[ - \tfrac{|V_\shortparallel|^2}{8 T} - \tfrac{(1 - q'') |\zeta_\shortparallel|^2}{2 T} \big] \\
			& \lesssim_{\rho, \u, T} (1 + |v-\u|)^{- \gamma} \exp \big( \tfrac{q |v - \u + V_\shortparallel|^2 + q |\zeta_\shortparallel - \frac{1}{2} V_\shortparallel|^2}{4 T}  - \tfrac{q'' |V_\shortparallel|^2}{8 T} - \tfrac{(1 - q'') |\zeta_\shortparallel - \frac{1}{2} V_\shortparallel|^2}{4 T} \big) \\
			& \lesssim_{\rho, \u, T} (1 + |v-\u|)^{- \gamma} \M^{- \frac{q}{2}} (v + V_\shortparallel) \exp \big( - \tfrac{q'' |V_\shortparallel|^2}{8 T} - \tfrac{(1 - q - q'') |\zeta_\shortparallel - \frac{1}{2} V_\shortparallel|^2}{4 T} \big) \,.
		\end{aligned}
	\end{equation}
	We then take $0 < q'' < 1 - q < 1$ for any fixed $q \in (0,1)$, i.e., $1 - q - q'' > 0$. From the definition of $W_2^\phi (v)$ in \eqref{W2-phi}, the bounds \eqref{W2-phi-4} and $\tfrac{\phi (v, V_\shortparallel) (1 + |v-\u|)^{- \gamma} }{ ( 1 + |v - \u| + |v - \u + V_\shortparallel| )^{1 - \gamma}} \leq \frac{1}{1 +r}$, it thereby infers that, for any fixed $0 < q < 1$,
	\begin{equation}\label{W2-phi-bnd}
		\begin{aligned}
			| \langle & \nu^{-1} \M^{- \frac{q}{2}} W_2^\phi , \M^{- \frac{q}{2}} f \rangle | \\
			\lesssim_{\rho, \u, T} & \iint_{\R^3 \times \R^3} \tfrac{\phi (v, V_\shortparallel) (1 + |v-\u|)^{- \gamma} }{|V_\shortparallel| ( 1 + |v - \u| + |v - \u + V_\shortparallel| )^{1 - \gamma}} |\M^{- \frac{q}{2}} (v) f (v) \cdot \M^{- \frac{q}{2}} (v + V_\shortparallel) f (v + V_\shortparallel)| \\
			& \qquad \qquad \times \exp \big( - \tfrac{q'' |V_\shortparallel|^2}{8 T} - \tfrac{(1 - q - q'') |\zeta_\shortparallel - \frac{1}{2} V_\shortparallel|^2}{4 T} \big) \d V_\shortparallel \d v \\
			\lesssim_{\rho, \u, T} & \tfrac{1}{1 + r} \iint_{\R^3 \times \R^3} |V_\shortparallel|^{-1} |\M^{- \frac{q}{2}} (v) f (v) \cdot \M^{- \frac{q}{2}} (v + V_\shortparallel) f (v + V_\shortparallel)| \exp \big( - \tfrac{q'' |V_\shortparallel|^2}{8 T} \big) \d V_\shortparallel \d v \\
			\lesssim_{\rho, \u, T} & \tfrac{1}{1 + r} \int_{\R^3} |\M^{- \frac{q}{2}} (v) f (v)|^2 \d v \cdot \int_{\R^3} |V_\shortparallel|^{-1} \exp \big( - \tfrac{q'' |V_\shortparallel|^2}{8 T} \big) \d V_\shortparallel \\
			\lesssim_{\rho, \u, T} & \tfrac{1}{1 + r} \| \M^{- \frac{q}{2}} f \|^2_{L^2} \,,
		\end{aligned}
	\end{equation}
	where the inequalities $\int_{\R^3} |V_\shortparallel|^{-1} \exp \big( - \tfrac{q'' |V_\shortparallel|^2}{8 T} \big) \d V_\shortparallel \leq 4 \pi + (\tfrac{8 \pi T}{q''})^\frac{3}{2} < \infty$ and $\iint_{\R^3 \times \R^3} |g (v) g (v + V_\shortparallel) w (V_\shortparallel)| \d V_\shortparallel \d v \leq \int_{\R^3} |g (v)|^2 \d v \cdot \int_{\R^3} |w (V_\shortparallel)| \d V_\shortparallel$ have also been employed. This completes the estimate for the quantity $\langle \nu^{-1} \M^{- \frac{q}{2}} W_2^\phi , \M^{- \frac{q}{2}} f \rangle$.
	
	\vspace*{2mm}
	
	\underline{Case 2(b-2). Estimates of $\langle \nu^{-1} \M^{- \frac{q}{2}} W_2^{1 -\phi} , \M^{- \frac{q}{2}} f \rangle$.}
	
	\vspace*{2mm}
	
	Observing that
	\begin{equation*}
		\begin{aligned}
			\Psi (v, V_\shortparallel) : = \tfrac{[1 - \phi (v, V_\shortparallel)] \nu^{- \frac{1}{2}} (v) \nu^{- \frac{1}{2}} (v + V_\shortparallel) \M^{-q} (v) }{(1 + |v - \u| + |v - \u + V_\shortparallel|)^{1 - \gamma}} \lesssim_{\rho, \u, T} (1 + r)^{- \gamma} \exp (\tfrac{2 q r^2}{T}) \,,
		\end{aligned}
	\end{equation*}
	one therefore derives from the definition of $W_2^{1 - \phi} (v)$ in \eqref{W2-1-phi} that
	\begin{equation}\label{W2-1-phi-bnd}
		\begin{aligned}
			| \langle & \nu^{-1} \M^{- \frac{q}{2}} W_2^{1 -\phi} , \M^{- \frac{q}{2}} f \rangle | \\
			& \leq \iint_{\R^3 \times \R^3} \Psi (v, V_\shortparallel) |V_\shortparallel|^{-1} |\nu^\frac{1}{2} (v) f (v) \cdot \nu^\frac{1}{2} (v + V_\shortparallel) f (v + V_\shortparallel)| \exp (- \tfrac{|V_\shortparallel|^2}{8 T}) \d V_\shortparallel \d v \\
			& \lesssim_{\rho, \u, T} (1 + r)^{- \gamma} \exp (\tfrac{2 q r^2}{T}) \int_{\R^3} |f (v)|^2 \nu (v) \d v \cdot \int_{\R^3} \int_{\R^3} |V_\shortparallel|^{-1} \exp (- \tfrac{|V_\shortparallel|^2}{8 T}) \d V_\shortparallel \\
			& \lesssim_{\rho, \u, T} (1 + r)^{- \gamma} \exp (\tfrac{2 q r^2}{T}) \| f \|^2_{L^2 (\nu)} \,.
		\end{aligned}
	\end{equation}
	This finishes the estimate for the quantity $\langle \nu^{-1} \M^{- \frac{q}{2}} W_2^{1 -\phi} , \M^{- \frac{q}{2}} f \rangle$.
	
	Consequently, collecting the estimates \eqref{W-W1-W2}, \eqref{W2-phi-bnd} and \eqref{W2-1-phi-bnd} reduces to
	\begin{equation}\label{W2-bnd}
		\begin{aligned}
			| \langle \nu^{-1} \M^{- \frac{q}{2}} W_2 , \M^{- \frac{q}{2}} f \rangle | \lesssim_{\rho, \u, T} \tfrac{1}{1 + r} \| \M^{- \frac{q}{2}} f \|^2_{L^2} + (1 + r)^{- \gamma} \exp (\tfrac{2 q r^2}{T}) \| f \|^2_{L^2 (\nu)}
		\end{aligned}
	\end{equation}
	for any fixed $0 < q < 1$ and $r > 1$. Furthermore, the bounds \eqref{K2chi-W1+W2}, \eqref{W1-bnd} and \eqref{W2-bnd} imply that, if $-3 < \gamma < 0$,
	\begin{equation*}
		\begin{aligned}
			| \langle \nu^{-1} \M^{- \frac{q}{2}} K_2^\chi f , \M^{- \frac{q}{2}} f \rangle | \lesssim_{\rho, \u, T} \tfrac{1}{1 + r} \| \M^{- \frac{q}{2}} f \|^2_{L^2} + (1 + r)^{- \gamma} \exp (\tfrac{2 q r^2}{T}) \| f \|^2_{L^2 (\nu)}
		\end{aligned}
	\end{equation*}
	holds for any fixed $0 < q < 1$ and $r > 1$. The proof of Lemma \ref{Lmm-K1K2-Hypo2} is completed.

\section{Weighted $L^2$ estimates of $K f$ for $- 3 < \gamma \leq - \frac{3}{2}$: Proof of Lemma \ref{Lmm-VSP}}\label{Sec:VSP}

In this section, the main goal is to derive the lower bound of the quantity $ \langle \langle v \rangle^{ 2 \gamma l } \M^{-q} \partial_v^\alpha ( \nu f ) , \partial_v^\alpha f \rangle $ and the upper bound of the quantity $\langle \langle v \rangle^{2 \gamma l} \M^{- q} \partial^\alpha_v (K f), g \rangle$, hence, to justify Lemma \ref{Lmm-VSP}.

	The bound \eqref{L2-fg1} in Lemma \ref{Lmm-VSP} can be verified by the similar arguments in Lemma 2 of \cite{Strain-Guo-2008-ARMA}. For simplicity, we omit the details of proof here. We then prove the inequality \eqref{L2-fg2}. Remark that the ideas to prove \eqref{L2-fg2} are similar to Lemma 2 of \cite{Strain-Guo-2008-ARMA}. However, in \cite{Strain-Guo-2008-ARMA}, the arbitrary $q \in (0,1)$ is restricted to a sufficiently small $q_0 > 0$, and the arbitrarily fixed constant state $(\rho, \u, T)$ is specialized as $(1,0,1)$. Then, for convenience of readers, we display the details of proving \eqref{L2-fg2} here.
	
	Since $K f = K_1 f - K_2 f$, we will estimate the $K_1 f$ part and $K_2 f$ part separately. Recalling the definition of $K_1 f$ in \eqref{K1} and the relation $\M (v) \M (v_1) = \M (v') \M (v_1') $, one has
	\begin{equation}\label{Star1}
		\begin{aligned}
			K_1 f (v) = & \int_{\R^3} |v_1 - v|^\gamma f (v_1) \M^\frac{1}{2} (v) \M^\frac{1}{2} (v_1) \d v_1 \\
			= & \int_{\R^3} |v_1|^\gamma f (v + v_1) \M^\frac{1}{2} (v) \M^\frac{1}{2} (v + v_1) \d v_1 \,,
		\end{aligned}
	\end{equation}
	which means that
	\begin{equation*}
		\begin{aligned}
			& \langle \langle v \rangle^{2 \gamma l} \M^{ - q } \partial_v^\alpha ( K_1 f ) , g \rangle \\
			= & \sum_{\alpha' \leq \alpha} C_\alpha^{\alpha'} \iint_{\R^3 \times \R^3} \langle v \rangle^{2 \gamma l} \M^{ - q } (v) \partial_v^{\alpha'} [ \M^\frac{1}{2} (v) \M^\frac{1}{2} ( v+ v_1) ] |v_1|^\gamma \partial_v^{ \alpha - \alpha' } f (v + v_1) g (v) \d v_1 \d v \,.
		\end{aligned}
	\end{equation*}
	Note that
	\begin{equation*}
		\begin{aligned}
			| \partial_v^{\alpha'} [ \M^\frac{1}{2} (v) \M^\frac{1}{2} ( v+ v_1) ] | \lesssim_{\rho, \u, T} \sum_{\alpha'' \leq \alpha'} \langle v \rangle^{|\alpha''|} \langle v + v_1 \rangle^{|\alpha' - \alpha''|} \M^\frac{1}{2} (v) \M^\frac{1}{2} (v + v_1) \,.
		\end{aligned}
	\end{equation*}
	It thereby infers that
	\begin{equation*}
		\begin{aligned}
			| \langle \langle v \rangle^{2 \gamma l} \M^{ - q } \partial_v^\alpha ( K_1 f ) , g \rangle | \lesssim & \sum_{ \substack{ \alpha' \leq \alpha \\ \alpha'' \leq \alpha' } } \iint_{\R^3 \times \R^3} \langle v \rangle^{2 \gamma l + |\alpha''|} |v_1|^\gamma | \partial_v^{ \alpha - \alpha' } f (v + v_1) | | g (v) | \\
			& \qquad \quad \times \langle v + v_1 \rangle^{|\alpha' - \alpha''|} \M^{ \frac{1}{2} - q } (v) \M^\frac{1}{2} (v + v_1) \d v_1 \d v \\
			= & \sum_{ \substack{ \alpha' \leq \alpha \\ \alpha'' \leq \alpha' } } \iint_{\R^3 \times \R^3} \tfrac{ \langle v \rangle^{ \gamma l + |\alpha''| - \frac{1}{2} \gamma } }{ \langle v_1 \rangle^{ \gamma ( l + \frac{1}{2} ) + |\alpha' - \alpha''| } } \big[ \M (v) \M (v_1) \big]^\frac{1-q}{2} |v_1 - v|^\gamma \\
			& \times | \langle v_1 \rangle^{ \gamma ( l + \frac{1}{2} ) } \M^{- \frac{q}{2}} \partial_x^{\alpha - \alpha'} f (v_1) | \cdot | \langle v \rangle^{ \gamma ( l + \frac{1}{2} ) } \M^{- \frac{q}{2}} g (v) | \d v_1 \d v \,.
		\end{aligned}
	\end{equation*}
	It is easy to see that
	\begin{equation*}
		\begin{aligned}
			\sum_{\alpha'' \leq \alpha'} \tfrac{ \langle v \rangle^{ \gamma l + |\alpha''| - \frac{1}{2} \gamma } }{ \langle v_1 \rangle^{ \gamma ( l + \frac{1}{2} ) + |\alpha' - \alpha''| } } \big[ \M (v) \M (v_1) \big]^\frac{1-q}{2} \lesssim \big[ \M (v) \M (v_1) \big]^\frac{1-q}{4}
		\end{aligned}
	\end{equation*}
	for $0 < q < 1$. Then
	\begin{equation}\label{Star2}
		\begin{aligned}
			| \langle \langle v \rangle^{2 \gamma l} \M^{ - q } \partial_v^\alpha ( K_1 f ) , g \rangle | \lesssim \sum_{\alpha' \leq \alpha} ( I_m^- + I_m^+ ) \,,
		\end{aligned}
	\end{equation}
	where
	\begin{equation*}
		\begin{aligned}
			I_m^\pm = \iint_{\R^3 \times \R^3} & |v_1 - v|^\gamma \big[ \M (v) \M (v_1) \big]^\frac{1-q}{4} \mathbf{1}_{\pm ( |v_1 - v| - m_0 ) \geq 0 } \\
			& \times | \langle v_1 \rangle^{ \gamma ( l + \frac{1}{2} ) } \M^{- \frac{q}{2}} \partial_v^{\alpha'} f (v_1) | \cdot | \langle v \rangle^{ \gamma ( l + \frac{1}{2} ) } \M^{- \frac{q}{2}} g (v) | \d v_1 \d v \,.
		\end{aligned}
	\end{equation*}
	
	For the quantity $I_m^-$, the H\"older inequality reduces to
	\begin{equation*}
		\begin{aligned}
			I_m^- \leq \sqrt{ I_{m1}^- I_{m2}^- } \,,
		\end{aligned}
	\end{equation*}
	where
	\begin{equation*}
		\begin{aligned}
			I_{m1}^- = \iint_{\R^3 \times \R^3} \mathbf{1}_{ |v_1 - v| \leq m_0 } |v_1 - v|^\gamma \big[ \M (v) \M (v_1) \big]^\frac{1-q}{4} | \langle v_1 \rangle^{ \gamma ( l + \frac{1}{2} ) } \M^{- \frac{q}{2}} \partial_v^{\alpha'} f (v_1) |^2 \d v_1 \d v \,,
		\end{aligned}
	\end{equation*}
	and
	\begin{equation*}
		\begin{aligned}
			I_{m2}^- = \iint_{\R^3 \times \R^3} \mathbf{1}_{ |v_1 - v| \leq m_0 } |v_1 - v|^\gamma \big[ \M (v) \M (v_1) \big]^\frac{1-q}{4} | \langle v \rangle^{ \gamma ( l + \frac{1}{2} ) } \M^{- \frac{q}{2}} g (v) |^2 \d v_1 \d v \,.
		\end{aligned}
	\end{equation*}
	Observe that
	\begin{equation*}
		\begin{aligned}
			\int_{\R^3} \mathbf{1}_{ |v_1 - v| \leq m_0 } |v_1 - v|^\gamma \M^\frac{1-q}{4} (v) \d v \lesssim m_0^{3 + \gamma} \,.
		\end{aligned}
	\end{equation*}
	Then the Fubini Theorem shows
	\begin{equation*}
		\begin{aligned}
			I_{m1}^- \lesssim & m_0^{3 + \gamma} \int_{\R^3} \M^\frac{1-q}{4} (v_1) | \langle v_1 \rangle^{ \gamma ( l + \frac{1}{2} ) } \M^{- \frac{q}{2}} \partial_v^{\alpha'} f (v_1) |^2 \d v_1 \\
			\lesssim & m_0^{3 + \gamma} \| \langle v \rangle^{ \gamma ( l + \frac{1}{2} ) } \M^{- \frac{q}{2}} \partial_v^{\alpha'} f \|^2_{L^2} \,.
		\end{aligned}
	\end{equation*}
	Similarly, one has
	\begin{equation*}
		\begin{aligned}
			I_{m2}^- \lesssim m_0^{3 + \gamma} \| \langle v \rangle^{ \gamma ( l + \frac{1}{2} ) } \M^{- \frac{q}{2}} g (v) \|^2_{L^2} \,.
		\end{aligned}
	\end{equation*}
	As a result, the quantity $I_m^-$ can be bounded by
	\begin{equation*}
		\begin{aligned}
			I_m^- \lesssim m_0^{3 + \gamma} \| \langle v \rangle^{ \gamma ( l + \frac{1}{2} ) } \M^{- \frac{q}{2}} \partial_v^{\alpha'} f \|_{L^2} \| \langle v \rangle^{ \gamma ( l + \frac{1}{2} ) } \M^{- \frac{q}{2}} g (v) \|_{L^2} \,.
		\end{aligned}
	\end{equation*}
	
	For the quantity $I_m^+$, it can be split as
	\begin{equation*}
		\begin{aligned}
			I_m^+ = I_m^{+, -} + I_m^{+, +} \,,
		\end{aligned}
	\end{equation*}
	where
	\begin{equation*}
		\begin{aligned}
			I_m^{+, \pm} = \iint_{\R^3 \times \R^3} & \mathbf{1}_{ |v_1 - v| \geq m_0 } |v_1 - v|^\gamma \big[ \M (v) \M (v_1) \big]^\frac{1-q}{4} \mathbf{1}_{\pm ( |v_1| - m' ) \geq 0 } \\
			& \times | \langle v_1 \rangle^{ \gamma ( l + \frac{1}{2} ) } \M^{- \frac{q}{2}} \partial_v^{\alpha'} f (v_1) | \cdot | \langle v \rangle^{ \gamma ( l + \frac{1}{2} ) } \M^{- \frac{q}{2}} g (v) | \d v_1 \d v \,.
		\end{aligned}
	\end{equation*}
	Here $m' > 0$ is sufficiently large to be determined later. Notice that
	\begin{equation*}
		\begin{aligned}
			\mathbf{1}_{ |v_1 - v| \geq m_0 } \mathbf{1}_{ |v_1| \geq m' } \big[ \M (v) \M (v_1) \big]^\frac{1-q}{4} \lesssim e^{- c' (m')^2}  \big[ \M (v) \M (v_1) \big]^\frac{1-q}{8} \,,
		\end{aligned}
	\end{equation*}
	where $c' = \frac{1-q}{32 T} > 0$. Then the quantity $I_m^{+, +}$ can be bounded by
	\begin{equation*}
		\begin{aligned}
			I_m^{+, +} \lesssim & e^{- c' (m')^2} \iint_{\R^3 \times \R^3} \mathbf{1}_{ |v_1 - v| \geq m_0 } \mathbf{1}_{ |v_1| \geq m' } |v_1 - v|^\gamma \big[ \M (v) \M (v_1) \big]^\frac{1-q}{8} \\
			& \qquad \qquad \times | \langle v_1 \rangle^{ \gamma ( l + \frac{1}{2} ) } \M^{- \frac{q}{2}} \partial_v^{\alpha'} f (v_1) | \cdot | \langle v \rangle^{ \gamma ( l + \frac{1}{2} ) } \M^{- \frac{q}{2}} g (v) | \d v_1 \d v \\
			\lesssim & e^{- c' (m')^2} \int_{\R^3} | \langle v \rangle^{ \gamma ( l + \frac{1}{2} ) } \M^{- \frac{q}{2}} g (v) | \M^\frac{1-q}{8} (v) \Big( \int_{\R^3} |v_1 - v|^{2 \gamma} \M^\frac{1-q}{4} (v_1) \d v_1 \Big)^\frac{1}{2} \\
			& \qquad \qquad \qquad \times \Big( \int_{\R^3} | \langle v_1 \rangle^{ \gamma ( l + \frac{1}{2} ) } \M^{- \frac{q}{2}} \partial_v^{\alpha'} f (v_1) |^2 \d v_1 \Big)^\frac{1}{2} \d v \\
			\lesssim & e^{- c' (m')^2} \| \langle v \rangle^{ \gamma ( l + \frac{1}{2} ) } \M^{- \frac{q}{2}} \partial_v^{\alpha'} f \|_{L^2} \int_{\R^3} \langle v \rangle^\gamma \M^\frac{1-q}{8} (v) | \langle v \rangle^{ \gamma ( l + \frac{1}{2} ) } \M^{- \frac{q}{2}} g (v) | \d v \\
			\lesssim & e^{- c' (m')^2} \| \langle v \rangle^{ \gamma ( l + \frac{1}{2} ) } \M^{- \frac{q}{2}} \partial_v^{\alpha'} f \|_{L^2} \| \langle v \rangle^{ \gamma ( l + \frac{1}{2} ) } \M^{- \frac{q}{2}} g \|_{L^2} \,,
		\end{aligned}
	\end{equation*}
	where the last second inequality is implied by the bound
	\begin{equation*}
		\begin{aligned}
			\int_{\R^3} |v_1 - v|^{2 \gamma} \M^\frac{1-q}{4} (v_1) \d v_1 \lesssim \langle v \rangle^{2 \gamma}
		\end{aligned}
	\end{equation*}
	derived from Lemma 2.1 of \cite{LS-2010-KRM}, and the last inequality is resulted from the fact
	$$\int_{\R^3} \langle v \rangle^{ 2 \gamma } \M^\frac{1-q}{4} (v) \d v \lesssim 1 \,.$$
	Together with the fact $\nu (v) \thicksim \langle v \rangle^\gamma$ and the compactness of $K_1$ restricted on $\{ |v_1 - v| \geq m_0 , |v_1| \leq m' \}$, the quantity $I_m^{+, -}$ can easily be bounded by
	\begin{equation*}
		\begin{aligned}
			I_m^{+, -} \leq C(m') \big( \varpi_0 \sum_{|\alpha'| = |\alpha|} \| \langle v \rangle^{ \gamma ( l + \frac{1}{2} ) } \M^{- \frac{q}{2}} \partial_v^{\alpha'} f \|_{L^2} + C_{\varpi_0} \| f \|_{L^2 (\nu)} \| \langle v \rangle^{ \gamma ( l + \frac{1}{2} ) } \M^{- \frac{q}{2}} g \|_{L^2}
		\end{aligned}
	\end{equation*}
	for sufficiently small $\varpi_0 > 0$ to be determined later. Then the quantity $I_m^+$ can be dominated by
	\begin{equation*}
		\begin{aligned}
			I_m^+ \lesssim \Big\{ ( e^{- c' (m')^2} + C (m') \varpi_0 ) \sum_{|\alpha'| = |\alpha|} \| \langle v \rangle^{ \gamma ( l + \frac{1}{2} ) } \M^{- \frac{q}{2}} \partial_v^{\alpha'} f \|_{L^2} \\
			+ C(m', \varpi_0) \| f \|_{L^2 (\nu)} \Big\} \| \langle v \rangle^{ \gamma ( l + \frac{1}{2} ) } \M^{- \frac{q}{2}} g \|_{L^2}
		\end{aligned}
	\end{equation*}
	for large $m' > 0$ and small $\varpi_0 > 0$ to be determined later. Collecting the estimates for the quantities $I_m^-$ and $I_m^+$ above, one has
	\begin{equation}\label{K1-bnd}
		\begin{aligned}
			| \langle \langle v \rangle^{2 \gamma l} \M^{ - q } \partial_v^\alpha ( K_1 & f ) , g \rangle | \lesssim \| \langle v \rangle^{ \gamma ( l + \frac{1}{2} ) } \M^{- \frac{q}{2}} g \|_{L^2} \Big\{ C(m', \varpi_0) \| f \|_{L^2 (\nu)}\\
			& + ( m_0^{\gamma + 3} + e^{- c' (m')^2} + C (m') \varpi_0 ) \sum_{|\alpha'| = |\alpha|} \| \langle v \rangle^{ \gamma ( l + \frac{1}{2} ) } \M^{- \frac{q}{2}} \partial_v^{\alpha'} f \|_{L^2} \Big\}
		\end{aligned}
	\end{equation}
	for any small $m_0, \varpi_0 > 0$ and large $m' > 0$.
	
	We then control the quantity $\langle \langle v \rangle^{2 \gamma l} \M^{ - q } \partial_v^\alpha ( K_2 f ) , g \rangle$. Let
	\begin{equation*}
		\begin{aligned}
			V = v_1 - v \,, \ V_\shortparallel = (V \cdot \omega) \omega \,,  V_\perp = V - V_\shortparallel \,, \omega = \tfrac{v' - v}{|v' - v|} \,,
		\end{aligned}
	\end{equation*}
	whose intuitively geometric relations are illustrated in Figure \ref{Fig1} before. As shown in Section 2 of \cite{Grad-1963}, one has
	\begin{equation*}
		\begin{aligned}
			\d \omega \d v_1 = \d \omega \d V = \tfrac{2 \d V_\perp \d V_\shortparallel}{|V_\shortparallel|^2} \,.
		\end{aligned}
	\end{equation*}
	Then the $K_2 f$ given in \eqref{K2} can be represented by
	\begin{equation*}
		\begin{aligned}
			K_2 f (v) = 4 \int_{\R^3} \int_{V_\perp \perp V_\shortparallel} |V_\shortparallel|^{-2} f (v+V_\shortparallel) \M^\frac{1}{2} (v+V_\perp) \M^\frac{1}{2} (v + V) b (\omega, V) \d V_\perp \d V_\shortparallel \,,
		\end{aligned}
	\end{equation*}
	which means that
	\begin{equation}\label{Star3}
		\begin{aligned}
			\partial_v^\alpha ( K_2 f ) = 4 \sum_{\alpha' \leq \alpha} C_\alpha^{\alpha'} \int_{\R^3} \int_{V_\perp \perp V_\shortparallel} & |V_\shortparallel|^{-2} \partial_v^{\alpha'} f (v+V_\shortparallel) \\
			& \times \partial_v^{\alpha - \alpha'} [ \M^\frac{1}{2} (v+V_\perp) \M^\frac{1}{2} (v + V) ] b (\omega, V) \d V_\perp \d V_\shortparallel \,.
		\end{aligned}
	\end{equation}
	Note that $| \partial_v^\alpha \M^\frac{1}{2} (v) | \lesssim \langle v \rangle^{|\alpha|} \M^\frac{1}{2} (v)$, which implies
	\begin{equation*}
		\begin{aligned}
			& | \partial_v^{\alpha - \alpha'} [ \M^\frac{1}{2} (v+V_\perp) \M^\frac{1}{2} (v + V) ] | \\
			\lesssim & \sum_{\alpha'' \leq \alpha - \alpha'} \langle v + V_\perp \rangle^{ |\alpha''| } \langle v + V \rangle^{ | \alpha - \alpha' - \alpha'' | } \M^\frac{1}{2} (v+V_\perp) \M^\frac{1}{2} (v + V) \,.
		\end{aligned}
	\end{equation*}
	Then one establishes
	\begin{equation}\label{K2-mpn}
		\begin{aligned}
			| \partial_v^\alpha ( K_2 f ) | \lesssim & \sum_{ \substack{ \alpha' \leq \alpha \\ \alpha'' \leq \alpha - \alpha' } } \int_{\R^3} \int_{V_\perp \perp V_\shortparallel} 2 |V_\shortparallel|^{-2} | \partial_v^{\alpha'} f (v+V_\shortparallel) | b (\omega, V) \\
			& \times \langle v + V_\perp \rangle^{ |\alpha''| } \langle v + V \rangle^{ | \alpha - \alpha' - \alpha'' | } \M^\frac{1}{2} (v+V_\perp) \M^\frac{1}{2} (v + V) \d V_\perp \d V_\shortparallel \\
			= & \sum_{ \substack{ \alpha' \leq \alpha \\ \alpha'' \leq \alpha - \alpha' } } \iint_{\R^3 \times \mathbb{S}^2} | \partial_v^{\alpha'} f (v') | \langle v_1' \rangle^{ |\alpha''| } \langle v_1 \rangle^{ | \alpha - \alpha' - \alpha'' | } \\
			& \qquad \qquad \qquad \qquad \qquad \times \M^\frac{1}{2} (v_1') \M^\frac{1}{2} (v_1) b (\omega, v_1 - v) \d \omega \d v_1 \\
			: = & \widetilde{K}_2^{m, -} + \widetilde{K}_2^{m, +} \,,
		\end{aligned}
	\end{equation}
	where
	\begin{equation*}
		\begin{aligned}
			\widetilde{K}_2^{m, \pm} = \sum_{ \substack{ \alpha' \leq \alpha \\ \alpha'' \leq \alpha - \alpha' } } \iint_{\R^3 \times \mathbb{S}^2} \mathbf{1}_{\pm ( |v_1 - v| - m_0 ) \geq 0 } | \partial_v^{\alpha'} f (v') | \langle v_1' \rangle^{ |\alpha''| } \langle v_1 \rangle^{ | \alpha - \alpha' - \alpha'' | } \\
			\times \M^\frac{1}{2} (v_1') \M^\frac{1}{2} (v_1) b (\omega, v_1 - v) \d \omega \d v_1 \,.
		\end{aligned}
	\end{equation*}
	Here $m_0 > 0$ is small to be determined later.
	
	We now estimate the quantity $\langle \langle v \rangle^{2 \gamma l} \M^{ - q } \widetilde{K}_2^{m, -} , g \rangle$. It is easy to see that
	\begin{equation*}
		\begin{aligned}
			& | \langle \langle v \rangle^{2 \gamma l} \M^{ - q } \widetilde{K}_2^{m, -} , g \rangle | \\
			\lesssim & \sum_{ \substack{ \alpha' \leq \alpha \\ \alpha'' \leq \alpha - \alpha' } } \iiint_{\R^3 \times \R^3 \times \mathbb{S}^2} \mathbf{1}_{|v_1 - v| \leq m_0} | \langle v \rangle^{ \gamma ( l + \frac{1}{2} ) } \M^{ - \frac{q}{2} } g (v) | \cdot | \langle v' \rangle^{ \gamma ( l + \frac{1}{2} ) } \M^{ - \frac{q}{2} } \partial_v^{\alpha'} f (v') | \\
			& \quad \times \tfrac{ \langle v \rangle^{\gamma ( l - \frac{1}{2} )} \M^{ - \frac{q}{2} } (v) }{ \langle v' \rangle^{\gamma ( l + \frac{1}{2} )} \M^{ - \frac{q}{2} } (v') } \langle v_1' \rangle^{ |\alpha''| } \langle v_1 \rangle^{ | \alpha - \alpha' - \alpha'' | } \M^\frac{1}{2} (v_1') \M^\frac{1}{2} (v_1) b ( v_1 - v , \omega ) \d \omega \d v_1 \d v \,.
		\end{aligned}
	\end{equation*}
	Observe that for $|v_1 - v| \leq m_0$,
	\begin{equation*}
		\begin{aligned}
			|v_1' - \u| = & | v_1 - [(v_1  - v) \cdot \omega] \omega - \u | = | v - \u + (v_1 - v) - [(v_1  - v) \cdot \omega] \omega | \\
			\geq & |v - \u| - | (v_1 - v) - [(v_1  - v) \cdot \omega] \omega | \geq |v - \u| - |v_1 - v| \geq |v - \u| - m_0 \,,
		\end{aligned}
	\end{equation*}
	and
	\begin{equation*}
		\begin{aligned}
			|v' - \u| = |v_1 - \u - (v_1 - v) + (v_1 - v) - [(v_1  - v) \cdot \omega] \omega| \geq |v_1 - \u| - |v_1 - v| \geq |v_1 - \u| - m_0 \,,
		\end{aligned}
	\end{equation*}
	where the relations \eqref{Colli-v} (or Figure \ref{Fig1}) have been used. Then for any $0 < \eta_0 < 1$,
	\begin{equation*}
		\begin{aligned}
			|v_1' - \u|^2 \geq (1 - \eta_0) |v - \u|^2 - \tfrac{m_0^2}{\eta_0 ( 1 - \eta_0 )} \,, \quad |v' - \u|^2 \geq (1 - \eta_0) |v_1 - \u|^2 - \tfrac{m_0^2}{\eta_0 ( 1 - \eta_0 )} \,,
		\end{aligned}
	\end{equation*}
	which tells us
	\begin{equation*}
		\begin{aligned}
			\M (v_1') \lesssim e^{ \frac{m_0^2}{2 T \eta_0 ( 1 - \eta_0 )} } \M^{1 - \eta_0} (v) \,, \quad \M (v') \lesssim e^{ \frac{m_0^2}{2 T \eta_0 ( 1 - \eta_0 )} } \M^{1 - \eta_0} (v_1) \,.
		\end{aligned}
	\end{equation*}
	Moreover, by \eqref{Colli-v}, one has
	\begin{equation*}
		\begin{aligned}
			\langle v' \rangle^{|\alpha''|} = \langle v + [ (v_1 - v) \cdot \omega ] \omega \rangle^{|\alpha''|} \lesssim \langle v \rangle^{|\alpha''|} \langle v_1 \rangle^{|\alpha''|}
		\end{aligned}
	\end{equation*}
	Furthermore, there holds
	\begin{equation*}
		\begin{aligned}
			\tfrac{ \langle v \rangle^{\gamma ( l - \frac{1}{2} )} }{ \langle v' \rangle^{\gamma ( l + \frac{1}{2} )} } \lesssim \langle v - v' \rangle^{|\gamma| ( |l| + \frac{1}{2} )} \lesssim \langle v \rangle^{|\gamma| ( |l| + \frac{1}{2} )} \langle v' \rangle^{|\gamma| ( |l| + \frac{1}{2} )} \lesssim \langle v \rangle^{|\gamma| ( 2 |l| + 1 )} \langle v_1 \rangle^{|\gamma| ( |l| + \frac{1}{2} )} \,.
		\end{aligned}
	\end{equation*}
	It thereby follows that
	\begin{equation*}
		\begin{aligned}
			& \tfrac{ \langle v \rangle^{\gamma ( l - \frac{1}{2} )} \M^{ - \frac{q}{2} } (v) }{ \langle v' \rangle^{\gamma ( l + \frac{1}{2} )} \M^{ - \frac{q}{2} } (v') } \langle v_1' \rangle^{ |\alpha''| } \langle v_1 \rangle^{ | \alpha - \alpha' - \alpha'' | } \M^\frac{1}{2} (v_1') \M^\frac{1}{2} (v_1) \\
			\lesssim & \langle v \rangle^{|\gamma| ( 2 |l| + 1 ) + |\alpha''|} \langle v_1 \rangle^{|\gamma| ( |l| + \frac{1}{2} ) + |\alpha - \alpha'|} \M^\frac{1 - q - \eta_0}{2} (v) \M^\frac{ 1 + q (1 - \eta_0) }{2} (v_1) \\
			\lesssim & [ \M (v) \M (v_1) ]^\frac{1 - q}{4}
		\end{aligned}
	\end{equation*}
	by taking $\eta_0 = 1 - q \in (0, 1)$. Then the H\"older inequality infers
	\begin{equation}\label{FgFf}
		\begin{aligned}
			| \langle \langle v \rangle^{2 \gamma l} \M^{ - q } \widetilde{K}_2^{m, -} , g \rangle | \lesssim \sum_{ \beta' \leq \beta } \sqrt{F_g F_f} \,,
		\end{aligned}
	\end{equation}
	where
	\begin{equation*}
		\begin{aligned}
			F_g = \iiint_{\R^3 \times \R^3 \times \mathbb{S}^2} \mathbf{1}_{|v_1 - v| \leq m_0} | \langle v \rangle^{\gamma ( l +  \frac{1}{2} )} \M^{ - \frac{q}{2} } g (v) |^2 [ \M (v) \M (v_1) ]^\frac{1 - q}{4} b ( v_1 - v , \omega ) \d \omega \d v_1 \d v \,,
		\end{aligned}
	\end{equation*}
	and
	\begin{equation*}
		\begin{aligned}
			F_f = \iiint_{\R^3 \times \R^3 \times \mathbb{S}^2} \mathbf{1}_{|v_1 - v| \leq m_0} | \langle v' \rangle^{\gamma ( l +  \frac{1}{2} )} \M^{ - \frac{q}{2} } \partial_v^{\alpha'} f (v') |^2 [ \M (v) \M (v_1) ]^\frac{1 - q}{4} b ( v_1 - v , \omega ) \d \omega \d v_1 \d v \,.
		\end{aligned}
	\end{equation*}
	
	Together with $\int_{\mathbb{S}^2} b ( v_1 - v , \omega ) \d \omega = |v_1 - v|^\gamma $ and $\int_{\R^3} \mathbf{1}_{|v_1 - v| \leq m} |v_1 - v|^\gamma \M^\frac{1-q}{4} (v_1) \d v_1 \lesssim m_0^{\gamma + 3}$, the symbol $F_g$ can be bounded by
	\begin{equation}\label{Fg-bnd}
		\begin{aligned}
			F_g \lesssim m_0^{\gamma + 3} \int_{\R^3} | \langle v \rangle^{\gamma ( l +  \frac{1}{2} )} \M^{ - \frac{q}{2} } g (v) |^2 \M^\frac{1-q}{4} (v) \d v \lesssim m_0^{\gamma + 3} \| \langle v \rangle^{\gamma ( l +  \frac{1}{2} )} \M^{ - \frac{q}{2} } g \|^2_{L^2} \,.
		\end{aligned}
	\end{equation}
	Observe that
	\begin{equation*}
		\begin{aligned}
			& \mathbf{1}_{|v_1 - v| \leq m_0} [ \M (v) \M (v_1) ]^\frac{1 - q}{4} b ( v_1 - v , \omega ) \d \omega \d v_1 \d v \\
			= & \mathbf{1}_{|v_1' - v'| \leq m_0} [ \M (v') \M (v_1') ]^\frac{1 - q}{4} b ( v_1' - v' , \omega ) \d \omega \d v_1' \d v'
		\end{aligned}
	\end{equation*}
	under the change of variables $(v_1, v) \mapsto (v_1', v')$. Then, combining with the estimate of \eqref{Fg-bnd}, one has
	\begin{equation}\label{Ff-bnd}
		\begin{aligned}
			F_f = & \iiint_{\R^3 \times \R^3 \times \mathbb{S}^2} \mathbf{1}_{|v_1 - v| \leq m_0} | \langle v \rangle^{\gamma ( l +  \frac{1}{2} )} \M^{ - \frac{q}{2} } \partial_v^{\alpha'} f (v) |^2 [ \M (v) \M (v_1) ]^\frac{1 - q}{4} b ( v_1 - v , \omega ) \d \omega \d v_1 \d v \\
			\lesssim & m_0^{\gamma + 3} \| \langle v \rangle^{\gamma ( l +  \frac{1}{2} )} \M^{ - \frac{q}{2} } \partial_v^{\alpha'} f \|^2_{L^2} \,.
		\end{aligned}
	\end{equation}
	It then follows from \eqref{FgFf}, \eqref{Fg-bnd} and \eqref{Ff-bnd} that
	\begin{equation}\label{K2-mn-bnd}
		\begin{aligned}
			| \langle \langle v \rangle^{2 \gamma l} \M^{ - q } \widetilde{K}_2^{m, -} , g \rangle | \lesssim m_0^{\gamma + 3} \sum_{ \alpha' \leq \alpha } \| \langle v \rangle^{\gamma ( l +  \frac{1}{2} )} \M^{ - \frac{q}{2} } \partial_v^{\alpha'} f \|_{L^2} \| \langle v \rangle^{\gamma ( l +  \frac{1}{2} )} \M^{ - \frac{q}{2} } g \|_{L^2} \,.
		\end{aligned}
	\end{equation}
	
	It remains to control the inner product $\langle \langle v \rangle^{2 \gamma l} \M^{ - q } \widetilde{K}_2^{m, +} , g \rangle$. We split
	\begin{equation*}
		\begin{aligned}
			\widetilde{K}_2^{m, +} = \widetilde{K}_{2,+}^{m, +} + \widetilde{K}_{2,-}^{m, +} \,,
		\end{aligned}
	\end{equation*}
	where
	\begin{equation*}
		\begin{aligned}
			\widetilde{K}_{2,\pm}^{m, +} = \sum_{ \substack{ \alpha' \leq \alpha \\ \alpha'' \leq \alpha - \alpha' } } \iint_{\R^3 \times \mathbb{S}^2} & \mathbf{1}_{|v_1 - v| \geq m_0} \mathbf{1}_{\pm ( 1 + |v| + |v'| - m' ) \geq 0} | \partial_v^{\alpha'} f (v') | \\
			& \times | \langle v_1' \rangle^{ |\alpha''| } \langle v_1 \rangle^{ | \alpha - \alpha' - \alpha'' | } \M^\frac{1}{2} (v_1') \M^\frac{1}{2} (v_1) b (\omega, v_1 - v) \d \omega \d v_1 \,.
		\end{aligned}
	\end{equation*}
	Here $m' > 0$ is sufficiently large to be determined later. Following the similar arguments in \eqref{W2-phi-bnd}, one easily has
	\begin{equation*}
		\begin{aligned}
			| \langle \langle v \rangle^{2 \gamma l} \M^{ - q } \widetilde{K}_{2, +}^{m, +} , g \rangle | \lesssim \tfrac{1}{1 + m'} \sum_{ \alpha' \leq \alpha } \| \langle v \rangle^{\gamma ( l +  \frac{1}{2} )} \M^{ - \frac{q}{2} } \partial_v^{\alpha'} f \|_{L^2} \| \langle v \rangle^{\gamma ( l +  \frac{1}{2} )} \M^{ - \frac{q}{2} } g \|_{L^2} \,.
		\end{aligned}
	\end{equation*}
	Moreover, from the similar compact argument of $\partial_\beta [ K_2^{1 - \Upsilon} g_1 ]$ in the end of proof of Lemma 2 of \cite{Strain-Guo-2008-ARMA}, one has
	\begin{equation*}
		\begin{aligned}
			| \langle \langle v \rangle^{2 \gamma l} & \M^{ - q } \widetilde{K}_{2, -}^{m, +} , g \rangle | \\
			& \lesssim \Big\{ \epsilon_0 \sum_{|\alpha'| = |\alpha|} \| \langle v \rangle^{\gamma ( l +  \frac{1}{2} )} \M^{ - \frac{q}{2} } \partial_v^{\alpha'} f \|_{L^2} + C ( \epsilon_0, m' ) \| f \|_{ L^2 (\nu) } \Big\} \| \langle v \rangle^{\gamma ( l +  \frac{1}{2} )} \M^{ - \frac{q}{2} } g \|_{L^2}
		\end{aligned}
	\end{equation*}
	for any small $\epsilon_0 > 0$ to be determined later. Then one establishes that
	\begin{equation}\label{K2-mp-bnd}
		\begin{aligned}
			& | \langle \langle v \rangle^{2 \gamma l} \M^{ - q } \widetilde{K}_2^{m, +} , g \rangle | \\
			\lesssim & \Big\{ ( \epsilon_0 + \tfrac{1}{1 + m'} ) \sum_{|\alpha'| \leq |\alpha|} \| \langle v \rangle^{\gamma ( l +  \frac{1}{2} )} \M^{ - \frac{q}{2} } \partial_v^{\alpha'} f \|_{L^2} + C ( \epsilon_0, m' ) \| f \|_{ L^2 (\nu) } \Big\} \| \langle v \rangle^{\gamma ( l +  \frac{1}{2} )} \M^{ - \frac{q}{2} } g \|_{L^2}
		\end{aligned}
	\end{equation}
	for small $\epsilon_0 > 0$ and large $m' > 0$. As a result, by \eqref{K2-mpn}, \eqref{K2-mn-bnd} and \eqref{K2-mp-bnd},
	\begin{equation}\label{K2-bnd}
		\begin{aligned}
			& | \langle \langle v \rangle^{2 \gamma l} \M^{ - q } \partial_v^\alpha ( K_2 f ) , g \rangle | \\
			\lesssim & ( m_0^{\gamma + 3} + \epsilon_0 + \tfrac{1}{1 + m'} ) \sum_{|\alpha'| \leq |\alpha|} \| \langle v \rangle^{\gamma ( l +  \frac{1}{2} )} \M^{ - \frac{q}{2} } \partial_v^{\alpha'} f \|_{L^2} \| \langle v \rangle^{\gamma ( l +  \frac{1}{2} )} \M^{ - \frac{q}{2} } g \|_{L^2} \\
			& + C ( \epsilon_0, m' ) \| f \|_{ L^2 (\nu) } \| \langle v \rangle^{\gamma ( l +  \frac{1}{2} )} \M^{ - \frac{q}{2} } g \|_{L^2}
		\end{aligned}
	\end{equation}
	for small $\epsilon_0, m_0 > 0$ and large $m' > 0$. For any $\eta \in (0, 1)$ fixed, we first take $\epsilon_0, m_0, m'$ such that
	\begin{equation*}
		\begin{aligned}
			\tfrac{1}{2} \eta = e^{- c' (m')^2} + m_0^{\gamma + 3} + \epsilon_0 + \tfrac{1}{1 + m'} \,.
		\end{aligned}
	\end{equation*}
	We then choose small $\varpi_0 > 0$ such that $ C (m') \varpi_0 = \frac{1}{2} \eta$. So,
	\begin{equation*}
		\begin{aligned}
			\eta = e^{- c' (m')^2} + m_0^{\gamma + 3} + \epsilon_0 + \tfrac{1}{1 + m'} + C (m') \varpi_0 \,.
		\end{aligned}
	\end{equation*}
	Then the bounds \eqref{K1-bnd} and \eqref{K2-bnd} conclude the inequality \eqref{L2-fg2}. The proof of Lemma \ref{Lmm-VSP} is finished.

\section{Estimate of weighted $L^2$ norm of $\partial_v^\alpha \partial_{t,x}^\beta f$: Proof of Lemma \ref{Lmm-SVP-Der}}\label{Sec:Lmm4.1}

In this section, the main goal is to control the norm $\| \langle v \rangle^{ \gamma (l + \frac{1}{2}) } \M^{ - \frac{q}{2} } \partial_v^\alpha \partial_{t,x}^\beta f \|_{L^2}$, i.e., give the proof of Lemma \ref{Lmm-SVP-Der}. The key is to control the commutator $[ \partial_{t,x}^\beta, \L ] f$.

Note that
\begin{equation}
	\begin{aligned}
		\partial_{t,x}^\beta \L f = \L \partial_{t,x}^\beta f + [ \partial_{t,x}^\beta , \L ] f = \partial_{t,x}^\beta g \,.
	\end{aligned}
\end{equation}
where $[X, Y] = XY - YX$ is the commutator operator. Then
\begin{equation}\label{Finish1}
	\begin{aligned}
		\langle \langle v \rangle^{  2 \gamma l} \M^{ - q } \partial_v^\alpha \partial_{t,x}^\beta \L f , \partial_v^\alpha \partial_{t,x}^\beta f \rangle = & \langle \langle v \rangle^{  2 \gamma l} \M^{ - q } \partial_v^\alpha  \L \partial_{t,x}^\beta f , \partial_v^\alpha \partial_{t,x}^\beta f \rangle \\
		& + \langle \langle v \rangle^{  2 \gamma l} \M^{ - q } \partial_v^\alpha [ \partial_{t,x}^\beta , \L ] f , \partial_v^\alpha \partial_{t,x}^\beta f \rangle \,.
	\end{aligned}
\end{equation}
Recall that $\L = \nu - K$. Lemma \ref{Lmm-VSP} therefore reduces to
\begin{equation}\label{Finish2}
	\begin{aligned}
		\langle \langle v \rangle^{  2 \gamma l} \M^{ - q } \partial_v^\alpha \L \partial_{t,x}^\beta f , \partial_v^\alpha \partial_{t,x}^\beta f \rangle \geq & ( 1 - 3 \eta ) \| \langle v \rangle^{ \gamma ( l + \frac{1}{2} ) } \M^{ - \frac{q}{2} } \partial_v^\alpha \partial_{t,x}^\beta f \|^2_{L^2} \\
		& - C_\eta \sum_{\alpha' < \alpha} \| \langle v \rangle^{ \gamma ( l + \frac{1}{2} ) } \M^{ - \frac{q}{2} } \partial_v^{\alpha'} \partial_{t,x}^\beta f \|^2_{L^2} - C_\eta \| \partial_{t,x}^\beta f \|^2_{L^2}
	\end{aligned}
\end{equation}
for any sufficiently small $\eta > 0$.

It remains to deal with the commutator part $ \langle \langle v \rangle^{  2 \gamma l} \M^{ - q } \partial_v^\alpha [ \partial_{t,x}^\beta , \L ] f , \partial_v^\alpha \partial_{t,x}^\beta f \rangle $. Observe that $\L = \nu - K_1 + K_2$. It thereby holds
\begin{equation}\label{Finish3}
	\begin{aligned}
		\langle \langle v \rangle^{  2 \gamma l} \M^{ - q } \partial_v^\alpha [ \partial_{t,x}^\beta , \L ] f , \partial_v^\alpha \partial_{t,x}^\beta f \rangle = & \langle \langle v \rangle^{  2 \gamma l} \M^{ - q } \partial_v^\alpha [ \partial_{t,x}^\beta , \nu ] f , \partial_v^\alpha \partial_{t,x}^\beta f \rangle \\
		& - \langle \langle v \rangle^{  2 \gamma l} \M^{ - q } \partial_v^\alpha [ \partial_{t,x}^\beta , K_1 ] f , \partial_v^\alpha \partial_{t,x}^\beta f \rangle \\
		& + \langle \langle v \rangle^{  2 \gamma l} \M^{ - q } \partial_v^\alpha [ \partial_{t,x}^\beta , K_2 ] f , \partial_v^\alpha \partial_{t,x}^\beta f \rangle \,.
	\end{aligned}
\end{equation}

\textbf{\em Step 1. Control of $ \langle \langle v \rangle^{  2 \gamma l} \M^{ - q } \partial_v^\alpha [ \partial_{t,x}^\beta , \nu ] f , \partial_v^\alpha \partial_{t,x}^\beta f \rangle $.} A straightforward calculation shows
\begin{equation*}
	\begin{aligned}
		\partial_v^\alpha [ \partial_{t,x}^\beta , \nu ] f = \sum_{ \alpha' \leq \alpha } \sum_{0 \neq \beta' \leq \beta} C_\alpha^{\alpha'} C_\beta^{\beta'} \partial_v^{ \alpha - \alpha' } \partial_{t,x}^{\beta'} \nu \partial_v^{\alpha'} \partial_{t,x}^{ \beta - \beta' } f \,.
	\end{aligned}
\end{equation*}
Due to $\nu = \int_{\R^3} |v_1 - v|^\gamma \M(v_1) \d v_1 = \int_{\R^3} |v_1|^\gamma \M ( v + v_1 ) \d v_1$, one has
\begin{equation*}
	\begin{aligned}
		\partial_v^{ \alpha - \alpha' } \partial_{t,x}^{\beta'} \nu = \int_{\R^3} |v_1|^\gamma \partial_v^{ \alpha - \alpha' } \partial_{t,x}^{\beta'} \M ( v + v_1 ) \d v_1 \,.
	\end{aligned}
\end{equation*}
It is easy to see
\begin{equation*}
	\begin{aligned}
		| \partial_v^{ \alpha - \alpha' } \partial_{t,x}^{\beta'} \M ( v + v_1 ) | \lesssim \mathfrak{e}_{|\beta'|} \langle v + v_1 \rangle^{ |\alpha - \alpha'| + 2 |\beta'| } \M (v + v_1) \lesssim \mathfrak{e}_m \M^\frac{1}{2} (v + v_1) \,,
	\end{aligned}
\end{equation*}
which infers
\begin{equation*}
	\begin{aligned}
		| \partial_v^{ \alpha - \alpha' } \partial_{t,x}^{\beta'} \nu | \lesssim \mathfrak{e}_m \int_{\R^3} |v_1|^\gamma \M^\frac{1}{2} ( v + v_1 ) \d v_1 = \mathfrak{e}_m \int_{\R^3} |v - v_1|^\gamma \M^\frac{1}{2} ( v_1 ) \d v_1 \lesssim \mathfrak{e}_m \nu (v) \,,
	\end{aligned}
\end{equation*}
where the last inequality is derived from Lemma 2.1 of \cite{LS-2010-KRM}. As a consequence,
\begin{equation}\label{Triangle1}
	\begin{aligned}
		| \partial_v^\alpha [ \partial_{t,x}^\beta , \nu ] f | \lesssim \mathfrak{e}_m \sum_{ \alpha' \leq \alpha } \sum_{0 \neq \beta' \leq \beta} \nu | \partial_v^{\alpha'} \partial_{t,x}^{ \beta - \beta' } f | \,,
	\end{aligned}
\end{equation}
which implies
\begin{equation}\label{Finish4}
	\begin{aligned}
		& | \langle \langle v \rangle^{  2 \gamma l} \M^{ - q } \partial_v^\alpha [ \partial_{t,x}^\beta , \nu ] f , \partial_v^\alpha \partial_{t,x}^\beta f \rangle | \\
		\lesssim & \mathfrak{e}_m \sum_{ \alpha' \leq \alpha } \sum_{0 \neq \beta' \leq \beta} \langle \langle v \rangle^{ 2 \gamma l } \M^{ - q } \nu | \partial_v^{\alpha'} \partial_{t,x}^{ \beta - \beta' } f | , | \partial_v^\alpha \partial_{t,x}^\beta f | \rangle \\
		\lesssim & \mathfrak{e}_m \sum_{ \alpha' \leq \alpha } \sum_{0 \neq \beta' \leq \beta} \| \langle v \rangle^{ \gamma ( l + \frac{1}{2} ) } \M^{ - \frac{q}{2} } \partial_v^{\alpha'} \partial_{t,x}^{ \beta - \beta' } f \|_{ L^2 } \| \langle v \rangle^{ \gamma ( l + \frac{1}{2} ) } \M^{ - \frac{q}{2} } \partial_v^\alpha \partial_{t,x}^\beta f \|_{ L^2 } \,,
	\end{aligned}
\end{equation}
Here the fact $\nu (v) \lesssim \mathfrak{e}_m \langle v \rangle^\gamma$ has been used.

\textbf{\em Step 2. Control of $ \langle \langle v \rangle^{  2 \gamma l} \M^{ - q } \partial_v^\alpha [ \partial_{t,x}^\beta , K_1 ] f , \partial_v^\alpha \partial_{t,x}^\beta f \rangle $.} By \eqref{Star1}, one easily has
\begin{equation}\label{Triangle2}
	\begin{aligned}
		\partial_v^\alpha [ \partial_{t,x}^\beta , K_1 ] f = \sum_{ \substack{ \alpha' \leq \alpha \\ 0 \neq \beta' \leq \beta } } C_\alpha^{\alpha'} C_\beta^{\beta'} \int_{\R^3} |v_1|^\gamma \partial_v^{ \alpha - \alpha' } \partial_{t,x}^{ \beta - \beta' } f ( t, x, v + v_1 ) \\
		\times \partial_v^{\alpha'} \partial_{t,x}^{ \beta' } [ \M^\frac{1}{2} (v) \M^\frac{1}{2} ( v + v_1 ) ] \d v_1 \,.
	\end{aligned}
\end{equation}
Note that
\begin{equation}\label{Triangle3}
	\begin{aligned}
		& | \partial_v^{\alpha'} \partial_{t,x}^{ \beta' } [ \M^\frac{1}{2} (v) \M^\frac{1}{2} ( v + v_1 ) ] | \\
		\lesssim & \mathfrak{e}_m^m \sum_{ \substack{ \alpha'' \leq \alpha' \\ \beta'' \leq \beta' } } \langle v \rangle^{ |\alpha''| + 2 |\beta''| } \langle v + v_1 \rangle^{ | \alpha' - \alpha'' | + 2 | \beta' - \beta'' | } \M^\frac{1}{2} (v) \M^\frac{1}{2} ( v + v_1 )
	\end{aligned}
\end{equation}
and
\begin{equation*}
	\begin{aligned}
		\sum_{ \substack{ \alpha'' \leq \alpha' \\ \beta'' \leq \beta' } } \frac{ \langle v \rangle^{ \gamma l + |\alpha''| + 2 |\beta''| - \frac{1}{2} \gamma } }{ \langle v_1 \rangle^{ \gamma ( l + \frac{1}{2} ) + |\alpha' - \alpha''| + 2 | \beta' - \beta'' | } } [ \M (v) \M (v_1) ]^\frac{1 - q}{2} \lesssim \mathfrak{e}_m [ \M (v) \M (v_1) ]^\frac{1 - q}{4} \,.
	\end{aligned}
\end{equation*}
Then the similar arguments in \eqref{Star2} show that
\begin{equation*}
	\begin{aligned}
		& | \langle \langle v \rangle^{  2 \gamma l} \M^{ - q } \partial_v^\alpha [ \partial_{t,x}^\beta , K_1 ] f , \partial_v^\alpha \partial_{t,x}^\beta f \rangle | \\
		\lesssim & \mathfrak{e}_m^m \sum_{ \substack{ \alpha' \leq \alpha \\ 0 \neq \beta' \leq \beta } } \iint_{\R^3 \times \mathbb{S}^2} |v - v_1|^\gamma [ \M (v) \M (v_1) ]^\frac{1 - q}{4} | \langle v_1 \rangle^{ \gamma ( l + \frac{1}{2} ) } \M^{ - \frac{q}{2} } \partial_v^{ \alpha - \alpha' } \partial_{t,x}^{ \beta - \beta' } f ( t, x, v_1 ) | \\
		& \qquad \qquad \qquad \qquad \qquad \qquad \qquad \qquad \quad \times | \langle v \rangle^{ \gamma ( l + \frac{1}{2} ) } \M^{ - \frac{q}{2} } \partial_v^\alpha \partial_{t,x}^\beta f ( t, x, v ) | \d v_1 \d v \,.
	\end{aligned}
\end{equation*}
The same arguments of \eqref{K1-bnd} then tell us
\begin{equation}\label{Finish5}
	\begin{aligned}
		& | \langle \langle v \rangle^{  2 \gamma l} \M^{ - q } \partial_v^\alpha [ \partial_{t,x}^\beta , K_1 ] f , \partial_v^\alpha \partial_{t,x}^\beta f \rangle | \\
		\lesssim & \mathfrak{e}_m^m \| \langle v \rangle^{ \gamma ( l + \frac{1}{2} ) } \M^{ - \frac{q}{2} } \partial_v^\alpha \partial_{t,x}^\beta f \|_{L^2} \sum_{ 0 \neq \beta' \leq \beta } \Big\{ C (m', \varpi_0) \| \partial_{t,x}^{\beta - \beta'} f \|_{ L^2 ( \nu ) } \\
		& + ( m_0^{\gamma + 3} + e^{ - c' (m')^2 } + C (m') \varpi_0 ) \sum_{|\alpha'| = |\alpha|} \| \langle v \rangle^{ \gamma ( l + \frac{1}{2} ) } \M^{ - \frac{q}{2} } \partial_v^{ \alpha' } \partial_{t,x}^{ \beta - \beta' } f \|_{ L^2 } \Big\}
	\end{aligned}
\end{equation}
for any small $m_0, \varpi_0 > 0$ and large $m' > 0$.

\textbf{\em Step 3. Control of $ \langle \langle v \rangle^{  2 \gamma l} \M^{ - q } \partial_v^\alpha [ \partial_{t,x}^\beta , K_2 ] f , \partial_v^\alpha \partial_{t,x}^\beta f \rangle $.} Following the similar derivations in \eqref{Star3}, one has
\begin{equation*}
	\begin{aligned}
		\partial_v^\alpha [ \partial_{t,x}^\beta , K_2 ] f = 4 \sum_{ \substack{ \alpha' \leq \alpha \\ 0 \neq \beta' \leq \beta } } C_\alpha^{\alpha'} C_\beta^{\beta'} \int_{\R^3} \int_{V_\perp \perp V_\shortparallel} | V_\shortparallel |^{-2} \partial_v^{ \alpha' } \partial_{t,x}^{ \beta - \beta' } f ( t, x, v + V_\shortparallel ) \\
		\times \partial_v^{ \alpha - \alpha' } \partial_{t,x}^{ \beta' } [ \M^\frac{1}{2} ( v + V_\perp ) \M^\frac{1}{2} ( v + V ) ] b ( \omega, V ) \d V_\perp \d V_\shortparallel \,.
	\end{aligned}
\end{equation*}
Observe that
\begin{equation*}
	\begin{aligned}
		& | \partial_v^{ \alpha - \alpha' } \partial_{t,x}^{ \beta' } [ \M^\frac{1}{2} ( v + V_\perp ) \M^\frac{1}{2} ( v + V ) ] | \\
		\lesssim & \mathfrak{e}_m^m \sum_{ \substack{ \alpha'' \leq \alpha - \alpha' \\ \beta'' \leq \beta' } } \langle v + V_\perp \rangle^{ |\alpha''| + 2 |\beta''| } \langle v + V \rangle^{ | \alpha - \alpha' - \alpha'' | + 2 |\beta' - \beta''| } \M^\frac{1}{2} ( v + V_\perp ) \M^\frac{1}{2} ( v + V ) \,.
	\end{aligned}
\end{equation*}
The same arguments in \eqref{K2-mpn} then indicate that
\begin{equation}\label{Triangle4}
	\begin{aligned}
		| \partial_v^\alpha [ \partial_{t,x}^\beta , K_2 ] f | \lesssim \mathfrak{e}_m^m \sum_{ \substack{ \alpha' \leq \alpha \\ \alpha'' \leq \alpha - \alpha' } } \sum_{ \substack{ 0 \neq \beta' \leq \beta \\ \beta'' \leq \beta } } \iint_{\R^3 \times \mathbb{S}^2} | \partial_v^{ \alpha' } \partial_{t,x}^{ \beta - \beta' } | \langle v_1 \rangle^{ | \alpha - \alpha' - \alpha'' | + 2 |\beta' - \beta''| } \\
		\times \langle v_1' \rangle^{ |\alpha''| + 2 | \beta'' | } \M^\frac{1}{2} (v_1') \M^\frac{1}{2} (v_1) b ( \omega, v_1 - v ) \d \omega \d v_1 \,.
	\end{aligned}
\end{equation}
It then follows from the same estimates in \eqref{K2-bnd} that
\begin{equation}\label{Finish6}
	\begin{aligned}
		& | \langle \langle v \rangle^{  2 \gamma l} \M^{ - q } \partial_v^\alpha [ \partial_{t,x}^\beta , K_2 ] f , \partial_v^\alpha \partial_{t,x}^\beta f \rangle | \\
		\lesssim & \mathfrak{e}_m^m \| \langle v \rangle^{ \gamma ( l + \frac{1}{2} ) } \M^{ - \frac{q}{2} } \partial_v^\alpha \partial_{t,x}^\beta f \|_{L^2} \sum_{ 0 \neq \beta' \leq \beta } \Big\{ C (\epsilon_0, m') \| \partial_{t,x}^{\beta - \beta'} f \|_{ L^2 ( \nu ) } \\
		& + ( m_0^{\gamma + 3} + \epsilon_0 + \tfrac{1}{m' + 1} ) \sum_{|\alpha'| \leq |\alpha|} \| \langle v \rangle^{ \gamma ( l + \frac{1}{2} ) } \M^{ - \frac{q}{2} } \partial_v^{ \alpha' } \partial_{t,x}^{ \beta - \beta' } f \|_{ L^2 } \Big\}
	\end{aligned}
\end{equation}
for small $m_0, \epsilon_0 > 0$ and large $m' > 0$.

\textbf{\em Step 4. Control of the quantity $ \sum_{ \beta' \leq \beta } \| \partial_{t,x}^{\beta'} f \|^2_{ L^2 ( \nu ) } $.} Without loss of generality, we only need to estimate the norm $\| \partial_{t,x}^\beta f \|^2_{ L^2 ( \nu ) } $. Since $ \L \partial_{t,x}^\beta f = \partial_{t,x}^\beta g - [ \partial_{t,x}^\beta , \L ] f $, the hypocoercivity \eqref{Hypo-L} shows that
\begin{equation}\label{Bracket1}
	\begin{aligned}
		\lambda \| \partial_{t,x}^\beta f \|^2_{ L^2 ( \nu ) } \leq \langle \L \partial_{t,x}^\beta f, \partial_{t,x}^\beta f \rangle = \langle \partial_{t,x}^\beta g , \partial_{t,x}^\beta g \rangle - \langle [ \partial_{t,x}^\beta , \L ] f , \partial_{t,x}^\beta f \rangle \,.
	\end{aligned}
\end{equation}
The H\"older inequality tells us
\begin{equation}\label{Bracket2}
	\begin{aligned}
		| \langle \partial_{t,x}^\beta g , \partial_{t,x}^\beta f \rangle | \leq \| \partial_{t,x}^\beta f \|_{ L^2 ( \nu ) } \| \nu^{ - \frac{1}{2} } \partial_{t,x}^\beta g \|_{L^2} \,.
	\end{aligned}
\end{equation}

It remains to estimate the quantity $ \langle [ \partial_{t,x}^\beta , \L ] f , \partial_{t,x}^\beta f \rangle $. Observe that
\begin{equation}\label{DbStar0}
	\begin{aligned}
		[ \partial_{t,x}^\beta , \L ] f = [ \partial_{t,x}^\beta , \nu ] f - [ \partial_{t,x}^\beta , K_1 ] f + [ \partial_{t,x}^\beta , K_2 ] f \,.
	\end{aligned}
\end{equation}
By \eqref{Triangle1} with $\alpha = 0$, one has
\begin{equation*}
	\begin{aligned}
		| [ \partial_{t,x}^\beta , \nu ] f | \lesssim \mathfrak{e}_m^m \sum_{ 0 \neq \beta' \leq \beta } \nu | \partial_{t,x}^{ \beta - \beta' } f | \,.
	\end{aligned}
\end{equation*}
Then
\begin{equation}\label{DbStar1}
	\begin{aligned}
		| \langle [ \partial_{t,x}^\beta , \nu ] f , \partial_{t,x}^\beta f \rangle | \lesssim \mathfrak{e}_m^m \sum_{ 0 \neq \beta' \leq \beta } \| \partial_{t,x}^{ \beta - \beta' } f \|_{ L^2 ( \nu ) } \| \partial_{t,x}^\beta f \|_{ L^2 ( \nu ) } \,.
	\end{aligned}
\end{equation}

By \eqref{Triangle2} with $\alpha = 0$, it infers
\begin{equation*}
	\begin{aligned}
		[ \partial_{t,x}^\beta , K_1 ] f = \sum_{ 0 \neq \beta' \leq \beta } C_\beta^{\beta'} \int_{\R^3} |v_1|^\gamma \partial_{t,x}^{ \beta - \beta' } f ( t, x, v + v_1 ) \partial_{t,x}^{ \beta' } [ \M^\frac{1}{2} (v) \M^\frac{1}{2} ( v + v_1 ) ] \d v_1 \,.
	\end{aligned}
\end{equation*}
The \eqref{Triangle3} then shows that for any small $\epsilon > 0$,
\begin{equation*}
	\begin{aligned}
		| [ \partial_{t,x}^\beta , K_1 ] f | \lesssim & \mathfrak{e}_m^m \sum_{ 0 \neq \beta' \leq \beta } \int_{\R^3} |v_1|^\gamma | \partial_{ t, x }^{ \beta - \beta' } f ( t, x, v + v_1 ) | [ \M (v) \M (v + v_1) ]^{ \frac{1}{2} - \epsilon } \d v_1 \\
		= & \mathfrak{e}_m^m \sum_{ 0 \neq \beta' \leq \beta } \underbrace{ \int_{\R^3} |v_1 - v|^\gamma | \partial_{ t, x }^{ \beta - \beta' } f ( t, x, v_1 ) | [ \M (v) \M (v_1) ]^{ \frac{1}{2} - \epsilon } \d v_1 }_{:= K_1^- ( | \partial_{ t, x }^{ \beta - \beta' } f | )} \,,
	\end{aligned}
\end{equation*}
where the operator $K_1^-$ is almost the same as the operator $K_1$ (just replacing the factor $[ \M (v) \M (v_1) ]^\frac{1}{2}$ in $K_1$ by $[ \M (v) \M (v_1) ]^{ \frac{1}{2} - \epsilon }$ for any small $\epsilon > 0$). Then by the boundedness arguments of $K_1$ (see \cite{Guo-CPAM06} for instance),
\begin{equation}\label{DbStar2}
	\begin{aligned}
		| \langle [ \partial_{t,x}^\beta , K_1 ] f , \partial_{t,x}^\beta f \rangle | \lesssim \mathfrak{e}_m^m \sum_{ 0 \neq \beta' \leq \beta } \| \partial_{t,x}^{ \beta - \beta' } f \|_{ L^2 ( \nu ) } \| \partial_{t,x}^\beta f \|_{ L^2 ( \nu ) } \,.
	\end{aligned}
\end{equation}

From the same arguments in \eqref{Triangle4} with $\alpha = 0$, it follows that for any small $\epsilon > 0$,
\begin{equation*}
	\begin{aligned}
		| [ \partial_{t,x}^\beta , K_2 ] f | \lesssim \mathfrak{e}_m^m \sum_{ 0 \neq \beta' \leq \beta } \underbrace{ \iint_{\R^3 \times \mathbb{S}^2} | \partial_{ t, x }^{ \beta - \beta' } f ( t, x, v' ) | [ \M (v_1') \M (v_1) ]^{ \frac{1}{2} - \epsilon } b ( \omega , v_1 - v ) \d \omega \d v_1 }_{:= K_2^- ( | \partial_{ t, x }^{ \beta - \beta' } f | )} \,,
	\end{aligned}
\end{equation*}
where the operator $K_2^-$ is almost the same as the operator $K_2$ (just replacing the factor $[ \M (v_1') \M (v_1) ]^\frac{1}{2}$ in $K_2$ by $[ \M (v_1') \M (v_1) ]^{ \frac{1}{2} - \epsilon }$ for any small $\epsilon > 0$). It thereby follows from the boundedness arguments of $K_2$ (see \cite{Guo-CPAM06} for instance) that
\begin{equation}\label{DbStar3}
	\begin{aligned}
		| \langle [ \partial_{t,x}^\beta , K_2 ] f , \partial_{t,x}^\beta f \rangle | \lesssim \mathfrak{e}_m^m \sum_{ 0 \neq \beta' \leq \beta } \| \partial_{t,x}^{ \beta - \beta' } f \|_{ L^2 ( \nu ) } \| \partial_{t,x}^\beta f \|_{ L^2 ( \nu ) } \,.
	\end{aligned}
\end{equation}
Therefore, the relations \eqref{DbStar0}, \eqref{DbStar1}, \eqref{DbStar2} and \eqref{DbStar3} imply that
\begin{equation}\label{Bracket3}
	\begin{aligned}
		| \langle [ \partial_{t,x}^\beta , \L ] f , \partial_{t,x}^\beta f \rangle | \lesssim \mathfrak{e}_m^m \sum_{ 0 \neq \beta' \leq \beta } \| \partial_{t,x}^{ \beta - \beta' } f \|_{ L^2 ( \nu ) } \| \partial_{t,x}^\beta f \|_{ L^2 ( \nu ) } \,.
	\end{aligned}
\end{equation}
As a result, \eqref{Bracket1}, \eqref{Bracket2} and \eqref{Bracket3} reduce to
\begin{equation*}
	\begin{aligned}
		\| \partial_{ t, x }^\beta f \|^2_{ L^2 (\nu) } \lesssim \mathfrak{e}_m^{2 m} \sum_{ 0 \neq \beta' \leq \beta } \| \partial_{t,x}^{ \beta - \beta' } f \|^2_{ L^2 ( \nu ) } + \| \nu^{ - \frac{1}{2} } \partial_{t,x}^\beta g \|^2_{L^2} \,,
	\end{aligned}
\end{equation*}
which, by employing the Induction Principle and combining with \eqref{L2-5}, results to
\begin{equation*}
	\begin{aligned}
		\| \partial_{ t, x }^\beta f \|^2_{ L^2 (\nu) } \lesssim & \mathfrak{e}_m^{2 m^2} \big( \| f \|^2_{ L^2 ( \nu ) } + \| \nu^{ - \frac{1}{2} } \partial_{t,x}^\beta g \|^2_{L^2} \big) \lesssim \mathfrak{e}_m^{2 m^2} \sum_{ \beta' \leq \beta } \| \nu^{ - \frac{1}{2} } \partial_{t,x}^{ \beta' } g \|^2_{L^2} \\
		\lesssim & \mathfrak{e}_m^{2 m^2} \sum_{ \beta' \leq \beta } \| \langle v \rangle^{ \gamma ( l - \frac{1}{2} ) } \M^{ - \frac{q}{2} } \partial_{t,x}^{ \beta' } g \|^2_{L^2} \,.
	\end{aligned}
\end{equation*}
Hence, we further have
\begin{equation}\label{Finish7}
	\begin{aligned}
		\sum_{ \beta' \leq \beta } \| \partial_{t,x}^{\beta'} f \|^2_{ L^2 ( \nu ) } \lesssim \mathfrak{e}_m^{2 m^2} \sum_{ \beta' \leq \beta } \| \langle v \rangle^{ \gamma ( l - \frac{1}{2} ) } \M^{ - \frac{q}{2} } \partial_{t,x}^{ \beta' } g \|^2_{L^2} \,.
	\end{aligned}
\end{equation}

\textbf{\em Step 5. To finish the estimate \eqref{L2-Dtxv} in Lemma \ref{Lmm-SVP-Der}.} Together with $\L f = g$ and the H\"older inequality, one has
\begin{equation}\label{Finish8}
	\begin{aligned}
		\langle \langle v \rangle^{  2 \gamma l} \M^{ - q } \partial_v^\alpha \partial_{t,x}^\beta \L f , \partial_v^\alpha \partial_{t,x}^\beta f \rangle \leq \| \langle v \rangle^{ \gamma ( l + \frac{1}{2} ) } \M^{ - \frac{q}{2} } \partial_v^\alpha \partial_{t,x}^\beta f \|_{ L^2 } \| \langle v \rangle^{ \gamma ( l - \frac{1}{2} ) } \M^{ - \frac{q}{2} } \partial_v^\alpha \partial_{t,x}^\beta g \|_{ L^2 } \,.
	\end{aligned}
\end{equation}
From collecting the bounds \eqref{Finish3}, \eqref{Finish4}, \eqref{Finish5} and \eqref{Finish6}, it follows
\begin{equation}\label{Finish9}
	\begin{aligned}
		& | \langle \langle v \rangle^{  2 \gamma l} \M^{ - q } \partial_v^\alpha [ \partial_{t,x}^\beta , \L ] f , \partial_v^\alpha \partial_{t,x}^\beta f \rangle | \\
		\lesssim & \mathfrak{e}_m^m \| \langle v \rangle^{ \gamma ( l + \frac{1}{2} ) } \M^{ - \frac{q}{2} } \partial_v^\alpha \partial_{t,x}^\beta f \|_{L^2} \\
		& \times \sum_{ 0 \neq \beta' \leq \beta } \Big\{ \| \partial_{t,x}^{\beta - \beta'} f \|_{ L^2 ( \nu ) } + \sum_{|\alpha'| \leq |\alpha|} \| \langle v \rangle^{ \gamma ( l + \frac{1}{2} ) } \M^{ - \frac{q}{2} } \partial_v^{ \alpha' } \partial_{t,x}^{ \beta - \beta' } f \|_{ L^2 } \Big\} \,.
	\end{aligned}
\end{equation}
It is further derived from \eqref{Finish1}, \eqref{Finish2} and \eqref{Finish9} that
\begin{equation}\label{Finish10}
	\begin{aligned}
		& \langle \langle v \rangle^{  2 \gamma l} \M^{ - q } \partial_v^\alpha \partial_{t,x}^\beta \L f , \partial_v^\alpha \partial_{t,x}^\beta f \rangle \\
		\geq & ( 1 - 4 \eta ) \| \langle v \rangle^{ \gamma ( l + \frac{1}{2} ) } \M^{ - \frac{q}{2} } \partial_v^\alpha \partial_{t,x}^\beta f \|^2_{L^2} \\
		& - C_\eta \sum_{\alpha' < \alpha} \| \langle v \rangle^{ \gamma ( l + \frac{1}{2} ) } \M^{ - \frac{q}{2} } \partial_v^{\alpha'} \partial_{t,x}^\beta f \|^2_{L^2} - C_\eta \sum_{\beta' \leq \beta} \| \partial_{t,x}^{ \beta' } f \|^2_{L^2} \\
		& - C_\eta \mathfrak{e}_m^{2 m} \sum_{ 0 \neq \beta' \leq \beta } \sum_{|\alpha'| \leq |\alpha|} \| \langle v \rangle^{ \gamma ( l + \frac{1}{2} ) } \M^{ - \frac{q}{2} } \partial_v^{ \alpha' } \partial_{t,x}^{ \beta - \beta' } f \|^2_{ L^2 }
	\end{aligned}
\end{equation}
for any sufficiently small $\eta > 0$. By taking $\eta = \frac{1}{8}$, the bounds \eqref{Finish7}, \eqref{Finish8} and \eqref{Finish10} conclude the estimate \eqref{L2-Dtxv}. Then the proof of Lemma \ref{Lmm-SVP-Der} is finished.

\vspace*{3mm}

\noindent{\bf Conflict of interest statement:} We declare that there are no conflicts of interest to this work.


\section*{Acknowledgments}
The author N. Jiang thanks for the valuable comments from Y. Guo after the first draft. This work was supported by National Key R\&D Program of China under the grant 2023YFA1010300. The author N. Jiang is supported by the grants from the National Natural Foundation of China under contract Nos. 11971360 and 11731008, and also supported by the Strategic Priority Research Program of Chinese Academy of Sciences, Grant No. XDA25010404. The author Y.-L. Luo is supported by grants from the National Natural Science Foundation of China under contract No. 12201220 and the Guang Dong Basic and Applied Basic Research Foundation under contract No. 2021A1515110210.

\bigskip

\bibliography{reference}

\end{document}